\documentclass[12pt,reqno]{amsbook}
\usepackage{amscd}
\usepackage{epsfig}
\usepackage{float}
\usepackage{url}
\newcommand{\Z}{{\mathbb Z}}
\newcommand{\Q}{{\mathbb Q}}
\newcommand{\C}{{\mathbb C}}
\newcommand{\R}{{\mathbb R}}
\renewcommand{\P}{{\mathbb P}}

\newcommand{\Mm}{{\mathbb M}}

\newcommand{\AAA}{{\mathcal A}}

\newcommand{\OO}{{\mathcal O}}

\newcommand{\RR}{{\mathcal R}}

\DeclareMathOperator{\id}{id}

\DeclareMathOperator{\pr}{pr}
\DeclareMathOperator{\rank}{rank}

\DeclareMathOperator{\tr}{tr}

\newcommand{\www}{\widetilde}
\newcommand{\oooo}{\overline}
\newcommand{\uuuu}{\underline}
\newcommand{\iiii}{\infty}
\newcommand{\nnn}{\nabla}
\newcommand{\whhh}{\widehat}
\newcommand{\orrr}{\overrightarrow}
\newcommand{\olll}{\overleftarrow}

\newcommand{\paa}{\partial}
\newcommand{\paaz}{{\partial_z}}
\newcommand{\zdz}{{z\partial_z}}
\newcommand{\xdx}{{x\partial_x}}
\newcommand{\tdt}{{t\partial_t}}

\newcommand{\bsp}{\left(\begin{smallmatrix}} 
\newcommand{\esp}{\end{smallmatrix}\right)} 
\newcommand{\MGcdot}{\,}
\newcommand{\MGcdots}{\cdots}
\newcommand{\MGtimes}{\times\,}   %for split formulas
\newcommand{\Mon}{\text{Mon}}   %for possible later use
\newcommand{\mon}{mon}   %for possible later use
\newcommand{\ini}{ini}  %for possible later use

\begin{document}

\theoremstyle{plain}
\newtheorem{lemma}{Lemma}[chapter]
\newtheorem{definition/lemma}[lemma]{Definition/Lemma}
\newtheorem{theorem}[lemma]{Theorem}
\newtheorem{proposition}[lemma]{Proposition}
\newtheorem{corollary}[lemma]{Corollary}
\newtheorem{conjecture}[lemma]{Conjecture}
\newtheorem{conjectures}[lemma]{Conjectures}

\theoremstyle{definition}
\newtheorem{definition}[lemma]{Definition}
\newtheorem{withouttitle}[lemma]{}
\newtheorem{remark}[lemma]{Remark}
\newtheorem{remarks}[lemma]{Remarks}
\newtheorem{example}[lemma]{Example}
\newtheorem{examples}[lemma]{Examples}
\newtheorem{notation}[lemma]{Notation}

\title[Painlev\'e III (0,0,4,-4),  connections, and movable poles]
{The Painlev\'e III equation\\  of type (0,0,4,-4),
\\ 
its associated 
vector bundles
with isomonodromic connections,
\\ 
and the geometry
\\ 
of the movable poles} 

\author{Martin A. Guest \and Claus Hertling}

\address{Martin A. Guest\\
Department of Mathematics,
Faculty of Science and Engineering,
Waseda University,
3-4-1 Okubo, Shinjuku, Tokyo 169-8555, Japan}
\email{martin@waseda.jp}

\address{Claus Hertling\\
Lehrstuhl f\"ur Mathematik VI, Universit\"at Mannheim, Seminargeb\"aude
A 5, 6, 68131 Mannheim, Germany}
\email{hertling@math.uni-mannheim.de}



\date{}

\begin{abstract}

The paper is about a Painlev\'e III equation 
of  type $P_{III}(D_6)$
and its relation to isomonodromic families of vector bundles on $\P^1$ with
meromorphic connections. 

The purpose of the paper is two-fold:  it offers a conceptual language for the geometrical objects underlying Painlev\'e equations,  and
it offers new results on a particular Painlev\'e III equation 
of  type $P_{III}(D_6)$, which we denote by $P_{III}(0,0,4,-4)$.  This is equivalent to the radial sine (or sinh) Gordon equation and, as such, it appears very widely in geometry and physics.  

The length of the paper is due to the systematic development
of the material in the language of vector bundles with meromorphic connections,
together with their additional structures which take care of all relevant symmetries.  Our motivation is to explain in a leisurely fashion the language of this general theory, by means of a concrete example.  

Complex multi-valued solutions on $\C^\ast$ are the natural context for most
of the paper, but in the last three chapters real solutions
on $\R_{>0}$ (with or without singularities) are addressed. The vector bundles appearing there can be regarded as mixing holomorphic and antiholomorphic
structures.  They are generalizations of variations of Hodge structures,
called TERP structures (for Twistor Extension Real Pairing) and are
related to $tt^*$ geometry,  and harmonic bundles. 
The paper gives a complete picture of semisimple rank 2 TERP structures.
General results on TERP structures can be applied to study the asymptotics
near $0$ and $\iiii$ of real solutions.

Finally,  results about the asymptotics of real solutions near $0$ and near $\iiii$
are combined with results on the global geometry of the moduli spaces of initial data and 
monodromy data. This leads to a new global picture of all 
zeros and poles of all real solutions  of $P_{III}(0,0,4,-4)$ on $\R_{>0}$.
\end{abstract}

\maketitle

\tableofcontents

\setcounter{chapter}{0}

\chapter{Introduction}\label{s1}
\setcounter{equation}{0}

The paper is about a Painlev\'e III equation 
of  type $P_{III}(D_6)$
and its relation to isomonodromic families of vector bundles on $\P^1$ with
meromorphic connections.  This equation is very classical, and indeed many of the results below can be extracted from the literature.  But we are motivated by more recent developments, of which this particular equation provides an illuminating and nontrivial example.  

The purpose of the paper is two-fold:  it offers a conceptual language for the geometrical objects underlying Painlev\'e equations,  and
it offers new results on a particular Painlev\'e III equation 
of  type $P_{III}(D_6)$, which we denote by $P_{III}(0,0,4,-4)$.  This is equivalent to the radial sine (or sinh) Gordon equation and, as such, it appears very widely in geometry and physics.  

The length of the paper is due to the systematic development
of the material in the language of vector bundles with meromorphic connections,
together with their additional structures which take care of all relevant symmetries.  Our motivation is to explain in a leisurely fashion the language of this general theory, by means of a concrete example.  

We emphasize that it is not necessary to read this paper from beginning to end;  different parts of the paper will be of interest to different readers.  However, we recommend that chapter \ref{s1} be read from beginning to end, in order to understand how to find those parts.  This chapter is a guide to the paper and a summary of results.

\noindent{\em Viewpoint of this paper}

The relevant geometrical objects are vector bundles with flat 
meromorphic connections.  By choosing local trivializations one obtains scalar (or matrix) differential equations.  While the sine-Gordon equation may itself look very natural, certain aspects of its solutions appear quite awkward if the underlying geometrical objects are ignored.  This holds even more so for the Painlev\'e equations, whose explicit formulations have historical but no other significance.

It is well known that the (nonlinear) Painlev\'e equations describe isomonodromic deformations of associated linear meromorphic differential equations, that is, families of linear equations whose monodromy data remains constant within the family.  This theory (like the theory of isospectral deformations)  is of great importance as it links nonlinear p.d.e.\ with algebraic methods of integrable systems theory.  It has far-reaching applications to differential geometry (harmonic maps and harmonic bundles), to variations of Hodge structure and their generalizations, and to quantum field theory in physics (most recently, mirror symmetry). However,  the foundation for all this is the underlying geometry of families of vector bundles with meromorphic connections on $\P^1$, rather than the linear equations themselves, and a description of this geometry forms our starting point.

We emphasize that two variables (or rather two spaces) are involved:  the \lq\lq time variable\rq\rq\ $x\in U\subseteq\C$ of the original Painlev\'e equation, and the deformation variable $z\in \P^1$ of the linear system.  Our vector bundles are rank 2 complex vector bundles
\[
\pi: E\to U\times\P^1
\]
whose restrictions $\pi\vert_{\{x\}\times \P^1}$ 
are holomorphic bundles on $\P^1$, or \lq\lq twistors\rq\rq.  

The first half of the article introduces the fundamental  moduli space of vector bundles with meromorphic connections on $\P^1$, then shows how it may be identified with two other spaces.  First we have the (slightly less canonical) space $M^{mon}$ of monodromy data, which arises through the Riemann-Hilbert correspondence.  Then we have  the (much less canonical) space $M^{ini}$ of initial conditions, which arises when local bases are chosen.  The familiar Painlev\'e and sine-Gordon equation appear only at this point.  

We obtain complex analytic isomorphisms
\[
M^{ini} \cong \text{moduli space of bundles} \cong M^{mon}
\]
which are highly nontrivial as they encode the solutions of the equations.  

The second half of the paper addresses in this way the complex multi-valued solutions of $P_{III}(0,0,4,-4)$
on the \lq\lq time variable space\rq\rq\ $\C^*$.  Of particular interest are the poles and zeros of these solutions and the behaviour  of the solutions at $x=0$ and $x=\iiii$.  So far all data is holomorphic.   

In the last four chapters  we consider real solutions (with singularities) 
on $\R_{>0}$. The vector bundles there mix holomorphic and antiholomorphic
structures.  They are generalizations of variations of Hodge structures,
called TERP structures (for Twistor Extension Real Pairing) and are
related to $tt^*$ geometry and harmonic bundles. In this language, the paper gives a complete picture of semisimple rank 2 TERP structures.
It also shows what general results on TERP structures tell us about the asymptotics of solutions
near $0$ and $\iiii$.

In the last chapter, results about the asymptotics near $0$ and $\iiii$
are combined with the global geometry of $M^{ini}$ and $M^{mon}$. This combination leads to a new global picture of all 
zeros and poles of all real solutions on $\R_{>0}$.  Very briefly, this can be stated as follows.   For real solutions $f$ of $P_{III}(0,0,4,-4)$, there are four types of behaviour at $x=0$ and four types of behaviour at $x=\iiii$:
\[
\renewcommand{\arraystretch}{1.3}
\begin{array}{ll}
\text{At $x=0$:}   &  \text{At $x=\iiii$:}   
\\
\text{(i)}\  f>0  &\text{(i)} \ f>0
\\
\text{(ii)}\  f<0    &\text{(ii)} \ f<0
\\
\text{(iii)} \ \text{$f$ has zeros} &\text{(iii)}\  \text{$f$ has zeros/poles}
\\
\text{(iv)} \ \text{$f^{-1}$ has zeros} &\text{(iv)} \ \text{$f^{-1}$ has zeros/poles}
\end{array}
\]
where the \lq\lq zeros\rq\rq\  converge to $0$ and the \lq\lq zeros/poles\rq\rq\   alternate and converge to $\iiii$.  Extrapolating naively to the whole of $\R_{>0}$, this gives 16 possibilities, but  (i) and (ii) cannot occur together for the same $f$, hence just 14 remain.  We prove that {\em precisely these 14 possibilities occur, and there are no other solutions.} That is, there are no solutions with {\it mixed zones} between the specified types at 0 and at $\iiii$. The 14 possibilities correspond to 14 explicit strata of the space $M^{mon}$.

We give most proofs in detail, and sometimes alternative proofs, without striving for the shortest one.  It would be possible to extract shorter proofs of our new results, depending on the background of the intended audience.   

Before giving a more detailed description of the contents, let us reiterate the motivation for our approach.  The main point is that the bundle point of view deals with intrinsic geometric objects.  These become progressively more concrete (but less canonical) as various choices are made.  The intrinsic approach leads to a clear understanding of the space of solutions, and of the global (qualitative) properties of the solutions themselves.  Moreover, it is an approach which could be adapted for other equations.   

We are aware that our approach has some disadvantages. 
The length of the paper and level of detail may be off-putting for the reader who is interested only in a specific aspect of the solutions. Despite this length, the paper is not entirely self-contained; at several key points we rely on references to the literature, notably for some of the asymptotic results.  Another omission is that we have not yet exploited Lie-theoretic language, nor the Hamiltonian or symplectic aspects.    

This project started when the authors attempted to reconcile the (algebraic/complex analytic) approach to variations of Hodge structures via meromorphic connections with the (differential geometric) approach to harmonic maps via loop groups.  The current article is written primarily from the first viewpoint; a supplementary article from the second viewpoint is in preparation.   However, both approaches are linked by the fundamental idea --- the Riemann-Hilbert correspondence --- that the monodromy data of the associated linear o.d.e.\ represents faithfully the solutions of the nonlinear equation.  In the first approach, this monodromy data, or rather the meromorphic connection from which it is derived, is explicitly at the forefront.  In the second approach, the monodromy data lies behind the so called generalized Weierstrass or DPW data, and $z\in\P^1$ is the spectral parameter.  

\noindent{\em Summary of this chapter}

In order to make the results in the paper approachable
and transparent, this chapter 1 is quite detailed.  Unlike the main body of the paper (chapters \ref{s2} -\ref{s18}), it starts in section \ref{s1.1} with the 
Painlev\'e III equations, and explains immediately and concretely the space $M^{ini}$
of initial conditions.  Although this is quite long, it is just a friendly introduction to essentially  well known facts on Painlev\'e III.
Section \ref{s1.2} gives, equally concretely, the space $M^{mon}$ of monodromy data (at this point, without explaining where it comes from).
Section \ref{s1.3} presents the main results on real solutions.  No special knowledge is required to understand these statements.

Section \ref{s1.4} is intended to orient the reader who is interested primarily in knowing where the above results can be found in the paper, and which parts of the paper 
might be most relevant for applications. 
The history of the subject is rich, and some
references and historical remarks are given in section \ref{s1.5}. 

In section \ref{s1.6} we summarize our treatment of vector bundles with meromorphic connections. The terminology is cumbersome,  but it is designed to be used, not just looked at.  It
reflects the properties of these bundles and the additional structures which underlie the symmetries of the equations.   Section \ref{s1.7} goes on to relate these bundles to the important geometrical concepts which motivated the entire project, and which are (for us) the main applications:  TERP structures, variations of Hodge structures, and tt* geometry.  

Finally, in section \ref{s1.8} we list some problems and questions that we leave for future consideration.

\section{A Painlev\'e equation and its space of initial data $M^{ini}$}\label{s1.1}

\noindent
The Painlev\'e III equation with parameters 
$(\alpha,\beta,\gamma,\delta)\in\C^4$
is the second order ordinary differential equation
\begin{equation}\label{1.1}
P_{III}(\alpha,\beta,\gamma,\delta):\ 
f_{xx}= \frac{f_x^2}{f}-\frac{1}{x}f_x+\frac{1}{x}(\alpha f^2+\beta)
+\gamma f^3+\delta\frac{1}{f}\,.
\end{equation}
This is a meromorphic o.d.e.\ whose coefficients have a simple pole at $x=0\in\C$, so we assume from now on that $x\in\C^\ast$.

It is a basic fact of o.d.e.\ theory that for any $x_0\in\C^*=\C\setminus\{0\}$ and any {\it regular initial datum}
$(f_0,\www f_0)\in\C^*\times \C$ there exists a unique holomorphic local solution
$f$ with $f(x_0)=f_0$ and $f_x(x_0)=\www f_0$. The \lq\lq Painlev\'e property\rq\rq\  
is the following nontrivial extension of this:

\begin{theorem}\label{t1.1} (Painlev\'e property)
Any local solution extends to a global multi-valued meromorphic function on $\C^*$.
\end{theorem}

Painlev\'e's proof was incomplete. The Painlev\'e property follows from the 
relation to isomonodromic families (and fundamental results on isomonodromic families) --- see \cite{FN80}, \cite{IN86}, \cite{FIKN06}, or theorem \ref{10.3} below.
Other proofs can be found, for example, in \cite{HL01}, \cite{GLSh02}.

The Painlev\'e property implies that analytic continuations of solutions with regular initial data give meromorphic functions on the universal covering of $\C^*$, i.e.\  their singularities are at worst poles (in fact they are all simple poles).  However the location of such poles varies drastically with the initial data.  In addition, there are other solutions corresponding to \lq\lq singular initial data\rq\rq, and the same applies to their poles as well.  The {\em global} description of all these solutions and their poles presents a nontrivial problem, which (as far as we know) has not so far been addressed in the existing literature.

If one chooses locally a logarithm 
\[
\varphi=2\log f,\  \text{i.e.\ } f=e^{\varphi/2},
\]
then $\varphi$ branches at the poles and zeros of $f$, but equation \eqref{1.1} simplifies to
\begin{equation}\label{1.2}
(\xdx)^2\varphi = 2x(\alpha e^{\varphi/2}+\beta e^{-\varphi/2})+
2x^2(\gamma e^\varphi + \delta e^{-\varphi})
\end{equation}
where $\partial_x$ means $d/dx$.

One sees immediately the following symmetries (here $k\in\Z$):
\begin{eqnarray}
\label{1.3}
\!\!\!\!\!\!\!\!\!
\begin{split}
&f\textup{ and }\varphi+4\pi i k\textup{ are solutions for }
(\alpha,\beta,\gamma,\delta) \iff \\
-&f\textup{ and }\varphi+2\pi i+4\pi i k\textup{ are solutions for }
(-\alpha,-\beta,\gamma,\delta) \iff \\
&f^{-1}\textup{ and }-\!\varphi+4\pi i k\textup{ are solutions for }
(-\beta,-\alpha,-\delta,-\gamma) \iff 
\\
-&f^{-1}\textup{ and }-\!\varphi+2\pi i+4\pi i k\textup{ are solutions for }
(\beta,\alpha,-\delta,-\gamma) . 
\end{split}
\end{eqnarray}

In \cite{OKSK06} four cases are identified:
\[
\begin{array}{lcc}
P_{III}(D_6)    &\quad&\gamma\delta\neq 0
\\
P_{III}(D_7)    & &\gamma=0,\alpha\delta\neq 0 \text{ or } \delta=0,\beta\gamma\neq 0
\\
P_{III}(D_8)  &  &  \gamma=0,\delta=0,\alpha\beta\neq 0
\\
P_{III}(Q)  & &  \alpha=0,\gamma=0  \text{ or } \beta=0,\delta=0
\end{array}
\]
The four parameters can be reduced to 2,1,0,1 parameters (respectively) by
rescaling $x$ and $f$.

In particular,  $P_{III}(D_6)$ can be reduced to $P_{III}(\alpha,\beta,4,-4)$ with 
$(\alpha,\beta)\in\C^2$.  
In this paper we consider only the Painlev\'e III equation $P_{III}(0,0,4,-4)$.
Then \eqref{1.2} becomes the radial sinh-Gordon equation
\begin{equation}\label{1.4}
(\xdx)^2\varphi = 16 x^2\sinh\varphi.
\end{equation}
All four symmetries \eqref{1.3} preserve the space of solutions of $P_{III}(0,0,4,-4)$.

The following lemma makes precise statements about the zeros and poles of solutions
of $P_{III}(0,0,4,-4)$ and gives a meaning to {\it singular initial data} at zeros or poles
of solutions. It generalizes to all $P_{III}(D_6)$ equations. The lemma is known, but for lack of a suitable reference we sketch the proof. 

\begin{lemma}\label{t1.2}
(a) Let 
\[
f(x)=a_1(x-x_0) + a_2(x-x_0)^2 +a_3(x-x_0)^3 +\MGcdots
\]
be a local solution of $P_{III}(0,0,4,-4)$ near $x_0\in\C^*$ with 
$f(x_0)=0$. Then
\begin{equation}\label{1.5}
a_1=\pm 2, \quad a_2=\tfrac12 a_1/x_0,
\end{equation}
there is no restriction on $a_3\in\C$,  and, for $n\ge 4$, $a_n$ is determined by $a_1,\dots,a_{n-1}$. 

(b) For any $\varepsilon_2=\pm 1$ and any $\widehat f_0\in\C$, there exists
a unique holomorphic solution $f$ of $P_{III}(0,0,4,-4)$ with $f(x_0)=0$ and
$\paa_xf(x_0)=-2\varepsilon_2$, $\paa_x^3 f(x_0)=\widehat f_0$.

(c) A solution can have only simple zeros and simple poles.
There are two types of zeros (according to whether $\paa_xf(x_0)=2$ or $-2$) and 
two types of poles (according to whether $\paa_x (f^{-1})(x_0)=2$ or $-2$).
There is a 1--1 correspondence between local solutions with zeros or poles at
$x_0\in\C^*$ and the set of singular initial data
\begin{equation}\label{1.6}
(\varepsilon_1,\varepsilon_2,\widehat f_0)\in\{\pm 1\}\times \{\pm 1\}\times\C
\end{equation}
given by
\begin{equation}\label{1.7}
f^{\varepsilon_1}(x_0)=0,\ \paa_x(f^{\varepsilon_1})(x_0)=-2\varepsilon_2,\ 
\paa_x^3(f^{\varepsilon_1})(x_0)=\widehat f_0.
\end{equation}
\end{lemma}

{\bf Proof:}
Part (a) is obtained by substituting the Taylor series into equation \eqref{1.1}.
In part (b) only the convergence of the resulting formal power series is nontrivial. 
In this paper it is proved indirectly via theorem \ref{t10.3} (b). 
Part (c) follows immediately from the symmetry $f\leftrightarrow f^{-1}$ in \eqref{1.3}.
\hfill$\Box$

For any $x_0\in\C^*$ the space of regular initial data is $\C^*\times \C$ 
with coordinates $(f_0,\www f_0)$ and the conditions
$f(x_0)=f_0$, $\paa_x f(x_0)=\www f_0$.
The space of singular initial data is given in \eqref{1.6} 
with the conditions \eqref{1.7}.
It is a rather easy exercise to find new coordinates such 
that the space of regular
initial data glues to any one component of the space of singular 
initial data to a chart
isomorphic to $\C\times\C$. This can be done simultaneously for all
$x_0\in\C^*$. We explain this next, and state the result in corollary \ref{t1.3}. 

Suppose that $x_0$ is a zero of $f$ (so $\varepsilon_1=1$) 
with $\varepsilon_2=1$ in \eqref{1.5}. Then for $x$ close to $x_0$,
\begin{eqnarray*}
f&=& (-2)(x-x_0) + \tfrac{1}{2}(-2/x_0)(x-x_0)^2 
+\tfrac{1}{6}\whhh f_0(x-x_0)^3+O((x-x_0)^4)\\
&=& (-2)(x-x_0)(1+\tfrac{1}{2}(1/x_0)(x-x_0))+O((x-x_0)^3)\\
&=& (-2)(x-x_0)(1-\tfrac14 (f/x_0))+O((x-x_0)^3),
\end{eqnarray*}
so
\[
(-2)(x-x_0)=f\MGcdot (1+\tfrac14{f}/{x_0})+O((x-x_0)^3).
\]
Thus 
\begin{eqnarray}
\paa_x f &=& (-2)+\tfrac{-2}{x_0}(x-x_0)+\tfrac{1}{2}\whhh  f_0(x-x_0)^2
+O((x-x_0)^3) \nonumber \\
&=& -2+\tfrac{-2}{x}(x-x_0)+(-2)(\tfrac{1}{x_0}-\tfrac{1}{x})(x-x_0)\nonumber \\
&& +\tfrac{1}{2}\whhh f_0(x-x_0)^2 + O((x-x_0)^3)\nonumber \\
&=& -2+\tfrac{-2}{x}(x-x_0)+\tfrac{-2}{x^2}(x-x_0)^2 
+\tfrac{1}{2}\whhh f_0(x-x_0)^2 + O((x-x_0)^3)\nonumber \\
&=& -2+\tfrac{f}{x}(1+\tfrac{f}{4x_0}) + \tfrac{-1}{2x^2}f^2 
+\tfrac{1}{8}\whhh f_0\MGcdot f^2 + O((x-x_0)^3)\nonumber \\
&=& -2+\tfrac{1}{x}f -\tfrac{1}{4x^2}f^2 + \tfrac{1}{8}\whhh f_0\MGcdot f^2
+O((x-x_0)^3).\label{1.8}
\end{eqnarray}
On the space of regular solutions $\C^*\times \C$ the coordinates 
$(f_0,\www f_0)$ can be replaced by the coordinates $(f_0,g_0)$ with
\begin{equation}\label{1.9}
\www f_0 =\frac{2f_0}{x}\MGcdot g_0
\end{equation}
or by the coordinates $(f_0,\www g_0)$ with
\begin{equation}\label{1.10}
g_0=-\frac{x}{f_0}+\frac{1}{2}+\frac{f_0}{2}\www g_0.
\end{equation}
Then \eqref{1.8} shows that $\www g_0$ extends to the set of 
singular initial data in lemma \ref{t1.2} (c) with 
$\varepsilon_1=\varepsilon_2=1$ and identifies the point
$(0,\www g_0)\in\{0\}\times \C$ with the point $(1,1,\whhh f_0)$
with
\begin{equation}\label{1.11}
\www g_0=\frac{-1}{4x_0}+\frac{x_0}{8}\whhh f_0.
\end{equation}

Now corollary \ref{t1.3} follows easily.
There and later the four pairs 
$(\varepsilon_1,\varepsilon_2)\in\{\pm 1\}\times\{\pm 1\}$
are enumerated by $k\in\{0,1,2,3\}$ via the bijection
\begin{equation}\label{1.12}
(1,1)\mapsto 0,\ (-1,1)\mapsto 1,\ (1,-1)\mapsto 2,\ (-1,-1)\mapsto 3.
\end{equation}

%The coordinates $(f_0,\www f_0)$ on $\C^*\times \C$ are replaced by
%\begin{eqnarray}\label{1.9}
%(f_0,g_0)\quad\text{with}\quad \www f_0=\frac{2f_0}{x_0}\MGcdot g_0.
%\end{eqnarray}
%The coordinate $\widehat f_0$ on the component 
%$(\varepsilon_1,\varepsilon_2)\times\C$
%of the set in \eqref{1.6} is replaced by
%\begin{eqnarray}\label{1.10}
%\www g_k\quad \text{with}\quad \varepsilon_2\widehat f_0 = 
%\frac{2}{x_0^2}+\frac{8}{x_0}\MGcdot \www g_k.
%\end{eqnarray}
%These inessential coordinate changes are made in order to have 
%prettier normal forms
%for the isomonodromic families in chapter \ref{s8}.

We may now introduce the first \lq\lq moduli space\rq\rq\ that will be of interest to us:
\[
M_{3FN}^{ini}=M_{3FN}^{reg}\cup M_{3FN}^{sing},
\quad
M_{3FN}^{sing}=\cup_{k=0}^3 M_{3FN}^{[k]}.
\]
The suffix $3FN$ will be explained in chapter \ref{s10} (see also section \ref{s1.5}).

\begin{corollary}\label{t1.3}
The set of regular and singular initial data for local solutions of 
$P_{III}(0,0,4,-4)$ for all $x_0\in\C^*$ is the 3-dimensional algebraic manifold 
$M_{3FN}^{ini}$ whose four affine charts $\C^*\times\C\times \C$ have
coordinates $(x_0,f_k,\www g_k)$ for $k\in\{0,1,2,3\}$.

Each chart consists of the set $M_{3FN}^{reg}$ of regular initial 
data and of one of the components $M_{3FN}^{[k]}$ of the set
$M_{3FN}^{sing}=\cup_{k=0}^3 M_{3FN}^{[k]}$ of singular initial data.
In any chart
$M_{3FN}^{reg}\cong \C^*\times\C^*\times\C$ and
$M_{3FN}^{[k]}=\C^*\times \{0\}\times\C$. 
The charts are related by
\begin{eqnarray}\label{1.13}
g_k&:=&-\frac{x_0}{f_k}+\frac{1}{2}+\frac{f_k}{2}\MGcdot \www g_k,\\ \label{1.14}
(f_0,g_0)&=&(f_1^{-1},-g_1)=(-f_2,g_2)=(-f_3^{-1},-g_3).
\end{eqnarray}
\end{corollary}

\eqref{1.14} comes from the symmetries in \eqref{1.3} and from 
\eqref{1.12} and \eqref{1.9}. 
\eqref{1.13} comes from \eqref{1.10}.

$M_{3FN}^{ini}$ is in fact the space of initial data of Okamoto \cite{Ok79}
for $P_{III}(0,0,4,-4)$.  He starts with a naive compactification of $M_{3FN}^{reg}$,
blows it up several times and then takes out a certain divisor.
The factor $f_0^2$ before $\www g_0$ in
\begin{eqnarray*}
\www f_0=\frac{2f_0}{x_0}\MGcdot g_0 = -2 + \frac{f_0}{x_0} + \frac{f_0^2}{2}\MGcdot \www g_0
\end{eqnarray*}
reflects the fact that the component 
$M_{3FN}^{[0]}$ is obtained by two subsequent blowing ups.
The complement of $M_{3FN}^{ini}$ in Okamoto's compactification is a divisor of type
$\www D_6$, which is the origin of the name $P_{III}(D_6)$. 
This compactification and divisor are not obvious from corollary \ref{t1.3}.

Descriptions of $M_{3FN}^{ini}$ by four affine charts (in fact, for the space of initial data
of all $P_{III}(D_6)$ equations) can be found in \cite{MMT99} and in \cite{Te07},
but they do not make explicit lemma \ref{t1.2} and the notion of singular initial data.

Up to here everything is essentially well known, but the following section contains some new results.

\section{The space of monodromy data $M^{mon}$}\label{s1.2}

\noindent
It is well known that Painlev\'e equations describe isomonodromic deformations of systems of linear differential equations. This permits a very different and very fruitful
point of view on $M_{3FN}^{ini}$ and the solutions of $P_{III}(0,0,4,-4)$.
It gives rise to a holomorphic isomorphism 
\[
\Phi_{3FN}:M^{ini}_{3FN}\to M_{3FN}^{mon}
\]
where $M_{3FN}^{mon}$ is the space of monodromy data.
It is treated in this section as a black box, but will be a major theme in the rest of this paper.

Let us introduce
\begin{eqnarray}\label{1.15}
\begin{split}
V^{mat}&:= \{(s,B)\in \C\times SL(2,\C)\, |\, 
B=
%\begin{pmatrix}
\bsp
b_1&b_2\\-b_2&b_1+sb_2
\esp
%\end{pmatrix}
\}
\\
\cong V^{mon}&:= \{(s,b_1,b_2)\in\C^3\, |\, b_1^2+b_2^2+sb_1b_2=1\}.
\end{split}
\end{eqnarray}
We have an automorphism     
\begin{eqnarray}\label{1.16}
\begin{split}
m_{[1]}^2:\C\times V^{mat}&\to \C\times V^{mat}\\
 (\beta,s,B) &\mapsto (\beta+4\pi i, s, \Mon_0^{mat}(s)^{-2}\MGcdot B)
\end{split}
\end{eqnarray}
where
\begin{eqnarray}\label{1.17}
\Mon_0^{mat}(s):= S^t\MGcdot S^{-1}=\begin{pmatrix}1&-s\\s&1-s^2\end{pmatrix},
\quad
S:=\begin{pmatrix}1&s\\0&1\end{pmatrix}.
\end{eqnarray}
(The significance of these matrices will be explained later.) 
Define
\begin{equation}\label{1.18}
M_{3FN}^{mon}:= \C\times V^{mat}/\langle m_{[1]}^2\rangle,
\end{equation}
where $\langle m_{[1]}^2\rangle$ is the group (isomorphic to $\Z$) generated by $m_{[1]}^2$.
While $V^{mat}$ and $\C\times V^{mat}$ are algebraic manifolds, 
$M_{3FN}^{mon}$ is only an analytic manifold.
But the fibres $M_{3FN}^{mon}(x)$ of the natural projection
\begin{eqnarray}\label{1.19}
pr_{3FN}^{mon}:M_{3FN}^{mon}\to \C^*,\quad
[(\beta,s,B)]\mapsto \tfrac12 e^{-\beta/2}=x
\end{eqnarray}
are algebraic, as any choice of $\beta$ with $\tfrac12 e^{-\beta/2}=x$
yields an isomorphism
\begin{eqnarray}\label{1.20}
pr_{mat}:M_{3FN}^{mon}(x)&\to& V^{mat},\\ 
{}[(\www \beta,s,\www B)]{}&\mapsto& \textup{ the pair }(s,B)\textup{ with }
{}[(\www \beta,s,\www B)]{}=[(\beta,s,B)].\nonumber
\end{eqnarray}
The trivial foliation on $\C\times V^{mat}$ over $\C$ with leaves
$\C\times \{(s,B)\}$ induces a foliation of $M_{3FN}^{mon}$ over $\C^*$.
A leaf through $[(\beta,s,B)]$ has finitely many branches if and only if 
$\Mon_0^{mat}(s)^2$ has finite order, and then the number of branches is
this order.
  
Theorem \ref{t1.4} below is a key result of this paper.  
It follows from theorem \ref{t10.3}, 
which in turn depends on the Riemann-Hilbert correspondence of chapter \ref{s2}.
Its origin and relation to results in the literature
are discussed in section \ref{s1.6}.

\begin{theorem}\label{t1.4}
The relation of $P_{III}(0,0,4,-4)$ to isomonodromic families gives rise to a
holomorphic isomorphism
\begin{eqnarray}\label{1.21}
\Phi_{3FN}:M_{3FN}^{ini}\to M_{3FN}^{mon}
\end{eqnarray}
such that any local solution of $P_{III}(0,0,4,-4)$ is mapped to a local leaf
in $M_{3FN}^{mon}$.
\end{theorem}

This means that the restriction of the meromorphic function $f_0$ on
$M_{3FN}^{ini}$ to one leaf in $M_{3FN}^{mon}$ (composed with the projection of this leaf
to $\C^*$) gives a multi-valued solution of $P_{III}(0,0,4,-4)$ and that all
solutions arise in this way. 
If a local solution $f$ is defined in a simply connected subset $U$ of $\C^*$
then, for any choice of a lift $U^\prime\subset \C$ of $U$ to the universal covering
$\C\to\C^*,\ \beta\mapsto \tfrac12 e^{-\beta/2}$, we obtain an element
 $(s,B)$ of $V^{mat}$.  Here $s$ is unique,
and $B$ depends on the choice of $U^\prime$. 

This $(s,B)$ encodes
the geometry of the solution $f$ in a very nice way. One instance is
the following obvious corollary.

\begin{corollary}\label{t1.5}
A global multi-valued solution has finitely many branches if and only if 
$\Mon_0^{mat}(s)^2$ has finite order, and then the number of branches is this order.
\end{corollary}

Less obvious are the following statements, which are amongst our main results:

(1) The pair $(s,B)$ and the asymptotics of $f$ for
$x\to 0$ determine one another (see section \ref{s1.4}).

(2) In the case of a real solution  of
$P_{III}(0,0,4,-4)$ on $\R_{>0}$, one can read off the sequence of zeros and poles of $f$
from the pair $(s,B)$ with $\beta\in\R$ (section
\ref{s1.3} and chapter \ref{s18}).

We shall deal with multi-valued solutions $f$ on $\C^*$ of $P_{III}(0,0,4,-4)$
where one branch over $\C-\R_{\leq 0}$ is distinguished.
The choice of $\beta$ with 
$\beta\in (-2\pi,2\pi)$ 
for this branch yields
a unique $(s,B)\in V^{mat}$, and the solution is denoted by
$f_{mult}(.,s,B)$. Then $V^{mat}$ parametrizes all multi-valued solutions
with a distinguished branch.

\section{Real solutions of $P_{III}(0,0,4,-4)$ on $\R_{>0}$}\label{s1.3}

\noindent
In this section we sketch the new results concerning real solutions  of $P_{III}(0,0,4,-4)$ on $\R_{>0}$.
Before this is formulated in theorem \ref{t1.6} below,
we discuss the subspaces $M_{3FN,\R}^{ini}$ and $M_{3FN,\R}^{mon}$ corresponding to initial conditions and monodromy data of real solutions.

$M_{3FN,\R}^{ini}$ is the real semialgebraic submanifold of $M_{3FN}^{ini}$
which is obtained from corollary \ref{t1.3} simply by replacing 
$x_0\in\C^*$ by $x_0\in\R_{>0}$ and the spaces $\C$, $\C^*$ by $\R$, $\R^*$.
Thus 
\begin{eqnarray}
M_{3FN,\R}^{ini} &=& M_{3FN,\R}^{reg}\cup
\left( \cup_{k=0}^3 M_{3FN,\R}^{[k]}\right)
=M_{3FN,\R}^{reg}\cup M_{3FN,\R}^{sing}, \nonumber
\\
M_{3FN,\R}^{reg}&\cong& \R_{>0}\times \R^*\times \R,\nonumber 
\\
M_{3FN,\R}^{[k]}&\cong& \R_{>0}\times \{0\}\times\R,\label{1.22}
\\
M_{3FN,\R}^{reg}\cup M_{3FN,\R}^{[k]}&\cong& \R_{>0}\times\R\times \R
\quad\textup{for any }k\in\{0,1,2,3\},\nonumber
\end{eqnarray}
and these last four spaces are charts which cover $M_{3FN,\R}^{ini}$.
The space $M_{3FN,\R}^{reg}$ has two connected components, one for $f(x_0)>0$ and one for $f(x_0)<0$. They are supplemented by the four
real hypersurfaces of singular initial data, which give the two types of zeros 
and poles.

\begin{figure}[h]
\begin{center}
\includegraphics[width=0.8\textwidth]{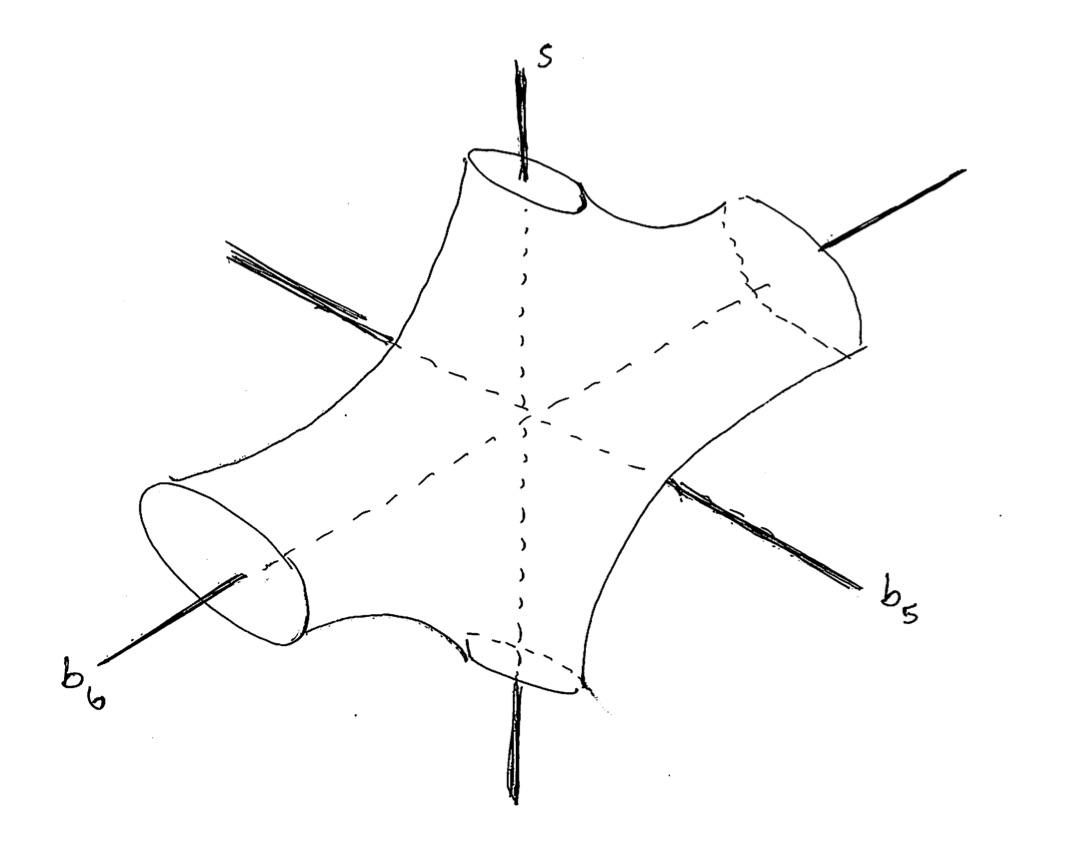} 
\end{center}
\caption{The space $V^{mat,\R}$.}\label{MGpic1}
\end{figure}

A real solution  of $P_{III}(0,0,4,-4)$ on $\R_{>0}$ is the restriction to 
$\R_{>0}$ of the distinguished branch of a multi-valued solution 
$f_{mult}(.,s,B)$. The space of pairs $(s,B)$ corresponding
to real solutions on $\R_{>0}$ is (see theorem \ref{t15.5} (a)) 
\begin{eqnarray}\label{1.23}
\begin{split}
V^{mat,\R}&:= \{(s,B)\in V^{mat}\, |\, (s,B)=(\oooo s,\oooo B^{-1})\}
\\
&= \{(s,B)\in V^{mat}\, |\, s,b_5,b_6\in\R\}
\\
\cong V^{mon,\R}&:= \{(s,b_5,b_6)\in\R^3\, |\, b_5^2+(\tfrac14 s^2-1)b_6^2=1\}
\end{split}
\end{eqnarray}
where 
$
b_5:=b_1+\tfrac12 s b_2, \quad b_6:=ib_2.
$
This surface is illustrated in Fig.\ \ref{MGpic1}.

The isomorphism $\Phi_{3FN}:M_{3FN}^{ini}\to M_{3FN}^{mon}$ restricts to
a real analytic isomorphism
\begin{equation}\label{1.25rev}
\Phi_{3FN,\R}:M_{3FN,\R}^{ini}\to M_{3FN,\R}^{mon}.
\end{equation}
Restricted to real numbers, the covering $\R\to\R_{>0},\ \beta\mapsto  \tfrac12 e^{-\beta/2}$ is bijective, so we have a natural isomorphism
\begin{equation}\label{1.24rev}
\psi_{mat}:M_{3FN,\R}^{mon} \stackrel{\cong}{\longrightarrow}
\R_{>0}\times V^{mat,\R}.
\end{equation}

At the end of the paper, in  remarks \ref{t18.5} and \ref{t18.7},
pictures are presented which represent the images in $M_{3FN,\R}^{mon}$
of the hypersurfaces $M_{3FN,\R}^{[k]}$.
These pictures (the one in remark \ref{t18.7} is conjectural) are consistent
with, and are derived from, theorem \ref{t1.6}. 

\begin{figure}[h]
\begin{center}
\includegraphics[width=1.0\textwidth]{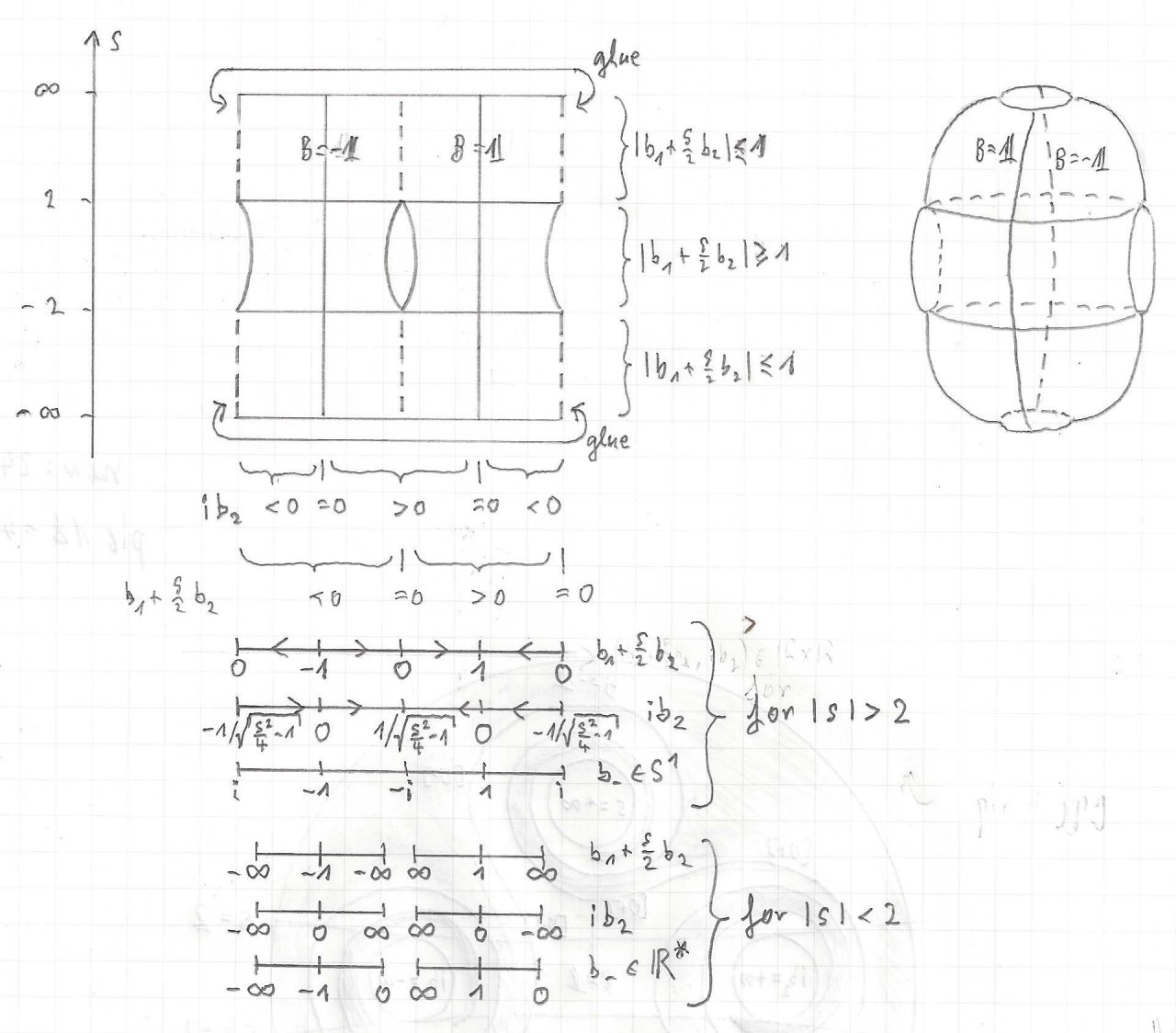} 
\end{center}
\caption{Stratification of the space $V^{mat,\R}$.}\label{pic1}
\end{figure}

The four strata
\begin{eqnarray}\label{1.26}
\begin{split}
\text{(i)}\quad  &|s|\leq 2, \ b_5\geq 1
\\
\text{(ii)}\quad  &|s|\leq 2, \ b_5\leq -1
\\
\text{(iii)}\quad  &s>2
\\
\text{(iv)}\quad  &s<-2
\end{split}
\end{eqnarray}
of $V^{mat,\R}$ will be important for the behaviour of solutions near 0.
The four strata
\begin{eqnarray}\label{1.27}
\begin{split}
\text{(i)}\quad  &B={\bf 1}_2
\\
\text{(ii)}\quad &B=-{\bf 1}_2
\\
\text{(iii)}\quad &b_6<0
\\
\text{(iv)}\quad &b_6>0
\end{split}
\end{eqnarray}
will be important for the behaviour of solutions near $\iiii$.
Taking all possible intersections we obtain 14 strata, as illustrated in
Fig.\ \ref{pic1}.
%%{\sc Later 2 picture of $V^{mat,\R}$}
%\includegraphics[width=1.0\textwidth]{pic1.jpg} 

The following notation allows a precise formulation of the behaviour of the solutions near $0$ and $\iiii$.
Let $f=f_{mult}(.,s,B)|_{\R_{>0}}$ be a real solution  of
$P_{III}(0,0,4,-4)$ on $\R_{>0}$. A zero $x_0\in\R_{>0}$ of $f$ with $\paa_x f(x_0)=\pm 2$
is denoted $[0\pm]$, a pole $x_1\in\R_{>0}$ of $f$ with $\paa_x(f^{-1})(x_1)=\pm 2$
is denoted $[\iiii\pm]$. 

If a $y_0\in\R_{>0}$ exists such that $f|_{(0,y_0]}$ is positive then $f$
and $f|_{(0,y_0]}$ have type $\olll{>0}$. 
If a $y_0\in\R_{>0}$ exists such that $f(y_0)<0$ and $f|_{(0,y_0]}$
has infinitely many zeros $x_1,x_2,x_3,\dots$ with $x_1>x_2>x_3>\dots$
and such that $x_k$ is of type $[0-]$ for odd $k$ and of type $[0+]$ for even $k$,
then $f$ and $f|_{(0,y_0]}$ have type $\olll{[0+][0-]}$.
The types $\olll{<0}$, $\olll{[0-][0+]}$, $\olll{[\iiii+][\iiii-]}$ and 
$\olll{[\iiii-][\iiii+]}$ for $f$ and $f|_{(0,y_0]}$ are defined analogously.

The types $\orrr{>0}$, $\orrr{<0}$, $\orrr{[0+][\iiii-]}$,
$\orrr{[\iiii-][0+]}$, $\orrr{[0-][\iiii+]}$ and $\orrr{[\iiii+][0-]}$
are defined analogously. Of course $f$ has type $\orrr{[0+][\iiii-]}$
if and only if it has type $\orrr{[\iiii-][0+]}$.
But for $f|_{[y_0,\iiii)}$ the smallest zero or pole is important.

Parts (a) and (b) of theorem \ref{t1.6} below give local information near $0$ and near $\iiii$.
Part (c) shows that there is no intermediate mixed zone:  the 
types near $0$ and near $\iiii$ determine the sequence of zeros and poles completely.

\begin{theorem}\label{t1.6} (Theorem \ref{t18.2} and theorem \ref{t18.4})
Fix a real solution $f=f_{mult}(.,s,B)|_{\R_{>0}}$  of $P_{III}(0,0,4,-4)$ on $\R_{>0}$
and its monodromy data $(s,B)\in V^{mat,\R}$.

(a) The type of $f$ near $\iiii$ depends only on the stratum in \eqref{1.27}
in which $(s,B)$ lies. The following table lists the types.
\begin{eqnarray}\label{1.28}
\begin{array}{c|c}
B={\bf 1}_2 & \orrr{>0}\\
B=-{\bf 1}_2 & \orrr{<0}\\ 
%ib_2<0 & \orrr{[\iiii+][0-]}\textup{ and }\orrr{[0-][\iiii+]}\\ 
b_6<0 & \orrr{[\iiii+][0-]}\textup{ and }\orrr{[0-][\iiii+]}\\ 
%ib_2>0 & \orrr{[\iiii-][0+]}\textup{ and }\orrr{[0+][\iiii-]}
b_6>0 & \orrr{[\iiii-][0+]}\textup{ and }\orrr{[0+][\iiii-]}
\end{array}
\end{eqnarray}

(b) The type of $f$ near $0$ depends only on the stratum in \eqref{1.26}
in which $(s,B)$ lies. The following table lists the types.
\begin{eqnarray}\label{1.29}
\begin{array}{c|c}
%|s|\leq 2,\ b_1+\frac{s}{2}b_2\geq 1 & \olll{>0}\\ 
|s|\leq 2,\ b_5\geq 1 & \olll{>0}\\ 
%|s|\leq 2\ \ b_1+\frac{s}{2}b_2\leq -1 & \olll{<0}\\
|s|\leq 2,\ b_5\leq -1 & \olll{<0}\\
s>2 & \olll{[0+][0-]}\textup{ and }\olll{[0-][0+]}\\ 
s<-2 & \olll{[\iiii+][\iiii-]}\textup{ and }\olll{[\iiii-][\iiii+]}\\
\end{array}
\end{eqnarray}

(c) There exists some $y_0\in\R_{>0}$ with $f(y_0)\neq 0$ such that $f|_{(0,y_0]}$
is of type (i), (ii), (iii), or (iv) near $0$ and $f|_{[y_0,\iiii)}$ is of type (i), (ii), (iii), or (iv) near $\iiii$.
Thus the sequence of zeros and/or poles on $\R_{>0}$ of $f$ is completely
determined by the
14 strata.
\end{theorem}

Parts (a) and (b) are known, but part (c) is new.  

For part (a) there are several sources. \cite{MTW77} studies the solutions
of type $\orrr{>0}$ or $\orrr{<0}$, but does not identify $(s,B)$.
\cite[ch.\  11]{IN86} obtains $B=\pm{\bf 1}_2$ for these solutions in the case 
$|s|\leq 2$, and it shows that all solutions with $B\neq\pm{\bf 1}_2$
and $|s|\leq 2$ have infinitely many zeros or poles near $\iiii$.
Probably the arguments in \cite[ch.\  11]{IN86} do not require $|s|\leq 2$.
Anyway, the same without the restriction to $|s|\leq 2$ follows from the 
characterization of nilpotent orbits of TERP structures by mixed TERP structures
in \cite[theorem 9.3]{HS07} and \cite[corollary 8.15]{Mo11b} (see section \ref{s1.7} and chapter \ref{s17}).

Part (b) follows from the asymptotic formulae for $x\to 0$ in \cite{Ni09}
which are rewritten in theorem \ref{t12.4}.

Part (c) follows from a combination of (a) and (b) with basic properties of the space $M_{3FN,\R}^{ini}$
and the isomorphism $\Phi_{3FN,\R}:M_{3FN,\R}^{ini}\to\R_{>0}\times V^{mat,\R}$.
One needs essentially only that $M_{3FN,\R}^{ini}$ has the four smooth 
hypersurfaces of singular initial data and that their images in 
$\R_{>0}\times V^{mat,\R}$ are transversal to the fibres of the projection
$pr_{mat}:\R_{>0}\times V^{mat,\R}\to V^{mat,\R}$, which follows from the
zeros and poles being simple.
It is also convenient to use the existence of the two smooth solutions
$f_{mult}(.,0,\pm{\bf 1}_2)=\pm 1$.

The argument for (c), which combines local asymptotic information near 0
and $\iiii$ with the global geometry of the moduli spaces, seems to be new
in the theory of Painlev\'e III equations.

\section{User's guide to the results}\label{s1.4}

\noindent
This section is for the reader who wishes to see concrete results on solutions of $P_{III}(0,0,4,-4)$ before (or instead of) the general theory of meromorphic connections and the proofs of those results.   
Such readers can restrict attention to chapters
\ref{s1}, \ref{s12}, \ref{s15} and \ref{s18}. 

Section \ref{s1.1} contains basic elementary facts about the equation $P_{III}(0,0,4,-4)$, as well as our formulation of the \lq\lq space of initial conditions\rq\rq.  
Section \ref{s1.2} gives the monodromy data (the derivation of which can be found in chapter  \ref{s2}).  

Theorem \ref{t12.4} in chapter \ref{s12} rewrites the asymptotic formulae for 
$x\to 0$ from  \cite{IN86}, \cite{FIKN06}, \cite{Ni09}.
They are remarkable for several reasons.

First, they yield a direct connection between 
solutions $f$
and their monodromy data $(s,B)$. 
The leading term is determined by and determines the pair $(s,B)$, so it determines
the solution. This is independent of the moduli space viewpoint, although it cannot be formulated as a global statement.  However, if one regards the leading term as an initial datum for a solution \lq\lq at $x=0$\rq\rq, then the space of such data can be identified with our space of monodromy data $V^{mat}$.

Second, the concrete form of the leading term is important.
Theorem \ref{t12.4} describes which leading terms are possible. 
Consider the stratification 
$V^{mat}= V^{mat,a}\cup V^{mat,b+}\cup V^{mat,b-}\cup V^{mat,c+}\cup V^{mat,c-}$, where
\begin{eqnarray}\label{1.30}
\begin{split}
%V^{mat,a}&:= \{(s,B)\in V^{mat}\, |\, s\in \C-[-2,2]\},  \\
V^{mat,a}&:= \{(s,B)\in V^{mat}\, |\, s\in \C-(\R_{\le-2}\cup\R_{\ge2})\},  \\
V^{mat,b\pm}&:= \{(s,B)\in V^{mat}\, |\,  \pm s\in \R_{>2}\},  \\
V^{mat,c\pm}&:= \{(s,B)\in V^{mat}\, |\, s=\pm 2\}. 
\end{split}
\end{eqnarray}
If $(s,B)\in V^{mat,a}$ then 
%$f_{mult}(.,s,B)$ 
$f$
has leading term
\begin{eqnarray}\label{1.31}
\frac{\Gamma(\frac{1}{2}-\alpha_-)}{\Gamma(\frac{1}{2}+\alpha_-)}\MGcdot b_-\MGcdot
\left(\frac{x}{2}\right)^{2\alpha_-}
\end{eqnarray}
where $\alpha_-$ is defined by
\begin{eqnarray}\label{1.32}
e^{-2\pi i\alpha_-}=\lambda_-
:=(1-\tfrac12 {s^2})-s\sqrt{\tfrac14{s^2}-1} \textup{ and }
\Re(\alpha_-)\in (-\tfrac{1}{2},\tfrac{1}{2})
\end{eqnarray}
and
\begin{eqnarray}\label{1.33}
b_- = (b_1+\tfrac12{s}b_2)+\sqrt{\tfrac14{s^2}-1}\ \MGcdot b_2. 
\end{eqnarray}
Here, $\lambda_-$ is an eigenvalue of $\Mon_0^{mat}(s)$, and $b_-$ is an eigenvalue  of $B$ with the
same eigenvector ($B$ and $\Mon_0^{mat}(s)$ commute). The function
$({x}/{2})^{2\alpha_-}$ is a multi-valued function with a distinguished branch.
\eqref{1.31} implies that on any multisector a neighbourhood of $0$ exists on which
%$f_{mult}(.,s,B)$ 
$f$
has no zeros or poles.

If $(s,B)\in V^{mat,b+}$ then 
%$f_{mult}(.,s,B)$ 
$f$
has leading term
\begin{eqnarray}\label{1.34}
-\frac{x}{t^{NI}}\MGcdot \sin\left(2t^{NI}\log \frac12{x} 
-2\arg\Gamma(1+it^{NI})+\delta^{NI}\right)
\end{eqnarray}
where
\begin{eqnarray}\label{1.35}
t^{NI} &:=& \frac{1}{2\pi} \log |\lambda_-|,\\
\delta^{NI}&\in & \C\textup{ with }\Re(\delta^{NI})\in [0,2\pi]\textup{ and }
e^{i\delta^{NI}}=b_- 
\label{1.36}
\end{eqnarray}
with $\lambda_-$ and $b_-$ as above.
In particular, 
%$f_{mult}(.,s,B)$ 
$f$
has in any sufficiently large multisector no poles near 0
and exactly one sequence of zeros converging to 0. They have all approximately the 
same argument $\frac{\log|b_-|}{2t^{NI}}$; see formulae \eqref{12.25} 
and \eqref{12.26}.

If $(s,B)\in V^{mat,c+}$ then 
%$f_{mult}(.,s,B)$ 
$f$
has  leading term
\begin{eqnarray}\label{1.37}
-2x\MGcdot (b_1+\tfrac12{s}b_2)\MGcdot \left( \log\frac{x}{2} - \frac{i\pi}{2}
(b_1+\tfrac12{s}b_2)b_2+\gamma_{Euler} \right) .
\end{eqnarray}
In particular, on any multisector a neighbourhood of $0$ exists on which 
%$f_{mult}(.,s,B)$ 
$f$
has no zeros or poles.

If $(s,B)\in V^{mat,b-\cup c-}$ then $f^{-1}$ has the monodromy data
$$(-s,\bsp 1&0\\0&-1\esp B\bsp 1&0\\0&-1\esp)
\in V^{mat,b+\cup c+},$$
and one can derive the leading term of $f$ from the leading term of $f^{-1}$.

The proof in \cite{Ni09} treats the three cases $(s,B)\in V^{mat,a}, V^{mat,b+}$ and
$V^{mat,c+}$ separately. In chapter \ref{s12} the formulae for 
$(s,B)\in V^{mat,b+}$ and for $(s,B)\in V^{mat,c+}$ are derived from the formula
(with the leading term {\it and} the next term) for $(s,B)\in V^{mat,a}$.
Also in other aspects the formulae are made more transparent here.

Chapter \ref{s13} offers a proof in the language of this paper.
But the heart of the proof is close to that in \cite{Ni09}.

Chapters \ref{s12} and \ref{s13} explain the shape of the formulae,
the only surprise being the constant factors, in particular 
$\Gamma(\frac{1}{2}-\alpha_-)/\Gamma(\frac{1}{2}+\alpha_-)$ in \eqref{1.31},
which come from Hankel functions.  

Chapter \ref{s18}  gives a complete picture of the real solutions on $\R_{>0}$ and their zeros and poles.  This is summarized in section \ref{s1.3} --- see theorem \ref{t1.6}.

The real solutions discussed so far correspond to the real solutions of the radial sinh-Gordon equation 
\[
(\xdx)^2\varphi = 16 x^2\sinh\varphi.
\]
Chapter \ref{s15} discusses also solutions on $\R_{>0}$ with values in $i\R_{>0}$ or
in $S^1$. They correspond to real solutions of the radial sinh-Gordon equation
\[
(\xdx)^2\varphi = -16 x^2\sinh\varphi
\]
with negative sign,
and the radial sine-Gordon equation
\[
(\xdx)^2\varphi = 16 x^2\sin\varphi.
\]
These solutions are much easier to deal with: they are all smooth on $\R_{>0}$, because the condition $\paa_x f(x_0)=\pm 2$ in lemma \ref{t1.2} is impossible
when the values are in $i\R_{>0}$ or in $S^1$.

\section{Related work on Painlev\'e III and meromorphic connections}\label{s1.5}

\noindent
This section summarizes some history, and some relations between this paper and the work of other authors. More information is given in chapter \ref{s9}.

The history of the relation between Painlev\'e III and meromorphic connections
is long. Already shortly after the discovery of the Painlev\'e equations
Garnier found a family of second order linear differential equations with
rational coefficients such that isomonodromic subfamilies are governed 
by solutions of $P_{III}(D_6)$. Okamoto \cite{Ok86} and Ohyama-Okumura \cite{OO06} 
considered this more systematically. 

At the beginning of the 1980's Flaschka-Newell \cite{FN80} and Jimbo-Miwa
\cite{JM81} wrote down first order linear systems of differential equations
in $2\times 2$ matrices whose isomonodromic families are governed
by solutions of $P_{III}(D_6)$. But the appearance of the solutions in the 
matrices was different, and the framework in \cite{FN80} works only for 
$P_{III}(0,0,4,-4)$.  We use the suffix 3FN in $M_{3FN}^{ini}$, $M_{3FN}^{mon}$ to indicate that we are working with this framework (the $3$ indicates the $III$ in $P_{III}(D_6)$).

Unpublished calculations of the second author show the following.
An isomonodromic family as in \cite{FN80} gives four solutions of $P_{III}(0,0,4,-4)$ which are related by
the symmetries in \eqref{1.3}.  On the other hand,  the recipe in \cite{JM81} gives 
for the same isomonodromic family one solution of
$P_{III}(0,4,4,-4)$. The four solutions are mapped to the one solution by the 
$4{:}1$ folding transformation in \cite{TOS05} and \cite{Wi04} which is called
$\psi^{[4]}_{III(D^{(1)}_6)}$ in \cite{TOS05}. 
With hindsight, this folding transformation could have been found in the early 1980's
by comparing \cite{FN80} and \cite{JM81}.

All of \cite{IN86}, \cite[ch.\  13-16]{FIKN06} and \cite{Ni09} work with the
framework in \cite{FN80}. We also work with this framework, but lift it to
vector bundles with additional structures.
This is important in order to obtain the full space $M_{3FN}^{mon}$ of monodromy data.
\cite{IN86} works with a Zariski open set.
\cite{Ni09}, \cite{FIKN06} treat also the complement, but keep it separate from the Zariski open set
(separatrix solutions and generic solutions). 

Lifting the framework to vector bundles 
is also important in order to write down the bundles which correspond to the 
singular initial data and to obtain the isomorphism 
$\Phi_{3FN}:M_{3FN}^{ini}\to M_{3FN}^{mon}$. 
Our $M_{3FN}^{mon}$ is the space of monodromy data of such bundles.

The isomorphism $\Phi_{3FN}:M_{3FN}^{ini}\to M_{3FN}^{mon}$ is a {\it Riemann-Hilbert
isomorphism} in the notation of \cite{FIKN06}, \cite{Mo11a} and other papers
and a {\it monodromy map} in the notation of \cite{Bo01}.
These references contain general results assuring the holomorphicity of such 
Riemann-Hilbert maps. But in order to obtain an isomorphism one has to be careful
about the objects included. 
In many situations, a set of vector bundles with meromorphic connections
and naively chosen conditions is not Hausdorff; one needs the right stability conditions.

The framework of \cite{JM81} for all $P_{III}(D_6)$ equations was taken up in the 
language of vector bundles with connections in \cite{PS09}. 
But there on both sides of the Riemann-Hilbert correspondence only Zariski open subset
were considered.
\cite{PT14} builds on \cite{PS09} and constructs a Riemann-Hilbert isomorphism.
It has some similarities, but also some differences from our isomorphism.
All reducible vector bundles with connections are completely reducible
in our case, and are not completely reducible in their case. 
They have more objects, whereas our objects (the $P_{3D6}$-TEJPA bundles) are richer.

\section{$P_{3D6}$-TEJPA bundles}\label{s1.6}

\noindent
Our fundamental ingredient is the concept of $P_{3D6}$ bundle, a
geometrical object which (on making certain choices) gives rise to the
$P_{III}(D_6)$ equations and their solutions.  In this section we
summarize how this concept will be developed and used.

The isomonodromic families for $P_{III}(0,0,4,-4)$ in \cite{FN80},
as well as the isomonodromic families for all $P_{III}(D_6)$ equations
in \cite{JM81}, are, in the terminology of chapters \ref{s2} and \ref{s4},
isomonodromic families of trace free $P_{3D6}$ bundles.
The appearance of the solutions of the $P_{III}$ equations in the normal
forms  differs in \cite{FN80}, \cite{JM81}; we shall work with the recipe in \cite{FN80}.
Such $P_{3D6}$ bundles have additional symmetries, giving rise to the notion of $P_{3D6}$-TEJPA bundles. 

A $P_{3D6}$ bundle is a $4$-tuple $(H,\nnn,u^1_0,u^1_\iiii)$ (definition \ref{t2.1}).
Here $H\to\P^1$ is a holomorphic vector bundle of rank 2 on $\P^1$, or {\em twistor}.
$\nnn$ is a meromorphic connection whose only poles are at $0$ and $\iiii$, and both have order $2$.
Both poles are semisimple. We denote by  $u^1_0\neq u^2_0$ the eigenvalues of
the endomorphism $[z\nnn_\zdz]:H_0\to H_0$ and by $u^1_\iiii\neq u^2_\iiii$ those of $[-\nnn_{\paa_z}]:H_\iiii\to H_\iiii$. 
The ordering of the eigenvalues is part of the data.

Chapter \ref{s2} introduces also $P_{3D6}$ monodromy tuples (intrinsic formulation of Stokes data) and 
$P_{3D6}$ numerical tuples (classical matrix formulation of Stokes data), and states a Riemann-Hilbert correspondence
between $P_{3D6}$ bundles and $P_{3D6}$ monodromy tuples (theorem \ref{t2.3}).
This goes back to Sibuya \cite{Si67}, \cite{Si90}.    
Modern versions of a Riemann-Hilbert correspondence with parameters which contain
theorem \ref{t2.3} as very special case are in \cite{Bo01}
and in \cite[theorem 4.3.1]{Mo11a}.
The latter is the most general version we know. See also remarks \ref{t2.4}.

Chapter \ref{s3} analyzes reducibility of $P_{3D6}$ bundles
 in terms of their monodromy tuples. This makes statements in \cite{PS09} more explicit.

Chapter \ref{s4} considers isomonodromic families of $P_{3D6}$ bundles.
They have four parameters, the eigenvalues $u^1_0, u^2_0, u^1_\iiii, u^2_\iiii$, 
but in fact they reduce to just one essential parameter.

A $P_{3D6}$ bundle is {\it trace free} if 
if its determinant bundle $\det(H,\nnn)$ with connection is the trivial
rank 1 bundle with trivial flat connection. Then
$u^2_0=-u^1_0, u^2_\iiii=-u^1_\iiii$.

Chapter \ref{s4} discusses 
a solution in \cite{Heu09} of the {\it inverse monodromy problem}, 
in the case of irreducible trace free $P_{3D6}$ bundles.
This problem asks whether, for a given monodromy group, there is an 
isomonodromic family of holomorphic vector bundles on $\P^1$ with 
meromorphic connections, such that each member has the given monodromy group
and such that a generic member is holomorphically trivial.  
By \cite{Heu09} (see theorem \ref{t4.2}), the generic members of any universal
isomonodromic family of trace free $P_{3D6}$ bundles are trivial holomorphic
bundles (they are {\it pure twistors}), and all others, which form a (possibly empty) hypersurface,  are isomorphic
to $\OO_{\P^1}(1)\oplus \OO_{\P^1}(-1)$ 
(they are {\it $(1,-1)$-twistors}; this notation is established in remark \ref{t4.1} (iv)).
For the case of trace free 
$P_{3D6}$-TEP bundles we shall given a short proof of this in theorem \ref{t8.2}
(a) and lemma \ref{t8.5}. 
A more involved proof can be found in \cite{Ni09}.

The $P_{3D6}$ bundles which are relevant to \cite{FN80} 
can be equipped with three types of additional structure:  a pairing $P$,
an automorphism $A$ and an automorphism $J$.
Chapter \ref{s6} introduces the pairing $P$ and $P_{3D6}$-TEP bundles
(T = {\it twistor} $\sim$ holomorphic vector bundle, 
E = {\it extension} $\sim$ meromorphic connection, P = pairing)
and studies their monodromy tuples and isomorphism classes.  A $P_{3D6}$ bundle can
be equipped with a pairing $P$ if and only if it is irreducible or completely reducible and the \lq\lq exponents of formal monodromy\rq\rq\ are all zero. 

In chapter \ref{s7}
we shall see that $P_{3D6}$-TEP bundles can be equipped with automorphisms $A$ and $J$ if and 
only if they are trace free. Then $A$ and $J$ are each unique up to the sign.
The choice of $A$ and $J$ can be considered as a marking 
which fixes these signs.
The choice distinguishes also two out of eight $4$-tuples of certain bases
of the monodromy tuple. The two $4$-tuples of bases differ only by a global sign.
They are useful for the moduli spaces and normals forms.
With the choice of $A$ and $J$,  the trace free $P_{3D6}$-TEP bundle 
becomes a $P_{3D6}$-TEJPA bundle.

Next we introduce moduli spaces $M_{3TJ}^{mon}$ (in chapter \ref{s7}) 
and $M_{3TJ}^{ini}$ (in chapter \ref{s8}), which are analogous to 
$M_{3FN}^{mon}$ and $M_{3FN}^{ini}$. Both can be identified with the moduli space
$M_{3TJ}$ of all isomorphism classes
of $P_{3D6}$-TEJPA bundles.  We obtain a natural analytic isomorphism
\begin{eqnarray}\label{1.38}
\Phi_{3TJ}: M_{3TJ}^{mon}\to M_{3TJ}^{ini}.
\end{eqnarray}
Like $\Phi_{3FN}$, 
it is a Riemann-Hilbert correspondence.  There are natural projections
\begin{equation}\label{1.39}
pr_{3TJ}^{mon} : M_{3TJ}^{mon}\to\C^*\times \C^*,\quad
pr_{3TJ}^{ini} : M_{3TJ}^{ini}\to\C^*\times \C^*
\end{equation}
given by
$(H,\nnn,u^1_0,u^1_\iiii,P,A,J)\to (u^1_0,u^1_\iiii)$.  The space $M_{3FN}^{ini}$ is defined at this point using a choice of normal form of the flat connection; 
only later, in theorem \ref{t10.3}, is $M_{3FN}^{ini}$ identified with the space of initial
data of the solutions of $P_{III}(0,0,4,-4)$.  The derivation of $P_{III}(0,0,4,-4)$ from the normal form of the connection is carried out in chapter \ref{s10}, after some preparation and general remarks on 
Painl\'eve equations in chapter \ref{s9}.

In fact we define $M_{3TJ}^{mon}$ and $M_{3TJ}^{ini}$ first, then construct
the spaces $M_{3FN}^{mon}$ and $M_{3FN}^{ini}$  (in chapter \ref{s10})
as pull-backs of $M_{3TJ}^{mon}$ and $M_{3TJ}^{ini}$ via the diagonal embedding
\begin{equation}\label{1.40}
c^{diag}:\C^*\to \C^*\times \C^*,\quad x\mapsto (x,x).
\end{equation}

Like $M_{3FN}^{mon}$, the space $M_{3TJ}^{mon}$ is the quotient of the
algebraic manifold $\C\times\C^*\times V^{mat}$ by a group action of
$\langle m_{[1]}\rangle\cong\Z$.  It inherits a foliation above $pr_{3TJ}^{mon}$
from this, and the fibres $M_{3TJ}^{mon}(u^1_0,u^1_\iiii)$ of $pr_{3TJ}^{mon}$
are isomorphic as algebraic manifolds to one another and to $V^{mat}$.

As in the case of $M_{3FN}^{ini}$, the fibres $M_{3TJ}^{ini}(u^1_0,u^1_\iiii)$ of 
$pr_{3TJ}^{ini}$ are algebraic surfaces with four affine charts, and $M_{3TJ}^{ini}$
is an algebraic manifold of dimension 4. 
It has three strata (unlike $M_{3FN}^{ini}$, which has five): the open stratum of normal forms of
$P_{3D6}$-TEJPA bundles which are pure twistors, and two smooth hypersurfaces
of $P_{3D6}$-TEJPA bundles which are $(1,-1)$-twistors. In the 1-parameter isomonodromic families of $P_{3D6}$-TEJPA bundles which
correspond to solutions of $P_{III}(0,0,4,-4)$, the $(1,-1)$-twistors correspond
to the zeros and poles of the solutions (theorem \ref{t10.3}).

Matrices which amount to normal forms (but without identifying $P$, $A$ and $J$) 
for the pure twistor $P_{3D6}$-TEJPA bundles
are also given in \cite{FN80}, \cite{IN86}, \cite{FIKN06}, \cite{Ni09}.
But the normal forms in theorem \ref{t8.2} for the $(1,-1)$ twistor $P_{3D6}$-TEJPA bundles are new.
The relations between our approach and \cite{IN86}, \cite{FIKN06}, \cite{Ni09} is discussed in detail in chapters \ref{s11}-\ref{s13}.

The three symmetries $R_1,R_2,R_3$ of $P_{3D6}$-TEJPA bundles introduced in chapter \ref{s7} are supplemented by further symmetries $R_4,R_5$ in chapter \ref{s14}; $R_5$ is a \lq\lq reality condition\rq\rq.  These are used in chapter \ref{s15} to study real solutions of 
$P_{III}(0,0,4,-4)$ on $\R_{>0}$.  Chapter \ref{s18} brings together the results obtained so far to describe the space of all real solutions.

\section{TERP(0) bundles, generalizations of Hodge structures}\label{s1.7}

\noindent
One motivation for the second author to study real solutions 
(possibly with zeros and/or poles) of $P_{III}(0,0,4,-4)$ on $\R_{>0}$  is that they
are equivalent to Euler orbits of semisimple rank 2 TERP(0) bundles.

A pure and polarized TERP(0) bundle generalizes a pure and polarized Hodge structure.
An isomonodromic family of them is called a pure and polarized TERP structure
and generalizes a variation of Hodge structures. TERP structures were defined in
\cite{He03} in order to give a framework for the data in \cite{CV91}, \cite{CV93}, \cite{Du93}.
The semisimple rank 2 case is the simplest essentially new case beyond variations of Hodge structures.

In \cite{CV91}, \cite{CV93}, \cite{Du93} it was noticed that this case is related to
real solutions  of $P_{III}(0,0,4,-4)$ on $\R_{>0}$. In \cite{CV91}, \cite{CV93}
some of the globally smooth and positive solutions (namely, some of the solutions
$f_{mult}(.,s,{\bf 1}_2)|_{\R_{>0}}$ with $|s|\leq 2$) are obtained from the 
data there.

The discussion in chapter \ref{s16} will show that semisimple rank 2 TERP structures
(not necessarily pure and polarized) correspond to 
solutions  of $P_{III}(0,0,4,-4)$ on $\R_{>0}$ with values in $\R$ or in $S^1$.
The solutions with values in $S^1$ are all globally smooth, and their TERP structures
are all pure, but not polarized.

A real solution $f$ on $\R_{>0}$ corresponds to an Euler orbit $\cup_{x>0}\text{TERP}_f(x)$,
a certain 1-parameter isomonodromic family of TERP(0) bundles.
$\text{TERP}_f(x)$ is a $(1,-1)$ twistor if and only if $x$ is a zero or a pole of $f$.
If $f(x)>0$ then $\text{TERP}_f(x)$ is pure and polarized;  if $f(x)<0$ then $\text{TERP}_f(x)$
is pure, but not polarized.

In \cite{Sch73} and \cite{CKS86} nilpotent orbits of Hodge structures are studied,
and a beautiful correspondence with polarized mixed Hodge structures is established.
In \cite{CV91}, \cite{CV93} a renormalization group flow leads to an 
{\it infrared limit} of the Euler orbits of those TERP(0) bundles which are studied
there implicitly. 

This motivated the definition in \cite{HS07} of a {\it nilpotent orbit}
of TERP(0) bundles (definition \ref{t17.3}): 
An Euler orbit $\cup_{x>0}G(x)$ of TERP(0) bundles (definition \ref{t17.1})
is a {\it nilpotent orbit} if for large $x$ the TERP(0) bundle $G(x)$ is pure and
polarized.

It also motivated the conjecture 9.2 in \cite{HS07} that (in the case when the 
formal decomposition of the pole at $0$ is valid without ramification) an Euler orbit
$\cup_{x>0}G(x)$ is a nilpotent orbit if and only if one 
(equivalently: any) TERP(0) bundle $G(x)$ is a {\it mixed TERP structure}. 
The general definition of a mixed TERP structure is given in \cite{HS07}, and
the semisimple case is in definition \ref{t17.6}. 
Roughly, $G(x)$ is a mixed TERP structure if the real structure and Stokes structure
are compatible and if the regular singular pieces of the formal decomposition
of the pole at $0$ induce certain polarized mixed Hodge structures.

The conjecture was proved in \cite{HS07} (the direction $\Leftarrow$ and the regular
singular case of $\Rightarrow$, building on \cite{Mo03}) 
and in \cite{Mo11b} (the general case of the more
difficult direction $\Rightarrow$, building on \cite{Mo11a}).
The semisimple case is formulated in theorem \ref{t17.9}.

In our context here it shows that a real solution $f_{mult}(.,s,B)|_{\R_{>0}}$ 
 of $P_{III}(0,0,4,-4)$ on $\R_{>0}$ is smooth and positive for large $x$ if and
only if $B={\bf 1}_2$ (corollary \ref{t17.10}). 

In \cite{HS07} the concept of {\it Sabbah orbit} is defined. An Euler orbit 
$\cup_{x>0}G(x)$ of TERP(0) bundles $G(x)$ is a Sabbah orbit if for small $x$
the TERP(0) bundle is pure and polarized. And \cite[theorem 7.3]{HS07} characterizes
Sabbah orbits by the property that they induce (in a different way than above) 
polarized mixed Hodge structures. This result gives a characterization of those
real solutions $f_{mult}(.,s,B)|_{\R_{>0}}$ which are smooth and positive for small $x$.
But we do not need it here, as the asymptotic formulae in \cite{Ni09} 
(see theorem \ref{t12.4} and section \ref{s1.4}) give information which is 
more precise and easier to handle.

Finally, a result in \cite{HS11} gives certain Euler orbits $\cup_{>0}G(x)$ 
of TERP(0) bundles such that all $G(x)$ are pure and polarized
(see theorem \ref{t17.11} for the semisimple case).
In our context here, this result applies and gives in fact all globally smooth
and positive solutions  of $P_{III}(0,0,4,-4)$ on $\R_{>0}$
(corollary \ref{t17.12}, remark \ref{t17.13} (i), theorem \ref{t18.2}).

Chapters \ref{s16} and \ref{s18} give by the correspondence with solutions
on $\R_{>0}$ (with values in $\R$ or in $S^1$) a rather complete picture of 
the semisimple rank 2 TERP structures and their Euler orbits.
The applications here of the general results on TERP structures illustrate what
these general results give and how they work.

\section{Open problems}\label{s1.8}

\noindent
We collect here some questions and conjectures.

(I) (Remarks \ref{t9.1} and \ref{t9.2} and section \ref{s1.5})
\cite{JM81} proposes a recipe which relates all solutions of $P_{III}(D_6)$
to isomonodromic families of trace free $P_{3D6}$ bundles (see also 
\cite[ch.\  5]{FIKN06}, \cite{PS09}, \cite{PT14}). 
A different recipe, which applies only for the solutions of $P_{III}(0,0,4,-4)$,
is proposed in \cite{FN80} (see also \cite{IN86}, \cite[ch.\  7-16]{FIKN06}, \cite{Ni09})
and this paper. The $P_{3D6}$ bundles in \cite{FN80} are related in \cite{JM81} 
to solutions of $P_{III}(0,4,4,-4)$. The second author made calculations which show
that these recipes are related by the $4{:}1$ folding transformation in \cite{Wi04}
and \cite{TOS05} which is called $\psi^{[4]}_{III(D_6^{(1)})}$ in \cite{TOS05}
and which maps almost all solutions of $P_{III}(0,0,4,-4)$ to almost all solutions
of $P_{III}(0,4,4,-4)$. 

It seems interesting to study what one obtains from the other two
folding transformations $\psi_{IV}^{[4]}$ and $\psi_{II}^{[2]}$ in \cite{TOS05}
and known recipes for $P_{IV}$ and $P_{II}$.

(II) (Remark \ref{t10.4} (i)) 
By \cite{Ok79} the algebraic surface  
$M_{3FN}^{ini}(x)$ of initial data at $x\in \C^*$ of solutions
of $P_{III}(0,0,4,-4)$ is the complement of a divisor $Y$ of type $\www D_6$ in a 
natural compactification $S$. Can the points in $Y$ be given a meaning in terms of
suitable generalizations of $P_{3D6}$ bundles?  If so, can this be applied
to a better understanding of the solutions of $P_{III}(0,0,4,-4)$?

(III) (Remark \ref{t10.4} (vi)) 
The multi-valued function $f_{mult}$ on $M_{3FN}^{ini}$ in remark \ref{t10.2}
unites all solutions of $P_{III}(0,0,4,-4)$ and depends holomorphically on
the parameters $(s,B)\in V^{mat}$. Can this dependence be controlled
(beyond the asymptotic formulae for small or large $x$)? 
In particular, are there differential equations governing the dependence on $s$ and $B$?
We do not see any, and wonder why not.

For constant $s$ and varying $B$, the monodromy of a family of $P_{3D6}$-TEP bundles
is constant, but the poles along $0$ and $\iiii$ of the covariant derivatives
along $B$ of sections have infinite order. Therefore it is not isomonodromic,
and it does not give differential equations in an obvious way.

(IV) (Conjecture \ref{t18.6} (a)) 
We expect that the hypersurfaces in $\R_{>0}\times V^{mat,a\cup c,\R}$ and in 
$\R_{>0}\times V^{mat,b\cup c,\R}$ which give the $(1,-1)$-twistors (the sheets of zeros and poles of solutions) are convex.

(V) Conjectures \ref{t18.6} (b) and (c) formulate our expectations on the behaviour
of the limits of the sheets of zeros and poles of solutions when $(s,B)$ approach
the holes $s=\pm\iiii$ or $ib_2=\pm\iiii$ in $V^{mat,\R}$. 

(VI) In chapter \ref{s12} we have analyzed and rewritten the asymptotic formulae
as $x\to 0$ for solutions $f_{mult}(.,s,B)$ from \cite{Ni09}.
Similar general formulae for $x\to\iiii$ are missing.
Special cases are in \cite{IN86}, for example Chapter 8 contains asymptotic formulae
as $x\to \iiii$ for solutions on $\R_{>0}$ which are smooth near $\iiii$.
Chapter 11 contains asymptotic formulae as $x\to\iiii$ for the zeros
and poles of real solutions on $\R_{>0}$. 
This should be extended to asymptotic formulae as $x\to\iiii$ 
for the multi-valued solutions $f_{mult}(.,s,B)$ on $\C^*$ and their zeros and poles.

(VII) The matrix 
$T(s)=\bsp 0&1\\-1&s\esp$
in \eqref{14.10} is equivalent to the parameter $s=\tr T$.
The parameters $s$ and $2b_5=\tr B$ in $V^{mat}$ have completely different roles, yet
there is a symmetry between them.
The defining equation for $V^{mat}$ is 
\begin{eqnarray}
0&=&4\det B-4=4b_1^2+4b_2^2+4sb_1b_2-4 \nonumber \\
&=& (2b_1+sb_2)^2-(s^2-4)b_2^2-4 \nonumber\\
&=& ((\tr B)^2-4) + (ib_2)^2 ((\tr T)^2)-4).\label{1.41}
\end{eqnarray}
The map
\begin{eqnarray}\label{1.42}
V^{mat}\dashrightarrow V^{mat},\quad (\tr T,\tr B,ib_2)\mapsto (\tr B,\tr T,\frac{1}{ib_2})
\end{eqnarray}
is birational and exchanges the traces of $T$ and of $B$, 
and it restricts to $V^{mat,\R}$.

Within $V^{mat,\R}$, either $|\tr T|\leq 2$ and $|\tr B|\geq 2 $ (in $V^{mat,a\cup c,\R}$) 
or $|\tr T|\geq 2$ and $|\tr B|\leq 2$ (in $V^{mat,b\cup c,\R}$).
We wonder whether there is more behind this symmetry.

Further,  we wonder about the apparent rotational symmetry by $\frac{\pi}{2}$ in the 
picture in remark \ref{t18.7}. Of course, because of the sizes of the spirals for
small $x_0$ and large $x_0$, one should invert small and large $x_0$.
Are the spirals, i.e.\ the four hypersurfaces in $\R_{>0}\times V^{mat,\R}$, exchanged
by a suitable automorphism which is roughly a rotation by $\frac{\pi}{2}$ in the
picture of $V^{mat,\R}$ and which exchanges small and large $x_0$?
If so, this might relate asymptotic formulae for $x\to 0$ with asymptotic
formulae for $x\to \iiii$.

\section{Acknowledgements}\label{s1.9}

\noindent
This paper was mainly written during several stays of the second author 
at Tokyo Metropolitan University. He thanks TMU for hospitality.
He also thanks Masa-Hiko Saito for fruitful discussions.
The first author thanks Alexander Its and Chang-Shou Lin for explanations
of related material.

The first author was partially supported by JSPS grants A25247005 and A21244004, and by Waseda University grant 2013B-083.  The second author was partially supported by 
DFG grant HE 2287/4-1 and JSPS grant S-11023.

\chapter{The Riemann-Hilbert correspondence for $P_{3D6}$ bundles}\label{s2}
\setcounter{equation}{0}

\noindent
This chapter will formulate the Riemann-Hilbert correspondence for those
holomorphic vector bundles on $\P^1$ with meromorphic connections
which are central for the Painl\'eve III($D_6$) equations.
Everything in this chapter is classical, though presented in the language of bundles. We shall give
references after theorem \ref{t2.3}. 

We define first $P_{3D6}$ bundles, then $P_{3D6}$ monodromy tuples,
then we state the Riemann-Hilbert correspondence in theorem 
\ref{t2.3}. After that, the correspondence
will be explained.  Finally we introduce $P_{3D6}$ numerical tuples, 
a more concrete version of monodromy tuples.

\begin{definition}\label{t2.1}
A {\it $P_{3D6}$ bundle} is a $4$-tuple $(H,\nabla,u^1_0,u^1_\iiii)$ consisting
of the following ingredients.
First, $H\to\P^1$ is a holomorphic vector bundle of rank 2 on $\P^1$ (a \lq\lq twistor\rq\rq),
and $\nabla$ is a meromorphic (hence flat) connection on $H$, which is
holomorphic on $H|_{\C^*}$ but has poles at $0$ and $\infty$
of order 2.  Writing $z$ for the coordinate on $\C\subset\P^1$,
the eigenvalues of the endomorphism $[z\nnn_\zdz]:H_0\to H_0$ (i.e.\ the coefficient of $z^{-2}$ in
a matrix representation of $\nabla_\paaz$) are assumed distinct, and will be denoted $u^1_0,u^2_0\in \C$.  In particular this endomorphism is semisimple.
Similarly, the eigenvalues of the
endomorphism $[-\nnn_\paaz]:H_\iiii\to H_\iiii$ are assumed distinct, and will be denoted
$u^1_\iiii,u^2_\iiii$.
\end{definition}

Distinguishing one eigenvalue at $0$ and one at $\infty$ is part of the definition
of $P_{3D6}$ bundle. For the structure of
$P_{3D6}$-bundles this choice is inessential, but it will facilitate our treatment of
normal forms of $P_{3D6}$ bundles and
families with varying eigenvalues.

Before defining $P_{3D6}$ monodromy tuples, it is necessary to introduce
sectors in $\C^*$ which will be used for the Stokes structures.  
They depend on $u^1_0\ne u^2_0, u^1_\iiii\ne u^2_\iiii\in \C$. At $z=0$ we define
\begin{eqnarray}\label{2.1}
\begin{split}
\zeta_0&:= i(u^1_0-u^2_0)/|u^1_0-u^2_0|\in S^1,\\ 
I^a_0&:=\{z\in S^1\, |\, \Re(\tfrac1z  (u^1_0-u^2_0))<0\},\\ 
I^b_0&:=-I^a_0=\{z\in S^1\, |\, \Re(\tfrac1z  (u^1_0-u^2_0))>0\},\\ 
I^+_0&:= I^a_0\cup I^b_0\cup\{\zeta_0\}=S^1-\{-\zeta_0\},\\ 
I^-_0&:=-I^+_0=I^a_0\cup I^b_0\cup\{-\zeta_0\}=S^1-\{\zeta_0\}.
\end{split}
\end{eqnarray}
At $z=\infty$ we use the
coordinate $\www z:=\frac{1}{z}$ on $\P^1-\{0\}$, so we define
\begin{eqnarray}\label{2.2}
\begin{split}
\zeta_\iiii&:= |u^1_\iiii-u^2_\iiii| / i(u^1_\iiii-u^2_\iiii)\in S^1,\\
I^a_\iiii&:=\{z\in S^1\, |\, \Re(z(u^1_\iiii-u^2_\iiii))<0\},\\ 
I^b_\iiii&:=-I^a_\iiii=\{z\in S^1\, |\, \Re(z(u^1_\iiii-u^2_\iiii))>0\},\\
I^+_\iiii&:= I^a_\iiii\cup I^b_\iiii\cup\{\zeta_\iiii\}=S^1-\{-\zeta_\iiii\},\\
I^-_\iiii&:=-I^+_\iiii=I^a_\iiii\cup I^b_\iiii\cup\{-\zeta_\iiii\}=S^1-\{\zeta_\iiii\}.
\end{split}
\end{eqnarray}
For any (connected and open) subset $I\subset S^1$, we denote by
\[
\whhh I:=\R_{>0}\MGcdot I:=\{z\in \C^*\, |\, z/|z|\in I\}
\]
the corresponding sector in $\C^*$.
The following pictures show the sectors 
for some choice of $u^1_0,u^2_0,u^1_\iiii,u^2_\iiii$ with
$u^1_0=u^1_\iiii=-u^2_0=-u^2_\iiii$.

%{\sc Later 2 pictures, the left for 0, the right for $\iiii$.}
\includegraphics[width=0.9\textwidth]{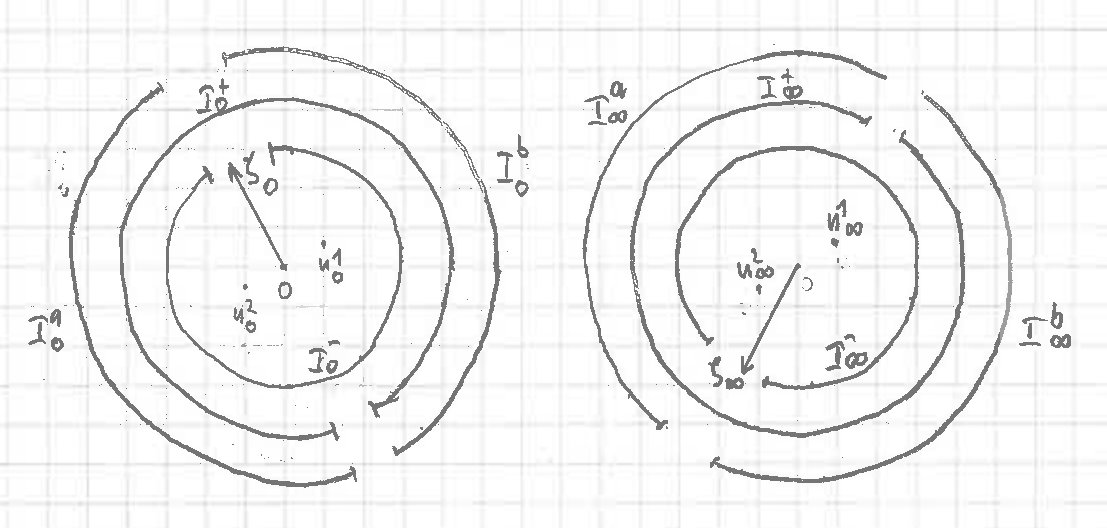} 

\begin{definition}\label{t2.2}
A {\it $P_{3D6}$ monodromy tuple} is a $17$-tuple
\begin{eqnarray}\label{2.3}
(u^j_0,u^j_\iiii,\alpha^j_0,\alpha^j_\iiii,L,L^{+j}_0,L^{-j}_0,
L^{+j}_\iiii,L^{-j}_\iiii (j=1,2)).
\end{eqnarray}
Here $u^j_0,u^j_\iiii,\alpha_0^j,\alpha_\iiii^j \in \C$, and 
$u^j_0,u^j_\iiii$
satisfy $u^1_0\neq u^2_0,u^1_\iiii\neq u^2_\iiii$; 
they will represent the eigenvalues of the \lq\lq pole parts\rq\rq\ of a connection.  
The various bundles will represent the decompositions at irregular 
singularities given by the classical formal theory, 
and $\alpha^j_0,\alpha^j_\iiii$ will represent the exponents of formal monodromy.  
More precisely, they are assumed to satisfy the following conditions. 
$L$ is a flat 
rank 2 vector bundle on $\C^*$.
$L^{\pm j}_0$ are flat rank 1 subbundles of $L|_{\whhh I^{\pm}_0}$, and 
$L^{\pm j}_\iiii$ are flat rank 1 subbundles of $L|_{\whhh I^{\pm}_\iiii}$,
with
\begin{eqnarray}\label{2.4}
L|_{\whhh I^{\pm}_0}=L^{\pm 1}_0\oplus L^{\pm 2}_0,\quad
L|_{\whhh I^{\pm}_\iiii}=L^{\pm 1}_\iiii\oplus L^{\pm 2}_\iiii,\\ \label{2.5}
L^{+1}_0|_{\whhh I^a_0}=L^{-1}_0|_{\whhh I^a_0},\quad
L^{+2}_0|_{\whhh I^b_0}=L^{-2}_0|_{\whhh I^b_0},\\ \label{2.6}
L^{+1}_\iiii|_{\whhh I^a_\iiii}=L^{-1}_\iiii|_{\whhh I^a_\iiii},\quad
L^{+2}_\iiii|_{\whhh I^b_\iiii}=L^{-2}_\iiii|_{\whhh I^b_\iiii}.
\end{eqnarray}
Denote (for $0$ and $\iiii$, this is abbreviated by $0/\iiii$) by 
$$t^{aij}_{0/\infty}: L^{+i}_{0/\iiii}|_{\whhh I^a_{0/\iiii}}
\to L^{-j}_{0/\iiii}|_{\whhh I^a_{0/\iiii}}$$
the projection induced by the restriction to $\whhh I^a_{0/\iiii}$ of the splitting 
$L|_{\whhh I^-_{0/\iiii}}=L^{-1}_{0/\iiii}\oplus L^{-2}_{0/\iiii}$,
and similarly $t^{bij}_{0/\iiii}$.
Then $t^{a11}_{0/\iiii}=\id, t^{a12}_{0/\iiii}=0$, thus
$t^{a22}_{0/\iiii}$ is an isomorphism,
and $t^{b22}_{0/\iiii}=\id, t^{b21}_{0/\iiii}=0$, thus
$t^{b11}_{0/\iiii}$ is an isomorphism. Glueing by $t^{ajj}_{0/\iiii}$
and $t^{bjj}_{0/\iiii}$ (for fixed $j$) gives flat rank 1 bundles 
$L^{+j}_{0,\infty}\cup_{t's} L^{-j}_{0/\iiii}$ on $\C^*$.
They and the numbers $\alpha_{0/\iiii}^j$ are linked by the
compatibility conditions
\begin{eqnarray}\label{2.7}
\begin{split}
\textup{eigenvalue of monodromy of }L^{+j}_0\cup_{t's} L^{-j}_0
%\textup{ is }
=
e^{-2\pi i\alpha^j_0},\\ 
\textup{eigenvalue of monodromy of }L^{+j}_\iiii\cup_{t's} L^{-j}_\iiii
%\textup{ is }
=e^{2\pi i\alpha^j_\iiii}.
\end{split}
\end{eqnarray}
\end{definition}

\begin{theorem}\label{t2.3} (Riemann-Hilbert correspondence)

(a)There is a canonical $1{:}1$ correspondence (described below) between 
$P_{3D6}$ bundles and $P_{3D6}$ monodromy tuples.

(b) It extends to a canonical $1{:}1$ correspondence between
holomorphic families of $P_{3D6}$ bundles and holomorphic
families of $P_{3D6}$ monodromy tuples (over the same complex base manifold).
\end{theorem}

\begin{remarks}\label{t2.4}
(i) More commonly, the correspondence (a) is formulated for meromorphic
bundles with meromorphic connections, not for holomorphic bundles,
e.g. \cite{Ma83a}. But in our semisimple case one can extract 
from \cite{Si90} the version for holomorphic bundles.
Below we shall follow \cite[II 5]{Sa02} for the treatment of holomorphic
bundles.

(ii) Several strong generalizations of the classical correspondence
for a meromorphic bundle with a meromorphic connection are given in 
\cite{Sa13}: for (not necessarily regular) holonomic D-modules
on curves \cite[theorem 5.14]{Sa13}, and for good meromorphic
connections in higher dimensions \cite[theorem 10.8 and theorem 12.16]{Sa13}.

(iii) It is also more common to formulate (a) alone,
i.e.\ for single bundles, rather than families.
Indeed, in the case of families, one has to be careful with the 
regular singular pieces which have to be encoded in the monodromy tuple.
But in our semisimple case this is not a problem and can also be
extracted from \cite{Si90}.
A remark on families is also made in \cite{Ma83a}, but for 
meromorphic bundles. 

(iv) A correspondence for families of holomorphic vector bundles is 
also important in \cite{FIKN06}, but there it is discussed in detail
only in an example. 

\cite{Bo01} gives a beautiful description of a Riemann-Hilbert correspondence
(he calls it the monodromy map) which covers theorem \ref{t2.3}.
He restricts to the semisimple case, and considers two versions,
one with and one without framings on both sides.
The version with framings has the advantage that all spaces are 
manifolds. The central statements are in \cite[ch.\  3,7]{Bo01}
after proposition 3.7 and after definition 7.3.

(v) The most general version of the Riemann-Hilbert correspondence
for holomorphic vector bundles with meromorphic connections that
we know is given in \cite[Theorem 4.3.1]{Mo11a}.  Theorem \ref{t2.3} is a very special case of this.
Mochizuki considers families of vector bundles on complex manifolds
of arbitrary dimension with meromorphic poles along normal crossing divisors
such that a formal decomposition exists locally
(without ramification) and such that the regular singular pieces
have logarithmic poles. Our $u^j_0,u^j_\iiii$ ($j=1,2$) correspond to the
{\it good sets of irregular values} there, and the rest of our 
monodromy data \eqref{2.3} is called the {\it full Stokes data.}
Theorem 4.3.1 in \cite{Mo11a} includes the possibility of 
subfamilies on which the family of connections extends to a flat connection
on $(z,t)$-space, where $t$ is the deformation parameter,
and which are thus isomonodromic subfamilies.
We did not include that in theorem \ref{t2.3}. We shall discuss
isomonodromic families in our situation separately in chapter \ref{s3}.

(vi) We are assuming here an obvious formulation of the
notion of a holomorphic family of $P_{3D6}$ bundles.
One has a holomorphic vector bundle $H^{(T)}$ on $\P^1\times T$, where $T$
is a (usually connected) complex manifold,
and a relative connection $\nnn^{(T)}$ which gives on each bundle
$H^{(t)}=H^{(T)}|_{\P^1\times \{t\}}$ for $t\in T$ a flat meromorphic
connection $\nnn^{(t)}$, such that $(H^{(t)},\nnn^{(t)})$ is a 
$P_{3D6}$ bundle. Then, automatically, the numbers $u^j_0(t),u^j_\iiii(t)$
($j=1,2$) vary holomorphically.

The sectors in \eqref{2.1}
and \eqref{2.2} also vary holomorphically. The notion of a holomorphic
family of $P_{3D6}$ monodromy tuples is now clear.
$L^{(T)}$ is a holomorphic vector bundle on $\C^*\times T$,
with a relative connection which gives on each bundle $L^{(t)}$, $t\in T$,
a flat structure. The flat subbundles $L^{\pm j}_0,L^{\pm j}_\iiii$
vary (within the varying sectors) holomorphically. The numbers
$\alpha^j_0,\alpha^j_\iiii$ depend holomorphically on $t\in T$.
\end{remarks}

{\bf From a $P_{3D6}$ bundle to a $P_{3D6}$ monodromy tuple:}
This will make one direction in theorem \ref{t2.3} (a) precise.
The following is well-known --- references are, for example, \cite{Ma83a}, \cite{Si90},  \cite{Sa02},  \cite{PS03}, \cite{FIKN06}.

Let $(H,\nnn,u^1_0,u^1_\iiii)$ be a $P_{3D6}$ bundle.
Then $u^j_0,u^j_\iiii$ are given by definition \ref{t2.1}, and
$L$ is the restriction $H|_{\C^*}$ with its flat structure.
The subbundles and the numbers $\alpha^j_0,\alpha^j_\iiii$ will
come from sectorial decompositions which \lq\lq lift\rq\rq\  the formal decompositions of Hukuhara, Turrittin and others
at the poles.

We discuss first the pole at $0$.
The pole at $\iiii$ will (with the new coordinate $\www z=\frac{1}{z}$)
be entirely analogous.
The formal decomposition at the pole at $0$ exists without ramification in the case of $P_{3D6}$ bundles,
 because the pole is semisimple with distinct 
eigenvalues (e.g. \cite[II 5.7 Theorem]{Sa02}). 
Usually the formal decomposition 
is written for meromorphic bundles. But again because the pole
is semisimple with distinct eigenvalues, 
it is valid also for holomorphic $P_{3D6}$ bundles
(e.g. \cite[II 5.8 Remark]{Sa02}). 

Writing $\C[[z]]$ for the algebra of formal power series and
$\C\{z\}$ for the subalgebra of convergent power series, one
can express the formal decomposition at $0$ as
\begin{eqnarray}\label{2.8}
\Psi^{for}_0:(\OO(H)_0,\nnn)\otimes_{\C\{z\}}\C[[z]]
\stackrel{\cong}{\to} \oplus_{j=1}^2
(\OO(H^j_0)_0,\nnn^j_0)\otimes_{\C\{z\}}\C[[z]]
\end{eqnarray}
where $\OO(H)$ is the sheaf of holomorphic sections of $H$,
$(H^j_0,\nnn^j_0)$ is a rank 1 holomorphic vector bundle on $\C$, with a meromorphic
(flat) connection with a pole at $0$ of order $\leq 2$ 
(2 if $u^j_0\neq 0$; 1 if $u^j_0=0$), with eigenvalue $u^j_0$ of the 
pole part $[z\nnn_\zdz]:(H^j_0)_0\to (H^j_0)_0$.
Then $(H^j_0,\nnn^j_0-d(\frac{-u^j_0}{z}))$ has a logarithmic pole at $0$.
That means that there is a unique $\alpha^j_0\in\C$ such that
$$\OO(H^j_0)=\OO_\C\MGcdot z^{\alpha^j_0}\MGcdot e^{-u^j_0/z}\MGcdot e^j_0$$
where $e^j_0$ is a multi-valued flat global section of $H^j_0|_{\C^*}$.
Then the monodromy of $(H^j_0,\nnn^j_0)$ has eigenvalue
$e^{-2\pi i\alpha^j_0}$, and 
\begin{eqnarray*}
(\nnn^j_0)_\zdz(z^{\alpha^j_0}\MGcdot e^{-u^j_0/z}\MGcdot e^j_0)
&=&\zdz(z^{\alpha^j_0}\MGcdot e^{-u^j_0/z})\MGcdot e^j_0\\
&=&(\alpha^j_0+{u^j_0}/{z})\MGcdot z^{\alpha^j_0}\MGcdot e^{-u^j_0/z}\MGcdot e^j_0.
\end{eqnarray*}

This formal decomposition lifts to a holomorphic 
decomposition in sectors which contain at most one of the 
two Stokes directions $\R_{>0}\MGcdot (\pm \zeta_0)$.  In sectors which contain exactly one of the two Stokes directions,
it lifts uniquely. So a unique lift of the formal decomposition exists in the sectors $\whhh I^\pm_0$, and they are maximal with this property.

To formulate the sectorial decomposition,
 we need the sheaf $\AAA$ on $S^1$ of holomorphic functions
in (intersections with neighbourhoods of $0$ of) sectors which have an 
asymptotic expansion in $\C[[z]]$. For the definition we refer to 
\cite{Ma83a} (it is instructive to compare the definitions in 
\cite{Sa02}, \cite{PS03}, \cite{Mo11a}).
For $f\in \AAA_\xi$, $\xi\in S^1$, we denote by $\whhh f\in\C[[z]]$ its
asymptotic expansion.

Then there exist two unique flat isomorphisms
\begin{eqnarray}\label{2.9}
\Psi^\pm_0:H|_{\whhh I^\pm_0}\stackrel{\cong}{\to}
H^1_0|_{\whhh I^\pm_0}\oplus H^2_0|_{\whhh I^\pm_0}
\end{eqnarray}
which
extend to isomorphisms
\begin{eqnarray}\label{2.10}
\Psi^\pm_0:(\OO(H)_0,\nnn)\otimes_{\C\{z\}}\AAA|_{I^\pm_0}
\stackrel{\cong}{\to} \oplus_{j=1}^2
(\OO(H^j_0)_0,\nnn^j_0)\otimes_{\C\{z\}}\AAA|_{I^\pm_0},
\end{eqnarray}
lifting $\Psi^{for}_0$.
The flat rank 1 subbundles
\begin{eqnarray}\label{2.11}
L^{\pm j}_0:=(\Psi^\pm_0)^{-1}(H^j_0|_{\whhh I^\pm_0})
\end{eqnarray}
satisfy 
$L|_{\whhh I^{\pm}_0}=L^{\pm 1}_0\oplus L^{\pm 2}_0$ 
and \eqref{2.5}.
Let us explain the part $L^{+1}_0|_{\whhh I^a_0}=L^{-1}_0|_{\whhh I^a_0}$ of 
\eqref{2.5}. If $e^{\pm 1}_0$ and $e^{\pm}_0$ are flat generating sections
of $L^{\pm 1}_0$ and $L^{\pm 2}_0$, then $z^{\alpha^1_0}e^{-u^1_0/z}e^{\pm 1}_0$
and $z^{\alpha^2_0}e^{-u^2_0}e^{\pm 2}_0$ are generating sections of
$\OO(H)|_0\otimes_{\C\{z\}}\AAA|_{I^\pm_0}$. 
If $e^{-1}_{0|I^a_0}=\varepsilon_1\MGcdot e^{+1}_{0|I^a_0}
+\varepsilon_2\MGcdot e^{+2}_{0|I^a_0}$ with $\varepsilon_2\neq 0$,
then also $z^{\alpha^1_0}e^{-u^1_0}e^2_0$ would be a section of 
$\OO(H)|_0\otimes_{\C\{z\}}\AAA|_{I^a_0}$. Then 
$z^{\alpha^1_0}e^{-u^1_0}$ would be in $z^{\alpha^2_0}e^{-u^2_0}\MGcdot \AAA|_{I^a_0}$.
But that is impossible as, for $z\in \whhh I^a_0$ close to 0,
$e^{-u_2/z}$ is exponentially much smaller than $e^{-u^1_0/z}$.

The discussion of the pole at $\iiii$ is entirely analogous,
with the new coordinate $\www z=\frac{1}{z}$ on $\P^1-\{0\}$.  We have
$\www z\partial_{\www z}=-\zdz$,  
$\www z^2\partial_{\www z}=-\paaz$, and 
\[
\OO(H^j_\iiii)=\OO_{\P^1-\{0\}}\MGcdot \www z^{\alpha^j_\iiii}\MGcdot
e^{-u^j_\iiii/\www z} = 
\OO_{\P^1-\{0\}}\MGcdot z^{-\alpha^j_\iiii}\MGcdot e^{-u^j_\iiii\MGcdot z}.\]
One has to observe that the monodromy with respect to $\www z$
is the inverse of the monodromy (with respect to $z$).
This is the reason for the different signs in the two compatibility
conditions in \eqref{2.7}.

All data of the $P_{3D6}$ monodromy tuple have now been defined.

{\bf From a $P_{3D6}$ monodromy tuple to a $P_{3D6}$ bundle:}
This will make the other direction in theorem \ref{t2.3} (a) precise.
It is as classical as the first direction.

Let $(u^j_0,u^j_\iiii,\alpha^j_0,\alpha^j_\iiii,L,L^{\pm j}_0,
L^{\pm j}_\iiii (j=1,2))$ be a $P_{3D6}$ monodromy tuple.
$L$ will be the restriction to $\C^*$ of a bundle $H\to\P^1$,
whose extensions to $0$ and $\iiii$ still have to be constructed.
We shall discuss the extension to $0$ (the extension to $\iiii$
will be analogous). Choose flat generating sections $e^{\pm j}_0$
of the flat rank 1 subbundles $L^{\pm j}_0$ of $L|_{\whhh I^\pm_0}$
with 
\begin{eqnarray}\label{2.12}
e^{-1}_0|_{\whhh I^a_0}=e^{+1}_0|_{\whhh I^a_0}\textup{ and }
e^{-2}_0|_{\whhh I^b_0}=e^{-2\pi i\alpha^2_0}\MGcdot e^{+2}_0|_{\whhh I^b_0}.
\end{eqnarray}
Consider the holomorphic sections
\begin{eqnarray}\label{2.13}
z^{\alpha^j_0}\MGcdot e^{-u^j_0/z}\MGcdot e^{\pm j}_0\quad\textup{ on }\whhh I^\pm_0.
\end{eqnarray}
Then there exist coefficients $a^{\pm ij}_0\in \AAA|_{I^\pm_0}$
with $\whhh a^{+ij}_0=\whhh a^{-ij}_0$ and
$\whhh a^{\pm ij}_0(0)=\delta_{ij}$ and such that for both $j=1,2$
the sections on the (intersections with a neighbourhood of $0$ of the) sectors 
$\whhh I^+_0$ and $\whhh I^-_0$
\begin{eqnarray}\label{2.14}
\sum_{i=1}^2 z^{\alpha^i_0}\MGcdot e^{-u^i_0/z}\MGcdot e^{\pm i}_0\MGcdot a^{\pm ij}_0
\end{eqnarray}
glue to a holomorphic section on (a neighbourhood of $0$ in) $\C^*$,
and the two sections form a basis. 
Here $z^{\alpha^i_0}$ is defined on $\whhh I^+_0$ 
using an arbitrary branch of $\log z$
on $\whhh I^+_0$.
It is defined on $\whhh I^-_0$ by 
extending the branch on $\whhh I^+_0$ counterclockwise 
(that means, over $\whhh I^a_0$)
to $\whhh I^-_0$. Then the extension $\OO(H)_0$ to 0
of $L=H|_{\C^*}$ is defined by these sections.
The sections and the coefficients are not at all unique, but the extension
is unique and has a pole of order 2 
at $0$ with eigenvalues $u^j_0$ of the pole part.

The same procedure works at $\iiii$, with the new coordinate $\www z={1}/{z}$.
This completes the construction of the $P_{3D6}$ bundle.\hfill$\Box$

\begin{remarks}\label{t2.5}
(i) Let $e^{\pm j}_0$ be as above. 
Then $\uuuu e^{\pm}_0:=(e^{\pm 1}_0,e^{\pm 2}_0)$ is a flat basis of 
$L|_{\whhh I^\pm_0}$. The construction above can be rephrased by saying
that there exist two matrices $A^\pm_0\in GL(2,\AAA|_{I^\pm_0})$ with 
$\whhh A^+_0=\whhh A^-_0$ and 
$\whhh A^\pm_0(0)={\bf 1}_2$ such that the two holomorphic
bases
\begin{eqnarray}\label{2.15}
\uuuu e^\pm_0\MGcdot \begin{pmatrix}z^{\alpha^1_0}e^{-u^1_0/z} & 0\\
0 & z^{\alpha^2_0}e^{-u^2_0/z}\end{pmatrix}\MGcdot A^\pm_0(z)
\end{eqnarray}
glue to a basis $\uuuu\varphi$ on (a neigborhood of $0$ in) $\C^*$,
with which the extension of $H|_{\C^*}$ to $0$ is defined.
%Then $\varphi_j(0)$ is an eigenvector of $[z\nnn_\zdz]:H_0to H_0$
%with eigenvalue $u^j_0$.

%MG %changed by CH
Here we use (for $n=m=2$) 
the convention that the matrix multiplication extends in the following 
way to a product of an $n$-tuple (= a row vector) with entries $v_1,....v_n$ 
in a $K$-vector space $V$ and a matrix $(a_{ij})\in M(n\times m,K)$,
$$ (v_1,...,v_n)(a_{ij})=(\sum_{i=1}^n a_{i1}v_i,...,\sum_{i=1}^n a_{im}v_i).$$
So, the product is an $m$-tuple with entries in $V$.

(ii) Given a basis $\uuuu\varphi$ of $\OO(H)_0$ the following
conditions are equivalent.
\begin{eqnarray}\label{2.16}
\textup{($\alpha$)} && [z\nnn_\zdz](\uuuu\varphi(0))=\uuuu\varphi(0)\MGcdot
\begin{pmatrix}u^1_0 & 0\\ 0 & u^2_0\end{pmatrix},\\ \label{2.17}
\textup{($\beta$)} && \nnn_\zdz\uuuu\varphi(z)
=\uuuu\varphi(z)\MGcdot \left[\frac{1}{z}\begin{pmatrix}u^1_0&0\\0&u^2_0\end{pmatrix}
+B(z)\right]
\\ \nonumber
&& \textup{with }B(z)\in M(2\times 2,\C\{z\}),\\ \label{2.18}
\textup{($\gamma$)} && 
\textup{flat generating sections }e^{\pm j}_0 \textup{ of }L^\pm_0
\textup{ and matrices } A^\pm_0\textup{ as}\\ \nonumber
&& \textup{in (i) exist such that the basis }\uuuu\varphi \textup{ takes the form \eqref{2.15} on }\whhh I^\pm_0.
\end{eqnarray}
The first equivalence is trivial, the second is at the heart of the story.
Here the correspondence 
\[
\uuuu\varphi(0)\longleftarrow \uuuu\varphi \longrightarrow \uuuu e^\pm_0
\]
is $1{:}\text{many}{:}1$. It induces a canonical $1{:}1$ correspondence
\begin{eqnarray}\label{2.19}
\uuuu\varphi(0) \longleftrightarrow \uuuu e^\pm_0
\end{eqnarray}
between bases $\uuuu\varphi(0)$ of $H_0$ of eigenvectors of $[z\nnn_\zdz]$
and bases $\uuuu e^\pm_0$ of $L|_{\whhh I^\pm_0}$ 
which respect the splitting 
$L|_{\whhh I^{\pm}_0}=L^{\pm 1}_0\oplus L^{\pm 2}_0$ 
and satisfy \eqref{2.12}.
\end{remarks}

The $P_{3D6}$ monodromy tuples are good for conceptual arguments, e.g.\ regarding irreducibility in chapter \ref{s3}, but less well suited for
calculations such as those in chapter \ref{s6}.
For this reason we shall construct {\it $P_{3D6}$ numerical tuples}
from $P_{3D6}$ monodromy tuples, although this involves making further choices.

Let a $P_{3D6}$ monodromy tuple be given. Choose bases
$\uuuu e^+_0=(e^{+1}_0,e^{+2}_0)$ and 
$\uuuu e^+_\iiii =(e^{+1}_\iiii,e^{+2}_\iiii)$ of flat generating sections
of the bundles $L^{+1}_0,L^{+2}_0,L^{+1}_\iiii,L^{+2}_\iiii$.
Then there are unique bases 
$\uuuu e^-_0=(e^{-1}_0,e^{-2}_0)$ and 
$\uuuu e^-_\iiii =(e^{-1}_\iiii,e^{-2}_\iiii)$ of flat generating sections
of the bundles $L^{-1}_0,L^{-2}_0,L^{-1}_\iiii,L^{-2}_\iiii$
such that \eqref{2.12} and the analogous condition at $\iiii$ hold.
Furthermore, there are unique numbers $s^a_0,s^b_0,s^a_\iiii,s^b_\iiii\in\C$
such that the matrices 
\begin{eqnarray}\label{2.20}
\begin{split}
S^a_0&:=
\begin{pmatrix}
1 & s^a_0\\ 0 & 1
\end{pmatrix},\quad
S^b_0:=
\begin{pmatrix} 
e^{-2\pi i\alpha^1_0} & 0 \\ 0 & e^{-2\pi i \alpha^2_0}
\end{pmatrix}
\MGcdot 
\begin{pmatrix}
1 & 0\\s^b_0 & 1
\end{pmatrix},
\\ 
S^a_\iiii&:=
\begin{pmatrix}
1 & s^a_\iiii\\ 0 & 1
\end{pmatrix},\ \ 
S^b_\iiii:=
\begin{pmatrix} 
e^{-2\pi i\alpha^1_\iiii} & 0 \\ 
0 & e^{-2\pi i \alpha^2_\iiii}
\end{pmatrix}
\MGcdot 
\begin{pmatrix}
1 & 0\\s^b_\iiii & 1
\end{pmatrix}
\end{split}
\end{eqnarray}
and the bases $\uuuu e^\pm_0,\uuuu e^\pm_\iiii$ satisfy 
\begin{eqnarray}\label{2.21}
\begin{split}
\uuuu e^-_0|_{\whhh I^a_0}=\uuuu e^+_0|_{\whhh I^a_0}\MGcdot S^a_0,
&\textup{ i.e.\ } \uuuu e^-_0(-z)=
\uuuu e^+_0(ze^{\pi i})\MGcdot S^a_0\textup{ for }z\in \whhh I^+_0,\\ 
\uuuu e^-_0|_{\whhh I^b_0}=\uuuu e^+_0|_{\whhh I^b_0}\MGcdot S^b_0,
&\textup{ i.e.\ } \uuuu e^-_0(-z)=
\uuuu e^+_0(ze^{-\pi i})\MGcdot S^b_0\textup{ for }z\in \whhh I^+_0,\\ 
\uuuu e^-_\iiii|_{\whhh I^a_\iiii}=\uuuu e^+_\iiii|_{\whhh I^a_\iiii}
\MGcdot S^a_\iiii, &\textup{ i.e.\ } \uuuu e^-_\iiii(-z)=
\uuuu e^+_\iiii(ze^{-\pi i})\MGcdot S^a_\iiii\textup{ for }z\in \whhh I^+_\iiii,
\\ 
\uuuu e^-_\iiii|_{\whhh I^b_\iiii}=\uuuu e^+_\iiii|_{\whhh I^b_\iiii}
\MGcdot S^b_\iiii, &\textup{ i.e.\ } \uuuu e^-_\iiii(-z)=
\uuuu e^+_\iiii(ze^{\pi i})\MGcdot S^b_\iiii
\textup{ for }z\in \whhh I^+_\iiii.
\end{split}
\end{eqnarray}
Here the meaning of $\uuuu e^+_0(ze^{\pi i})$ (or $\uuuu e^+_0(ze^{-\pi i})$)
is that the basis $\uuuu e^+_0$ of flat sections is extended (flatly)
counterclockwise (or clockwise) to $\whhh I^-_0$, and similarly for 
$e^+_\infty(ze^{\pm \pi i})$.

Writing
\[
\text{$\Mon=$  monodromy of $L$ }
\]
we have
\begin{eqnarray}\label{2.22}
\begin{split}
\Mon(\uuuu e^+_0)(z) &= \uuuu e^+_0(ze^{2\pi i}) 
= \uuuu e^-_0(-ze^{\pi i})\MGcdot (S^a_0)^{-1}\\
&=\uuuu e^+_0(z)\MGcdot S^b_0\MGcdot (S^a_0)^{-1},\\ 
\Mon(\uuuu e^-_\iiii)(z) &= \uuuu e^-_\iiii(ze^{2\pi i}) 
= \uuuu e^+_\iiii(-ze^{\pi i})\MGcdot S^a_\iiii\\ 
&=\uuuu e^-_\iiii(z)\MGcdot (S^b_\iiii)^{-1}\MGcdot S^a_\iiii.
\end{split}
\end{eqnarray}

Let us choose a number $\beta\in \C^*$ such that 
\begin{eqnarray}\label{2.23}
e^{-\beta}=(u^1_0-u^2_0)(u^1_\iiii-u^2_\iiii).
\end{eqnarray}
Then 
\begin{eqnarray*}
e^\beta\MGcdot (-\zeta_0)\MGcdot \R_{>0} &=& 
\zeta_\iiii\MGcdot \R_{>0},\\
e^{\beta}\MGcdot \whhh I^+_0&=&\whhh I^-_\iiii,\\
e^\beta\MGcdot \whhh I^a_0 &=& \whhh I^a_\iiii.
\end{eqnarray*}
For an arbitrary point $e^\xi\in \whhh I^+_0$, this $\beta$ 
encodes a (homotopy class) 
of a path from this point $e^\xi$ to the point $e^{\xi+\beta}$,
namely a path whose lift to the universal covering $\exp:\C\to\C^*$
with starting point $\xi$ ends at $\xi+\beta$.
The homotopy class of the path is unique.
The path (or its homotopy class) is called $[\beta]$.

Using this path, the connection matrix 
\[
B(\beta)=\begin{pmatrix}b_1 & b_2\\ b_3 & b_4\end{pmatrix}\in GL(2,\C)
\]
can be defined by
\begin{equation}\label{2.24}
\uuuu e^-_\iiii = (\uuuu e^+_0\textup{ extended flatly along }[\beta]
\textup{ from }\whhh I^+_0\textup{ to } \whhh I^-_\iiii)\MGcdot B(\beta).
\end{equation}
From this definition we have
\begin{eqnarray}\label{2.25}
b_1b_4-b_2b_3\neq 0
\end{eqnarray}
and
\begin{eqnarray*}
&&(\uuuu e^+_0\textup{ extended along }[\beta]
\textup{ to } \whhh I^-_\iiii)\\
&=&\Mon\left[(\uuuu e^+_0\textup{ extended along }[\beta-2\pi i]
\textup{ to } \whhh I^-_\iiii)\right].
\end{eqnarray*}
This implies
\begin{eqnarray*}
&&(\uuuu e^+_0\textup{ extended along }[\beta]
\textup{ to } \whhh I^-_\iiii)\MGcdot B(\beta-2\pi i)\\
&=&\Mon\left[(\uuuu e^-_\iiii)\right]
=\uuuu e^-_\iiii\MGcdot (S^b_\iiii)^{-1}\MGcdot S^a_\iiii\\
&=&(\uuuu e^+_0\textup{ extended along }[\beta]
\textup{ to } \whhh I^-_\iiii)\MGcdot B(\beta)
\MGcdot (S^b_\iiii)^{-1}\MGcdot S^a_\iiii\\
\textup{and also }
&=&\Mon\left[(\uuuu e^+_0\textup{ extended along }[\beta]
\textup{ to } \whhh I^-_\iiii)\MGcdot B(\beta)\right]\\
&=&(\uuuu e^+_0\textup{ extended along }[\beta]
\textup{ to } \whhh I^-_\iiii)\MGcdot S^b_0\MGcdot (S^a_0)^{-1}\MGcdot B(\beta).
\end{eqnarray*}
This shows
\begin{eqnarray}\label{2.26}
B(\beta)\MGcdot (S^b_\iiii)^{-1}\MGcdot S^a_\iiii
&=&S^b_0\MGcdot (S^a_0)^{-1}\MGcdot B(\beta)\\ \label{2.27}
&=& B(\beta -2\pi i).
\end{eqnarray}
The group $\Z$ acts on the pairs $(\beta,B(\beta))$; its generator $[1]$
acts via
\begin{eqnarray}\label{2.28}
[1]\cdot(\beta,B(\beta))&=&(\beta+2\pi i,B(\beta+2\pi i))
\\ \nonumber
&=&(\beta+2\pi i,(S^b_0\MGcdot (S^a_0)^{-1})^{-1}\MGcdot B(\beta)).
\end{eqnarray}

The bases $\uuuu e^\pm_0$ and $\uuuu e^\pm_\iiii$ can be changed by diagonal
base changes
\begin{eqnarray*}
\www{\uuuu e}^\pm_0=\uuuu e^\pm_0\MGcdot
\begin{pmatrix}\lambda_1&0\\0&\lambda_2\end{pmatrix},\quad
\www{\uuuu e}^\pm_\iiii=\uuuu e^\pm_\iiii\MGcdot
\begin{pmatrix}\lambda_3&0\\0&\lambda_4\end{pmatrix}
\end{eqnarray*}
with $(\lambda_1,\lambda_2,\lambda_3,\lambda_4)\in(\C^*)^4$.
Then  $u^j_0,u^j_\iiii,\alpha^j_0,\alpha^j_\iiii,\beta$ are unchanged, but
\begin{eqnarray}\label{2.29}
\begin{split}
\www S^{a/b}_0&=\begin{pmatrix}\lambda_1&0\\0&\lambda_2\end{pmatrix}^{-1}
\MGcdot S^{a/b}_0\MGcdot\begin{pmatrix}\lambda_1&0\\0&\lambda_2\end{pmatrix},
\\
\www S^{a/b}_\iiii&=\begin{pmatrix}\lambda_3&0\\0&\lambda_4\end{pmatrix}^{-1}
\MGcdot S^{a/b}_\iiii\MGcdot\begin{pmatrix}\lambda_3&0\\0&\lambda_4\end{pmatrix},
\\
\www B(\beta)&=\begin{pmatrix}\lambda_1&0\\0&\lambda_2\end{pmatrix}^{-1}
\MGcdot B(\beta)\MGcdot\begin{pmatrix}\lambda_3&0\\0&\lambda_4\end{pmatrix},
\\
(\www s^a_0,  & \www  s^b_0,  \www s^a_\iiii,\www s^b_\iiii,
\www b_1,\www b_2,\www b_3,\www b_4) 
\\ 
&=\left(\tfrac{\lambda_2}{\lambda_1}s^a_0,\tfrac{\lambda_1}{\lambda_2}s^b_0, \tfrac{\lambda_4}{\lambda_3}s^a_\iiii,\tfrac{\lambda_3}{\lambda_4}s^b_\iiii,
\tfrac{\lambda_3}{\lambda_1}b_1,\tfrac{\lambda_4}{\lambda_1}b_2,
\tfrac{\lambda_3}{\lambda_2}b_3,\tfrac{\lambda_4}{\lambda_2}b_4\right).
\end{split}
\end{eqnarray}
The actions of $\Z$ in \eqref{2.28} and of $(\C^*)^4$ in \eqref{2.29}
commute. 

\begin{definition}\label{t2.6}
A {\it$ P_{3D6}$ numerical tuple} is a $17$-tuple of complex numbers
\begin{eqnarray}\label{2.30}
(u^j_0,u^j_\iiii,\alpha^j_0,\alpha^j_\iiii (j=1,2),
s^a_0,s^b_0,s^a_\iiii,s^b_\iiii,\beta,b_1,b_2,b_3,b_4)
\end{eqnarray}
such that $u^1_0\neq u^2_0,u^1_\iiii\neq u^2_\iiii$, and
\eqref{2.23}, \eqref{2.25} and \eqref{2.26} are satisfied.
\end{definition}

The group $(\C^*)^4\times \Z$ acts via \eqref{2.28} and \eqref{2.29}
on the set of $P_{3D6}$ numerical tuples.

\begin{lemma}\label{t2.7}
(a) There is a canonical $1{:}1$ correspondence between $P_{3D6}$ monodromy tuples
and $(\C^*)^4\times \Z$-orbits of $P_{3D6}$ numerical tuples.

(b) Any holomorphic family of $P_{3D6}$ numerical tuples induces a 
holomorphic family of $P_{3D6}$ bundles.
Any holomorphic family of $P_{3D6}$ bundles can locally be represented
by a holomorphic family of $P_{3D6}$ numerical tuples.
\end{lemma}

{\bf Proof:} Starting from a $P_{3D6}$ monodromy tuple, 
the $(\C^*)^4\times \Z$-orbit of $P_{3D6}$ numerical tuples has been 
constructed above. Going backwards through the construction, one can 
construct from a $P_{3D6}$ numerical tuple a $P_{3D6}$ monodromy tuple.
The statements regarding holomorphic families are obvious.
\hfill$\Box$

\begin{remarks}\label{t2.8}
(i) Let us fix the $9$-tuple
$(u^j_0,u^j_\iiii,\alpha^j_0,\alpha^j_\iiii (j=1,2),\beta)$. The quotient set
\[
\{(s^a_0,s^b_0,s^a_\iiii,s^b_\iiii,b_1,b_2,b_3,b_4)\}/(\C^*)^4
\]
of those $8$-tuples which, together with the fixed part, form a $P_{3D6}$ numerical
tuple, parametrizes the set of isomorphism classes of $P_{3D6}$ bundles
with the given fixed part.

In the case of trace free $P_{3D6}$ bundles
(see remark \ref{t4.1} (ii) and (iii))  
this quotient was considered in \cite{PS09}.
The categorical quotient was calculated there. 
See remark \ref{t3.3} (i) and (ii).

(ii) In a $P_{3D6}$ monodromy tuple, the link between the data at $0$ 
and at $\iiii$ is built in, but the compatibility condition 
\eqref{2.7} has to be formulated. In a $P_{3D6}$ numerical tuple,
the compatibility condition \eqref{2.7} is built in, but the link
between the data at $0$  and at $\iiii$ has to be formulated --- it is the
connection matrix $B(\beta)$.
\end{remarks}

\chapter{(Ir)Reducibility}\label{s3}
\setcounter{equation}{0}

\noindent
A pair $(H,\nnn)$, where $H\to M$ is a holomorphic vector bundle on a 
complex manifold $M$, and $\nnn$ is a (flat) meromorphic connection, 
is said to be reducible
if there exists a subbundle $G\subset H$ with $0<\rank G<\rank H$ which
is (at all nonsingular points of the connection) a flat subbundle.
Such a $G$ will simply be called a flat subbundle. 
A pair $(H,\nnn)$ is completely reducible if it decomposes into
a sum of flat rank 1 subbundles.

The (ir)reducibility of a $P_{3D6}$ bundle is easily characterized
in terms of properties of its $P_{3D6}$ monodromy tuple,
and the reducible $P_{3D6}$ bundles can be described explicitly.
This is the content of lemma \ref{t3.1} and corollary \ref{t3.2} below.
They are essentially contained in \cite{PS09}, though not in this detail.

\begin{lemma}\label{t3.1} (cf.\  \cite{PS09})
(a) A $P_{3D6}$ bundle is reducible but not completely reducible
if and only if exactly one of the following four conditions $(C_{ij})$,
$i,j=1,2,$ holds. 
%\begin{list}{}{}
%\item[$(C_{ij})$:]
%$\alpha^i_0+\alpha^j_\iiii\in\Z$, 
%$L^{+i}_0$ and $L^{-i}_0$ glue to a (flat rank 1) subbundle $L^i_0$ of $L$,
%$L^{+j}_\iiii$ and $ L^{-j}_\iiii$ glue to a (flat rank 1) subbundle $L^j_\iiii$
%of $L$, and $L^i_0=L^j_\iiii$.
%\end{list}
\begin{list}{}{}
\item
\text{Condition $(C_{ij})$:}
\item
\quad
$\alpha^i_0+\alpha^j_\iiii\in\Z$
\item
\quad
$L^{+i}_0$ and $L^{-i}_0$ glue to a (flat, rank 1) subbundle $L^i_0$ of $L$
\item
\quad
$L^{+j}_\iiii$ and $ L^{-j}_\iiii$ glue to a (flat, rank 1) subbundle $L^j_\iiii$
of $L$
\item
\quad
$L^i_0=L^j_\iiii$.
\end{list}
If $(C_{ij})$ holds, let $e^{ij}$ be a multi-valued flat global generating
section of $L^i_0=L^j_\iiii$. Then the sections
\[
\text{
$\sigma^{ij}:= z^{\alpha^i_0}e^{-u^i_0/z-u^j_\iiii z}\MGcdot e^{ij}$
and $z^{-\alpha^j_\iiii-\alpha^i_0}\MGcdot \sigma^{ij}$
}
\] 
generate the restriction
to $\C$ and the restriction to $\P^1-\{0\}$ of a flat subbundle $W$ of
$(H,\nnn)$, and $\deg W=-\alpha^i_0-\alpha^j_\iiii$.

(b) A $P_{3D6}$ bundle is completely reducible if and only if condition 
$(C_{ij})$ (for some $i,j$) and also condition $(C_{3-i,3-j})$ hold.
\end{lemma}

{\bf Proof:}
If $W\subset H$ is a flat rank 1 subbundle, then 
$W|_{\C^*}\subset H|_{\C^*}=L$ must be compatible with all four splittings
$L|_{\whhh I^\pm_{0/\iiii}}=L^{\pm 1}_{0/\iiii}\oplus L^{\pm 2}_{0/\iiii}$
of $L$ (on $\whhh I^\pm_0$ and $\whhh I^\pm_\iiii$).
Then $W|_{\C^*}$ must  coincide with one subbundle of the 
splitting on each of the four sets
$\whhh I^\pm_0$ and $\whhh I^\pm_\iiii$. This shows that one of the four conditions $(C_{ij})$ holds.

Conversely, if $(C_{ij})$ holds, then a flat rank 1 subbundle $W$ 
can be constructed as in part (a) of the lemma.

If a second condition holds, it must be $(C_{3-i,3-j})$, and
then one can construct a second flat rank 1 subbundle, so that 
$(H,\nnn)$ is completely reducible. 

Conversely, if $(H,\nnn)$ is completely reducible, two conditions
$(C_{ij})$ and $(C_{3-i,3-j})$ hold. This establishes (a) and (b).
\hfill $\Box$

\begin{corollary}\label{t3.2} (cf.\  \cite{PS09})
Fix $(u^j_0,u^j_\iiii,\alpha^j_0,\alpha^j_\iiii(j=1,2),\beta)$
with $e^{-\beta}=(u^1_0-u^2_0)(u^1_\iiii-u^2_\iiii)$.
Consider the set $M(u^j_0,u^j_\iiii,\alpha^j_0,\alpha^j_\iiii)$
of isomorphism classes of $P_{3D6}$ bundles whose $P_{3D6}$ monodromy tuples
contain $u^j_0,u^j_\iiii,\alpha^j_0,\alpha^j_\iiii$. 

(a) Then
\begin{eqnarray}\label{3.1}
\deg H=\deg \det H &=& -\alpha^1_0-\alpha^2_0-\alpha^1_\iiii-\alpha^2_\iiii\in\Z,\\
\label{3.2}
\textup{so }\alpha^1_0+\alpha^1_\iiii\in\Z &\iff&
\alpha^2_0+\alpha^2_\iiii\in\Z,\\ \label{3.3}
\alpha^1_0+\alpha^2_\iiii\in\Z &\iff&
\alpha^2_0+\alpha^1_\iiii\in\Z.
\end{eqnarray}

(b) If neither \eqref{3.2} nor \eqref{3.3}
hold, then $M(u^j_0,u^j_\iiii,\alpha^j_0,\alpha^j_\iiii)$ does not 
contain reducible $P_{3D6}$ bundles.

(c) If conditions \eqref{3.2} or conditions \eqref{3.3} hold, $M(u^j_0,u^j_\iiii,\alpha^j_0,\alpha^j_\iiii)$
contains three families of reducible bundles, and these families map
to one point in the categorical quotient 
\begin{equation}\label{3.4}
\{P_{3D6}\textup{ numerical tuples containing }
(u^j_0,u^j_\iiii,\alpha^j_0,\alpha^j_\iiii,\beta)\}/\!/(\C^*)^4
\end{equation}
where $(\C^*)^4$
acts as in \eqref{2.29}. 
One family is parametrized by a point; it is a completely reducible
$P_{3D6}$ bundle. The other two families are parametrized by $\P^1$
and contain reducible, but not completely reducible, $P_{3D6}$ bundles.

(d) If both conditions \eqref{3.2} and \eqref{3.3} hold,
the categorical quotient has two special points, and above each of them
there are three families as in (c).
\end{corollary}

{\bf Proof:}
(a) By definition, $\deg H=\deg\det H$. The bundle $\det(H,\nnn)$
is a product of a factor $e^{-(u^1_0+u^2_0)/z-(u^1_\iiii+u^2_\iiii)z}
\times\text{\it (trivial flat bundle)}$ 
and a bundle with logarithmic poles at $0$ and $\iiii$, with 
eigenvalues $\alpha^1_0+\alpha^2_0$ and $-(\alpha^1_\iiii+\alpha^2_\iiii)$
of $[\nnn_\zdz]$ on $(\det H)_0$ and $(\det H)_\iiii$. This establishes \eqref{3.1}.
Now \eqref{3.2} and \eqref{3.3} follow trivially.

(b) If none of the conditions in \eqref{3.2} and \eqref{3.3}
hold, then none of the conditions $(C_{ij})$ hold
for a bundle in 
$M(u^j_0,u^j_\iiii,\alpha^j_0,\alpha^j_\iiii)$. Hence, by lemma \ref{t3.1},
none of these bundles can be reducible.

(c)+(d) First we consider condition $(C_{11})$. A $P_{3D6}$ monodromy
tuple satisfies $(C_{11})$ if and only if any associated $P_{3D6}$ numerical tuple
satisfies $s^b_0=s^b_\iiii=b_3=0$. As $u^j_0$ and $u^j_\iiii$ are fixed, we can also
fix some $\beta$ with \eqref{2.23} and restrict to the connection matrix
$B=B(\beta)$. Then the set of isomorphism classes of $P_{3D6}$ monodromy tuples
corresponds to the quotient (as a set) of the $5$-tuples $(s^a_0,s^a_\iiii,b_1,b_2,b_4)$
by the action of $(\C^*)^4$ as in \eqref{2.29} where 
$(\lambda_1,\lambda_2,\lambda_3,\lambda_4)$ maps 
$(s^a_0,s^a_\iiii,b_1,b_2,b_4)$ to
\begin{equation}\label{3.5}
\left(\tfrac{\lambda_2}{\lambda_1}s^a_0,\tfrac{\lambda_4}{\lambda_3}s^a_\iiii,
\tfrac{\lambda_3}{\lambda_1}b_1,\tfrac{\lambda_4}{\lambda_1}b_2,
\tfrac{\lambda_4}{\lambda_2}b_4\right).
\end{equation}
As $b_1\neq 0$ and $b_4\neq 0$, they can be normalized to $b_1=b_4=1$.
Then the action reduces to the action of $(\C^*)^2$ on triples $(s^a_0,s^a_\iiii,b_2)$,
where $(\lambda_1,\lambda_2)$ maps this triple to 
\[
\left(\tfrac{\lambda_2}{\lambda_1}s^a_0,\tfrac{\lambda_2}{\lambda_1}s^a_\iiii,
\tfrac{\lambda_2}{\lambda_1}b_2\right).
\]
Equation \eqref{2.26} connects $s^a_0,s^a_\iiii$ and $b_2$ linearly.
It is equivalent to
\begin{align*}
\begin{pmatrix}1&b_2\\0&1\end{pmatrix}
&\begin{pmatrix}e^{2\pi i\alpha^1_\iiii} & 0\\0&e^{2\pi i \alpha^2_\iiii}\end{pmatrix}
\begin{pmatrix}1&s^a_\iiii\\ 0 &1\end{pmatrix}
\\
&=
\begin{pmatrix}e^{-2\pi i\alpha^1_0} & 0\\0&e^{-2\pi i \alpha^2_0}\end{pmatrix}
\begin{pmatrix}1&-s^a_\iiii\\ 0 &1\end{pmatrix}
\begin{pmatrix}1&b_2\\0&1\end{pmatrix}
\end{align*}
and thus equivalent to \eqref{3.2} and 
$s^a_\iiii = -s^a_0 + (1-e^{2\pi i(\alpha^1_0-\alpha^2_0)})b_2.$
This means that the isomorphism classes are parametrized by the orbits of the standard $\C^*$ 
action on $\C^2$ with coordinates $(s^a_0,b_2)$. The point $(s^a_0,b_2)=(0,0)$
corresponds to the completely reducible bundle which satisfies $(C_{11})$ and
$(C_{22})$. The other orbits form a $\P^1$ and correspond to the reducible
but not completely reducible $P_{3D6}$ bundles which satisfy only $(C_{11})$.
Obviously these orbits cannot be separated from the orbit $\{(0,0)\}$
by invariant functions  and are thus mapped to the same point in the
categorical quotient in \eqref{3.4}.

The $P_{3D6}$ bundles which satisfy $(C_{22})$ are represented in exactly the 
same way by $\C^2$ with the standard $\C^*$ action. The orbit $\{(0,0)\}$
corresponds again to the completely reducible bundle which satisfies $(C_{11})$ and
$(C_{22})$, and the other orbits form a $\P^1$ and correspond to the reducible
but not completely reducible $P_{3D6}$ bundles which satisfy only $(C_{22})$.
Therefore these are also mapped to the same point in the categorical
quotient in \eqref{3.4}.

The $P_{3D6}$ bundles which satisfy either $(C_{12})$ or $(C_{21})$ or both,
behave in exactly the same way. But they are mapped to another point in the
categorical quotient in \eqref{3.4}, as their $P_{3D6}$ numerical tuples
satisfy $s^a_0=s^a_\iiii=0$ and either $b_4=0$ or $b_1=0$ or both.
This proves (c) and (d). \hfill$\Box$

\begin{remarks}\label{t3.3}
(i) In the case of trace free $P_{3D6}$ bundles, the 
categorical quotient in \eqref{3.4} had been calculated in \cite{PS09}.
It is an affine cubic surface in $\C^3$. 
If neither \eqref{3.2} nor \eqref{3.3} hold, it is smooth and a geometric
quotient and equips $M(u^j_0,u^j_\iiii,\alpha^j_0,\alpha^j_\iiii)$
with a canonical structure of an affine algebraic manifold.

If either \eqref{3.2} or \eqref{3.3} hold, it is not a geometric quotient.
Then it has an $A_1$ singularity, and above the singular point,
one has the three families of reducible bundles.
Then $M(u^j_0,u^j_\iiii,\alpha^j_0,\alpha^j_\iiii)$ with the natural
quotient topology is not Hausdorff.

One obtains a moduli space if one either restricts to the completely
reducible $P_{3D6}$ bundle or to one of the two families of reducible,
but not completely reducible $P_{3D6}$ bundles. In the latter two cases,
the moduli space is obtained from the categorical quotient by
blowing up the $A_1$ singularity with an exceptional divisor $\P^1$.

If \eqref{3.2} and \eqref{3.3} hold, the categorical quotient
has two $A_1$ singularities and above each of the two singular points
one has three families of reducible bundles. Again 
$M(u^j_0,u^j_\iiii,\alpha^j_0,\alpha^j_\iiii)$ with the quotient
topology is not Hausdorff.
Again one can obtain a moduli space by choosing at either point
one of the three families of reducible bundles.

(ii) The eigenvalues of the monodromy of a reducible bundle are the 
eigenvalues $e^{-2\pi i\alpha^1_0}$ and $e^{-2\pi i\alpha^2_0}$ of the formal monodromy.

If at least one of the conditions \eqref{3.2} or \eqref{3.3} holds, then
the eigenvalues $e^{-2\pi i\alpha^1_0}$ and $e^{-2\pi i\alpha^2_0}$ coincide
if and only if \eqref{3.2} and \eqref{3.3} hold both. In this case a family 
parametrized by $\P^1$ of reducible but not completely reducible $P_{3D6}$ bundles,
splits into a
family parametrized by $\C$, where the monodromy has a $2\times 2$ Jordan block
and a single bundle with semisimple monodromy.

(iii) The case of interest to us is the case with 
$\alpha^j_0=\alpha^j_\iiii=0$. There the categorical quotient has
two $A_1$-singularities. But we shall consider at each of the two points
only the corresponding completely reducible connection. Therefore we obtain
a geometric quotient. See theorem \ref{t7.6} (b).
\end{remarks}

\chapter{Isomonodromic families}\label{s4}
\setcounter{equation}{0}

\noindent
Isomonodromic families of $P_{3D6}$ bundles
can be characterized easily in terms of the $P_{3D6}$ monodromy
tuples. We shall see that a universal isomonodromic family in four parameters
exists, but that only one parameter is essential.

A result of Heu \cite{Heu09} applies and shows that
generic members of isomonodromic families of irreducible $P_{3D6}$ bundles are
trivial as holomorphic vector bundles, and that the other members (if such 
exist) have sheaves of holomorphic sections isomorphic to 
$\OO_{\P^1}(1)\oplus \OO_{\P^1}(-1)$.  This implies solvability of
the inverse monodromy problem --- which asks whether some given monodromy data can be realized by a pure twistor with meromorphic connection --- in this case.

Consider a holomorphic family 
$(H^{(t)},\nnn^{(t)},u^{1,(t)}_0,u^{1,(t)}_\iiii)_{t\in T}$
of $P_{3D6}$ bundles over a complex base manifold $T$
and the corresponding family of $P_{3D6}$ monodromy tuples
$(u^{j,(t)}_0,u^{j,(t)}_\iiii,\alpha^{j,(t)}_0,\alpha^{j,(t)}_\iiii,
L^{(t)},L^{\pm j,(t)}_0,L^{\pm j,(t)}_\iiii)_{t\in T}$.
Suppose that all $L^{(t)}$ have isomorphic monodromy.
Then one can locally (over $S\subset T$, say) 
extend the family of connections on the bundles $L^{(t)}$
to a flat connection on $\cup_{t\in S}L^{(t)}$.
But in general the covariant derivatives of holomorphic sections
in the $T$ direction 
will have essential singularities along
$\{0\}\times S$ and $\{\iiii\}\times S$.

The family $\cup_{t\in S}H^{(t)}$ with extended flat
connection is considered {\it isomonodromic} only if this extended connection
has poles of Poincar\'e rank $\leq$ 1 along $\{0\}\times S$ and
$\{\iiii\}\times S$.  In other words
the covariant derivatives of
holomorphic sections by logarithmic vector fields along $\{0\}\times S$
and $\{\iiii\}\times S$ have at most simple poles along
$\{0\}\times S$ and $\{\iiii\}\times S$.

It is well known (see \cite{Ma83b}) that this holds if the splittings 
$$L|_{\whhh I^{\pm,(t)}_{0/\iiii}}=L^{\pm 1,(t)}_{0/\iiii}\oplus
L^{\pm 2,(t)}_{0,\iiii}$$
vary in a flat way and if the $\alpha^{j,(t)}_0,\alpha^{j,(t)}_\iiii$ are constant.
Thus, the only freedom is the holomorphic
variation of  $u^{j,(t)}_0,u^{j,(t)}_\iiii$.
This gives four parameters for a universal isomonodromic 
family over the universal covering $U^{univ}$ of 
\[
U:=\{(u^1_0,u^2_0,u^1_\iiii,u^2_\iiii)\in\C^4\, |\,
u^1_0\neq u^2_0,u^1_\iiii\neq u^2_\iiii\},
\]
with covering map $c^{univ}:U^{univ}\to U$.

However,  three of the four parameters are inessential in the following sense.
Define $c^{uc*}$ by
\begin{eqnarray}\label{4.1}
c^{uc*}:U\to\C^*,\quad (u^1_0,u^2_0,u^1_\iiii,u^2_\iiii)
\mapsto (u^1_0-u^2_0)(u^1_\iiii-u^2_\iiii).
\end{eqnarray}
Then the $P_{3D6}$ bundles in one fibre of $c^{uc*}$ can be obtained
by the following operations from one another.
If $(\lambda_0,\lambda_\iiii)\in\C^2$ then one can define for any 
$P_{3D6}$ bundle the following isomonodromic family in two parameters
\begin{eqnarray}
&&(H^{(\lambda_0,\lambda_\iiii)},\nnn^{(\lambda_0,\lambda_\iiii)},
u^{j,(\lambda_0,\lambda_\iiii)}_0,u^{j,(\lambda_0,\lambda_\iiii)}_\iiii)
\label{4.2}\\ \nonumber
\textup{ with}&&\OO(H^{(\lambda_0,\lambda_\iiii)}):=
e^{\lambda_0/z\ +\ \lambda_\iiii z}\MGcdot \OO(H),\\ \nonumber
&&\nnn^{(\lambda_0,\lambda_\iiii)}:=\nnn,\\  \nonumber
\textup{ thus}&& u^{j,(\lambda_0,\lambda_\iiii)}_0=u^j_0-\lambda_0,\ 
u^{j,(\lambda_0,\lambda_\iiii)}_\iiii=u^j_\iiii-\lambda_\iiii.
\end{eqnarray}
If $\lambda_1\in\C^*$ and $m_{\lambda_1}:\C\to\C,\ z\to \lambda_1\MGcdot z$,
then one can define the following isomonodromic family in one
parameter,
\begin{eqnarray}
(H^{(\lambda_1)},\nnn^{(\lambda_1)}):=m_{\lambda_1}^* (H,\nnn),
\label{4.3} \\ \nonumber
\textup{ thus}\quad u^{j,(\lambda_1)}_0={u^j_0}/{\lambda_1},\ 
u^{j,(\lambda_1)}_\iiii=\lambda_1\MGcdot u^j_\iiii.
\end{eqnarray}
Combining these operations, one obtains from one $P_{3D6}$ bundle
an isomonodromic family with trivial transversal monodromy 
over one fibre of $c^{uc*}$.
This isomonodromic family is considered {\it inessential}.

The restriction of the universal family over $U^{univ}$ to one fibre of
$c^{uc*}\circ c^{univ}$ pulls down to a family
over the fibre of $c^{uc*}$. The whole universal family pulls down to a family
over the fibre product $U^{beta}$, which is defined by $c^{uc*}$ and $e$,
where $e:\C\to \C^*,\zeta\mapsto e^\zeta$, is the universal
covering of $\C^*$:
\begin{eqnarray}\label{4.4}
\begin{CD}
U^{beta} @>c^{uc}>> \C\\
@Vc^{beta}VV @VVeV\\
U @>>c^{uc*}> \C^*.
\end{CD}
\end{eqnarray}

The covering $c^{beta}:U^{beta}\to U$ is a smaller covering of $U$ 
than the universal covering $U^{univ}\to U$.
Finally, by the operations above, all bundles in one $c^{uc}$-fibre 
of $U^{beta}$ are isomorphic as holomorphic vector bundles
because the isomonodromic family in one $c^{uc}$-fibre is inessential.

\begin{remarks}\label{t4.1}
(i) 
%MG
%The properties {\it irreducibility} and {\it complete reducibility} 
%and {\it reducibility, but not complete reducibility,} 
The properties of being
\begin{align*}
&\text{irreducible}\\
&\text{completely reducible}\\
&\text{reducible, but not completely reducible}
\end{align*}
are independent of the values $u^j_0,u^j_\iiii$.
Therefore all bundles in one isomonodromic family have
the same property. Furthermore, if they are all
reducible, a flat rank 1 subbundle of one $P_{3D6}$ bundle extends 
at least locally to a flat family of flat subbundles (on $\C^*\times T$).

If one and hence all bundles in an  isomonodromic family are completely reducible, they are
all isomorphic as holomorphic bundles, because the family splits locally
into two holomorphic families of rank 1 subbundles, and the degrees
of these rank 1 subbundles are constant.

(ii) Definition: Let $(H\to\P^1,\nnn)$ be a holomorphic vector bundle
with a meromorphic connection. It is called {\it trace free} if its 
determinant bundle $\det(H,\nnn)$ with connection is the trivial bundle
with trivial flat connection.

(iii) The description of the determinant bundle $\det(H,\nnn)$ of a
$P_{3D6}$ bundle in the proof of corollary \ref{t3.2} (a)
shows that the $P_{3D6}$ bundle is trace free if and only if
\begin{eqnarray}\label{4.5}
\sum_{j=1}^2 u^j_0=\sum_{j=1}^2u^j_\iiii=
\sum_{j=1}^2\alpha^j_0=\sum_{j=1}^2\alpha^j_\iiii=0.
\end{eqnarray}

(iv) The following terminology will be useful here and in chapter
\ref{s7}. Given a holomorphic vector bundle $H\to \P^1$ of some rank 
$n\ge1$, by Birkhoff-Grothendieck there exist unique integers $k_1,\dots,k_n$
with $k_1\geq \dots \geq k_n$ such that 
\[
\OO(H)\cong \OO_{\P^1}(k_1)\oplus \dots\oplus \OO_{\P^1}(k_n).
\]
The vector bundle $H$ is then called a $(k_1,\dots,k_n)$-twistor.
A $(0,\dots,0)$-twistor is also called {\it pure}. 
\end{remarks}

\begin{theorem}\label{t4.2} \cite{Heu09}
Let $(H,\nnn,u^1_0,u^1_\iiii)$ be a trace free irreducible $P_{3D6}$ bundle.
Then there exists a discrete (possibly empty) subset $\Sigma\subset\C$ such that,
in the universal isomonodromic family above $U^{beta}$,
all bundles in $(c^{uc})^{-1}(\C-\Sigma)$ are pure (i.e.
trivial as holomorphic vector bundles), 
and the bundles in $(c^{uc})^{-1}(\Sigma)$ are $(1,-1)$-twistors.
\end{theorem}

{\bf Proof:}
%Because of remark \ref{t4.1} (ii), we can restrict to the case
%of trace free (and irreducible) bundles. 
The restriction of the isomonodromic
family over $U^{beta}$ to the 2-dimensional subset $U^{beta,tr}\subset U^{beta}$ 
of trace free $P_{3D6}$ bundles still has 
two parameters: the parameter $\lambda_1$ from above is inessential,
and the parameter in the base of the map $c^{uc}:U^{beta,tr}\to\C$
is essential.

In a much more general situation of trace free rank 2 bundles with meromorphic
connections on closed complex curves, in \cite{Heu10} universal
isomonodromic families are constructed. In \cite{Heu09} the subsets of non-pure
twistors are estimated. In our situation the result of \cite{Heu09} is that
the subset
\begin{eqnarray}\label{4.6}
T_k&:=&\{t\in U^{beta,tr}\, |\, \kappa (H^{(t)})\leq k\}\quad 
\textup{with}\\  \nonumber
\kappa(H^{(t)})&:=&\min(-2\deg W^{(t)}\, |\, 
W^{(t)}\subset H^{(t)}\textup{ is a rank 1 subbundle})
\end{eqnarray}
is a closed analytic subset of codimension at least $-1-k$. 
Because of the one inessential parameter in the 2-dimensional
manifold $U^{beta,tr}$, $T_k$ consists of fibres of 
%MG
$c^{uc}:U^{beta,tr}\to\C$.
Therefore it is either of codimension $0$ or $1$ or is empty. Thus
\[
T_{-4}=\{t\in U^{beta,tr}\, |\, H^{(t)}\textup{ is an }(l,-l)\textup{-twistor
for some }l\geq 2\}
\]
has codimension $\geq -1-(-4)$, so it is empty.  Similarly,
\[
T_{-2}=\{t\in U^{beta,tr}\, |\, H^{(t)}\textup{ is an }
(1,-1)\textup{-twistor}\}
\]
has codimension $\geq -1-(-2)=1$ or is empty. So it is empty, or 
$T_{-2}=(c^{uc})^{-1}(\Sigma)$ for a discrete set $\Sigma\subset \C$.
\hfill$\Box$

\begin{remarks}\label{t4.3}
(i) For the case of $P_{3D6}$-TEP bundles (definition \ref{t6.1}),
which is our main interest, we shall reprove theorem \ref{t4.2}
in chapter \ref{s8} (in theorem \ref{t8.2} (a) and in lemma \ref{t8.5}) 
by elementary calculations using normal forms of such bundles. 
For this case theorem \ref{t4.2} is also proved in \cite{Ni09}.
However,  Heu's result is much more general, and the elementary
proof in chapter \ref{s8} seems simpler.

(ii) Theorem \ref{t4.2} will allow us to restrict to special cases
of the following favorable situation. 

Let $V\to\P^1\times T$ be a holomorphic rank 2 vector bundle on 
$\P^1$ times a (connected complex analytic) manifold $T$, such that any
bundle $V^{(t)}:=V|_{\P^1\times\{t\}}$ for $t\in T$ is either pure or
a $(1,-1)$-twistor.

First,  the set
$$\Theta:=\{t\in T\, |\, V^{(t)}\textup{ is a }(1,-1)\textup{-twistor}\}$$
is either empty or an analytic hypersurface or $\Theta=T$
\cite[I 5.3 Theorem]{Sa02}.

Second, (with $\pi:\P^1\times T\to T$ the projection)
\begin{eqnarray}\label{4.7}
&&\pi_*\OO(V)\textup{ is locally free of rank 2 and}
\\ \label{4.8}
&&\pi_*\OO(V)|_t:=\pi_*\OO(V)_t/{\bf m}_{t,T}
\cong \pi_*\OO(V^{(t)})\textup{ for any }t\in T,
\end{eqnarray}
which implies that any global section of some $V^{(t)}$ extends
(non-uniquely) to a global section of $V|_{\P^1\times S}$ for some
small neighbourhood $S\subset T$ of $t$.

Moreover, \eqref{4.7} and \eqref{4.8} imply that for any small open set $S\subset T$
one can choose two global sections on $V|_{\P^1\times S}$ such that 
they restrict for any $t\in S$ to a basis of the 2-dimensional 
space $\pi_*\OO(V^{(t)})$.

The proof of \eqref{4.7} and \eqref{4.8} is an application
of fundamental results on coherent sheaves in \cite[ch.\  III]{BS76}
and basic facts on the cohomology groups $H^l(\P^1,\OO_{\P^1}(k))$.
It is well known \cite[I 2]{Sa02} that 
\begin{align*}
\dim H^0(\P^1,\OO_{\P^1}(k))&=k+1 \ \textup{ for }k\geq -1,
\\
\dim H^l(\P^1,\OO_{\P^1}(k))&= 0 \ \textup{ for }l\geq 1, k\geq -1.
\end{align*}
This implies that
\begin{align*}
\dim H^0(\P^1,\OO(V^{(t)}))&=2,
\\
\dim H^l(\P^1,\OO(V^{(t)}))&=0\textup{ for }l\geq 1.
\end{align*}
The latter implies that
$R^l\pi_*\OO(V)=0$ for $l\geq 1$ \cite[III Cor. 3.11]{BS76}.
Now the base change theorem \cite[III Cor. 3.5]{BS76} applies
and gives \eqref{4.8}.
Finally, $\dim H^0(\P^1,\OO(V^{(t)}))=2$ for all $t\in T$
and \cite[III Lemma 1.6]{BS76} give \eqref{4.7}.

(iii) All statements and arguments in (ii) generalize to holomorphic
vector bundles $V\to\P^1\times T$ of rank $n\geq 2$ such that for any
$t\in T$ the bundle $V^{(t)}$ is a $(k_1,\dots,k_n)$-twistor 
with $k_n\geq -1$ and $\sum_{j=1}^n k_j=0$.
\end{remarks}

\chapter{Useful formulae: three $2\times 2$ matrices}\label{s5}
\setcounter{equation}{0}

\noindent
In this chapter we collect together some elementary formulae concerning three matrices which appear in the monodromy data: a Stokes matrix $S$, a monodromy matrix $\Mon$, and a connection matrix $B$.  

Let us fix $s\in\C$ and define two matrices
\begin{eqnarray}\label{5.1}
S:=\begin{pmatrix}1&s\\0&1\end{pmatrix},\ 
\Mon_0^{mat}:=S^t\MGcdot S^{-1}=\begin{pmatrix}1&-s\\s&1-s^2\end{pmatrix}
\in SL(2,\C).
\end{eqnarray}
In this chapter (and later) the notation $\sqrt{\frac14{s^2}-1}$ indicates the square root with argument in $[0,\pi)$ if $s\neq \pm 2$
(and $0$ if $s=\pm 2$).

First we consider the case $s=\pm 2$. Then $\Mon_0^{mat}$ has 
a $2\times 2$ Jordan block with eigenvalues $\lambda_+=\lambda_-=-1$.
A suitable basis of $M(2\times 1,\C)$ is 
\begin{equation}\label{5.2}
v_1=\begin{pmatrix}1\\ {s}/{2}\end{pmatrix},\ 
v_2=\begin{pmatrix}0\\ 1\end{pmatrix}
\end{equation}
with
$\Mon_0^{mat}\MGcdot v_1=-v_1$, 
$\Mon_0^{mat}\MGcdot v_2=-v_2-sv_1$.
We define $\alpha_\pm := \textup{sign}(s)\MGcdot (\mp \tfrac{1}{2})$.
Then $e^{-2\pi i \alpha_\pm}=-1=\lambda_\pm$.

Now we consider the case $s\in\C-\{\pm 2\}$. 
Then $\Mon_0^{mat}$ is semisimple with two distinct eigenvalues
\begin{eqnarray}\label{5.3}
\lambda_\pm =(1-\tfrac12{s^2})\pm s\sqrt{\tfrac14{s^2}-1}.
\end{eqnarray}
Thus
$\lambda_++\lambda_-=2-s^2$, $\lambda_+\lambda_-=1$.
For eigenvectors we take
\begin{eqnarray}\label{5.4}
v_\pm =
\begin{pmatrix}
\vphantom{\dfrac12}
1
\\
\mp\sqrt{\frac14{s^2}-1}+\frac12{s}
\end{pmatrix}
\end{eqnarray}
with $\Mon_0^{mat}\MGcdot v_\pm=\lambda_\pm\MGcdot v_\pm$.

For $s\in \C-(\R_{\geq 2}\cup \R_{\leq -2})$ there exist unique
$\alpha_\pm\in\C$ with
\begin{eqnarray}\label{5.5}
\Re(\alpha_\pm)\in(-\tfrac{1}{2},\tfrac{1}{2})
\textup{ and }e^{-2\pi i\alpha_\pm}=\lambda_\pm.
\end{eqnarray}
For $s\in \R_{>2}\cup\R_{<-2}$ there exist unique
$\alpha_\pm$ with
\begin{eqnarray}\label{5.6}
e^{-2\pi i \alpha_\pm}=\lambda_\pm
\textup{ and }\Re(\alpha_\pm)=\textup{sign}(s)\MGcdot(\mp\tfrac{1}{2}).
\end{eqnarray}

\begin{lemma}\label{t5.1}
For any $s\in \C$:
\newline
(a) 
\begin{eqnarray}\label{5.7}
\lambda_+(s)\MGcdot\lambda_-(s)=1,\quad \alpha_+(s)+\alpha_-(s)=0,
\\ \label{5.8}
\lambda_\pm(s)=\lambda_\mp(-s),\quad 
\alpha_\pm(s)=\alpha_\mp(-s),
\\ \label{5.9}
\lambda_+,\lambda_-\in S^1\iff s\in [-2,2]
\iff \alpha_\pm\in[-\tfrac{1}{2},\tfrac{1}{2}], 
\\ \label{5.10}
|\lambda_+|<1\iff \Im(\alpha_+)<0\iff \Im(s)>0\textup{ or }s\in \R_{>2},
\\
\label{5.11}
|\lambda_+|>1\iff \Im(\alpha_+)>0\iff \Im(s)<0\textup{ or }s\in \R_{<-2},
\\
\lambda_+,\lambda_-\in\R_{>0}\iff s\in i\R,
\iff\alpha_+,\alpha_-\in i\R,\label{5.12}
\\ 
\Im(\lambda_+),\Im(\lambda_-)>0\iff \Re(\alpha_+)\in(-\tfrac{1}{2},0)
\nonumber \\ 
\iff\Re(s)>0, s\notin\R_{\geq 2}.\label{5.13}
\end{eqnarray}
The map
\begin{eqnarray}\nonumber
\C&\to& \{z\in\C\, |\, \!-\!\tfrac{1}{2}\!<\!\Re(z)\!<\!\tfrac{1}{2}\}
\!\cup\! (-\tfrac{1}{2}+i\R_{\leq 0})\!\cup\! (\tfrac{1}{2}+i\R_{\geq 0})
\\
s&\mapsto& \alpha_+(s)\label{5.14}
\end{eqnarray}
is bijective. It is continuous in each of the three regions
$\C-(\R_{>2}\cup\R_{<-2})$, $\{z\in \C\, |\, \Im(z)>0\}\cup\R_{\geq 2}$,
$\{z\in\C\, |\, \Im(z)<0\}\cup \R_{\leq -2}$.
The maps $s\mapsto\sqrt{\tfrac14{s^2}-1}$ and $s\mapsto \lambda_+(s)$
are also continuous in these three regions.
\newline
(b) 
\begin{eqnarray}\label{5.15}
&&\begin{split}
(-i)e^{-\pi i\alpha_\pm} &= -\sqrt{\tfrac14{s^2}-1}\pm \tfrac12{s}
\textup{ for all }s\in\C,\\ 
\textup{and }&= \tfrac12 {(1+\lambda_\pm)}/{\sqrt{\tfrac14{s^2}-1}}
\textup{ for }s\neq\pm 2,
\end{split}
\\ \label{5.16}
&&\cos(\pi\alpha_\pm)= (-i)\sqrt{\tfrac14{s^2}-1},\ 
\sin(\pi\alpha_\pm)=\mp\tfrac12{s}.
\end{eqnarray}
\newline
(c) In the special case $s\in\R_{>2}$, define
\begin{eqnarray}\label{5.17}
t^{NI}=\frac{1}{2\pi}\log|\lambda_-|=\frac{1}{\pi}\log\sqrt{|\lambda_-|}
\in\R_{>0}.
\end{eqnarray}
Then
\begin{eqnarray}\label{5.18}
\alpha_\pm &=& \mp\left(\tfrac{1}{2}+it^{NI}\right),\\ \label{5.19}
\sqrt{|\lambda_\pm|}&=& |e^{-\pi i\alpha_\pm}| 
= |e^{\mp\pi t^{NI}}|= \mp\sqrt{\tfrac14{s^2}-1}+\tfrac{s}{2},
\\ \label{5.20}
\cosh(\pi t^{NI})&=&\tfrac12{s},\ \sinh(\pi t^{NI})=\sqrt{\tfrac14{s^2}-1}.
\end{eqnarray}
\end{lemma}

{\bf Proof:}
(a) \eqref{5.7}, \eqref{5.8} and \eqref{5.12} are obvious. For \eqref{5.9} observe that
$\lambda_+\MGcdot\lambda_-=1$. Now \eqref{5.10} and \eqref{5.11} follow
easily as follows.  They can be verified at some special values, and then
\eqref{5.9} and the bijectivity and continuity of the map in \eqref{5.14}
give \eqref{5.10} and \eqref{5.11} for all values. Similarly, \eqref{5.13}
follows from \eqref{5.12}.
The map \eqref{5.14} is bijective because the map 
$\C-\{\pm 2\}\to\C^*,\ s\mapsto \lambda_+(s)$, is injective.
The continuity in the first region is trivial. The continuity in the other
two regions follows easily from definitions \eqref{5.5} and \eqref{5.6}.
The continuity of the maps $s\mapsto\sqrt{\tfrac14{s^2}-1}$ and 
$s\mapsto \lambda_+(s)$ in the three regions follows from
the choice $\arg\sqrt{\tfrac14{s^2}-1}\in[0,\pi)$.

(b) The case $s=\pm 2$ is trivial. In the case $s\neq \pm 2$, we have
\begin{eqnarray}\label{5.21}
\left(\frac{1+\lambda_\pm}{2\sqrt{\tfrac14{s^2}-1}}\right)^2
&=& \frac{1+2\lambda_\pm+\lambda_\pm^2}{s^2-4}
=\frac{(4-s^2)\lambda_\pm}{s^2-4}\\ \nonumber
&=&-\lambda_\pm = -e^{-2\pi i \alpha_\pm}.
\end{eqnarray}
This and \eqref{5.3} give the formulae in \eqref{5.15} up to a sign.
At $s=0$ ($\lambda_\pm =1, \alpha_\pm=0$) the sign is correct.
Because of the continuity of the three maps in (a) in the three regions,
the sign in \eqref{5.15} is correct for all $s\neq \pm 2$.
\eqref{5.16} follows from \eqref{5.15} and from $\alpha_++\alpha_-=0$.

(c) These are easy consequences of part (b).
\hfill$\Box$

\begin{lemma}\label{t5.2}
(a) For any matrix $B\in SL(2,\C)$ the following equivalences hold.
\begin{eqnarray}
B\MGcdot S=S\MGcdot (B^t)^{-1}\
&\iff& S^t\MGcdot B^t=B^{-1}\MGcdot S^t\nonumber\\
&\iff& B\MGcdot S^t=S^t\MGcdot (B^t)^{-1}\label{5.22}\\ 
&\iff& B=\begin{pmatrix}b_1&b_2\\ -b_2&b_1+sb_2\end{pmatrix}.\nonumber
\end{eqnarray}
If $s\neq 0$ then these conditions are equivalent to 
\begin{eqnarray}\label{5.23}
\Mon_0^{mat}\MGcdot B= B\MGcdot \Mon_0^{mat}.
\end{eqnarray}

(b) If the conditions in (a) hold, and $s\neq \pm 2$, then
\begin{eqnarray}\label{5.24}
\begin{split}
B\MGcdot v_\pm &= b_\pm\MGcdot v_\pm\\ 
\textup{with }b_\pm &= 
(b_1+\tfrac12{s}b_2)\mp\sqrt{\tfrac14{s^2}-1}\MGcdot b_2
= b_1\mp i e^{-\pi i \alpha_\pm}\MGcdot b_2
\end{split}
\end{eqnarray}
and $b_+\MGcdot b_-=\det B=1$. If $s=\pm 2$, then
\begin{eqnarray}\label{5.25}
\begin{split}
B\MGcdot v_1&=\www b_1\MGcdot v_1,\ B\MGcdot v_2=\www b_1\MGcdot v_2+b_2\MGcdot v_1\\
\textup{with }\www b_1&= b_1+\tfrac12{s}b_2\in\{\pm 1\}. 
\end{split}
\end{eqnarray}
If $s\neq \pm 2$, $B$ as in (a) is determined by its eigenvalue
$b_-$, and any value $b_-\in\C^*$ is realizable by a matrix $B$ as in (a). 
If $s=\pm 2$, $B$ is determined by its eigenvalue
$\www b_1\in\{\pm 1\}$ and by $b_2$, and any pair 
$(\www b_1,b_2)\in \{\pm 1\}\times\C$ is realizable by a matrix as in (a).
\end{lemma}

The proof consists of elementary calculations and is omitted.
The eigenvalues $b_\pm$ of $B$ will be important in chapter \ref{s12}.

The following matrix $T(s)$ is a square root of $-\Mon_0^{mat}(s)$.
Its appearance looks surprising. Lemma \ref{t5.3} formulates its properties.
We shall use it in the proofs of lemma \ref{t13.1} (a) and theorem \ref{t14.1}
(a).

\begin{lemma}\label{t5.3}
Fix $s\in \C$. The matrix
\begin{eqnarray}\label{5.26}
T(s):=\begin{pmatrix}0&1\\-1&s\end{pmatrix}
\end{eqnarray}
satisfies
\begin{eqnarray}\label{5.27}
T(s)^{-1}&=&\begin{pmatrix}s&-1\\1&0\end{pmatrix},\ 
T(s)^2=-\Mon_0^{mat}(s),\\ 
T(s)\MGcdot B&=&B\MGcdot T(s)\quad\textup{for }B\textup{ as in }\eqref{5.22}.
\label{5.28}
\end{eqnarray}
If $s\neq \pm 2$, $T(s)$ is semisimple with eigenvectors $v_\pm$
and eigenvalues $\mp i e^{-\pi i\alpha_\pm}$, i.e.\ 
\begin{eqnarray}\label{5.29}
T(s)(v_+\  v_-)=(v_+\  v_-)\begin{pmatrix}-ie^{-\pi i\alpha_+} & 0 \\
0 & ie^{-\pi i\alpha_-}\end{pmatrix}.
\end{eqnarray}
If $s=\pm 2$, $T(s)$ has the eigenvalue $\tfrac12{s}$ and a 
$2\times 2$ Jordan block and satisfies
\begin{eqnarray}\label{5.30}
T(s)(v_1\ v_2)=(v_1\ v_2)\begin{pmatrix}\tfrac12{s} & 1 \\ 0 &\tfrac12{s}
\end{pmatrix}.
\end{eqnarray}
\end{lemma}
The proof consists of elementary calculations and is omitted.

\chapter{$P_{3D6}$-TEP bundles}\label{s6}
\setcounter{equation}{0}

\noindent
In this paper we are interested in the Painlev\'e III($D_6$) equation
of type $(\alpha,\beta,\gamma,\delta)=(0,0,4,-4)$
and in the isomonodromic families of $P_{3D6}$ bundles which are 
associated to it in \cite{FN80}, \cite{IN86}, \cite{FIKN06}, \cite{Ni09}.
These $P_{3D6}$ bundles are special and can be equipped with rich 
additional structure. We shall develop this structure in two steps
in chapters \ref{s6} and \ref{s7}. The most important
part is the TEP structure below. In chapter \ref{s7} it will be further
enriched to a TEJPA structure. 
Isomonodromic families of $P_{3D6}$-TEJPA bundles will correspond to
solutions of the  equation $P_{III}(0,0,4,-4)$ (theorem \ref{t10.3}).

First, some notation will be fixed.
The following holomorphic or antiholomorphic involutions of $\P^1$
will be used here and later:
\begin{eqnarray}\label{6.1}
j:\P^1\to\P^1,&& z\mapsto -z,\\ \label{6.2}
\rho_c:\P^1\to\P^1,&& z\mapsto \tfrac{1}{c}z
\quad\textup{for some }c\in\C^*,\\ \label{6.3}
\gamma:\P^1\to\P^1,&& z\mapsto {1}/{\oooo z},\\ \label{6.4}
\sigma:\P^1\to\P^1,&& z\mapsto -{1}/{\oooo z}.
\end{eqnarray}
Apart from $\rho_c$, these were also used in \cite{He03}, \cite{HS07}, \cite{HS11},
\cite{Mo11b}, \cite{Sa02}, \cite{Sa05b}.

\begin{definition}\label{t6.1}
(a) A TEP bundle is a holomorphic
vector bundle $H\to\P^1$ of rank $n\ge1$,
with a (flat) meromorphic connection $\nnn$ and a $\C$-bilinear pairing $P$.  
The connection is holomorphic on $\C^*$ and the pole at $0$ has order $\leq 2$.
The pairing
\begin{eqnarray}\label{6.5}
\textup{(pointwise)}&&P:H_z\times H_{j(z)}\to \C, \quad\textup{for all }z\in\P^1,
\\
\textup{(for sections)}&&P:\OO(H)\times j^*\OO(H)\to \OO_{\P^1},
\quad\OO_{\P^1}\textup{-linear}\nonumber
\end{eqnarray}
is symmetric, nondegenerate, and
flat on $H|_{\C^*}$.  
Symmetric means that
\begin{eqnarray*}
P(a(z),b(-z))=P(b(-z),a(z))
\end{eqnarray*}
for $a(z)\in H_z,b(-z)\in H_{-z}$,
and flatness means
\begin{eqnarray*}
\zdz P(a(z),b(-z))=P(\nnn_\zdz a(z),b(-z))+P(a(z),\nnn_\zdz b(-z))
\end{eqnarray*}
for $a\in \Gamma(U,\OO(H)),b\in \Gamma(j(U),\OO(H))$.

(b) A $P_{3D6}$-TEP bundle is a $P_{3D6}$ bundle which is also a TEP bundle.
\end{definition}

\begin{remarks}\label{t6.2}
(i) The definition of TEP bundles here differs from the definition of TEP 
structures in \cite{He03} and \cite{HS07}. The TEP structures there are
restrictions to $\C$ of TEP bundles here. 

(ii) The name TEP is intended to be self-referencing (cf.\ \cite{He03}, \cite{HS07}): 
T $=$ twistor $=$ holomorphic vector bundle on $\P^1$, 
E $=$ {\it extension} $=$ a meromorphic connection, P $=$ pairing.
By adding letters, one obtains richer structures (see chapters \ref{s7}
and \ref{s16}).

(iii) In formula \eqref{6.10} and later, we write the matrix of the pairing 
$P$ for a pair $(v_1,v_2)$ of vectors in $H_z$ and a pair $(v_3,v_4)$ of 
vectors in $H_{-z}$ as
$$P((v_1\ v_2)^t,(v_3\ v_4))=P(\begin{pmatrix}v_1\\v_2\end{pmatrix},(v_3\ v_4))=
\begin{pmatrix}P(v_1,v_3)& P(v_1,v_4)\\P(v_2,v_3)&P(v_2,v_4)\end{pmatrix}.$$
\end{remarks}

\begin{theorem}\label{t6.3}
(a) A $P_{3D6}$ bundle can be enriched to a $P_{3D6}$-TEP bundle if and only if
it is 
\begin{eqnarray}\label{6.6}
&&\text{irreducible or completely reducible, and}
\\
\label{6.7}
&&\alpha^1_0=\alpha^2_0=\alpha^1_\iiii=\alpha^2_\iiii=0.
\end{eqnarray}

(b) All TEP structures on a given $P_{3D6}$ bundle are isomorphic.

(c) Only two completely reducible $P_{3D6}$ bundles satisfying \eqref{6.7} exist.
In both cases, each flat rank 1 subbundle can be equipped with a 
TEP structure, which is unique up to rescaling.
The orthogonal sums of these TEP structures are the only TEP structures
on the $P_{3D6}$ bundle.
Therefore, such TEP structures are parametrized by $\C^*\times \C^*$.

(d) The TEP structures on an irreducible $P_{3D6}$ bundle with \eqref{6.7}
differ only by rescaling. Therefore they are parametrized by $\C^*$.

(e) Any $P_{3D6}$-TEP bundle has bases 
$\uuuu e^\pm_0=(e^{\pm 1}_0,e^{\pm 2}_0)$ and
$\uuuu e^\pm_\iiii=(e^{\pm 1}_\iiii,e^{\pm 2}_\iiii)$ such that
\begin{eqnarray}\label{6.8}
&&e^{\pm j}_0\textup{ and }e^{\pm j}_\iiii 
\textup{ are flat bases of }
L^{\pm j}_0\textup{ and }L^{\pm j}_\iiii,\\ \label{6.9}
&&e^{+1}_{0/\iiii}|_{\whhh I^a_{0/\iiii}} 
=e^{-1}_{0/\iiii}|_{\whhh I^a_{0/\iiii}}\textup{ \ and \ } 
e^{+2}_{0/\iiii}|_{\whhh I^b_{0/\iiii}}
=e^{-2}_{0/\iiii}|_{\whhh I^b_{0/\iiii}} 
\end{eqnarray}
(this is \eqref{2.12} and its analogue at $\iiii$ in the case of \eqref{6.7})
and 
\begin{eqnarray}\label{6.10} 
&&P((\uuuu e^\pm_0)^t(z),\uuuu e^\mp_0(-z))
=P((\uuuu e^\pm_\iiii)^t(z),\uuuu e^\mp_\iiii(-z))={\bf 1}_2
%=\begin{pmatrix}1&0\\0&1\end{pmatrix}
,
\\ \label{6.11}
&&\det B(\beta)=1 \textup{ for }\beta\textup{ and }B(\beta)
\textup{ as in \eqref{2.23},\eqref{2.24}.}
\end{eqnarray}

If $\uuuu e^\pm_0,\uuuu e^\pm_\iiii$ is a $4$-tuple of such bases,
then all others are obtained by modifying signs.
Altogether there are eight $4$-tuples of such bases, namely:
\begin{eqnarray}\label{6.12}
\varepsilon_0(e^{\pm 1}_0,\varepsilon_1 e^{\pm 2}_0),\ 
\varepsilon_0\varepsilon_2(e^{\pm 1}_\iiii,\varepsilon_1 e^{\pm 2}_\iiii)
\textup{ for }\varepsilon_0,\varepsilon_1,\varepsilon_2\in\{\pm 1\}.
\end{eqnarray}
For any $4$-tuple of such bases, the $P_{3D6}$ numerical tuple satisfies
\begin{eqnarray}\label{6.13}
\begin{split}
s^a_0&= s^b_0=-s^a_\iiii=-s^b_\iiii=:s, \\ 
B(\beta)&= \begin{pmatrix}b_1&b_2\\-b_2&b_1+sb_2\end{pmatrix}
\quad (\textup{with }\det B(\beta)=1).
\end{split}
\end{eqnarray}

(f) The automorphism group of a $P_{3D6}$-TEP bundle is
\begin{eqnarray}\label{6.14}
\begin{split}
\textup{Aut}(H,\nnn,P)&=\{\pm \id\}
\quad\textup{(irreducible case)}\\ 
\textup{Aut}(H,\nnn,P)&=\{\pm \id,\pm G\}
\quad\textup{(completely reducible case)}
\end{split}
\end{eqnarray}
with
\begin{eqnarray}\label{6.15}
&&G:\uuuu e^\pm_0\mapsto (e^{\pm 1}_0,-e^{\pm 2}_0),\
\uuuu e^\pm_\iiii\mapsto \varepsilon_2(e^{\pm 1}_\iiii,-e^{\pm 2}_\iiii)
\textup{ \ and}\\ \nonumber
&&\varepsilon_2=1\textup{ in the case }
(C_{11}),(C_{22})\ (\textup{i.e.\ }b_2=b_3=0),\\ \nonumber
&&\varepsilon_2=-1\textup{ in the case }
(C_{12}),(C_{21})\ (\textup{i.e.\ }b_1=b_4=0).
\end{eqnarray}

(g) Any $P_{3D6}$ numerical tuple of an irreducible $P_{3D6}$ bundle
with \eqref{6.7} satisfies either
$s^a_0\MGcdot s^b_0=s^a_\iiii\MGcdot s^b_\iiii\neq 0$ or
$s^a_0=s^b_0=s^a_\iiii=s^b_\iiii=0$.
\end{theorem}

{\bf Proof:}
We shall give the proofs in this order:
(a)($\Rightarrow$),  (g), (e), (a)($\Leftarrow$), (c), (d), (b), (f).

(a)($\Rightarrow$): Let $(H,\nnn,u^1_0,u^1_\iiii,P)$ be a $P_{3D6}$-TEP bundle.
First we discuss the pole at $0$ (the pole at $\iiii$ will be analogous).
Let $\uuuu e^\pm_0=(e^{\pm 1}_0,e^{\pm 2}_0)$ be bases of
$L|_{\whhh I^\pm_0}$ consisting of flat generating sections
of $L^{\pm 1}_0,L^{\pm 2}_0$ with \eqref{2.12}.
Choose one branch of $\log z$ on $\whhh I^+_0$ and extend it counterclockwise
to $\whhh I^-_0$. We shall write $(ze^{\pi i})^{\alpha^j_0}$ instead of 
$(-z)^{\alpha^j_0}$ to indicate this extension.

Choose a basis $\uuuu\varphi$ of $\OO(H)_0$ as in remark \ref{t2.5} (ii),
with $A^\pm_0\in GL(2,\AAA_{I^\pm_0})$, 
$\whhh A^+_0=\whhh A^-_0$ and $\whhh A^\pm_0(0)={\bf 1}_2.$
Then, for $z\in \whhh I^+_0$, the matrix 
$P((\varphi(z))^t,\varphi(-z))\in GL(2,\C\{z\})$ is 
\begin{eqnarray}\label{6.16}
&&(A^+_0(z))^t\MGcdot \begin{pmatrix}\gamma^{11}(z) &\gamma^{12}(z)\\
\gamma^{21}(z)&\gamma^{22}(z)\end{pmatrix}\MGcdot A^-_0(-z)\\
\nonumber
\textup{with}&& \gamma^{ij}(z)=z^{\alpha^i_0+\alpha^j_0}e^{\pi i\alpha^j_0}
e^{-(u^i_0-u^j_0)/z}\MGcdot P(e^{+ i}_0(z),e^{- j}_0(-z)).
\end{eqnarray}
Because $\whhh I^+_0$ contains the Stokes line $\R_{>0}\MGcdot \zeta_0$,
both functions $e^{-(u^1_0-u^2_0)/z}$ and $e^{-(u^2_0-u^1_0)/z}$
are unbounded on $\whhh I^+_0$. Therefore
\begin{eqnarray}\label{6.17}
\begin{split}
P(e^{+1}_0,e^{-2}_0)&=0=P(e^{+2}_0,e^{-1}_0),\\ 
P(e^{+1}_0,e^{-1}_0)&\neq 0,\ P(e^{+2}_0,e^{-2}_0)\neq 0.
\end{split}
\end{eqnarray}
Therefore the splittings
$L|_{\whhh I^+_0}=L^{\pm 1}_0\oplus L^{\pm 2}_0$ 
are dual to one another with respect to the pairing.
Thus $L^{+1}_0$ and $L^{-1}_0$ glue to a rank 1 subbundle of $L$ if and only
if $L^{+2}_0$ and $L^{-2}_0$ glue to a rank 1 subbundle of $L$.
Since $z^{2\alpha^i_0}$ must be holomorphic and nonvanishing at $0$, 
$\alpha^1_0=\alpha^2_0=0$. 

Everything so far holds analogously at $\iiii$.
This establishes \eqref{6.7}.
If $(H,\nnn)$ is reducible, then one condition $(C_{ij})$ from lemma
\ref{t3.1} (a) holds. But because of the duality of the subbundles, at $0$ 
as well as at $\iiii$, then also $(C_{3-i,3-j})$ holds,
and $(H,\nnn)$ is completely reducible by lemma \ref{t3.1} (b).
(a)($\Rightarrow$) is proved.

(g) Let $(H,\nnn)$ be an irreducible $P_{3D6}$ bundle satisfying \eqref{6.7}.
Choose bases $\uuuu e^\pm_0$ and $\uuuu e^\pm_\iiii$ as in the construction
of $P_{3D6}$ numerical tuples from $P_{3D6}$ monodromy tuples in chapter \ref{s2}.
Because of \eqref{2.22},
the monodromy matrices with respect to the bases $\uuuu e^+_0$ and 
$\uuuu e^-_\iiii$ are 
\begin{eqnarray}\label{6.18}
S^b_0(S^a_0)^{-1} \!=\!
\begin{pmatrix}1&-s^a_0\\s^b_0&1\!-\!s^a_0s^b_0\end{pmatrix},
(S^b_\iiii)^{-1}S^a_\iiii \!=\! 
\begin{pmatrix} 1 & s^a_\iiii\\
-s^b_\iiii& 1\!-\!s^a_\iiii s^b_\iiii\end{pmatrix}.
\end{eqnarray}
Therefore
\begin{eqnarray}\label{6.19}
s^a_0s^b_0=2-\tr(\Mon)=s^a_\iiii s^b_\iiii.
\end{eqnarray}
If $s^a_0=s^b_0=0$ (or $s^a_\iiii=s^b_\iiii=0$) then $\Mon=\id$ and 
also $s^a_\iiii=s^b_\iiii=0$ (or $s^a_0=s^b_0=0$).

Suppose $s^a_0\neq 0$ and $s^b_0=0$ (the other case $s^a_0=0$
and $s^b_0\neq 0$ is analogous). Then the monodromy is unipotent with 
a $2\times 2$ Jordan block.
As $s^b_0=0$, $L^{+1}_0$ and $L^{-1}_0$ glue to a flat rank 1 subbundle
of $L$. This must be the bundle of eigenspaces of the monodromy.

Since $s^a_\iiii s^b_\iiii=s^a_0 s^b_0=0$, the same argument applies
at $\iiii$. For some $k\in \{1,2\}$, $L^{+k}_\iiii$ and $L^{-k}_\iiii$ must glue
to the bundle of eigenspaces of the monodromy. Thus 
$(H,\nnn)$ is reducible, a contradiction.

(e) First, let us ignore \eqref{6.11}. Then
\eqref{6.17} and the analogous condition at $\iiii$ show the existence
of bases $\uuuu e^\pm_0$ and $\uuuu e^\pm_\iiii$ with 
\eqref{6.8}, \eqref{6.9}, and \eqref{6.10}.
It is also clear that each of
\[
(e^{+1}_0,e^{-1}_0),\ (e^{+2}_0,e^{-2}_0),\ 
(e^{+1}_\iiii,e^{-1}_\iiii),\ (e^{+2}_\iiii,e^{-2}_\iiii)
\]
is unique up to a sign. The following calculation, which will also be useful
for (a)($\Leftarrow$), will show that $\det B(\beta)=\pm 1$.
Recall \eqref{2.26} and \eqref{2.21} and the meaning of 
$\uuuu e^+_0(ze^{\pi i})$ explained after \eqref{2.21}.
For $z\in\whhh I^+_0$ and $\www z=z\MGcdot e^\beta\in \whhh I^-_\iiii$
\begin{eqnarray}
{\bf 1}_2 &=& P((\uuuu e^+_\iiii(-\www z))^t,\uuuu e^-_\iiii(\www z))
\nonumber \\ \nonumber
&=& P((\uuuu e^-_\iiii(\www ze^{\pi i})\MGcdot (S^a_\iiii)^{-1})^t,
\uuuu e^-_\iiii(\www z))\\ \nonumber
&=& ((S^a_\iiii)^{-1})^t\MGcdot 
P((\uuuu e^-_\iiii(\www ze^{\pi i}))^t,\uuuu e^-_\iiii(\www z))\\ \nonumber
&=& ((S^a_\iiii)^{-1})^t\MGcdot 
P((\uuuu e^+_0(ze^{\pi i})\MGcdot B(\beta))^t,\uuuu e^+_0(z)\MGcdot B(\beta))\\ 
\nonumber
&=& ((S^a_\iiii)^{-1})^t\MGcdot (B(\beta))^t\MGcdot 
P((\uuuu e^+_0(ze^{\pi i}))^t,\uuuu e^-_0(-ze^{\pi i}))\MGcdot
(S^b_0)^{-1}\MGcdot B(\beta) \\ \label{6.20} 
&=& ((S^a_\iiii)^{-1})^t\MGcdot (B(\beta))^t\MGcdot {\bf 1}_2\MGcdot
(S^b_0)^{-1}\MGcdot B(\beta).
\end{eqnarray}
This shows that $\det B(\beta)=\pm 1$.
If $\det B(\beta)=-1$, one can replace $(e^{+2}_\iiii,e^{-2}_\iiii)$
by $(-1)(e^{+2}_\iiii,e^{-2}_\iiii)$.
The new $4$-tuple of bases (or the old if $\det B(\beta)$ was already $1$)
satisfies \eqref{6.10} and \eqref{6.11}.
The possible sign changes in \eqref{6.12} are also clear now.

It remains to prove \eqref{6.13}. 
If $z\in \whhh I^a_0$ then $-z\in \whhh I^b_0$, and \eqref{2.21} and 
\eqref{6.10} show 
\begin{eqnarray*}
{\bf 1}_2&=& P((e^-_0)^t(z),e^+_0(-z)) = 
(S^a_0)^t\MGcdot P((e^+_0)^t(z),e^-_0(-z))\MGcdot (S^b_0)^{-1}\\
&=& (S^a_0)^t\MGcdot (S^b_0)^{-1}
=\begin{pmatrix}1 & 0 \\ s^a_0 & 1\end{pmatrix}
\begin{pmatrix}1 & 0 \\-s^b_0 & 1\end{pmatrix}
=\begin{pmatrix} 1 & 0 \\ s^a_0-s^b_0 & 1 \end{pmatrix},
\end{eqnarray*}
thus $s^a_0=s^b_0$.
Analogously one obtains $s^a_\iiii=s^b_\iiii$.
Define $s:=s^a_0=s^b_0$ and 
\[
S:=\begin{pmatrix}1 & s \\ 0 & 1\end{pmatrix}=S^a_0=(S^b_0)^t.
\]
Because of (g), $s^2=s^a_0\MGcdot s^b_0=s^a_\iiii\MGcdot s^b_\iiii$,
thus $ s^a_\iiii=s$ or $s^a_\iiii=-s$. 

If $s^a_\iiii=s\neq 0$ then \eqref{6.20} leads to 
\begin{eqnarray*}
{\bf 1}_2&=&(S^{-1})^t\MGcdot B(\beta)^t\MGcdot {\bf 1}_2 \MGcdot (S^t)^{-1}\MGcdot
B(\beta),\\
B(\beta)^{-1}\MGcdot S^t &=& (S^{-1})^t\MGcdot B(\beta)^t,
\end{eqnarray*}
and, with $B(\beta)=\bsp b_1 & b_2 \\ b_3 & b_4\esp$,
\begin{eqnarray*}
\begin{pmatrix}b_4-sb_2 & -b_2 \\ -b_3+sb_1 & b_1\end{pmatrix}
&=& \begin{pmatrix}b_4 & -b_2\\-b_3 & b_1\end{pmatrix}
\begin{pmatrix}1 & 0\\s & 1\end{pmatrix}\\
&=&\begin{pmatrix}1 & 0\\-s & 1\end{pmatrix}
\begin{pmatrix} b_1 & b_3\\ b_2 & b_4\end{pmatrix}
= \begin{pmatrix}b_1 & b_3 \\ b_2-sb_1 & b_4-sb_3\end{pmatrix},
\end{eqnarray*}
so $b_3=-b_2,\ sb_1=-sb_1,\ b_1=0, -sb_2=-sb_3,\ b_2=b_3=0, b_4=0,$
a contradiction. Thus $s^a_\iiii=-s$.
Now \eqref{6.20} gives 
$$B(\beta)^{-1}\MGcdot S^t = S^t\MGcdot B(\beta),$$
which is one of the equivalent conditions in lemma \ref{t5.2} (a).
This establishes \eqref{6.13}.

(a)($\Leftarrow$):
The 1st step is the construction of bases $\uuuu e^\pm_0,\uuuu e^\pm_\iiii$
which satisfy \eqref{6.8}, \eqref{6.9}, \eqref{6.11} and \eqref{6.13}.
The 2nd step is the definition of $P$ on
$L|_{\whhh I^+_0}\times L|_{\whhh I^-_0}$.
The 3rd step is to show that $P$ extends to $H$ and defines a TEP structure
and that the bases from the 1st step satisfy also \eqref{6.10}.

{\bf 1st step:}
Because of (g) there are three cases:
\begin{list}{}{}
\item[1st case:] 
$(H,\nnn)$ is irreducible and $s^a_0s^b_0=s^a_\iiii s^b_\iiii\neq 0.$
\item[2nd case:]
$(H,\nnn)$ is irreducible and $s^a_0=s^b_0=s^a_\iiii=s^b_\iiii=0.$
\item[3rd case:]
$(H,\nnn)$ is completely reducible (and $s^a_0=s^b_0=s^a_\iiii=s^b_\iiii=0$).
\end{list}

{\bf 1st case:}
Choose an $s\in\C^*$ with $s^2=2-\tr \Mon$. It is unique up to a sign.
The action of $(\C^*)^4$ in \eqref{2.29} on the $P_{3D6}$ numerical tuples
shows that we can choose flat generating sections $e^{\pm j}_0,e^{\pm j}_\iiii$
of $L^{\pm j}_0,L^{\pm j}_\iiii$ with \eqref{6.8} and \eqref{6.9}
and with
\begin{eqnarray}\label{6.21}
&&s=s^a_0=s^b_0=-s^a_\iiii=-s^b_\iiii\\ \label{6.22}
\textup{and}&&\det B(\beta)=1 \textup{ for some (or any) }\beta.
\end{eqnarray}
Then
\begin{eqnarray}\label{6.23}
S^a_0=\begin{pmatrix}1&s\\0&1\end{pmatrix}=:S,\ 
S^b_0=S^t,\ S^a_\iiii = S^{-1},\ S^b_\iiii=(S^t)^{-1}.
\end{eqnarray}
The monodromy matrices of $\uuuu e^+_0$ and $\uuuu e^-_\iiii$ are both
\begin{eqnarray}\label{6.24}
S^b_0(S^a_0)^{-1}=S^tS^{-1}=\Mon_0^{mat}:=
\begin{pmatrix}1&-s\\s& 1-s^2\end{pmatrix}
=(S^b_\iiii)^{-1}S^a_\iiii.
\end{eqnarray}
Condition \eqref{2.26} for $B(\beta)$ becomes
$$\Mon_0^{mat} \MGcdot B(\beta) = B(\beta)\MGcdot \Mon_0^{mat}.$$
Lemma \ref{t5.2} (a) gives the second half of \eqref{6.13},
$$B(\beta) = \begin{pmatrix} b_1&b_2\\ -b_2&b_1+sb_2\end{pmatrix}.$$

{\bf 2nd case:}
Here and in the 3rd case $\Mon=\id$, condition \eqref{2.26} is empty,
and the subbundles $L^{+j}_0$ and $L^{-j}_0$ and $L^{+j}_\iiii$ and $L^{-j}_\iiii$
glue to rank 1 subbundles $L^j_0$ and $L^j_\iiii$ ($j=1,2$) of $L$.

Choose generating flat sections $e^j_0$ and $e^j_\iiii$ of $L^j_0$ and 
$L^j_\iiii$. Then $\uuuu e_0=(e^1_0,e^2_0)$ and 
$\uuuu e_\iiii=(e^1_\iiii,e^2_\iiii)$ are flat global bases of $L$.
In the irreducible case
\begin{eqnarray}\label{6.25}
\uuuu e_\iiii = \uuuu e_0\MGcdot B\textup{ with }
B=\begin{pmatrix}b_1&b_2\\b_3&b_4\end{pmatrix}
\textup{ and all }b_i\neq 0
\end{eqnarray}
(here $B$ is independent of the choice of $\beta$).
We seek $(\lambda_1,\lambda_2,\lambda_3,\lambda_4)\in (\C^*)^4$
such that 
\begin{eqnarray*}
\www B:=\begin{pmatrix}\lambda_1^{-1}&0\\0&\lambda_2^{-1}\end{pmatrix}
B\begin{pmatrix}\lambda_3&0\\0&\lambda_4\end{pmatrix}
=\begin{pmatrix}\frac{\lambda_3}{\lambda_1}b_1 &\frac{\lambda_4}{\lambda_1}b_2\\
\frac{\lambda_3}{\lambda_2}b_3& \frac{\lambda_4}{\lambda_2}b_4\end{pmatrix}
=\begin{pmatrix}\www b_1& \www b_2\\ \www b_3& \www b_4\end{pmatrix}
\end{eqnarray*}
satisfies
$\www b_4=\www b_1$, $\www b_3=-\www b_2$, $\det \www B=1$, i.e.,
\begin{eqnarray}\label{6.26}
\tfrac{\lambda_4}{\lambda_2}b_4=\tfrac{\lambda_3}{\lambda_1}b_1,\ 
\tfrac{\lambda_3}{\lambda_2}b_3=-\tfrac{\lambda_4}{\lambda_1}b_2,\ 
1=\tfrac{\lambda_3\lambda_4}{\lambda_1\lambda_2}\det B.
\end{eqnarray}
This is equivalent to
\begin{eqnarray*}
\left(\frac{\lambda_2}{\lambda_1}\right)^2=-\frac{b_4b_3}{b_1b_2},\ 
\frac{\lambda_3}{\lambda_4}= \frac{\lambda_1}{\lambda_2}\MGcdot
\frac{b_4}{b_1},\ \lambda_3\lambda_4=\lambda_1\lambda_2\det B^{-1}.
\end{eqnarray*}
After making an arbitrary choice of $\lambda_1$, we see that $\lambda_2$ exists and is unique
up to a sign, and then $\lambda_3$ and $\lambda_4$
exist and are unique up to a common sign.

{\bf 3rd case:}
Generating sections $e^j_0$ and $e^j_\iiii$ as in the 2nd case are chosen.
In the completely reducible case
\begin{eqnarray}\label{6.27}
\begin{split}
\uuuu e_\iiii =\uuuu e_0\MGcdot B \textup{ with }B=
\begin{pmatrix}b_1&b_2\\b_3&b_4\end{pmatrix}\textup{ and} \\ 
\textup{either }b_2=b_3=0\textup{ or }b_1=b_4=0.
\end{split}
\end{eqnarray}
In both cases \eqref{6.26} is solvable. In the case $b_2=b_3=0$, we have
$\www B=\varepsilon\MGcdot {\bf 1}_2$ where $\varepsilon\in\{\pm 1\}$
has to be chosen.
Then one just needs ${\lambda_4}/{\lambda_2}=\varepsilon\MGcdot b_4^{-1}$,
${\lambda_3}/{\lambda_1}=\varepsilon\MGcdot b_1^{-1}$.
In the case $b_1=b_4=0$, we have
$\www B=\varepsilon\MGcdot \bsp 0&1\\-1&0\esp$ 
where $\varepsilon\in\{\pm 1\}$ has to be chosen.
Then one just needs ${\lambda_4}/{\lambda_1}=\varepsilon\MGcdot b_2^{-1}$,
${\lambda_3}/{\lambda_2}=-\varepsilon\MGcdot b_3^{-1}$.

{\bf 2nd step:}
In the 1st step bases $\uuuu e^\pm_0,\uuuu e^\pm_\iiii$ were constructed
which satisfy \eqref{6.8}, \eqref{6.9}, \eqref{6.11}
and \eqref{6.13}. Define $P$ on $L|_{\whhh I^+_0}\times L|_{\whhh I^-_0}$ by
\begin{eqnarray*}
P((\uuuu e^+_0)^t(z),\uuuu e^-_0(-z)):={\bf 1}_2.
\end{eqnarray*}

{\bf 3rd step:}
If $z\in \whhh I^a_0$ then $-z\in \whhh I^b_0$, and \eqref{2.21} gives
\begin{eqnarray}\label{6.28}
\begin{split}
P((\uuuu e^-_0)^t(z),\uuuu e^+_0(-z))
&= (S^a_0)^t\MGcdot P((\uuuu e^+_0)^t(z),e^-_0(-z))\MGcdot (S^b_0)^{-1} \\ 
&= S^t\MGcdot {\bf 1}_2\MGcdot (S^t)^{-1}={\bf 1}_2.
\end{split}
\end{eqnarray}
If $z\in \whhh I^b_0$ then $-z\in \whhh I^a_0$, and \eqref{2.21} gives
\begin{eqnarray}\label{6.29}
\begin{split}
P((\uuuu e^-_0)^t(z),\uuuu e^+_0(-z))
&= (S^b_0)^t\MGcdot P((\uuuu e^+_0)^t(z),e^-_0(-z))\MGcdot (S^a_0)^{-1} \\ 
&= S\MGcdot {\bf 1}_2\MGcdot S^{-1}={\bf 1}_2.
\end{split}
\end{eqnarray}
Therefore $P$ is well defined on $H|_{\C^*}$, nondegenerate and symmetric,
and satisfies the first half of \eqref{6.10}.
The calculation in \eqref{6.20} 
(without the first equality ${\bf 1}_2=...$) gives
\[
P((\uuuu e^+_\iiii)^t(-\www z),\uuuu e^-_\iiii(\www z))
=S^t\MGcdot B(\beta)^t\MGcdot {\bf 1}_2\MGcdot (S^t)^{-1}\MGcdot B(\beta).
\]
Because of \eqref{6.13} and lemma \ref{t4.2} (a) this is ${\bf 1}_2$.
With the symmetry of $P$ this establishes the second half of \eqref{6.10}.

It remains to show that $P$ extends in a nondegenerate way to $0$ and $\iiii$.
Choose a basis $\uuuu\varphi$ of $\OO(H)_0$ as in the proof of 
(a)($\Rightarrow$).
The matrix $P((\uuuu\varphi(z))^t,\varphi(-z))$ in \eqref{6.16}
for $z\in \whhh I^+_0$ becomes
\[
(A^+_0(z))^t\MGcdot {\bf 1}_2\MGcdot A^-_0(-z).
\]
With $\whhh A^+_0=\whhh A^-_0$ and $\whhh A^\pm_0(0)={\bf 1}_2$ 
this proves the claim.
The extension to $\iiii$ is shown analogously.
(a)($\Leftarrow$) is proved.

(c) In the notation of the 2nd and 3rd case in the proof of (a)($\Leftarrow$),
the two completely reducible bundles are the bundles with
\begin{eqnarray*}
\textup{either}&& b_2=b_3=0,\ L^1_0=L^1_\iiii,\ L^2_0=L^2_\iiii,
\textup{ i.e. }(C_{11}),(C_{22}),\\
\textup{or}&& b_1=b_4=0,\ L^1_0=L^2_\iiii,\ L^2_0=L^1_\iiii,
\textup{ i.e. }(C_{12}),(C_{21}).\\
\end{eqnarray*}
In both cases, (e) shows that any TEP structure splits into orthogonal
TEP structures on the two flat rank 1 subbundles.
It is clear that they can be rescaled independently.

(d) One has to consider separately the 1st case and the 2nd case in the proof
of (a)($\Leftarrow$) and see in both cases that the possible choices of
bases $\uuuu e^\pm_0$ and $\uuuu e^\pm_\iiii$ differ only by the signs 
in \eqref{6.11} and a common scalar (which absorbs the sign $\varepsilon_0$).
In both cases this is easy.

(b) This follows from (c) and (d).

(f) An automorphism of a $P_{3D6}$-TEP bundle maps a $4$-tuple 
$\uuuu e^\pm_0,\uuuu e^\pm_\iiii$ of bases in (e) to a $4$-tuple
$\www{\uuuu e}^\pm_0,\www{\uuuu e}^\pm_\iiii$ which must be one of the
eight $4$-tuples in \eqref{6.12}.
Then 
\[
s=\www s=\varepsilon_1\MGcdot s,\ b_1=\www b_1=\varepsilon_2\MGcdot b_1,\ 
b_2=\www b_2=\varepsilon_1\varepsilon_2\MGcdot b_2.
\]
Conversely, any  $(\varepsilon_0,\varepsilon_1,\varepsilon_2)
\in(\{\pm 1\})^3$ with these properties induces an automorphism.
In the irreducible case $\varepsilon_1=\varepsilon_2=1$ and 
$\varepsilon_0\in\{\pm 1\}$.
In the completely reducible case $s=0,b_2=0$: 
$\varepsilon_2=1$, $\varepsilon_0,\varepsilon_1\in\{\pm 1\}$.
In the completely reducible case $s=0,b_1=0$:
$\varepsilon_1\varepsilon_2=1$, $\varepsilon_0,\varepsilon_1=\varepsilon_2
\in\{\pm 1\}$.
\hfill $\Box$

\chapter[$P_{3D6}$-TEJPA bundles and their monodromy tuples ]
{$P_{3D6}$-TEJPA bundles and moduli spaces of their 
monodromy tuples}\label{s7}

\setcounter{equation}{0}

\noindent
We are concerned with $P_{3D6}$ bundles which satisfy \eqref{6.6} and \eqref{6.7}.
They can be equipped with TEP structures, and the TEP structure is unique
up to isomorphism.
Therefore from now on we consider $P_{3D6}$-TEP bundles.
As explained in the previous chapter, they possess eight distinguished $4$-tuples of bases  $\uuuu e^\pm_0,\uuuu e^\pm_\iiii$ (theorem \ref{t6.3} (e)).

This is good, but not good enough.
We want to distinguish two of the $4$-tuples, which differ only by
a global sign. This may be regarded as a \lq\lq marking\rq\rq\  or \lq\lq framing\rq\rq.  It can be expressed very elegantly in the 
choice of two isomorphisms $A$ and $J$, which are unique up to signs.
The choice of the signs corresponds to the marking.
The isomorphisms $A,J$ express certain symmetries of  $P_{3D6}$-TEP bundles (together with a reality condition, $J$ is related to the R in TERP structures, whose relation to real solutions of Painlev\'e III (0,0,4,-4) will be developed in chapters
\ref{s16} and \ref{s17}).

We recall from \eqref{6.1} and \eqref{6.2} the involutions 
\begin{eqnarray*}
j:\P^1&\to&\P^1,\ z\mapsto -z,\\ 
\rho_c:\P^1&\to&\P^1,\ z\mapsto \tfrac{1}{c} z
\end{eqnarray*}
of $\P^1$.

\begin{definition}\label{t7.1}
(a) A TEPA bundle is a TEP bundle together with a $\C$-linear  isomorphism
\begin{eqnarray}\label{7.1}
\textup{(pointwise)}&&A:H_z\to H_{j(z)},\quad \textup{for all }z\in\P^1,\\
\nonumber 
\textup{(for sections)}&&A:\OO(H)\to j^*\OO(H),
\quad \OO_{\P^1}\textup{-linear}
\end{eqnarray}
which is flat on $H|_{\C^*}$ and which satisfies
\begin{eqnarray}\label{7.2}
A^2&=& -\id\quad\textup{and}\\ \label{7.3}
P(A\,a,A\,b)&=& P(a,b) \text{ for all $a,b\in \OO(H)$. }
\end{eqnarray}

(b) A TEJPA bundle is a TEPA bundle together with some $c\in\C^*$
and a $\C$-linear isomorphism
\begin{eqnarray}\label{7.4}
\textup{(pointwise)}&&J:H_z\to H_{\rho_c(z)},\quad\textup{for all }z\in\P^1,\\
\nonumber 
\textup{(for sections)}&&J:\OO(H)\to \rho_c^*\OO(H),
\quad \OO_{\P^1}\textup{-linear}
\end{eqnarray}
which is flat on $H|_{\C^*}$ and which satisfies
\begin{eqnarray}\label{7.5}
J^2&=& \id,\\ \label{7.6}
P(J\,a,J\,b)&=& P(a,b)\quad\textup{for all $a,b\in \OO(H)$, and}\\ \label{7.7}
A\circ J&=& -J\circ A.
\end{eqnarray}

(c) A $P_{3D6}$-TEPA bundle is a $P_{3D6}$ bundle which is also
a TEPA bundle. A $P_{3D6}$-TEJPA bundle is a $P_{3D6}$ bundle which 
is, with $c={u^1_\iiii}/{u^1_0}$, also a TEJPA bundle.
\end{definition}

\begin{remarks}\label{t7.2}
(i) TEJP bundles and $P_{3D6}$-TEJP bundles are defined in the obvious way.
One can also extract from $\www J:=A\circ J$ in a TEJPA bundle
the properties of $\www J$ and define TE$\www J$P bundles 
and $P_{3D6}$-TE$\www J$P bundles. We shall not need these,
but we make a comment on them in remark \ref{t7.4} (i) and (ii).

(ii) If a $P_{3D6}$ bundle is also a TEJPA bundle then automatically
$c=\pm{u^1_\iiii}/{u^1_0}$.
The choice $c={u^1_\iiii}/{u^1_0}$ in definition \ref{t7.1} (c)
is good enough and avoids a pointless discussion of both cases.

(iii) The existence of $A$ is probably specific to
$P_{3D6}$ bundles. On the other hand, the cousin of $J$, 
the R in TERP structures, exists in great generality.

(iv) If $(H,\nnn,P)$ is a TEP bundle 
then also $j^*(H,\nnn,P)$ is a TEP bundle. 
In the case of a TEPA bundle, $A$ is an isomorphism of TEP bundles.

(v) Let $(H,\nnn,P,A,J)$ be a TEJPA bundle. 
Then the pole at $\iiii$ has order $\leq 2$,
$\rho_c^*(H,\nnn,P)$ is a TEP bundle, 
and $J$ is an isomorphism of TEP bundles.
\end{remarks}

\begin{theorem}\label{t7.3}
(a) A $P_{3D6}$-TEP bundle can be enriched to a $P_{3D6}$-TEPA bundle
if and only if 
\begin{eqnarray}\label{7.8}
u^1_0+u^2_0=0=u^1_\iiii+u^2_\iiii.
\end{eqnarray}
The isomorphism $A$ is then unique, up to a sign.

(b) A $P_{3D6}$-TEPA bundle can be enriched to a $P_{3D6}$-TEJPA bundle.
The isomorphism $J$ is unique up to a sign.

It follows that a $P_{3D6}$-TEP bundle can be enriched in four ways to a
$P_{3D6}$-TEJPA bundle. Given
$(H,\nnn,u^1_0,u^1_\iiii,P,A,J)$, these are
\begin{eqnarray}\label{7.9}
(H,\nnn,u^1_0,u^1_\iiii,P,\varepsilon_1\MGcdot A,\varepsilon_2\MGcdot J)
\quad\textup{with }\varepsilon_1,\varepsilon_2\in\{\pm 1\}.
\end{eqnarray}

(c) Let $(H,\nnn,u^1_0,u^1_\iiii,P,A,J)$ be a $P_{3D6}$-TEJPA bundle.
It has a $4$-tuple of bases $\uuuu e^\pm_0,\uuuu e^\pm_\iiii$, unique up to a global sign,  with \eqref{6.8},
\eqref{6.9} and 
\begin{eqnarray}\label{7.10}
P((\uuuu e^+_0)^t,\uuuu e^-_0)\!=\!{\bf 1}_2, 
A(\uuuu e^+_0)\!=\!\uuuu e^-_0 \!\begin{pmatrix}0&-1\\1&0\end{pmatrix}, 
J(\uuuu e^+_0)\!=\!\uuuu e^+_\iiii \!\begin{pmatrix}1&0\\0&-1\end{pmatrix}.
\end{eqnarray}
This is one of the eight $4$-tuples in theorem \ref{t6.3} (e).
It satisfies \eqref{6.10}, \eqref{6.11}, \eqref{6.13} and 
\begin{eqnarray}\label{7.11}
A(\uuuu e^\pm_0(z))\!=\!
\uuuu e^{\mp}_0(-z)\begin{pmatrix}0&\!-1\!\\1&0\end{pmatrix},
A(\uuuu e^\pm_\iiii(z))\!=\!
\uuuu e^{\mp}_\iiii(-z)\begin{pmatrix}0&\!-1\!\\1&0\end{pmatrix},\\
\label{7.12}
J(\uuuu e^\pm_0(z))\!=\!
\uuuu e^{\pm}_\iiii(\rho_c(z))\begin{pmatrix}1&0\\0&\!-1\!\end{pmatrix},
J(\uuuu e^\pm_\iiii(z))\!=\!
\uuuu e^{\pm}_0(\rho_c(z))\begin{pmatrix}1&0\\0&\!-1\!\end{pmatrix}.
\end{eqnarray}

(d) If the data in (c) and a $P_{3D6}$-TEJPA bundle in \eqref{7.9}
for some $\varepsilon_1,\varepsilon_2\in\{\pm 1\}$ are given,
the $4$-tuple of bases as in (c) for this $P_{3D6}$-TEJPA bundle
is 
\begin{eqnarray}\label{7.13}
\varepsilon_0(e^{\pm 1}_0,\varepsilon_1 e^{\pm 2}_0),
\varepsilon_0\varepsilon_2(e^{\pm 1}_\iiii,\varepsilon_1 e^{\pm 2}_\iiii)
\quad\textup{with arbitrary }\varepsilon_0\in\{\pm 1\}.
\end{eqnarray}
Thus the four possible enrichments to $P_{3D6}$-TEJPA bundles of a 
$P_{3D6}$-TEP bundle correspond to the four equivalence classes
which are obtained from the eight $4$-tuples in theorem \ref{t6.3} (e)
modulo a global sign.

(e) If the $P_{3D6}$ bundle which underlies a $P_{3D6}$-TEJPA bundle
is completely reducible, then the automorphism 
$G\in\textup{Aut}(H,\nnn,P)$ from theorem \ref{t6.3} (f) maps 
$A$ and $J$ to
\begin{eqnarray}\label{7.14}
\begin{split}
-A\textup{ and }J\textup{ in the case }(C_{11}), (C_{22}),\\ 
-A\textup{ and }-J\textup{ in the case }(C_{12}), (C_{21}).
\end{split}
\end{eqnarray}
\end{theorem}

{\bf Proof:}
(a) Let $(H,\nnn,u^1_0,u^1_\iiii,P)$ be a $P_{3D6}$-TEP bundle.
First we shall prove (a)($\Rightarrow$), then (a)($\Leftarrow$),
then the uniqueness of $A$ up to sign.

(a)($\Rightarrow$):
The bundle $j^*(H,\nnn,P)$ is almost a $P_{3D6}$-TEP bundle.
All that is missing is the choice of one of the eigenvalues
$\{-u^1_0,-u^2_0\}$ at $0$, and one of the eigenvalues 
$\{-u^1_\iiii,-u^2_\iiii\}$ at $\iiii$.
If an isomorphism $A:(H,\nnn,P)\to j^*(H,\nnn,P)$ exists
then $\{u^1_0,u^2_0\}=\{-u^1_0,-u^2_0\}$ and 
$\{u^1_\iiii,u^2_\iiii\}=\{-u^1_\iiii,-u^2_\iiii\}.$
Since $u^1_0\neq u^2_0$ and $u^1_\iiii\neq u^2_\iiii$,
this implies $u^1_0+u^2_0=0$ and $u^1_\iiii+u^2_\iiii=0$,
hence \eqref{7.8}.

(a)($\Leftarrow$):
Suppose that \eqref{7.8} holds.
Then $(j^*H,j^*\nnn,u^1_0,u^1_\iiii,j^*P)$ is a $P_{3D6}$-TEP bundle.
Let us use a tilde to indicate the data of its $P_{3D6}$ monodromy tuple:
\begin{eqnarray}\label{7.15}
\begin{split}
I^\pm_0=j^*I^\mp_0,\ I^\pm_\iiii=j^*I^\mp_\iiii,\ 
I^{a/b}_0=j^*I^{b/a}_0,\ I^{a/b}_\iiii=j^*I^{b/a}_\iiii,\\ 
\www L=j^*L,\ \www L^{\pm k}_0=j^* L^{\mp 3-k}_0,\ 
\www L^{\pm k}_\iiii = j^* L^{\mp 3-k}_\iiii.
\end{split}
\end{eqnarray}
Let $\uuuu e^\pm_0, \uuuu e^\pm_\iiii$ be one of the eight $4$-tuples
of bases in theorem \ref{t6.3} (e)
for the $P_{3D6}$-TEP bundle $(H,\nnn,u^1_0,u^1_\iiii,P)$.

{\bf Claim:}
The $4$-tuple of bases
\begin{eqnarray}\label{7.16}
\www{\uuuu e}^\pm_0(z):=\uuuu e^\mp_0(-z)\begin{pmatrix}0&-1\\1&0\end{pmatrix},\ 
\www{\uuuu e}^\pm_\iiii(z):=\uuuu e^\mp_\iiii(-z)
\begin{pmatrix}0&-1\\1&0\end{pmatrix}
\end{eqnarray}
is one of the eight $4$-tuples in theorem \ref{6.3} (e) 
for the $P_{3D6}$-TEP bundle $(j^*H,j^*\nnn,u^1_0,u^1_\iiii,j^*P)$, and 
\begin{eqnarray}\label{7.17}
\www s = s,\quad \www B(\beta)=B(\beta).
\end{eqnarray}

The claim implies that there exists an isomorphism 
$A:(H,\nnn,P)\to j^*(H,\nnn,P)$ with
\begin{eqnarray}\label{7.18}
A(\uuuu e^\pm_0(z))=\uuuu e^\mp_0(\!-\!z)\begin{pmatrix}0&\!-\!1\\1&0\end{pmatrix},\ 
A(\uuuu e^\pm_\iiii(z)):=\uuuu e^\mp_\iiii(-z)
\begin{pmatrix}0&\!-\!1\\1&0\end{pmatrix}.
\end{eqnarray}
Then \eqref{6.10} for the new and old $4$-tuples shows that
$P(A\,a,A\,b)=P(a,b)$, and $A^2=-\id$ follows from \eqref{7.18}.
It remains to prove the claim.

{\bf Proof of the claim:}
\eqref{6.9} for $\www{\uuuu e}^\pm_0$ follows from
\begin{eqnarray*}
\www{\uuuu e}^{+1}_0|_{I^a_0}= j^* \uuuu e^{-2}_0|_{j^*I^b_0} 
= j^* \uuuu e^{+2}_0|_{j^*I^b_0} = \www{\uuuu e}^{-1}_0|_{I^a_0},\\
\www{\uuuu e}^{+2}_0|_{I^b_0}= -j^* \uuuu e^{-1}_0|_{j^*I^a_0} 
= -j^* \uuuu e^{+1}_0|_{j^*I^a_0} = \www{\uuuu e}^{-2}_0|_{I^b_0}
\end{eqnarray*}
and the analogue at $\iiii$ is similar.
Next, \eqref{6.10} for $\www{\uuuu e}^\pm_0$ follows from
\begin{eqnarray*}
P((\www{\uuuu e}^\pm_0)^t(z),\www{\uuuu e}^\mp_0(-z))
&=&\begin{pmatrix}0&-1\\1&0\end{pmatrix}^t
P((\uuuu e^\mp_0)^t(-z),\uuuu e^\pm_0(z))
\begin{pmatrix}0&-1\\1&0\end{pmatrix}\\
&=&\begin{pmatrix}0&-1\\1&0\end{pmatrix}^t
{\bf 1}_2
\begin{pmatrix}0&-1\\1&0\end{pmatrix}
={\bf 1}_2
\end{eqnarray*}
and the analogue at $\iiii$ is similar.
The following calculation shows $\www B(\beta)=B(\beta)$
(and thus $\det \www B(\beta)=1$, which is \eqref{6.11}).
Choose $z\in \whhh I^a_0$ and $\beta$ as in \eqref{2.23} and define 
$y:=ze^\beta$. Then $z$ and $y$ satisfy
$$-z\in \whhh I^b_0,\ y\in \whhh I^a_\iiii,\ 
-y=-ze^\beta\in \whhh I^b_\iiii.$$
Recall \eqref{2.21}:
\begin{eqnarray*}
&&(\www{\uuuu e}^+_0(-z)\textup{ extended along }[\beta]\textup{ to }(-y))
\MGcdot \www B(\beta) \\
&=& \www{\uuuu e}^-_\iiii(-y) 
= \uuuu e^+_\iiii(y)\MGcdot\begin{pmatrix}0&-1\\1&0\end{pmatrix}\\
&=& \uuuu e^-_\iiii(y)\MGcdot (S^a_\iiii)^{-1}
\begin{pmatrix}0&-1\\1&0\end{pmatrix}\\
&=& (\uuuu e^+_0(z)\textup{ extended along }[\beta]\textup{ to }y)
\MGcdot B(\beta)(S^a_\iiii)^{-1}
\begin{pmatrix}0&-1\\1&0\end{pmatrix}\\
&=& (\uuuu e^-_0(z)\textup{ extended along }[\beta]\textup{ to }y)
\MGcdot (S^a_0)^{-1}B(\beta)(S^a_\iiii)^{-1}
\begin{pmatrix}0&-1\\1&0\end{pmatrix}\\
&=& (\www{\uuuu e}^+_0(-z)\textup{ extended along }[\beta]\textup{ to }(-y))
\MGcdot \\ 
&&\quad\MGtimes \begin{pmatrix}0&-1\\1&0\end{pmatrix}^{-1}
(S^a_0)^{-1}B(\beta)(S^a_\iiii)^{-1}
\begin{pmatrix}0&-1\\1&0\end{pmatrix}\\
&=& (\www{\uuuu e}^+_0(-z)\textup{ extended along }[\beta]\textup{ to }(-y))
\MGcdot B(\beta).
\end{eqnarray*}
The last equality uses \eqref{6.13} and lemma \ref{t5.2} (a).
The facts proved up to now show that $\www{\uuuu e}^\pm_0,\www{\uuuu e}^\pm_\iiii$
is one of the eight $4$-tuples from theorem \ref{6.3} (e) for
$(j^*H,j^*\nnn,u^1_0,u^1_\iiii,j^*P)$. Therefore
$\www s^a_0=\www s^b_0=-\www s^a_\iiii=-\www s^b_\iiii=:\www s$ holds.
Finally $\www s=s$ follows with 
$z\in \whhh I^a_0$ and $-z\in \whhh I^b_0$ from the following calculation:
\begin{eqnarray*}
\www{\uuuu e}^-_0(-z)&=& \uuuu e^+_0(z)\begin{pmatrix}0&-1\\1&0\end{pmatrix}
= \uuuu e^-_0(z)(S^a_0)^{-1}\begin{pmatrix}0&-1\\1&0\end{pmatrix}\\
&=& \www{\uuuu e}^+_0(-z)\begin{pmatrix}0&-1\\1&0\end{pmatrix}^{-1}
(S^a_0)^{-1}\begin{pmatrix}0&-1\\1&0\end{pmatrix}\\
&=& \www{\uuuu e}^+_0(-z)\begin{pmatrix}1&0\\s&1\end{pmatrix}.
\end{eqnarray*}
The claim is proved. \hfill $(\Box)$

{\bf The uniqueness of $A$ up to a sign:}
Suppose $A^{(1)}$ is a second isomorphism which gives a TEPA bundle.
Then $A^{(1)}\circ A\in\textup{Aut}(H,\nnn,P)$.
If $(H,\nnn)$ is irreducible then $\textup{Aut}(H,\nnn,P)=\{\pm \id\}$
by theorem \ref{t6.3} (f), and $A^{(1)}=\pm A^{-1}=\mp A.$
If $(H,\nnn)$ is completely reducible then
$\textup{Aut}(H,\nnn,P)=\{\pm \id,\pm G\}$ with $G$ as in theorem 
\ref{t6.3} (f). But then 
\[
A\circ G=-G\circ A,\quad\textup{thus }
(G\circ A^{-1})^2=(G\circ A)^2=\id\neq-\id,
\]
hence $A^{(1)}\circ A=\pm G$ is impossible.

(b) The proof is similar to the proof of (a).
First the existence of $J$ is shown, then the uniqueness up to a sign.
Observe that
\[
c=\frac{u^1_\iiii}{u^1_0}=\frac{-u^1_\iiii}{-u^1_0}=\frac{u^2_\iiii}{u^2_0},
\quad \frac{u^j_0}{z}=u^j_\iiii\MGcdot\rho_c(z).
\]
Therefore $(\rho_c^*H,\rho_c^*\nnn,u^1_0,u^1_\iiii,\rho_c^*P)$ 
is a $P_{3D6}$-TEP bundle.
The data of its $P_{3D6}$ monodromy tuple are denoted with a tilde.
\begin{eqnarray}\label{7.19}
\begin{split}
I^\pm_0=\rho^*_cI^\pm_\iiii,\ I^\pm_\iiii=\rho_c^*I^\pm_0,\ 
I^{a/b}_0=\rho_c^*I^{a/b}_\iiii,\ I^{a/b}_\iiii=\rho_c^*I^{a/b}_0,\\ 
\www L=j^*L,\ \www L^{\pm j}_0=\rho_c^* L^{\pm j}_\iiii,\ 
\www L^{\pm j}_\iiii = \rho_c^* L^{\pm j}_0.
\end{split}
\end{eqnarray}
Let $\uuuu e^\pm_0, \uuuu e^\pm_\iiii$ be one of the eight $4$-tuples
of bases in theorem \ref{t6.3} (e)
for the $P_{3D6}$-TEP bundle $(H,\nnn,u^1_0,u^1_\iiii,P)$,
which satisfies additionally \eqref{7.18}.

{\bf Claim:}
The $4$-tuple of bases
\begin{eqnarray}\label{7.20}
\www{\uuuu e}^\pm_0(z):=\uuuu e^\pm_\iiii(\rho_c(z))
\begin{pmatrix}1&0\\0&-1\end{pmatrix},\ 
\www{\uuuu e}^\pm_\iiii(z):=\uuuu e^\pm_0(\rho_c(z))
\begin{pmatrix}1&0\\0&-1\end{pmatrix}
\end{eqnarray}
is one of the eight $4$-tuples in theorem \ref{t6.3} (e) 
for the $P_{3D6}$-TEP bundle
$(\rho_c^*H,\rho_c^*\nnn,u^1_0,u^1_\iiii,\rho_c^*P)$, and 
\begin{eqnarray}\label{7.21}
\www s = s,\quad \www B(\beta)=B(\beta).
\end{eqnarray}

{\bf Proof of the claim:}
\eqref{6.9} for $\www{\uuuu e}^\pm_0$ follows from
\begin{eqnarray*}
\www{\uuuu e}^{+1}_0|_{I^a_0}= \rho_c^* \uuuu e^{+1}_\iiii|_{\rho_c^*I^a_\iiii} 
= \rho_c^* \uuuu e^{-1}_\iiii|_{\rho_c^*I^a_\iiii} 
= \www{\uuuu e}^{-1}_0|_{I^a_0},\\
\www{\uuuu e}^{+2}_0|_{I^b_0}= -\rho_c^* \uuuu e^{+2}_\iiii|_{\rho_c^*I^b_\iiii} 
= -\rho_c^* \uuuu e^{-2}_\iiii|_{\rho_c^*I^b_\iiii} 
= \www{\uuuu e}^{-2}_0|_{I^b_0}
\end{eqnarray*}
and the analogue at $\iiii$ is similar.
Next, \eqref{6.10} for $\www{\uuuu e}^\pm_0$ follows from
\begin{align*}
P((\www{\uuuu e}^\pm_0)^t(z),\www{\uuuu e}^\mp_0(-z))
&\!=\!\begin{pmatrix}1&0\\0&-1\end{pmatrix}^t
\! P((\uuuu e^\pm_\iiii)^t(\rho_c(z)),\uuuu e^\mp_\iiii(-\rho_c(z)))
\! \begin{pmatrix}1&0\\0&-1\end{pmatrix}
\\
&=\begin{pmatrix}1&0\\0&-1\end{pmatrix}^t
{\bf 1}_2
\begin{pmatrix}1&0\\0&-1\end{pmatrix}
={\bf 1}_2
\end{align*}
and the analogue at $\iiii$ is similar.
The following calculation shows $\www B(\beta)=B(\beta)$
(and thus $\det \www B(\beta)=1$, which is \eqref{6.11}).
Choose $y\in \whhh I^a_0$ and $\beta$ as in \eqref{2.23}
and define $z:=ye^\beta$. Then 
$$z\in \whhh I^a_\iiii,\ \rho_c(y)=\rho_c(z)e^\beta,\ 
\rho_c(z)\in \whhh I^a_0,\ 
\rho_c(y)\in \whhh I^a_\iiii.
$$ 
Recall \eqref{2.21}.
The fifth equality in the following calculation uses the fact that
$\rho_c$ inverts the path $[\beta]$. 
The last equality uses \eqref{6.13} and lemma \ref{t5.2} (a).

\begin{eqnarray*}
&&(\www{\uuuu e}^+_0(\rho_c(z))
\textup{ extended along }[\beta]\textup{ to }\rho_c(y))
\MGcdot \www B(\beta) \\
&=& \www{\uuuu e}^-_\iiii(\rho_c(y)) 
= \uuuu e^-_0(y)\MGcdot\begin{pmatrix}1&0\\0&-1\end{pmatrix}\\
&=& \uuuu e^+_0(y)\MGcdot S^a_0
\begin{pmatrix}1&0\\0&-1\end{pmatrix}\\
&=& (\uuuu e^-_\iiii(z)\textup{ extended along }[-\beta]\textup{ to }y)
\MGcdot B(\beta)^{-1}S^a_0
\begin{pmatrix}1&0\\0&-1\end{pmatrix}\\
&=& (\uuuu e^+_\iiii(z)\textup{ extended along }[-\beta]\textup{ to }y)
\MGcdot S^a_\iiii B(\beta)^{-1}S^a_0
\begin{pmatrix}1&0\\0&-1\end{pmatrix}\\
&=& (\rho_c^*(\uuuu e^+_\iiii(z))\textup{ extended along }[\beta]
\textup{ to }\rho_c(y))
\MGcdot \\ 
&& %\begin{pmatrix}1&0\\0&-1\end{pmatrix}
\quad\MGtimes
S^a_\iiii B(\beta)^{-1} S^a_0
\begin{pmatrix}1&0\\0&-1\end{pmatrix}\\
&=& (\www{\uuuu e}^+_0(\rho_c(z))\textup{ extended along }[\beta]
\textup{ to }\rho_c(y))
\MGcdot \\
&& \quad\MGtimes \begin{pmatrix}1&0\\0&-1\end{pmatrix}
S^a_\iiii B(\beta)^{-1} S^a_0
\begin{pmatrix}1&0\\0&-1\end{pmatrix}\\
&=& (\www{\uuuu e}^+_0(\rho_c(z))\textup{ extended along }[\beta]
\textup{ to }\rho_c(y)) \MGcdot B(\beta).
\end{eqnarray*}

The facts proved up to now show that $\www{\uuuu e}^\pm_0,\www{\uuuu e}^\pm_\iiii$
is one of the eight $4$-tuples from theorem \ref{6.3} (e) for
$(\rho_c^*H,\rho_c^*\nnn,u^1_0,u^1_\iiii,\rho_c^*P)$. Therefore
$\www s^a_0=\www s^b_0=-\www s^a_\iiii=-\www s^b_\iiii=:\www s$ holds.
Finally $\www s=s$ follows with 
$z\in \whhh I^a_0$ and $\rho_c(z)\in \whhh I^a_\iiii$ from the calculation
\begin{eqnarray*}
\www{\uuuu e}^-_\iiii(\rho_c(z))&=& 
\uuuu e^-_0(z)\begin{pmatrix}1&0\\0&-1\end{pmatrix}
= \uuuu e^+_0(z)S^a_0\begin{pmatrix}1&0\\0&-1\end{pmatrix}\\
&=& \www{\uuuu e}^+_\iiii(\rho_c(z))\begin{pmatrix}1&0\\0&-1\end{pmatrix}^{-1}
S^a_0\begin{pmatrix}1&0\\0&-1\end{pmatrix}\\
&=& \www{\uuuu e}^+_\iiii(\rho_c(z))\begin{pmatrix}1&-s\\0&1\end{pmatrix}.
\end{eqnarray*}
The claim is proved. \hfill $(\Box)$

The claim implies that an isomorphism 
$J:(H,\nnn,P)\to \rho_c^*(H,\nnn,P)$ exists with
\begin{eqnarray}\label{7.22}
\! J(\uuuu e^\pm_0(z))\!=\!\uuuu e^\pm_\iiii(\rho_c(z))
\!\! \begin{pmatrix}1&0\\0&-1\end{pmatrix},\ 
\! J(\uuuu e^\pm_\iiii(z))\!=\!\uuuu e^\pm_0(\rho_c(z))
\!\! \begin{pmatrix}1&0\\0&-1\end{pmatrix}.
\end{eqnarray}
\eqref{6.10} for the new and the old $4$-tuple shows
$P(J\,a,J\,b)=P(a,b)$. \eqref{7.22} implies $J^2=\id$,
and $A\circ J=-J\circ A$ follows from \eqref{7.22} and \eqref{7.18}.

{\bf The uniqueness of $J$ up to a sign:}
Suppose $J^{(1)}$ is a second isomorphism which gives a TEJPA structure.
Then $J^{(1)}\circ J\in\textup{Aut}(H,\nnn,P)$.
If $(H,\nnn)$ is irreducible then $\textup{Aut}(H,\nnn,P)=\{\pm \id\}$
by theorem \ref{t6.3} (f), and $J^{(1)}=\pm J^{-1}=\pm J.$
If $(H,\nnn)$ is completely reducible then
$\textup{Aut}(H,\nnn,P)=\{\pm \id,\pm G)$ with $G$ as in theorem 
\ref{t6.3} (f). But then $A\circ G=-G\circ A$.
Then $J^{(1)}=\pm G\circ J^{-1}=\pm G\circ J$ is impossible because
$$A\circ (G\circ J)=-G\circ A\circ J=(G\circ J)\circ A,
\textup{ but }A\circ J^{(1)}=-J^{(1)}\circ A.$$

(c) The existence of a $4$-tuple with all properties in (c) follows
from the constructions in (a) and (b).
$\varepsilon_0(e^\pm_1,\varepsilon_1 e^\pm_2)$ is determined up to the 
signs $\varepsilon_0,\varepsilon_1\in\{\pm 1\}$ by 
\eqref{6.9} and $P((\uuuu e^+)^t,\uuuu e^-)={\bf 1}_2$.
Then $A(\uuuu e^+_0)=\uuuu e^-_0\MGcdot
\bsp0&-1\\1&0\esp$
fixes the sign $\varepsilon_1=1$.
$J(\uuuu e^+_0)=\uuuu e^+_\iiii\MGcdot
\bsp 1&0\\0&-1\esp$
determines $\uuuu e^+_\iiii$. \eqref{6.9} determines $\uuuu e^-_\iiii$.

(d) This follows from the uniqueness in (c), from \eqref{7.9}, and
from \eqref{7.10}.

(e) This follows from (d) and theorem \ref{t6.3} (f).
\hfill$\Box$

\begin{remarks}\label{t7.4}
(i) Given a $P_{3D6}$-TEP bundle, a $J$ which enriches it to a 
$P_{3D6}$-TEJP bundle is unique up to a sign in the irreducible cases
and in the completely reducible case $(C_{12}),(C_{21})$,
because $\textup{Aut}(H,\nnn,P)=\{\pm \id\}$ in the irreducible cases
and $G\circ J=-J\circ G$ in the case $(C_{12}),(C_{21})$,
so then $(G\circ J)^2=-\id$ and $G\circ J$ is not a candidate for another $J$.
But in the completely reducible case $(C_{11}),(C_{21})$, we have $G\circ J=J\circ G$
by \eqref{7.13}, and $\pm J,\pm G\circ J$ give four possible enrichments
to a $P_{3D6}$-TEJP bundle.

(ii) This apparent asymmetry between the two completely reducible cases
is resolved by considering also $\www J$ as in remark \ref{t7.2} (i).
But the asymmetry between the reducible and the irreducible cases is not resolvable.
It leads to the two $A_1$-singularities in the moduli space
$M_{3T}^{mon}$ in theorem \ref{t7.6}
\end{remarks}

Now we turn to moduli spaces of $P_{3D6}$-TEJPA bundles and 
$P_{3D6}$-TEP bundles. Using the distinguished $4$-tuples of bases
in theorem \ref{t7.3} (c) and their $P_{3D6}$ numerical tuples, we shall see that
these moduli spaces are quotients of algebraic manifolds and come equipped with foliations, in a way which makes the isomonodromic families transparent.

We introduce now the fundamental moduli spaces
\begin{align*}
M_{3TJ}&=\text{\{isomorphism classes of $P_{3D6}$-TEJPA bundles\}},
\\
M_{3T}&=\text{\{isomorphism classes of $P_{3D6}$-TEP bundles\}},
\end{align*}
and we denote by 
$pr^u$ the natural projections
\begin{eqnarray}\label{7.23}
\begin{split}
pr^u:M_{3TJ}\to\C^*\times\C^*,
(H,\nnn,u^1_0,u^1_\iiii,P,A,J)\mapsto (u^1_0,u^1_\iiii),\\ 
pr^u:M_{3T}\to\C^*\times\C^*,
(H,\nnn,u^1_0,u^1_\iiii,P)\mapsto (u^1_0,u^1_\iiii),
\end{split}
\end{eqnarray}
with fibres $M_{3TJ}(u^1_0,u^1_\iiii)$, $M_{3T}(u^1_0,u^1_\iiii)$, respectively.

Consider the covering
\begin{eqnarray}\label{7.24}
c^{path}:\C\times\C^*\to\C^*\times\C^*,\ 
(\beta,u^1_0)\mapsto(u^1_0,\tfrac14{e^{-\beta}}/u^1_0).
\end{eqnarray}
It is the covering in which the triples $(\beta,u^1_0,u^1_\iiii)$ with
$e^{-\beta}=4u^1_0u^1_\iiii$ live, which is \eqref{2.23}
in the case $u^2_0=-u^1_0,u^2_\iiii=-u^1_\iiii$.
Recall for $s\in\C$ the definition in \eqref{5.1} of the matrices
\[
S(s)=\begin{pmatrix}1&s\\0&1\end{pmatrix},\quad
\Mon_0^{mat}(s)=S^t\MGcdot S^{-1}=\begin{pmatrix}1&-s\\s&1-s^2\end{pmatrix}.
\]
Define the affine algebraic surface
\begin{eqnarray}\label{7.25}
V^{mon}:=\{(s,b_1,b_2)\in\C^3\, |\, b_1^2+b_2^2+sb_1b_2=1\}
\end{eqnarray}
and the canonically isomorphic space
\begin{eqnarray}\label{7.26}
V^{mat}=\{(s,B)\in\C\times SL(2,\C)\, |\,
B=
\bsp
b_1&b_2\\-b_2&b_1+sb_2
\esp\}.
\end{eqnarray}

\begin{theorem}\label{t7.5}
(a) $V^{mon}$ and $V^{mat}$ are smooth.

(b) The map
\begin{eqnarray}\label{7.27}
\Phi^{path}:(c^{path})^*M_{3TJ}&\to& \C\times \C^*\times V^{mat}\\
((H,\nnn,u^1_0,u^1_\iiii,P,A,J),\beta)&\mapsto& 
(\beta,u^1_0,s,B(\beta))\nonumber
\end{eqnarray}
is a bijection. Here $s$ and $B(\beta)$ are the data from \eqref{6.13}
which are associated to the distinguished $4$-tuple of bases
$\uuuu e^\pm_0,\uuuu e^\pm_\iiii$ in theorem \ref{t7.3} (c)
for the $P_{3D6}$-TEJPA bundle $(H,\nnn,u^1_0,u^1_\iiii,P,A,J)$.

$\Phi^{path}$ equips the set $(c^{path})^*M_{3TJ}$ with the structure of an
affine algebraic manifold, which is denoted by
$((c^{path})^*M_{3TJ})^{mon}$.

(c) $\Phi^{path}$ maps the isomorphism classes of $P_{3D6}$-TEJPA bundles
in $(c^{path})^*M_{3TJ}$ to the orbits of the action of the group $\Z$
on $\C\times \C^*\times V^{mat}$, whose generator $[1]$ acts by the
algebraic automorphism
\begin{eqnarray}\label{7.28}
m_{[1]}:(\beta,u^1_0,s,B)\mapsto (\beta+2\pi i,u^1_0,s,(\Mon_0^{mat})^{-1}\MGcdot B).
\end{eqnarray}
The quotient $\C\times \C^*\times V^{mat}/\langle m_{[1]}\rangle$
is an analytic manifold.
$\Phi^{path}$ induces a bijection
\begin{eqnarray}\label{7.29}
\Phi_{3TJ}:M_{3TJ}\to \C\times\C^*\times V^{mat}/\langle m_{[1]}\rangle.
\end{eqnarray}
This equips $M_{3TJ}$ with the structure of an analytic manifold,
which is denoted by $M_{3TJ}^{mon}$.
The projection $pr^u:M_{3TJ}^{mon}\to\C^*\times \C^*$ is an
analytic morphism, the fibres are
algebraic manifolds $M_{3TJ}^{mon}(u^1_0,u^1_\iiii)$,
isomorphic to $V^{mat}$.
Any choice of $\beta$ (with \eqref{2.23}) induces an isomorphism
$M_{3TJ}^{mon}(u^1_0,u^1_\iiii)\to V^{mat}$.

(d) The trivial foliation on $\C\times\C^*\times V^{mat}$ with
leaves $\C\times\C^*\times \{(s,B)\}$ induces foliations on
$\C\times\C^*\times V^{mat}/\langle m_{[1]}\rangle$
and on $M_{3TJ}^{mon}$.

A family of $P_{3D6}$-TEJPA bundles over a (complex analytic) base manifold $T$
is isomonodromic if and only if the induced holomorphic map
$T\to M_{3TJ}^{mon}$ takes values in one leaf.

(e) A leaf in $M_{3TJ}^{mon}$ 
contains only finitely many branches over $\C^*\times\C^*$ if and only
if $\Mon_0^{mat}$ has finite order.
\end{theorem}

{\bf Proof:}
(a) The three partial derivatives of $g:=b_1^2+b_2^2+sb_1b_2$,
$$\frac{\paa g}{\paa b_1}=2b_1+sb_2,\ \frac{\paa g}{\paa b_2}=2b_2+sb_1,\ 
\frac{\paa g}{\paa s}=-b_1b_2,
$$
do not vanish simultaneously on $V^{mon}$.

(b) The map $\Phi^{path}$ is injective: the data
$(\beta,u^1_0,s,B(\beta))$ contain a $P_{3D6}$ numerical tuple of the
$P_{3D6}$ bundle and determine $P,A$ and $J$ by \eqref{7.10},
and they contain $\beta$.

The map $\Phi^{path}$ is surjective: setting $\alpha^j_0=\alpha^j_\iiii=0$,
the data $(\beta,u^1_0,s,B)$ satisfy condition \eqref{2.26} of a
$P_{3D6}$ numerical tuple by lemma \ref{t5.2} (a), so they define
a $P_{3D6}$-bundle. In order to see that \eqref{7.10} defines a 
$P_{3D6}$-TEJPA bundle one has to go through the constructions
of $P,A$ and $J$ in the proofs of theorem \ref{t6.3} (a) 
and theorem \ref{t7.3} (a) and (b).

(c) The first statement including \eqref{7.28} follows from \eqref{2.28}.
The quotient $\C\times\C^*\times V^{mat}$ is not an algebraic manifold,
because $\langle m_{[1]}\rangle\cong\Z$ does not have finite order.
But it is an analytic manifold. The rest is clear.

(d) A family of $P_{3D6}$-TEJPA bundles is isomonodromic if the 
underlying family of $P_{3D6}$-TEP bundles is isomonodromic and if
$P,A$ and $J$ are flat (in the parameters).
Then the (up to a global sign) distinguished $4$-tuple of bases
$\uuuu e^\pm_0,\uuuu e^\pm_\iiii$ varies flatly.
As this $4$-tuple defines $s$ and $B(\beta)$ and the underlying family
of $P_{3D6}$ bundles is isomonodromic, $s$ and $B(\beta)$ are constant.

Conversely, a family over one leaf is an isomonodromic family of 
$P_{3D6}$ bundles by the discussion in chapter \ref{s4}, and
the bases $\uuuu e^\pm_0,\uuuu e^\pm_\iiii$ used for its construction
vary flatly. By \eqref{7.10} $P,A$ and $J$ are flat.

(e) This is obvious.\hfill$\Box$

Define algebraic automorphisms $R_1,R_2$ and $R_3$ of 
$\C\times\C^*\times V^{mon}$ by
\begin{eqnarray}
R_1:(\beta,u^1_0,s,b_1,b_2)&\mapsto& (\beta,u^1_0,-s,b_1,-b_2)\nonumber \\ 
R_2:(\beta,u^1_0,s,b_1,b_2)&\mapsto& (\beta,u^1_0,s,-b_1,-b_2)\label{7.30}\\ 
\nonumber
R_3=R_1\circ R_2:(\beta,u^1_0,s,b_1,b_2)&\mapsto& (\beta,u^1_0,-s,-b_1,b_2).
\end{eqnarray}
Then $G^{mon}:=\{\id,R_1,R_2,R_3\}$ is a group isomorphic to $\Z_2\times \Z_2$.
The action on $\C\times \C^*\times V^{mat}$ commutes with
$m_{[1]}$ and respects the fibres of the projection to $\C\times\C^*$.
Denote the induced automorphisms on
$\C\times\C^*\times V^{mat}/\langle m_{[1]}\rangle$, 
$V^{mat}$, $((c^{path})^*M_{3TJ})^{mon}$, $M_{3TJ}^{mon}$ and on the
fibres $M_{3TJ}^{mon}(u^1_0,u^1_\iiii)$ of $pr^u$ also by
$R_1,R_2$ and $R_3$.

\begin{theorem}\label{t7.6}
(a) Consider the action of $G^{mon}$ on $V^{mon}$.
The orbits of the four points $(0,\pm 1,0)$ and $(0,0,\pm 1)$ have length 2,
all other orbits have length 4. The quotient is an affine algebraic
variety with two $A_1$-singularities at the orbits of $(0,\pm 1,0)$
and $(0,0,\pm 1)$ and smooth elsewhere.
With $(y_1,y_2,y_3)=(s^2,b_1^2,b_2^2)$ it is
\begin{eqnarray}\label{7.31}
V^{mon}/G^{mon}\cong \{(y_1,y_2,y_3)\in\C^3\, |\, 
y_1y_2y_3\!-\!(y_2\!+\!y_3\!-\!1)^2\!=\!0\}.
\end{eqnarray}

(b) The set of isomorphism classes of $P_{3D6}$-TEP bundles $M_{3T}$
is $M_{3TJ}/G^{mon}$. The natural bijections
\begin{eqnarray*}
\Phi^{path}_{3T}:(c^{path})^*M_{3T}&\to & 
\C\times \C^*\times V^{mon}/G^{mon},\\
\Phi_{3T}:M_{3T}&\to & 
\C\times \C^*\times V^{mon}/G^{mon}\times \langle m_{[1]}\rangle
\end{eqnarray*}
induces on $(c^{path})^* M_{3T}$ the structure of an algebraic variety,
denoted by $((c^{path})^* M_{3T})^{mon}$,
and on $M_{3T}$ the structure of an analytic variety, denoted by $M_{3T}^{mon}$,
with 
\begin{eqnarray}\label{7.32}
M_{3T}^{mon}= M_{3TJ}^{mon}/G^{mon}.
\end{eqnarray}
The two completely reducible $P_{3D6}$-TEP bundles in
$M_{3T}(u^1_0,u^1_\iiii)$ lie at the two $A_1$-singularities
of $M_{3T}^{mon}(u^1_0,u^1_\iiii)=M_{3TJ}^{mon}(u^1_0,u^1_\iiii)/G^{mon}$.
\end{theorem}

{\bf Proof:}
(a) As $G^{mon}$ is finite and $V^{mon}$ is an affine algebraic variety, 
the quotient is a geometric quotient and an affine algebraic variety.
The statement on the lengths of the orbits is easily checked.
It implies immediately that the quotient of the smooth $V^{mon}$
has two $A_1$-singularities and is smooth elsewhere.

It remains to calculate equation(s) for the quotient. With
$(y_1,y_2,y_3,y_4)=(s^2,b_1^2,b_2^2,sb_1b_2)$ one obtains
\begin{eqnarray*}
\textup{Spec}(V^{mon})&=& \C[s,b_1,b_2]/(b_1^2+b_2^2+sb_1b_2-1),\\
\textup{Spec}(V^{mon}/G^{mon}) &=& (\textup{Spec}(V^{mon}))^{G^{mon}}\\
&=& (\textup{image of }\C[y_1,y_2,y_3,y_4]\textup{ in Spec}(V^{mon}))\\
&=& \C[y_1,y_2,y_3,y_4]/(y_2+y_3+y_4-1,y_1y_2y_3-y_4^2)\\
&=& \C[y_1,y_2,y_3]/(y_1y_2y_3-(y_2+y_3-1)^2).
\end{eqnarray*}

(b) Consider a $P_{3D6}$-TEJPA bundle $(H,\nnn,u^1_0,u^1_\iiii,P,A,J)$
and, for some $\beta$, its image $(\beta,u^1_0,s,b_1,b_2)$ in 
$\C\times \C^*\times V^{mon}$ under $\Phi^{path}$.
By theorem \ref{t7.3} (b) there are altogether four enrichments of the
underlying $P_{3D6}$-TEP bundle to $P_{3D6}$-TEJPA bundles, they differ
only by signs $\varepsilon_1,\varepsilon_2\in\{\pm 1\}$ and
are listed in \eqref{7.9}.
Their distinguished $4$-tuples of bases from theorem \ref{7.3} (c)
are listed in \eqref{7.13}. From this one reads off
that their images under $\Phi^{path}$ in $\C\times\C^*\times V^{mon}$
are $(\beta,u^1_0,\varepsilon_1 s,\varepsilon_2 b_1,
\varepsilon_1\varepsilon_2 b_2)$.
This shows that $M_{3T}=M_{3TJ}/G^{mon}$.

The four enrichments of a $P_{3D6}$-TEP bundle with irreducible 
$P_{3D6}$-bundle are non-isomorphic because
$\textup{Aut}(H,\nnn,P)=\{\pm \id\}=\textup{Aut}(H,\nnn,P,A,J)$.

In the completely reducible cases the four enrichments form pairs 
of isomorphic enrichments by theorem \ref{t7.3} (e).
Therefore in each space $M_{3T}^{mon}(u^1_0,u^1_\iiii)$ the 
completely reducible $P_{3D6}$-TEP bundles lie at the
two $A_1$-singularities.
\hfill$\Box$

\chapter[Normal forms of $P_{3D6}$-TEJPA bundles and moduli spaces]
{Normal forms of $P_{3D6}$-TEJPA bundles and their 
moduli spaces}\label{s8}
\setcounter{equation}{0}

\noindent
Any $P_{3D6}$-TEJPA bundle is pure or a $(1,-1)$-twistor
(see remark \ref{t4.1} (iv) for these notions). 
This follows from theorem \ref{t4.2} and theorem \ref{t6.3} (a).
In this chapter we shall give an elementary independent proof.
But our main purpose is to give normal forms.
They allow us to classify $P_{3D6}$-TEJPA bundles in a new way.
This induces a new structure as algebraic manifold $M_{3TJ}^{ini}$
on the set $M_{3TJ}$ of isomorphism classes of 
$P_{3D6}$-TEJPA bundles and a stratification into an 
open submanifold of pure bundles and two codimension 1
submanifolds of $(1,-1)$-twistors, though the two codimension 1 submanifolds intersect each subspace 
$M_{3TJ}^{ini}(u^1_0,u^1_\iiii)$ in four codimension 1
submanifolds. $M_{3TJ}^{ini}(u^1_0,u^1_\iiii)$ has four natural charts
isomorphic to $\C^2$, each containing all pure bundles and one of the four
families of $(1,-1)$-twistors.

Theorem \ref{t8.2} will consider $M_{3TJ}(u^1_0,u^1_\iiii)$ for a
fixed pair $(u^1_0,u^1_\iiii)$. Theorem \ref{t8.4} will
consider $M_{3TJ}$. The proofs use repeatedly some observations
of a general nature. They are collected in the following remarks.

\begin{remarks}\label{t8.1}
(i) The proof of the existence and uniqueness of the normal forms in theorem
\ref{t8.2} is not so difficult. They are largely determined by the
properties of the values at $0$ and $\iiii$ of the sections with 
which they are defined. And these properties can easily be rewritten
using the correspondence in \eqref{2.19}, in terms of properties of flat
generating sections of the bundles $L^\pm_0,L^\pm_\iiii$.
This is made precise in (ii).

(ii) Consider a $P_{3D6}$-TEJPA bundle.
Let $\uuuu e_0^\pm$ and $\uuuu e_\iiii^\pm$ be the (up to a global
sign, unique) $4$-tuple of bases from theorem \ref{t7.3} (c).
By the correspondence in \eqref{2.19} they correspond to bases
$\uuuu v_0$ of $H_0$ and $\uuuu v_\iiii$ of $H_\iiii$. These bases satisfy
\begin{equation}\label{8.1}
[z\nnn_\zdz]\uuuu v_0=\uuuu v_0\MGcdot u^1_0
\begin{pmatrix}1&0\\0&\!-\!1\end{pmatrix},
[\!-\!\nnn_{\paa_z}]\uuuu v_\iiii=\uuuu v_\iiii\MGcdot u^1_\iiii
\begin{pmatrix}1&0\\0&\!-\!1\end{pmatrix},
\end{equation}
\begin{equation}\label{8.2}
P(\uuuu v_0^t,\uuuu v_0)=P(\uuuu v_\iiii^t,\uuuu v_\iiii)=
{\bf 1}_2,
\end{equation}
\begin{equation}\label{8.3}
A(\uuuu v_0)=\uuuu v_0\MGcdot 
\begin{pmatrix}0&\!-\!1\\1&0\end{pmatrix},\quad
A(\uuuu v_\iiii)=\uuuu v_\iiii\MGcdot 
\begin{pmatrix}0&\!-\!1\\1&0\end{pmatrix},
\end{equation}
\begin{equation}\label{8.4}
J(\uuuu v_0)=\uuuu v_\iiii\MGcdot 
\begin{pmatrix}1&0\\0&\!-\!1\end{pmatrix},\quad 
J(\uuuu v_\iiii)=\uuuu v_0\MGcdot 
\begin{pmatrix}1&0\\0&\!-\!1\end{pmatrix}.
\end{equation}
If any basis $\uuuu{\www v}_0$ of $H_0$ satisfies the same
equations \eqref{8.1} to \eqref{8.3} as $\uuuu v_0$, then
it coincides up to a global sign with $\uuuu v_0$.
An analogous statement holds for a basis $\uuuu{\www v}_\iiii$ of $H_\iiii$.

(iii) Consider the matrix
\begin{eqnarray}\label{8.5}
C:=\begin{pmatrix}1&1\\-i&i\end{pmatrix},\quad\textup{hence}\quad
C^{-1}=\frac{1}{2}\begin{pmatrix}1&i\\1&-i\end{pmatrix}.
\end{eqnarray}
It satisfies
\begin{eqnarray}\label{8.6}
\begin{split}
C^{-1}\begin{pmatrix}1&0\\0&-1\end{pmatrix}C
&=&\begin{pmatrix}0&1\\1&0\end{pmatrix},\\
C^t \MGcdot {\bf 1}_2 \MGcdot C
&=&\begin{pmatrix}0&2\\2&0\end{pmatrix},\\ 
C^{-1}\begin{pmatrix}0&-1\\1&0\end{pmatrix}C
&=&\begin{pmatrix}i&0\\0&-i\end{pmatrix}.
\end{split}
\end{eqnarray}
Therefore the bases $\uuuu w_0:=\uuuu v_0C$ and $\uuuu w_\iiii:=\uuuu v_\iiii C$
satisfy
\begin{equation}\label{8.7}
[z\nnn_\zdz]\uuuu w_0=\uuuu w_0\MGcdot u^1_0
\begin{pmatrix}0&1\\1&0\end{pmatrix},
[-\nnn_{\paa_z}]\uuuu w_\iiii=\uuuu w_\iiii\MGcdot u^1_\iiii
\begin{pmatrix}0&1\\1&0\end{pmatrix},
\end{equation}
\begin{equation}\label{8.8}
P(\uuuu w_0^t,\uuuu w_0)=P(\uuuu w_\iiii^t,\uuuu w_\iiii)=
\begin{pmatrix}0&2\\2&0\end{pmatrix},
\end{equation}
\begin{equation}\label{8.9}
A(\uuuu w_0)=\uuuu w_0\MGcdot 
\begin{pmatrix}i&0\\0&-i\end{pmatrix},\quad
A(\uuuu w_\iiii)=\uuuu w_\iiii\MGcdot 
\begin{pmatrix}i&0\\0&-i\end{pmatrix},
\end{equation}
\begin{equation}\label{8.10}
J(\uuuu w_0)=\uuuu w_\iiii\MGcdot 
\begin{pmatrix}0&1\\1&0\end{pmatrix},\quad 
J(\uuuu w_\iiii)=\uuuu w_0\MGcdot 
\begin{pmatrix}0&1\\1&0\end{pmatrix}.
\end{equation}
If any basis $\uuuu{\www w}_0$ of $H_0$ satisfies the same
equations \eqref{8.7} to \eqref{8.9} as $\uuuu w_0$, then
it coincides up to a global sign with $\uuuu w_0$.
An analogous statement holds for a basis $\uuuu{\www w}_\iiii$ of $H_\iiii$.

(iv) Let $(H,\nnn,u^1_0,u^1_\iiii,P,A,J)$ be a $P_{3D6}$-TEJPA bundle.
$A$ and $J$ act on the space $\Gamma(\P^1,\OO(H))$ of global sections
of $H$. Since $A^2=-\id$ and $J^2=\id$, they act semisimply
with eigenvalues in $\{\pm i\}$ and $\{\pm 1\}$.
Since $AJ=-JA$, they exchange eigenspaces.
Therefore the four eigenspaces all have the same dimension, which
is half the dimension of $\Gamma(\P^1,\OO(H)$. Thus this dimension
is even, and $H$ is a $(k,-k)$-twistor with $k=0$ or $k>0$ and odd
($\deg H=0$ by \eqref{3.1} and \eqref{6.7}).
Theorem \ref{t8.2} (a) will show $k\in\{0,1\}$.

If $\sigma\in\OO(H)|_{\C^*}$ and 
$A(\sigma(z))=\varepsilon_1\MGcdot i \MGcdot\sigma(-z)$ with some 
$\varepsilon_1\in \{\pm 1\}$, then for $k\in\Z$
\begin{eqnarray}\label{8.11}
A(z^k\sigma(z))=(-1)^k\MGcdot\varepsilon_1\MGcdot i\MGcdot(-z)^k\sigma(-z).
\end{eqnarray}

(v) Denote by $\OO_{\P^1}(a,b)$ for $a,b\in\Z$ the sheaf of holomorphic
functions on $\P^1$ with poles of order $\leq a$ at $0$ and poles
of order $\leq b$ at $\iiii$.

Suppose that a $P_{3D6}$-TEJPA bundle 
$(H,\nnn,u^1_0,u^1_\iiii,P,A,J)$ is a $(k,-k)$-twistor for some $k>0$.
One can choose sections $\sigma_1\in \Gamma(\P^1,\OO(H))$ and
$\sigma_2\in \Gamma(\C,\OO(H))$ such that
$$\OO(H)= \OO_{\P^1}(0,k)\MGcdot\sigma_1\oplus \OO_{\P^1}(0,-k)\MGcdot \sigma_2.$$
Then $(\sigma_1(0),\sigma_2(0))$ is a basis of $H_0$,
and for $a,b\in\Z$
\begin{eqnarray*}
\Gamma(\P^1,\OO_{\P^1}(a,b)\MGcdot\OO(H))
&=& \bigl(\C z^{-a}\sigma_1\oplus \C z^{-a+1}\sigma_1\oplus \dots
\oplus \C z^{k+b}\sigma_1\bigr)\\ \nonumber
&\oplus& \bigl(\C z^{-a}\sigma_2\oplus \C z^{-a+1}\sigma_2\dots
\oplus \C z^{-k+b}\sigma_2\bigr)
\end{eqnarray*}
(the second sum is $0$ if $-a>-k+b$). The section $\sigma_1$ is unique
up to a scalar. Therefore it is an eigenvector of $A$, and 
\[
A(\sigma_1(z))=\varepsilon_1\MGcdot i\MGcdot \sigma_1(-z)
\quad\textup{for some unique }\varepsilon_1\in\{\pm 1\}.
\]
The section $\sigma_2$ is not unique. But it can also be chosen as an eigenvector 
of $A$, with
\[
A(\sigma_2(z))=-\varepsilon_1\MGcdot i \MGcdot \sigma_2(-z).
\]
The reason for this is that $\Gamma(\P^1,\OO_{\P^1}(0,k)\MGcdot \OO(H))$
is the sum of the two $k+1$-dimensional eigenspaces of $A$ with 
eigenvalues $\pm\varepsilon_1\MGcdot i$,
\begin{eqnarray*}
&&\C\sigma_1\oplus \C z^2\sigma_1\oplus \dots\oplus \C z^{2k}\sigma_1\\
\textup{and}&&
\C z\sigma_1\oplus\C z^3\sigma_1\oplus \dots\oplus \C z^{2k-1}\sigma_1
\oplus\www\sigma_2
\end{eqnarray*}
for some section $\www \sigma_2$. One can choose $\sigma_2=\www\sigma_2$.

(vi) Suppose that $(H,\nnn,u^1_0,u^1_\iiii,P,A,J)$ is a $P_{3D6}$-TEJPA bundle,
and $\uuuu\sigma$ is basis of $\OO(H)_0$ with
$$A(\uuuu\sigma(z))=\uuuu\sigma(-z)\MGcdot
\begin{pmatrix}i&0\\0&-i\end{pmatrix}.$$
Then \eqref{8.11} shows that 
the flatness $\nnn_\zdz \circ A=A\circ\nnn_\zdz$ of $A$ is equivalent to
\begin{eqnarray}\label{8.12}
\nnn_\zdz \uuuu\sigma(z) &=& \uuuu\sigma(z)\MGcdot\frac{1}{z}
\begin{pmatrix}a(z)&b(z)\\ c(z)&d(z)\end{pmatrix}\\
\textup{with}&&a(z),d(z)\in z\C\{z^2\},b(z),c(z)\in\C\{z^2\}.\nonumber
\end{eqnarray}
Because of $P(A\,a,A\,b)=P(a,b)$ and $A^2=-\id$, the eigenspaces of $A$ are
isotropic with respect to $P$. This isotropy, the symmetry of $P$,
$P(A\,a,A\,b)=P(a,b)$ and $A(\uuuu\sigma(z))=\uuuu\sigma(-z)\MGcdot
\begin{pmatrix}i&0\\0&-i\end{pmatrix}$
show that
\begin{eqnarray}\label{8.13}
P(\uuuu\sigma(z)^t,\uuuu\sigma(-z))=e(z)\begin{pmatrix}0&1\\1&0\end{pmatrix}
\textup{ with }e(z)\in\C\{z^2\}.
\end{eqnarray}
The flatness of $P$ is equivalent to
\begin{eqnarray}\label{8.14}
\begin{split}
\zdz e(z)\begin{pmatrix}0&1\\1&0\end{pmatrix}
&=P(\nnn_\zdz\uuuu\sigma(z)^t,\uuuu\sigma(-z))
+P(\uuuu\sigma(z)^t,\nnn_\zdz\uuuu\sigma(-z)) \\ 
&=\frac{1}{z}\begin{pmatrix}a(z)&c(z)\\b(z)&d(z)\end{pmatrix}
e(z)\begin{pmatrix}0&1\\1&0\end{pmatrix}\\ 
&\quad\quad+\frac{1}{-z}e(z)\begin{pmatrix}0&1\\1&0\end{pmatrix}
\begin{pmatrix}a(-z)&b(-z)\\c(-z)&d(-z)\end{pmatrix}\\ 
&= \frac{1}{z}e(z)(a(z)+d(z))\begin{pmatrix}0&1\\1&0\end{pmatrix}.
\end{split}
\end{eqnarray}
Therefore $a(z)+d(z)=0$ if and only if $\zdz e(z)=0$.
\end{remarks}

\begin{theorem}\label{t8.2}
(a) Any $P_{3D6}$-TEJPA bundle is pure or a $(1,-1)$-twistor.

(b) Any pure $P_{3D6}$-TEJPA bundle has four normal forms 
(any of which can be chosen), and these are listed
in \eqref{8.15} to \eqref{8.18}. They are indexed by 
$(\varepsilon_1,\varepsilon_2)\in\{\pm 1\}^2$ or by $k\in\{0,1,2,3\}$,
the correspondence between $(\varepsilon_1,\varepsilon_2)$ and $k$ being given by
\begin{eqnarray*}
&&\begin{array}{c|c|c|c|c}(\varepsilon_1,
\varepsilon_2)&(1,1)&(-1,1)&(1,-1)&(-1,-1)\\
\hline k&0&1&2&3\end{array}
\end{eqnarray*}
In each case, the basis $\uuuu\sigma_k$ of $\Gamma(\P^1,\OO(H))$ is unique, 
and $(f_k,g_k)\in\C^*\times \C$ below is also unique.
\begin{eqnarray}\label{8.15}
&&\nnn_\zdz \uuuu\sigma_k
\!=\! \uuuu\sigma_k\MGcdot\left[
\frac{u^1_0}{z}\begin{pmatrix}0&1\\1&0\end{pmatrix}
\!-\!  g_k\!\begin{pmatrix}1&0\\0&-1\end{pmatrix}
\!-\! z\MGcdot u^1_\iiii\!\begin{pmatrix}0&f^{-2}_k\\f^2_k&0\end{pmatrix}\right], 
\hspace*{1cm}\\ \label{8.16}
&&P((\uuuu\sigma_k)^t(z),\uuuu\sigma_k(-z))
=\begin{pmatrix}0&2\\2&0\end{pmatrix},\\ \label{8.17} 
&&A(\uuuu\sigma_k(z))=\uuuu\sigma_k(-z)\MGcdot \varepsilon_1\MGcdot
\begin{pmatrix}i&0\\0&-i\end{pmatrix},\\ \label{8.18}
&&J(\uuuu\sigma_k(z))=\uuuu\sigma_k(\rho_c(z))\MGcdot \varepsilon_2\MGcdot
\begin{pmatrix}0&f^{-1}_k\\f_k&0\end{pmatrix} .
\end{eqnarray}
The bases are related by
\begin{eqnarray}\label{8.19}
\uuuu\sigma_0=\uuuu\sigma_2,\ \uuuu\sigma_1=\uuuu\sigma_3=\uuuu\sigma_0
\MGcdot\begin{pmatrix}0&1\\1&0\end{pmatrix}.
\end{eqnarray}

(c) Formulae \eqref{8.15} to \eqref{8.18} define for any $k\in\{0,1,2,3\}$ and 
for any $(f_k,g_k)\in\C^*\times\C$
a pure $P_{3D6}$-TEJPA bundle. Therefore the set 
of pure $P_{3D6}$-TEJPA bundles has four natural charts with
coordinates $(u^1_0,u^1_\iiii,f_k,g_k)$, and each chart is isomorphic
to $\C^*\times\C^*\times\C^*\times \C$. With this structure as an 
affine algebraic manifold the set is called $M_{3TJ}^{reg}$.
Obviously the restriction of $pr^u$ to $M_{3TJ}^{reg}$,
$$pr^u:M_{3TJ}^{reg}\to\C^*\times \C^*, P_{3D6}\textup{-TEJPA bundle}\mapsto
(u^1_0,u^1_\iiii),$$
is an algebraic morphism. The fibres are called $M_{3TJ}^{reg}(u^1_0,u^1_\iiii)$.
The coordinates of the four charts are related by 
\begin{eqnarray}
(f_1,g_1)&=&(f_0^{-1},-g_0),\nonumber \\ \label{8.20}
(f_2,g_2)&=&(-f_0,g_0),\\ \nonumber (f_3,g_3)&=&(-f_0^{-1},-g_0).
\end{eqnarray}

(d) Fix a pair $(u^1_0,u^1_\iiii)\in\C^*\times\C^*$ and a square root
$\sqrt{c}$ of $c={u^1_\iiii}/{u^1_0}$.
Each $P_{3D6}$-TEJPA bundle which is a $(1,-1)$-twistor has a unique
normal form given by \eqref{8.21} to \eqref{8.24}.
Here, $k\in\{0,1,2,3\}$, $(\varepsilon_1,\varepsilon_2)\in\{\pm 1\}^2$
($k$ and $(\varepsilon_1,\varepsilon_2)$ are related as in (b)),
$\www g_k\in\C$ and the basis $\uuuu\psi$ are unique.
$\uuuu\psi$ is a basis of $H|_\C$ such that $(z\psi_1,z^{-1}\psi_2)$
is a basis of $H|_{\P^1-\{0\}}$.
The pair $[\varepsilon_1,\varepsilon_2\sqrt{c}]$ is called the type
of the $P_{3D6}$-TEJPA bundle. The proof will give a
different characterization of this type.
\begin{eqnarray}\label{8.21}
&&\nnn_\zdz \uuuu\psi
=\uuuu\psi\MGcdot\left[
\frac{u^1_0}{z}\begin{pmatrix}0&1\\1&0\end{pmatrix}
-\frac{1}{2}\begin{pmatrix}1&0\\0&-1\end{pmatrix}\right. \\ \nonumber
&&\hspace*{2cm}
-\www g_k\MGcdot \sqrt{c}\MGcdot  z\begin{pmatrix}0&1\\0&0\end{pmatrix}
+ \left. u^1_\iiii\MGcdot c\MGcdot z^3\begin{pmatrix}0&1\\0&0\end{pmatrix}\right] \\
\label{8.22}
&&P((\uuuu\psi)^t(z),\uuuu\psi(-z))
=\begin{pmatrix}0&2\\2&0\end{pmatrix},\\ \label{8.23} 
&&A(\uuuu\psi(z))=\uuuu\psi(-z)\MGcdot \varepsilon_1\MGcdot
\begin{pmatrix}i&0\\0&-i\end{pmatrix},\\ \label{8.24}
&&J(\uuuu\psi(z))=\uuuu\psi(\rho_c(z))\MGcdot \varepsilon_2\MGcdot
\begin{pmatrix}-\frac{1}{\sqrt{c}}z^{-1}&0\\0&\sqrt{c}z\end{pmatrix} .
\end{eqnarray}

(e) Fix a pair $(u^1_0,u^1_\iiii)\in\C^*\times\C^*$ and a square root
$\sqrt{c}$ of $c={u^1_\iiii}/{u^1_0}$.
Formulae \eqref{8.21} to \eqref{8.24} define for any $k\in\{0,1,2,3\}$ and 
the corresponding $(\varepsilon_1,\varepsilon_2)\in\{\pm 1\}^2$ and 
for any $\www g_k\in\C$ a $P_{3D6}$-TEJPA bundle which is a $(1,-1)$-twistor.

Therefore the set of $P_{3D6}$-TEJPA bundles in $M_{3TJ}(u^1_0,u^1_\iiii)$ 
which are $(1,-1)$-twistors has four components, one for each type,
and each component has a coordinate $\www g_k$ and is, with this coordinate,
isomorphic to $\C$. With this structure as an 
affine algebraic manifold the set is called 
$M_{3TJ}^{[\varepsilon_1,\varepsilon_2\sqrt{c}]}(u^1_0,u^1_\iiii)$.
The union of the four sets is called $M_{3TJ}^{sing}(u^1_0,u^1_\iiii)$.

(f) Fix a pair $(u^1_0,u^1_\iiii)\in\C^*\times\C^*$, a square root
$\sqrt{c}$ of $c={u^1_\iiii}/{u^1_0}$, a $k\in\{0,1,2,3\}$
and the corresponding pair $(\varepsilon_1,\varepsilon_2)\in\{\pm 1\}^2$.
Each $P_{3D6}$-TEJPA bundle in $M_{3TJ}^{reg}(u^1_0,u^1_\iiii)
\cup M_{3TJ}^{[\varepsilon_1,\varepsilon_2\sqrt{c}]}(u^1_0,u^1_\iiii)$
has a unique normal form listed in \eqref{8.25} to \eqref{8.28}.

Here, $(f_k,\www g_k)\in\C\times\C$ and the basis $\uuuu\varphi_k$ 
are unique.
$\uuuu\varphi_k$ is a basis of $H|_\C$ such that 
$(z\varphi_{k,1}-\frac{f_k}{\sqrt{c}}\varphi_{k,2},z^{-1}\varphi_{k,2})$
is a basis of $H|_{\P^1-\{0\}}$.

\begin{eqnarray}\label{8.25}
\begin{split}
&\nnn_\zdz \uuuu\varphi_k
=\uuuu\varphi_k\MGcdot\left[
\frac{u^1_0}{z}\begin{pmatrix}0&1\\1&0\end{pmatrix}
-\frac{1}{2}(1+\www g_k f_k)\begin{pmatrix}1&0\\0&-1\end{pmatrix}
\right. \\ 
&\hspace*{1cm}
-\www g_k\MGcdot\sqrt{c}\MGcdot z\begin{pmatrix}0&1\\0&0\end{pmatrix}
-u^1_\iiii\MGcdot f_k^2\MGcdot z\begin{pmatrix}0&0\\1&0\end{pmatrix}
\\ 
&\hspace*{1cm} 
+u^1_\iiii\MGcdot\sqrt{c}\MGcdot f_k\MGcdot z^2\begin{pmatrix}1&0\\0&-1\end{pmatrix}
+\left. u^1_\iiii\MGcdot c\MGcdot z^3\begin{pmatrix}0&1\\0&0\end{pmatrix}\right]
\end{split}
\end{eqnarray}
\begin{eqnarray}\label{8.26}
&&P((\uuuu\varphi_k)^t(z),\uuuu\varphi_k(-z))
=\begin{pmatrix}0&2\\2&0\end{pmatrix},\\ \label{8.27} 
&&A(\uuuu\varphi_k(z))=\uuuu\varphi_k(-z)\MGcdot \varepsilon_1\MGcdot
\begin{pmatrix}i&0\\0&-i\end{pmatrix},\\ \label{8.28}
&&J(\uuuu\varphi_k(z))=\uuuu\varphi_k(\rho_c(z))\MGcdot \varepsilon_2\MGcdot
\begin{pmatrix}-\frac{1}{\sqrt{c}}z^{-1}&0\\f_k&\sqrt{c}z\end{pmatrix} .
\end{eqnarray}

(g) In the situation of (f), $\C^2$ with coordinates $(f_k,\www g_k)$
is a chart for $M_{3TJ}^{reg}(u^1_0,u^1_\iiii)
\cup M_{3TJ}^{[\varepsilon_1,\varepsilon_2\sqrt{c}]}(u^1_0,u^1_\iiii)$.
Here the points $(0,\www g_k)$ correspond to the points in 
$M_{3TJ}^{[\varepsilon_1,\varepsilon_2\sqrt{c}]}(u^1_0,u^1_\iiii)$
with coordinate $\www g_k$. The points $(f_k,\www g_k)\in\C^*\times\C$
correspond to the points
\begin{eqnarray}\label{8.29}
(f_k,g_k)
=(f_k,-u^1_0\frac{\sqrt{c}}{f_k}+\frac{1}{2}+\frac{f_k}{2}\www g_k)
\end{eqnarray}
in $M_{3TJ}^{reg}(u^1_0,u^1_\iiii)$ with coordinates $(f_k,g_k)$ from (b)
and (c).
For $(0,\www g_k)$, the bases in (f) and (d) are related by
$\uuuu\psi=\uuuu\varphi_k$. This means that the normal form in (d) is the restriction
to $f_k=0$ of the normal form in (f).
For $(f_k,\www g_k)$ the bases in (f) and (b) are related by 
\begin{eqnarray}\label{8.30}
\uuuu\varphi_k =\uuuu\sigma_k\MGcdot
\begin{pmatrix}1&\frac{\sqrt{c}}{f_k}\MGcdot z\\0&1\end{pmatrix}
\end{eqnarray}

Therefore $M_{3TJ}(u^1_0,u^1_\iiii)$ is covered by four affine charts
isomorphic to $\C^2$ and with coordinates $(f_k, \www g_k)$.
Each chart contains $M_{3TJ}^{reg}(u^1_0,u^1_\iiii)$ and one of the 
four components $M_{3TJ}^{[\varepsilon_1,\varepsilon_2\sqrt{c}]}(u^1_0,u^1_\iiii)$.
With the induced structure as an algebraic
manifold it is called $M_{3TJ}^{ini}(u^1_0,u^1_\iiii)$.
The coordinate changes between the charts can be read off from
\eqref{8.29} and \eqref{8.20}.
\end{theorem}

{\bf Proof:}
(a) Let $(H,\nnn,u^1_0,u^1_\iiii,P,A,J)$ be a $P_{3D6}$-TEJPA bundle
which is a $(k,-k)$-twistor for some $k>0$.
We shall show $k=1$.

{\bf 1st proof:}
If $(H,\nnn)$ is irreducible, this follows immediately from theorem \ref{t4.2}.
If $(H,\nnn)$ is completely reducible, 
then the last statement in lemma \ref{3.1} (a) shows
that $H$ is a $(0,0)$-twistor, so this case is impossible as $k>0$.

{\bf 2nd proof:} We give a 2nd proof.
First, because the following proof is elementary and independent
of theorem \ref{t4.2}. Second, because the proof is a preparation
of the proof of (d).

Choose sections $\sigma_1$ and $\sigma_2$ as in remark \ref{t8.1} (v),
that means with
\begin{eqnarray}\label{8.31}
\OO(H)=\OO_{\P^1}(0,k)\MGcdot \sigma_1\oplus\OO_{\P^1}(0,-k)\MGcdot \sigma_2,\\
\label{8.32}
A(\uuuu \sigma(z))=\uuuu\sigma(-z)\MGcdot
\varepsilon_1\begin{pmatrix}i&0\\0&-i\end{pmatrix}
\textup{ for some }\varepsilon_1\in\{\pm 1\}.
\end{eqnarray}
The section $\sigma_1$ is unique up to rescaling.
The section $\sigma_2$ is unique up to adding an element of 
$\oplus_{l=0}^{k-1}\C z^{2l+1}\sigma_1$ and rescaling.
The sections $\sigma_1$ and $\sigma_2$ can be chosen such that
\begin{eqnarray}\label{8.33}
\begin{split}
\uuuu\sigma(0)&=\uuuu w_0\quad\textup{ if }\varepsilon_1=1,\\
\uuuu\sigma(0)&=\uuuu w_0\MGcdot\begin{pmatrix}0&1\\1&0\end{pmatrix}
\quad \textup{ if }\varepsilon_1=-1
\end{split}
\end{eqnarray}
($\uuuu w_0$ is defined in remark \ref{t8.1} (iii)).
This determines $\sigma_1$ uniquely, but $\sigma_2$ only up to 
adding an element as above.
$\uuuu\sigma$ is a basis of $H|_\C$, and \eqref{8.7} and \eqref{8.12} give
\begin{eqnarray}\label{8.34}
\begin{split}
\nnn_\zdz \uuuu\sigma(z) &=
\uuuu\sigma\left[\frac{u^1_0}{z}\begin{pmatrix}0&1\\1&0\end{pmatrix}
+\begin{pmatrix}a&b \\ \www c&d\end{pmatrix}\right]\\
\textup{for some }&a,b\in\C\{z^2\},\www c,d\in z\MGcdot \C\{z^2\}.
\end{split}
\end{eqnarray}
$(z^k\sigma_1,z^{-k}\sigma_2)$ is a basis of $H|_{\P^1-\{0\}}$ with
\begin{eqnarray}\label{8.35}
\begin{split}
\nnn_\zdz(z^k\sigma_1,z^{-k}\sigma_2)&=(z^k\sigma_1,z^{-k}\sigma_2)
\MGcdot\left[\frac{u^1_0}{z}\begin{pmatrix}0&z^{-2k}\\z^{2k}&0\end{pmatrix}
\right. \\
&+\left. \begin{pmatrix}a&bz^{-2k}\\ \www c z^{2k}&d\end{pmatrix}
+\begin{pmatrix}k&0\\0&-k\end{pmatrix}\right] . 
\end{split}
\end{eqnarray}
But $(H,\nnn)$ has a pole of order 2 at $\iiii$. Therefore
$$u^1_0z^{2k-1}+\www c(z)z^{2k}\in\C\oplus \C\MGcdot z.$$
This is only possible if $k=1$ (and then $\www c(z)=0$).
This finishes the 2nd proof of (a).

(b) It is obvious that the normal forms for $k=1,2,3$ are obtained 
from the normal form for $k=0$ by the base change in \eqref{8.19}
and the change of the scalars in \eqref{8.20}.
Therefore it is sufficient to prove existence and uniqueness of the
normal form for $k=0$.

First we assume existence and prove uniqueness. Formulae
\eqref{8.15} to \eqref{8.17} and the last statement of remark \ref{t8.1} (iii)
show $\uuuu\sigma(0)=\pm\uuuu w_0$. As $H$ is pure,
the basis $\uuuu\sigma$ is uniquely determined by $\uuuu\sigma(0)=\pm\uuuu w_0$.
Then also $(f_0,g_0)$ is unique.

Now we show existence. Choose that basis $\uuuu\sigma$ of 
$\Gamma(\P^1,\OO(H))$ which satisfies $\uuuu\sigma(0)=\uuuu w_0$.
As $P$ is holomorphic on $\OO(H)\times j^*\OO(H)$, it takes constant
values on sections of $\Gamma(\P^1,\OO(H))$, and $A$ maps global
sections to global sections. Therefore \eqref{8.8} and \eqref{8.9}
imply \eqref{8.16} and \eqref{8.17} with $\varepsilon_1=1$.

By remark \ref{t8.1} (iv), $J$ exchanges the eigenspaces of $A$
in $\Gamma(\P^1,\OO(H))$. Therefore, and because of $J^2=\id$, 
a scalar $f_0\in\C^*$ exists such that \eqref{8.18} 
holds with $\varepsilon_2=1$.
Comparison with \eqref{8.10} shows
\[
\uuuu\sigma(\iiii)=\uuuu w_\iiii\MGcdot
\begin{pmatrix}f_0&0\\0&f_0^{-1}\end{pmatrix}.
\]
Comparison with \eqref{8.7} shows
\[
[-\nnn_{\paa_z}]\uuuu\sigma(\iiii)
=\uuuu\sigma(\iiii)\MGcdot u^1_\iiii\MGcdot
\begin{pmatrix}0&f_0^{-2}\\f_0^2&0\end{pmatrix}.
\]
This equation together with \eqref{8.7} for $\uuuu\sigma(0)=\uuuu w_0$ 
and \eqref{8.12} shows
\begin{eqnarray*}
\nnn_\zdz \uuuu\sigma_0
=\uuuu\sigma_0\MGcdot\left[
\frac{u^1_0}{z}\begin{pmatrix}0&1\\1&0\end{pmatrix}
-\begin{pmatrix}a&0\\0&d\end{pmatrix}
-z\MGcdot u^1_\iiii\begin{pmatrix}0&f^{-2}_0\\f^2_0&0\end{pmatrix}\right]
\end{eqnarray*}
for some $a,d\in\C$. Now \eqref{8.16} and \eqref{8.14}
show $a+d=0$. One sets $g_0:=a$.
Existence of the normal form for $k=0$ is proved.

(c) Everything is clear as soon as the first statement is proved.
As above, it is sufficient to prove the first statement for the normal
form with $k=0$. So we have to show that \eqref{8.15} to \eqref{8.18}
define in the case $k=0$ for any $(f_0,g_0)\in\C^*\times \C$ a
$P_{3D6}$-TEJPA bundle.

The conditions 
\begin{align*}
A^2&=-\id,\ J^2=\id,\ AJ=-JA,\ 
\\
&P(A\,a,A\,b)=P(a,b),\ P(J\,a,J\,b)=P(a,b)
\end{align*}
are obviously satisfied, and $P$ is symmetric and nondegenerate.
It remains to show that $P,A$ and $J$ are flat.

$A$ is flat because \eqref{8.12} holds with
$a(z)=-d(z)=-g_0\MGcdot z$, $b(z)=u^1_0-z^2\MGcdot u^1_\iiii f^{-2}_0$,
$c(z)=u^1_0-z^2\MGcdot u^1_\iiii f^2_0$.
And $P$ is flat because \eqref{8.14} holds with $e(z)=2$ and
$a(z)+d(z)=0$. Finally, the flatness of $J$ is equivalent to
\[
\nnn_\zdz J(\uuuu\sigma(z)) = 
J(\nnn_\zdz\uuuu\sigma(z)),
\]
which holds because of the following calculations.
The first one uses $\zdz=-\rho_c(z)\paa_{\rho_c(z)}.$
\begin{eqnarray*}
&&\nnn_\zdz J(\uuuu\sigma(z))\\
&=& \nnn_\zdz \sigma(\rho_c(z))\MGcdot 
\begin{pmatrix}0&f_0^{-1}\\f_0&0\end{pmatrix}
= -\nnn_{\rho_c(z)}\paa_{\rho_c(z)} \sigma(\rho_c(z))\MGcdot 
\begin{pmatrix}0&f_0^{-1}\\f_0&0\end{pmatrix}\\
&=& -\uuuu\sigma(\rho_c(z)) 
\left[\frac{u^1_0}{\rho_c(z)}\begin{pmatrix}0&1\\1&0\end{pmatrix}
-g_0\begin{pmatrix}1&0\\0&-1\end{pmatrix}\right. \\
&&\hspace*{3cm}\left. -\rho_c(z)\MGcdot u^1_\iiii\begin{pmatrix}0&f^{-2}_0\\f^2_0&0\end{pmatrix}\right]
\begin{pmatrix}0&f_0^{-1}\\f_0&0\end{pmatrix}\\
&=& \uuuu\sigma(\rho_c(z)) 
\left[-z\MGcdot u^1_\iiii\begin{pmatrix}f_0&0\\0&f_0^{-1}\end{pmatrix}
-g_0\begin{pmatrix}0&-f_0^{-1}\\f_0&0\end{pmatrix}
+\frac{u^1_0}{z}\begin{pmatrix}f_0^{-1}&0\\0&f_0\end{pmatrix}\right]
\end{eqnarray*}
and 
\begin{eqnarray*}
&&J(\nnn_\zdz \uuuu\sigma(z))\\
&=& \sigma(\rho_c(z)) 
\begin{pmatrix}0&f_0^{-1}\\f_0&0\end{pmatrix}
\left[\frac{u^1_0}{z}\begin{pmatrix}0&1\\1&0\end{pmatrix}
-g_0\begin{pmatrix}1&0\\0&-1\end{pmatrix}
-z\MGcdot u^1_\iiii\begin{pmatrix}0&f^{-2}_0\\f^2_0&0\end{pmatrix}\right]\\
&=& \uuuu\sigma(\rho_c(z)) 
\left[\frac{u^1_0}{z}\begin{pmatrix}f_0^{-1}&0\\0&f_0\end{pmatrix}
-g_0\begin{pmatrix}0&-f_0^{-1}\\f_0&0\end{pmatrix}
-z\MGcdot u^1_\iiii\begin{pmatrix}f_0&0\\0&f^{-1}_0\end{pmatrix}\right] .
\end{eqnarray*}

(d) First, existence of the normal form will be shown, then uniqueness.

Choose $\sigma_1$ and $\sigma_2$ as in the 2nd proof of part (a).
Then \eqref{8.31} to \eqref{8.35} hold with $k=1$.
First we shall refine the choice of $\sigma_2$ and then call the
new basis $\uuuu\psi$.

$J$ acts on 
$$\Gamma(\P^1,\OO(H))=\C \sigma_1\oplus\C z\sigma_1,$$ 
it maps
$$\Gamma(\P^1,\OO_{\P^1}(0,1)\MGcdot\OO(H))=
\C\sigma_1\oplus\C z\sigma_1\oplus \C z^2\sigma_1\oplus \C \sigma_2$$
to 
$$\Gamma(\P^1,\OO_{\P^1}(1,0)\MGcdot\OO(H))=
\C z^{-1}\sigma_1\oplus \C\sigma_1\oplus\C z\sigma_1\oplus 
\C z^{-1}\sigma_2,$$
and it exchanges in any case the eigenspaces of $A$. Therefore
\begin{eqnarray*}
J(\uuuu\sigma(z)) &=& \uuuu\sigma(\rho_c(z))\MGcdot
\begin{pmatrix}\gamma_1 z^{-1}&\gamma_2\\0&\gamma_3 z\end{pmatrix}\\
\textup{with}&& \gamma_1,\gamma_3\in\C^*,\gamma_2\in\C,
\gamma_1\MGcdot\gamma_3=-1.
\end{eqnarray*}
$J^2=\id$ shows $\gamma_1^2c=1={\gamma_3^2}/{c}$.
We had chosen a fixed root $\sqrt{c}$. 
Therefore there is a unique $\varepsilon_2\in\{\pm 1\}$ with
\[
(\gamma_1,\gamma_3)=\varepsilon_2\MGcdot(-\frac{1}{\sqrt{c}},\sqrt{c}).
\]
Now we choose the new basis of $H|_\C$ 
\[
\uuuu\psi:=\uuuu\sigma\MGcdot
\begin{pmatrix}1&  \tfrac12 \varepsilon_2\MGcdot {\sqrt{c}\gamma_2}z\\ 
0 &1\end{pmatrix}.
\]
It satisfies \eqref{8.24}.

As $(z\sigma_1)(0)=0$ it also satisfies $\uuuu\psi(0)=\uuuu w_0$.
Comparison with \eqref{8.10} shows
\begin{eqnarray*}
(z\psi_1,z^{-1}\psi_2)(\iiii)\varepsilon_2
\begin{pmatrix}-\sqrt{c}&0\\0&\frac{1}{\sqrt{c}}\end{pmatrix}
=J(\uuuu\psi(0))=J(\uuuu w_0)=\uuuu w_\iiii
\begin{pmatrix}0&1\\1&0\end{pmatrix},\\
\textup{so}\quad 
(z\psi_1,z^{-1}\psi_2)(\iiii)=
\uuuu w_\iiii\MGcdot \varepsilon_2
\begin{pmatrix}0&\sqrt{c}\\ -\frac{1}{\sqrt{c}}&0\end{pmatrix}.
\end{eqnarray*}
Comparison with \eqref{8.7} shows
\begin{eqnarray}\label{8.36}
[-\nnn_{\paa_z}](z\psi_1,z^{-1}\psi_2)(\iiii)
= (z\psi_1,z^{-1}\psi_2)(\iiii)
\begin{pmatrix}0&-u^1_\iiii\MGcdot c\\ -u^1_0&0\end{pmatrix}.
\end{eqnarray}

The equations \eqref{8.31} to \eqref{8.35} hold with $k=1$ also for 
the new basis $\uuuu\psi$ instead of the basis $\uuuu\sigma$. 
As $(H,\nnn)$ has a pole of order 2 at $\iiii$, \eqref{8.35} shows
$$\www c=0,\ a,d\in\C,\ b=b_1z+b_3z^3\textup{ with }b_1,b_3\in\C.$$
Comparison of \eqref{8.35} and \eqref{8.36} gives
$b_3=u^1_\iiii\MGcdot c.$
This gives for the moment the following approximation of \eqref{8.21},
\begin{eqnarray}\label{8.37}
&&\nnn_\zdz \uuuu\psi
=\uuuu\psi\MGcdot\left[
\frac{u^1_0}{z}\begin{pmatrix}0&1\\1&0\end{pmatrix}
+\begin{pmatrix}a&0\\0&d\end{pmatrix}\right. \\ \nonumber
&&\hspace*{2cm}+b_1\MGcdot z\begin{pmatrix}0&1\\0&0\end{pmatrix}
+ \left. u^1_\iiii\MGcdot c\MGcdot z^3\begin{pmatrix}0&1\\0&0\end{pmatrix}\right] .
\end{eqnarray}

Because of \eqref{8.23} for $\uuuu\psi$ and because of remark \ref{t8.1} (vi),
$\psi_1$ and $\psi_2$ are isotropic with respect to $P$.
Because of
$$P(\psi_1(z),\psi_2(-z))=-P(z\psi_1(z),(-z)^{-1}\psi_2(-z))$$
the function $P(\psi_1(z),\psi_2(-z))$ is also holomorphic at $\iiii$,
thus constant. Together with $\uuuu\psi(0)=\uuuu w_0$ and \eqref{8.8}
this establishes \eqref{8.22}. Now \eqref{8.14} shows $a+d=0$.

It remains to show $a=-\frac{1}{2}.$
We use the flatness of $J$ for this. The matrices $M_1$ and $M_2$ in
\[
J(\uuuu\psi(z))\MGcdot M_1=\nnn_\zdz J(\uuuu\psi(z))
=J(\nnn_\zdz\uuuu\psi(z))=J(\uuuu\psi(z))\MGcdot M_2
\]
coincide. The constant part of $M_1$ can be read off from \eqref{8.37}
and \eqref{8.24} and is, using $\nnn_\zdz=-\nnn_{\rho_c(z)}\paa_{\rho_c(z)}$,
$\bsp -a-1&0\\0&-d+1\esp$.
The constant part of $M_2$ can be read off \eqref{8.37} and \eqref{8.24}
and is $\bsp a&0\\0&d\esp$.
This shows $a=-\frac{1}{2}$.
Existence of the normal form is proved.

It remains to prove the uniqueness of the normal form.
\eqref{8.31} and \eqref{8.32} show that $\varepsilon_1$ is unique.
\eqref{8.21} to \eqref{8.23} and remark \ref{t8.1} (iii) show
\begin{eqnarray*}
\uuuu\psi(0)&=& \uuuu w_0\quad \textup{ if }\varepsilon_1=1,\\
\uuuu\psi(0)&=& \uuuu w_0\MGcdot \begin{pmatrix}0&1\\1&0\end{pmatrix}
\quad\textup{ if }\varepsilon_1=-1.
\end{eqnarray*}
This and \eqref{8.23} and the condition that $(z\psi_1,z^{-1}\psi_2)$ 
is a basis of $H|_{\P^1-\{0\}}$ determines $\psi_1$ uniquely
and $\psi_2$ up to addition of an element of $\C z\psi_1$.
Now \eqref{8.24} and the argument in the proof of the existence
determine $\varepsilon_2$ and show that also $\psi_2$ is unique.

(e) As in the proof of (c), everything is clear as soon as the
first statement is shown, i.e., that \eqref{8.21} to \eqref{8.24}
define a $P_{3D6}$-TEJPA bundle for any $k,(\varepsilon_1,\varepsilon_2)$
(compatible as in (b)) and any $\www g_k\in\C$.

One can check this in a similar way to the proof of (c).
Only the flatness of $J$ requires a lengthy calculation.
Alternatively one can refer to part (f). There the formulae in (d) turn
up as the special case $f_k=0$. 

(f) For the pure $P_{3D6}$-TEJPA bundles in (b), one can 
calculate in a straightforward way their data $\nnn,P,A$ and $J$
with respect to the new basis $\uuuu\varphi_k$ (on $H|_\C$)
and the new scalars $(f_k,\www g_k)$, starting from the normal
forms in (b) with respect to $\uuuu\sigma_k$ and $(f_k,g_k)$,
and using the formulae \eqref{8.30} and \eqref{8.29}.
We omit the calculation. It leads to \eqref{8.25} -- \eqref{8.28}.

This calculation and the uniqueness of the normal forms in (b) 
show the existence and uniqueness of the normal forms in (f) for pure 
$P_{3D6}$-TEJPA bundles. 
Part (d) and the fact that (d) and (f) for $f_k=0$ coincide,
show the existence and uniqueness of the normal form in (f) for 
non-pure $P_{3D6}$-TEJPA bundles.

Then (e) becomes obvious. As the formulae
\eqref{8.25} to \eqref{8.28} give $P_{3D6}$-TEJPA bundles for 
all $(f_k,\www g_k)\in\C^*\times \C$, they give $P_{3D6}$-TEJPA bundles
also for all values $(0,\www g_k)\in\{0\}\times\C$.

(g) This is also clear now.
\hfill$\Box$

\begin{remark}\label{t8.3}
The group $G^{mon}=\{\id,R_1,R_2,R_3\}$ was introduced before theorem \ref{t7.6}.
It acts on $M^{mon}_{3TJ}$, and the orbits are the sets of (four or two) 
$P_{3D6}$-TEJPA bundles whose underlying $P_{3D6}$-TEP bundles are
isomorphic. Theorem \ref{t8.2} (b)+(c) shows that it acts also on the 
algebraic manifold $M^{reg}_{3TJ}$, via
\begin{eqnarray}
R_1:(u^1_0,u^1_\iiii,f_0,g_0)&\mapsto&(u^1_0,u^1_\iiii,f_0^{-1},-g_0)\\ 
\nonumber
R_2:(u^1_0,u^1_\iiii,f_0,g_0)&\mapsto&(u^1_0,u^1_\iiii,-f_0,g_0)\\ \label{8.38}
R_3:(u^1_0,u^1_\iiii,f_0,g_0)&\mapsto&(u^1_0,u^1_\iiii,-f_0^{-1},-g_0).
\nonumber
\end{eqnarray}
In other words, for any $P_{3D6}$-TEJPA bundle $T\in M^{reg}_{3TJ}$,
\begin{eqnarray}\label{8.39}
(f_k,g_k)(T)=(f_0,g_0)(R_k(T)),\quad (f_k,g_k)(R_k(T))=(f_0,g_0)(T).
\end{eqnarray}
The action of $G^{mon}$ extends to an action on $M^{ini}_{3TJ}(u^1_0,u^1_\iiii)$.
For $T\in M^{reg}_{3TJ}(u^1_0,u^1_\iiii)\cup M^{[+,+\sqrt{c}]}(u^1_0,u^1_\iiii)$ 
the image $R_k(T)$
is in $M^{reg}_{3TJ}(u^1_0,u^1_\iiii)\cup 
M^{[\varepsilon_1,\varepsilon_2\sqrt{c}]}(u^1_0,u^1_\iiii)$
with the $(\varepsilon_1,\varepsilon_2)$ which corresponds to $k$
as in theorem \ref{8.2} (b), and
\begin{eqnarray}\label{8.40}
(f_k,\www g_k)(R_k(T))=(f_0,\www g_0)(T).
\end{eqnarray}
\end{remark}

In theorem \ref{t8.2} (d) -- (g), only $P_{3D6}$-TEJPA bundles
with fixed $(u^1_0,u^1_\iiii)$ were considered, because of the 
appearance of $\sqrt{c}$. Now we shall consider the 
full space $M_{3TJ}$, and the pull-back by a $2{:}1$-covering, which
will resolve the ambiguity of the square root $\sqrt{c}$.
Consider the covering
\begin{eqnarray}\label{8.41}
c^{2{:}1}:\C^*\times\C^*&\to& \C^*\times \C^*\\ \nonumber
(x,y)&\mapsto& ({x}/{y},xy)=(u^1_0,u^1_\iiii).
\end{eqnarray}
The covering $c^{path}$ from \eqref{7.24} factorizes through it via
\begin{eqnarray}\label{8.42}
c^{path}:
\C\times\C^*\to \C^*\times \C^*\stackrel{c^{2{:}1}}{\to} \C^*\times \C^*,\\
(\beta, u^1_0)\mapsto(\tfrac12{e^{-\beta/2}},\tfrac12{e^{-\beta/2}}/{u^1_0}),
(x,y)\mapsto ({x}/{y},xy).\nonumber
\end{eqnarray}

The following theorem is a direct consequence of theorem \ref{t8.2}.

\begin{theorem}\label{t8.4}
(a) The set $(c^{2{:}1})^*M_{3TJ}$ carries a natural structure as an algebraic
manifold, with which it is called $((c^{2{:}1})^*M_{3TJ})^{ini}$.
It is obtained by glueing the four affine charts in theorem \ref{t8.2} (g)
for the spaces $M^{ini}_{3TJ}(u^1_0,u^1_\iiii)$ to four affine charts for
$(c^{2{:}1})^*M_{3TJ}$. 
Here $y=\sqrt{c}$ with $y$ as in \eqref{8.41}.
Each chart is isomorphic to $\C^*\times\C^*\times \C\times\C$
with coordinates $(x,y,f_k,\www g_k)$.
The set of non-pure $P_{3D6}$-TEJPA bundles consists of the four 
hyperplanes, one in each chart, with $f_k=0$, $k=0,1,2,3$.
Here they can be denoted by $((c^{2{:}1})^*M_{3TJ})^{[k]}$.

(b) The set $M_{3TJ}$ inherits from $((c^{2{:}1})^*M_{3TJ})^{ini}$
the structure of an algebraic manifold, with which it is called
$M^{ini}_{3TJ}$. 
Then $((c^{2{:}1})^*M_{3TJ})^{ini}=(c^{2{:}1})^*M_{3TJ}^{ini}$.
But a description of $M_{3TJ}^{ini}$ in charts is not obvious.

The set of non-pure $P_{3D6}$-TEJPA bundles consists only of 
two hypersurfaces. They can be denoted $M^{[\varepsilon_1]}_{3TJ}$.
Each of them intersects each fibre $M^{ini}_{3TJ}(u^1_0,u^1_\iiii)$
in the two hypersurfaces 
$M^{[\varepsilon_1,\varepsilon_2\sqrt{c}]}_{3TJ}(u^1_0,u^1_\iiii)$.
\end{theorem}

Theorem \ref{t4.2} makes a stronger statement than theorem \ref{t8.2} (a).
It gives a solution of the inverse monodromy problem for trace free
$P_{3D6}$ bundles.
We state the special case of $P_{3D6}$-TEJPA bundles 
(or $P_{3D6}$-TEP bundles) in the following lemma again,
and we provide an elementary proof, independent of the proofs in  
\cite{Heu09} and \cite{Ni09}.

\begin{lemma}\label{t8.5}
No isomonodromic family in $M_{3TJ}$ is contained in 
$M^{sing}_{3TJ}=\cup_{\varepsilon_1=\pm 1}M^{[\varepsilon_1]}_{3TJ}$.
The intersection of any isomonodromic family with $M^{sing}_{3TJ}$ 
is empty or a hypersurface in the isomonodromic family.
\end{lemma}

{\bf Proof:}
The first statement implies the second, because $M^{mon}_{3TJ}$
and $M^{ini}_{3TJ}$ give the same complex analytic manifold, and
the isomonodromic families as well as $M^{sing}_{3TJ}$ are analytic
subvarieties of this complex manifold.

We assume that an isomonodromic family is contained in 
$M^{[k]}_{3TJ}\subset (c^{2{:}1})^*M^{ini}_{3TJ}$ for some $k\in\{0,1,2,3\}$.
Then the normal form in theorem \ref{t8.2} (d) for 
$u^1_0=u^1_\iiii=:x\in U$ with $U\subset\C^*$ simply connected
must give for some holomorphic function $\www g_k\in\OO_U$
an isomonodromic family with basis of sections $\uuuu\psi(z,x)$, $x\in U$.
Because of \eqref{8.21} to \eqref{8.23}, the bases
$\uuuu e^\pm_0(z,x)$ which correspond by \eqref{2.19}
to $\uuuu\psi(0)\MGcdot C^{-1}$, are flat. The basis 
$\uuuu\sigma:= \uuuu\psi\MGcdot C^{-1}$
can be written as in \eqref{2.15}, 
\begin{eqnarray*}
\uuuu\psi(z,x)&=&\uuuu e^\pm_0(z,x)\MGcdot
\begin{pmatrix}e^{-x/z}&0\\0&e^{x/z}\end{pmatrix}A_0^\pm(z,x)C\\
\textup{with}&& \whhh A_0^+=\whhh A_0^-\textup{ and }
\whhh{A}_0^\pm(0,x)={\bf 1}_2\textup{ for all }x.
\end{eqnarray*}
This formula and the flatness of $\uuuu e^\pm_0$ give
\begin{eqnarray*}
\nnn_{x\partial_x}\uuuu\sigma(z,x) 
=\uuuu\sigma(z,x)\left[-\frac{x}{z}\begin{pmatrix}0&1\\1&0\end{pmatrix}
+\begin{pmatrix}a&b\\ \www c&d\end{pmatrix}\right]
\end{eqnarray*}
with $a,b,\www c,d\in\OO_U[z]$. Now one proceeds as in the proof of 
theorem \ref{t8.2} (d). One uses \eqref{8.21} to \eqref{8.23} 
and arguments as in remark \ref{t8.1} (vi)
and derives that $\www c=0$, $a=-d\in\OO_U$, $b=b_1z-xz^3$  with $b_1\in\OO_U$ and
\begin{eqnarray*}
\nnn_{x\paa_x}\uuuu\sigma(z,x)=
\uuuu\sigma(z,x)\left[-\frac{x}{z}\begin{pmatrix}0&1\\1&0\end{pmatrix}
+a\begin{pmatrix}1&0\\0&-1\end{pmatrix}
+\begin{pmatrix}0&b_1z-xz^3\\0&0\end{pmatrix}\right] .
\end{eqnarray*}
Now one can calculate
\[
(\nnn_{x\paa_x}\nnn_\zdz-\nnn_\zdz\nnn_{x\paa_x})\uuuu\sigma(z,x).
\]
The calculation shows that it cannot be $0$, although it must be $0$.
Because of this contradiction, the assumption above is wrong.
\hfill$\Box$.

\chapter{Generalities on the Painlev\'e equations}\label{s9}
\setcounter{equation}{0}

\noindent
The Painlev\'e equations are six families $P_I$, $P_{II}$, $P_{III}$, 
$P_{IV}$, $P_V$, $P_{VI}$ of second order differential equations
in  
\[
U=\C, \ \  \C^\ast=\C-\{0\}, \text{ or } \C-\{0,1\}
\]
of the form
\[
f_{xx}=R(x,f,f_x)
\]
where $R$ is holomorphic for $x\in U$ and rational in $f$ and $f_x$.
The following table lists the number of essential parameters in $R$
and the subset $U$.
\[
\renewcommand{\arraystretch}{1.3}
\begin{array}{c|c|c|c|c|c|c}
 & P_{VI} & \ P_V\  & \ P_{IV}\  & \ P_{III}\  & \ \ P_{II}\ \  & \ \ P_I\ \  
 \\ 
 \hline
\text{parameters} & 4 & 3 & 2 & 2 & 1 & 0 
\\ 
\hline
U & \C\!-\!\{0,\!1\} & \C^* & \C & \C^* & \C & \C 
\end{array}
\]
The equations are distinguished by the following {\em Painlev\'e property}:
Any local holomorphic solution extends to a global multi-valued  meromorphic solution in $U$ (single-valued if $U=\C$).
This means that solutions branch only at 
$0,1$ for $P_{VI}$, $0$ for $P_V$ and $P_{III}$, and on $U$ the only singularities are poles.
The positions of the poles in $U$ 
depend on the solution, so they are
{\em movable}. In other words, the Painlev\'e equations are distinguished
by the property that the only movable singularities of solutions are poles.
For proofs of the Painlev\'e property which do not use 
isomonodromic families of connections, see \cite{GLSh02} and references
there. In theorem \ref{10.3} we shall rewrite in our language the proof in 
\cite{FN80} for $P_{III}(0,0,4,-4)$ which uses isomonodromic families.

Between 1895 and 1910, Painlev\'e and Gambier discovered 50 such
families of differential equations, which included $P_I$ to $P_{VI}$.
But the general solutions of the other 44 equations could be reduced to
those of $P_I$ to $P_{VI}$, to rational or elliptic functions, or to solutions
of linear or first order differential equations.

In this paper we are interested in a special case of the Painlev\'e III equations.
A priori, there are four parameters $(\alpha,\beta,\gamma,\delta)\in\C^4$:
the equation $P_{III}(\alpha,\beta,\gamma,\delta)$ is 
\begin{eqnarray}\label{9.1}
f_{xx}&=& \frac{f_x^2}{f}-\frac{1}{x}f_x+\frac{1}{x}(\alpha f^2+\beta)
+\gamma f^3+\delta\frac{1}{f}.
\end{eqnarray}
Let us choose locally a logarithm $\varphi=2\log f$ of $f$, i.e.\  $f=e^{\varphi/2}$. Then $\varphi$ branches at the poles and zeros of $f$. 
\eqref{9.1} becomes 
\begin{eqnarray}\label{9.2}
(\xdx)^2\varphi &=& 2x(\alpha e^{\varphi/2}+\beta e^{-\varphi/2})+
2x^2(\gamma e^\varphi + \delta e^{-\varphi})
\end{eqnarray}
which is more symmetric than \eqref{9.1}. 

One sees that
\begin{eqnarray}
\label{9.3}
\!\!\!\!\!\!\!\!\!
\begin{split}
&f\textup{ and }\varphi+4\pi i k\textup{ are solutions for }
(\alpha,\beta,\gamma,\delta) \iff \\
-&f\textup{ and }\varphi+2\pi i+4\pi i k\textup{ are solutions for }
(-\alpha,-\beta,\gamma,\delta) \iff \\
&f^{-1}\textup{ and }-\!\varphi+4\pi i k\textup{ are solutions for }
(-\beta,-\alpha,-\delta,-\gamma) \iff 
\\
-&f^{-1}\textup{ and }-\!\varphi+2\pi i+4\pi i k\textup{ are solutions for }
(\beta,\alpha,-\delta,-\gamma) . 
\end{split}
\end{eqnarray}
Rescaling $x$ by $r\in\C^*$ and $f$ by $s\in\C^*$ shows
\begin{align}
f&\textup{ and }\varphi \textup{ are solutions for }
(\alpha,\beta,\gamma,\delta)\iff  \label{9.4}
\\
&s\MGcdot f(r\MGcdot x)\textup{ and }\varphi(r\MGcdot x)\!+\!2\log \!s 
\textup{ are solutions for }(\tfrac rs\alpha,rs\beta,\tfrac{r^2}{s^{2}}\gamma,r^2s^2\delta).\nonumber
\end{align}
In \cite{OKSK06} the equations \eqref{9.1} are split into four cases:
\[
\begin{array}{lcc}
P_{III}(D_6)    &\quad&\gamma\delta\neq 0
\\
P_{III}(D_7)    & &\gamma=0,\alpha\delta\neq 0 \text{ or } \delta=0,\beta\gamma\neq 0
\\
P_{III}(D_8)  &  &  \gamma=0,\delta=0,\alpha\beta\neq 0
\\
P_{III}(Q)  & &  \alpha=0,\gamma=0  \text{ or } \beta=0,\delta=0
\end{array}
\]
Using the symmetries \eqref{9.3} and the rescalings \eqref{9.4},
these four cases can be reduced to the following:
\[
\begin{array}{lcc}
P_{III}(D_6)    &\quad&(\alpha,\beta,4,-4)
\\
P_{III}(D_7)    & &(2,\beta,0,-4)
\\
P_{III}(D_8)  &  &  (4,-4,0,0)
\\
P_{III}(Q)  & &  (0,1,0,\delta),(0,0,0,1),(0,0,0,0)
\end{array}
\]

All cases for $P_{III}(Q)$ are solvable by quadratures, so they are usually
ignored.
The cases $P_{III}(D_7)$ and $P_{III}(D_8)$ are not treated in 
\cite{Ok79}, \cite{Ok86}, \cite{FN80}, \cite{JM81}, but in
\cite{OKSK06}, \cite{OO06}, \cite{PS09}.
Within the cases $P_{III}(D_6)(\alpha,\beta,4,-4)$, the case
$P_{III}(0,0,4,-4)$ lies at the center of the parameter space, as it is the fixed
point of the group of the four symmetries in \eqref{9.3}. This group
acts on the space of its (global multi-valued meromorphic) solutions.
In this paper we are interested in this case. Then \eqref{9.2}
is the radial sinh-Gordon equation
\begin{eqnarray}\label{9.5}
(\xdx)^2\varphi &=& 16 x^2\sinh\varphi.
\end{eqnarray}

Chapter \ref{s10} will connect 
two objects associated to $P_{III}(0,0,4,-4)$: 
the {\em spaces of initial conditions} of Okamoto \cite{Ok79}
and certain {\em isomonodromic families of meromorphic connections}.
The interrelations between these objects have been studied thoroughly
in the case of $P_{VI}$ by M.-H. Saito and others. But for $P_I$
to $P_V$ only first steps have been done in \cite{PS09}.
We shall treat the case $P_{III}(0,0,4,-4)$ in chapter \ref{s10}.
It will also be the basis for the later chapters.

The rational function $R(x,f,f_x)$ in $f$ and $f_x$ in the Painlev\'e equations
contains only the inverses indicated in the following table:
\[
\renewcommand{\arraystretch}{1.3}
\begin{array}{c|c|c|c|c|c|c}
 & P_{VI} & \ P_V\  & \ P_{IV}\  & \ P_{III}\  &  P_{II}  &  P_I 
\\ 
\hline
 & \!f^ {\!-1},(f\!-\!1)^{\!-1},(f\!-\!x)^{\!-1}\! & \!f^{\!-1},(f\!-\!1)^{\!-1}\! & f^{\!-1} & f^{\!-1} &
 & 
\\ 
\hline
U^\prime & \C-\{0,1,x_0\} & \C-\{0,1\} & \C^* & \C^* & \C & \C 
\end{array}
\]
Therefore at any point $x_0\in U$, any pair $(f(x_0),f_x(x_0))\in U^\prime\times \C$
with $U^\prime$ as shown, 
determines a unique local holomorphic
solution $f$, and, by the Painlev\'e property, a unique global multi-valued
meromorphic solution. Thus $U^\prime\times\C$ is a naive space of initial
conditions at $x_0$ for the solutions of the Painlev\'e equation.
But it misses the initial conditions for the solutions which are singular at $x_0$.
Okamoto \cite{Ok79} gave complete spaces of initial conditions for 
all $x_0\in U$ for most of the Painlev\'e equations. He missed some cases,
including $P_{III}(D_7)$ and $P_{III}(D_8)$. These two cases were treated
in \cite{OKSK06}.

In chapter \ref{s10} we shall describe the spaces of initial conditions for
the case $P_{III}(0,0,4,-4)$ by four charts such that each extends the naive
space $U^\prime\times \C=\C^*\times \C$ by initial data for singular solutions.
This is not what Okamoto did. He started with the compactification 
$\P^2$ of $U^\prime\times \C$, blew it up several times in appropriate
ways and threw out some superfluous hypersurfaces \cite{Ok79}.
Descriptions by charts are provided in 
\cite{ShT97}, \cite{MMT99}, \cite{Mat97}, \cite{NTY02},  \cite{Te07}  but they all build on \cite{Ok79}
and do not interpret the charts as extensions of $U^\prime\times\C$ by
initial data for singular solutions.

The relations between the Painlev\'e equations and isomonodromic families
of meromorphic connections have a long history.
Soon after the discovery of the Painlev\'e equations, Fuchs and Garnier
found second order linear differential equations with rational coefficients
whose isomonodromic families are governed by equations of type
$P_{VI}$ (Fuchs) or $P_I$ to $P_V$ (Garnier).

Their work was continued in \cite{Ok86} and (almost) completed in
\cite{OO06}. \cite{Ok86} gives six types of second order linear differential
equations, for generic members of the six Painlev\'e equations.
\cite{OO06} adds four types. The ten types together cover all Painlev\'e 
equations which cannot be solved by quadratures.
They can be classified by the orders of the poles
of second order linear differential equations.
The following table is essentially taken from \cite{PS09} and is equivalent
to data in \cite{OO06}.
\begin{eqnarray*}
\renewcommand{\arraystretch}{1.3}
\begin{array}{l|l|l}
 & \textup{orders of poles} & \text{ parameters}\\ \hline
P_{VI} & 1\ \ 1\ \ 1\ \ 1 & 4 \\
P_V & 1\ \ 1\ \ 2 & 3 \\
P_{V,deg} & 1\ \ 1\ \ \frac{3}{2} & 2 \\
P_{III}(D_6) & 2\ \ 2 & 2 \\
P_{III}(D_7) & \frac{3}{2}\ \ 2 & 1 \\
P_{III}(D_8) & \frac{3}{2}\ \ \frac{3}{2} & 0 \\
P_{IV} & 1\ \ 3 & 2 \\
P_{II} (\sim P_{34}) & 1\ \ \frac{5}{2} & 1\\
P_{II} & 4 & 1 \\
P_I & \frac{7}{2} & 0
\end{array}
\end{eqnarray*}

Remarkably, there are two types for $P_{II}$, namely type $P_{II}[1,\frac{5}{2}]$
and type $P_{II}[4]$.
In \cite{Ok86} the six types $P_{VI}$, $P_{V}$, $P_{III}(D_6)$, $P_{IV}$, 
$P_{II}[4]$, $P_I$ of the above ten are given.

Second order linear differential equations can be rewritten as first
order linear systems in $2\times 2$ matrices, and these are equivalent
to trivial holomorphic vector bundles of rank 2 on $\P^1$ with 
meromorphic connections. Such data have been associated to the Painlev\'e
equations in \cite{FN80}, \cite{JM81}, \cite{IN86}, \cite{FIKN06}, 
\cite{PS09}, \cite{PT14}. 
The vector bundle point of view is used in \cite{PS09} and \cite{PT14}.

\cite{JM81} gives six types of such data. We expect that they are equivalent
to the six types in \cite{Ok86}. \cite[ch.\  5]{FIKN06} takes up 
the six types in \cite{JM81}.
\cite{FN80} considers only two types, a type equivalent to
the type $P_{II}[1,\frac{5}{2}]$
and a new type for $P_{III}(0,0,4,-4)$. The type in \cite{FN80} 
equivalent to the type $P_{II}[1,\frac{5}{2}]$ is the two-fold branched
cover (branched at $0$ and $\infty$ in $\P^1$) with poles of order 1 (at $\infty$)
and 4 (at $0$), such that the formal decomposition of Hukuhara
and Turrittin exists at the pole of order 4.
\cite{IN86} and \cite[ch.\  7-16]{FIKN06} take up both types
in \cite{FN80} and connect Stokes data and the central connection matrix
with the asymptotic behaviour near $0$ and $\infty$ in $U$ of solutions of 
$P_{II}$ and $P_{III}(0,0,4,-4)$.

\begin{remark}\label{t9.1}
The relation between the type for $P_{III}(0,0,4,-4)$ in \cite{FN80}
and the type for $P_{III}(D_6)$ in \cite{JM81} is quite remarkable
and seems to have been unnoticed.
We claim that the following points (i) and (ii) are true. 
We intend to prove this elsewhere.

(i) The isomonodromic families for $P_{III}(0,0,4,-4)$ in \cite{FN80} coincide
essentially with the isomonodromic families for $P_{III}(0,4,4,-4)$ in \cite{JM81},
but the solutions of the Painlev\'e equations are built into the 
isomonodromic familes in different ways.

(ii) One isomonodromic family gives rise to four solutions of $P_{III}(0,0,4,-4)$
(from the symmetries in \eqref{9.3}),
and to one solution of $P_{III}(0,4,4,-4)$. These solutions are connected
by the $4{:}1$ folding transformation in \cite{TOS05} and \cite{Wi04}
which is called $\psi_{III(D_6^{(1)})}^{[4]}$ in \cite{TOS05} 
and which exists only for the solutions
of $P_{III}(0,0,4,-4)$ and $P_{III}(0,4,4,-4)$.
\end{remark}

\begin{remarks}\label{t9.2}
(i) \cite{IN86}, \cite[ch.\  7 -12]{FIKN06}, \cite{Ni09} and this paper 
work with the isomonodromic families for 
$P_{III}(0,0,4,-4)$ in \cite{FN80}. Theorem \ref{t10.3} makes the 
relation between these isomonodromic families and the solutions of
$P_{III}(0,0,4,-4)$ precise.

(ii) This relation for $P_{III}(0,0,4,-4)$ is not indicated in the table
above of the above ten types, and it is not considered in \cite{OO06}, \cite{PS09}.
Therefore we regard it as a new type. It is in this sense that
\cite{OO06} \lq\lq almost\rq\rq\ completes the work of Fuchs and Garnier.

(iii) If the isomonodromic family for $P_{III}(D_6)$ in \cite{JM81}
and the isomonodromic family for $P_{III}(0,0,4,-4)$ in \cite{FN80} 
had been compared earlier, the $4{:}1$ folding transformation 
$\psi_{III(D_6^{(1)})}^{[4]}$ 
in \cite{TOS05}, \cite{Wi04} might have been found earlier.
Anyway, this comparison provides the isomonodromic interpretation
of $\psi_{III(D_6^{(1)})}^{[4]}$ asked for in \cite[ch 6 0.]{FIKN06}
(bottom of page 221).

(iv) In \cite{TOS05} two other folding transformations $\psi_{IV}^{[3]}$ and $\psi_{II}^{[2]}$
are given, and according to \cite[ch.\  6]{FIKN06} 
$\psi_{IV}^{[3]}$ is due to Okamoto (1986), and 
$\psi_{II}^{[2]}$ is due to Gambier (1910).
They map $P_{IV}(0,-\frac{2}{9})$ to $P_{IV}(1,0)$ and 
$P_{II}(0)$ to $P_{II}(-\frac{1}{2})$.
It seems interesting to study the isomonodromic families for 
$P_{IV}(0,-\frac{2}{9})$ and $P_{II}(0)$ which one obtains 
via pull-back with $\psi_{IV}^{[3]}$ and $\psi_{II}^{[2]}$ 
from the isomonodromic families in \cite{JM81} for 
$P_{IV}(1,0)$ and $P_{II}(-\frac{1}{2})$.
We believe that the relations of these isomonodromic families
with $P_{IV}(0,-\frac{2}{9})$ and $P_{II}(0)$ have not yet been considered.
They are not considered in \cite{OO06} and \cite{PS09}.
\end{remarks}

\begin{remark}\label{t9.3}
\cite{OO06} gives more than the 10 types above.
It discusses different ways to classify the Painlev\'e equations, which
amount to different numbers: 5, 6, 8, 10, 14.

6 is the original number of families of Painlev\'e equations.

5 is obtained by unifying $P_{II}$ and $P_I$ in one family of equations.

By the scaling transformations which rescale $x$ and $f$, these are 
separated into 14 types of equations. For example $P_{III}$ is separated
into $P_{III}(D_6)$, $P_{III}(D_7)$, $P_{III}(D_8)$, $P_{III}(Q)$.

4 of these 14 types can be solved by quadrature. The remaining 10 cases
are given in the table above and are related to 10 types of second order
linear differential equations.

But $P_{II}$ and $P_{34}$ are equivalent, and $P_{III}(D_6)$ and $P_{V,deg}$
are equivalent. That leads to 8 cases.
\end{remark}

\chapter{Solutions of the Painlev\'e equation $P_{III}(0,0,4,-4)$}\label{s10}
\setcounter{equation}{0}

\noindent
Now we come to the relation between $P_{3D6}$-TEJPA bundles and the
Painlev\'e equation $P_{III}(0,0,4,-4)$. 
First, recall the covering $c^{2{:}1}$ from \eqref{8.41}
and the covering $c^{path}$ from \eqref{7.24} which factorizes through
$c^{2{:}1}$:
\begin{eqnarray}\label{10.1}
&c^{path}:\C\times\C^*\to\C^*\times\C^*\stackrel{c^{2{:}1}}{\to}\C^*\times \C^*,&\\
&(\beta,u^1_0)\mapsto (\tfrac12{e^{-\beta/2}},\frac{e^{-\beta/2}}{2u^1_0}),
(x,y)\stackrel{c^{2{:}1}}\mapsto 
({x}/{y},xy)=(u^1_0,\frac{e^{-\beta}}{4 u^1_0})=(u^1_0,u^1_\iiii).&
\nonumber
\end{eqnarray}
Recall also from \eqref{4.3} that within $(c^{2{:}1})^*M_{3TJ}^{mon}$,
the families of $P_{3D6}$-TEJPA bundles where only $y$ varies
(but not $(x,s,B)$) are {\it inessential isomonodromic families}
with trivial transversal monodromy.
The next lemma combines this observation with the data in theorem \ref{t7.5}
and theorem \ref{t8.4}.

\begin{lemma}\label{t10.1}
(a) Define 
\begin{eqnarray}\label{10.2}
\begin{split}
M_{3FN}&:= ((c^{2{:}1})^*M_{3TJ})_{|y=1}, \\ 
M_{3FN}(x)&:=((c^{2{:}1})^*M_{3TJ})(x,1)=M_{3TJ}(x,x), \\
M_{3FN}^{mon}&:= ((c^{2{:}1})^*M_{3TJ}^{mon})_{|y=1}, \\ 
M_{3FN}^{mon}(x)&:=((c^{2{:}1})^*M_{3TJ}^{mon})(x,1)= M_{3TJ}^{mon}(x,x),\\
M_{3FN}^{ini}&:= ((c^{2{:}1})^*M_{3TJ}^{ini})_{|y=1}, \\ 
M_{3FN}^{ini}(x) &:=((c^{2{:}1})^*M_{3TJ}^{ini})(x,1)=M_{3TJ}^{ini}(x,x).
\end{split}
\end{eqnarray}
(Here, 3FN stands for the Flaschka-Newell version of $P_{III}$, which is exactly
$P_{III}(0,0,4,-4)$.)
Then $M_{3FN}^{mon}$ and $M_{3FN}^{ini}$ give $M_{3FN}$ the same structure as an
analytic manifold, $M_{3FN}^{ini}$ is an algebraic manifold, and the fibres
$M_{3FN}(x)$ of the projection $pr_{3FN}:M_{3FN}\to \C^*$ carry two algebraic structures
$M_{3FN}^{mon}(x)$ and $M_{3FN}^{ini}(x)$.
There are canonical isomorphisms
\begin{eqnarray}
\begin{split}
(c^{2{:}1})^* M_{3TJ}^{mon} \!\cong\! M_{3FN}^{mon}\!\times\!\C^*,\
M_{3TJ}^{mon}({x}/{y},xy)\!\cong\! M_{3FN}^{mon}(x)\ \forall\ y\in\C^*,\\
(c^{2{:}1})^* M_{3TJ}^{ini} \!\cong\! M_{3FN}^{ini}\!\times\!\C^*,\ 
M_{3TJ}^{ini}({x}/{y},xy)\!\cong\! M_{3FN}^{ini}(x)\ \forall\ y\in\C^*,\label{10.3}
\end{split}
\end{eqnarray}
where $\C^*$ is equipped with the coordinate $y$.
$M_{3FN}$ can be seen as the set of isomorphism classes of inessential
isomonodromic families (over $\C^*$ with coordinate $y$ and with trivial
transversal monodromy) of $P_{3D6}$-TEJPA bundles.

(b) $M_{3FN}^{mon}$ is the quotient of the algebraic manifold $\C\times V^{mat}$
by the action of $\langle m_{[1]}^2\rangle$. The automorphism $m_{[1]}^2$ is algebraic,
but the group $\langle m_{[1]}^2\rangle$ is isomorphic to $\Z$, therefore the quotient
is an analytic manifold, and $M_{3FN}^{mon}$ is described by the isomorphism
\begin{eqnarray}\label{10.4}
\Phi^{mon}:M_{3FN}^{mon}&\to& \C\times V^{mat}/\langle m_{[1]}^2\rangle,\\
(H,\nnn,x,x,P,A,J)&\mapsto& [(\beta,s,B(\beta))]\textup{ for }
\beta\textup{ with }\tfrac12{e^{-\beta/2}}=x.\nonumber
\end{eqnarray}
Here $s$ and $B(\beta)$ are associated to $(H,\nnn,x,x,P,A,J)$ as 
in theorem \ref{t7.5} (b).
And $m_{[1]}^2$ acts on 
\begin{eqnarray}\label{10.5}
\C\times V^{mat}= \{(\beta,s,B)\in \C\times\C\times SL(2,\C)\, |\, 
B=\begin{pmatrix}b_1&b_2\\-b_2& b_1+sb_2\end{pmatrix}\}
\end{eqnarray}
by the algebraic automorphism
\begin{eqnarray}\label{10.6}
m_{[1]}^2:(\beta,s,B)\mapsto (\beta+4\pi i,s,(\Mon_0^{mat})^{-2}\MGcdot B).
\end{eqnarray}
Any choice of $\beta$ with $\tfrac12{e^{-\beta/2}}=x$ induces an isomorphism
$M_{3FN}^{mon}(x)\to V^{mat}$.

The trivial foliation on $\C\times V^{mat}$ with leaves $\C\times\{(s,B)\}$
induces foliations on $\C\times V^{mat}/\langle m_{[1]}^2\rangle$ and on
$M_{3FN}^{mon}$.
It lifts by the first isomorphism in \eqref{10.3} to the foliation 
on $(c^{2{:}1})^*M_{3TJ}^{mon}$ from theorem \ref{t7.5} (c).
The leaves are the maximal isomonodromic families within $M_{3FN}$.

(c) $M_{3FN}^{ini}$ is given by four natural charts. The charts have coordinates
$(x,f_k,\www g_k)$ for $k=0,1,2,3$ and are isomorphic to 
$\C^*\times\C\times\C$. Each chart intersects each other chart in 
$\C^*\times\C^*\times\C$. The coordinate changes are given by \eqref{8.20}
and 
\begin{eqnarray}\label{10.7}
(x,f_k,g_k)&=& (x,f_k,-\frac{x}{f_k}+\frac{1}{2}+\frac{f_k}{2}\MGcdot\www g_k).
\end{eqnarray}
The intersection of all charts is (with the induced algebraic structure)
$M_{3FN}^{reg}$ and is in each chart $\C^*\times\C^*\times\C$.
$M_{3FN}^{reg}$ is the set of pure twistors in $M_{3FN}$.
The $k$-th chart ($k=0,1,2,3$) consists of $M_{3FN}^{reg}$ and 
a smooth hypersurface $M_{3FN}^{[k]}$ isomorphic to 
$\C^*\times\{0\}\times \C$, which is the set of $(1,-1)$-twistors in this chart.
\end{lemma}

{\bf Proof:}
(a) This follows immediately from theorem \ref{t7.5} and theorem \ref{t8.4}
and the observation above on the inessential isomonodromic families
of $P_{3D6}$-TEJPA bundles.

(b) This follows from theorem \ref{t7.5} and \eqref{10.1} and $y=1$, which give
\begin{eqnarray}\label{10.8}
u^1_\iiii=u^1_0=\tfrac12{e^{-\beta/2}}=x.
\end{eqnarray}

(c) It remains to see that the charts and coordinates in theorem \ref{t8.4}
are compatible with the isomorphisms in the second line of \eqref{10.3},
in other words a $P_{3D6}$-TEJPA bundle $T\in M_{3TJ}^{ini}(x,x)$,
and its pull-back $m_y^*T\in M_{3TJ}^{ini}(\frac{x}{y},xy)$ 
with $m_y:\C\to\C,z\mapsto yz$, have the same
coordinates, which are $(f_0,g_0)$ if $T\in M_{3TJ}^{reg}(x,x)$ and which are 
$\www g_k$ if $T\in M_{3TJ}^{[k]}(x,x)$. 
But this follows from inspection of 
\eqref{8.15}-\eqref{8.18} and \eqref{8.21}-\eqref{8.24}:
Writing $u^1_0={x}/{y},u^1_\iiii=xy,\sqrt{c}=y,$ in all these formulae
$z$ turns up only as part of the product $yz$.
\hfill $\Box$

\begin{remarks}\label{t10.2}
For any leaf in $M_{3FN}^{mon}$, the restrictions of $f_0,g_0,\www g_0$
to the leaf will turn out (see theorem \ref{t10.3} (a) and (b))
to be multi-valued meromorphic functions on the underlying
punctured plane $\C^*$ with coordinate $x$.
In order to deal with these multi-valued functions properly, and in order
to be more concrete, we fix some elementary facts and notation.

(i) A leaf in $M_{3FN}^{mon}$ with a distinguished branch means a leaf
where one branch over $\C-\R_{\leq 0}$ is distinguished.
Leaves with distinguished branches are naturally parameterized by $V^{mat}$:
There is the natural isomorphism $M_{3FN}^{mon}(1)\cong V^{mat}$
from lemma \ref{t10.1} with $\beta=0,x=1$.
Any leaf is mapped to the intersection point 
of its distinguished branch with $M_{3FN}^{mon}(1)$.

(ii) A (meromorphic) function $h$ on a leaf is a multi-valued (meromorphic) 
function in $x\in\C^*$.
Distinguishing a branch of the leaf distinguishes a branch of the 
multi-valued function. A multi-valued function on $\C^*$ with a distinguished
branch over $\C-\R_{\leq 0}$ corresponds to a single-valued holomorphic function
on $\C$ which gives for $\xi\in\C$ with $\Im(\xi)\in(-\pi,\pi)$ the
value at $e^\xi=x$ of the distinguished branch of the multi-valued function.

(iii) In this way, the function $f_0$ on $M_{3FN}^{ini}$
induces a family parameterized by $V^{mat}$ of multi-valued functions
$f_{mult}(.,s,B)$ in $x\in\C^*$ with distinguished branches
and the corresponding single-valued functions $f_{univ}(.,s,B)$ in $\xi\in\C$
with $e^\xi=x$. More concretely, $f_{univ}$ is the meromorphic function
\begin{eqnarray}\label{10.9}
&&f_{univ}:\C\times V^{mat}\to\C,\\
&&f_{univ}=f_0\circ (\Phi_{mon}^{-1}\circ pr_{m_{[1]}}),\nonumber\\
&&\textup{with }pr_{m_{[1]}}:\C\times V^{mat}\to\C\times V^{mat}/
\langle m_{[1]}^2\rangle.\nonumber
\end{eqnarray}
Similarly, $g_0$ and $\www g_0$ give $g_{mult},g_{univ}$ and 
$\www g_{mult},\www g_{univ}$.

(iv)  The coordinates $\beta$ and $\xi$ are related by
\begin{eqnarray}\label{10.10}
\tfrac12{e^{-\beta/2}}=x=e^\xi,\quad -\beta/2=\xi+\log 2.
\end{eqnarray}
With respect to the coordinates $(\xi,s,B)$ on $\C\times V^{mat}$,
the action of $m_{[1]}$ on $\C\times V^{mat}$ takes the form
\begin{eqnarray}\label{10.11}
m_{[1]}:(\xi,s,B)\to(\xi-i\pi,s,\Mon_0^{mat}(s)^{-1}\MGcdot B)
\end{eqnarray}
Because of \eqref{7.28}, $f_{univ},g_{univ}$ and $\www g_{univ}$
are invariant under this action on $\C\times V^{mat}$. 
If one leaf with a distinguished branch
has the parameter $(s,B)$, the same leaf where the next branch after
one anticlockwise turn around $0$ is distinguished, has the parameter
$(s,\Mon_0^{mat}(s)^{-2}\MGcdot B)$. 

(v) The following properties are equivalent because of (iv):
\begin{list}{}{}
\item[$(\alpha)$]
A multi-valued function $f_{mult}(.,s,B)$ has only finitely many branches.
\item[$(\beta)$]
The leaf in $M_{3FN}^{mon}$ with distinguished branch and with parameter $(s,B)$
has finitely many branches. 
\item[$(\gamma)$]
The matrix $\Mon_0^{mat}(s)$ has finite order.
Equivalently: its eigenvalue $\lambda_+(s)$ is a root of unity and $\lambda_+(s)\neq -1$.
Equivalently: the value $\alpha_+(s)$ is in $(-\frac{1}{2},\frac{1}{2})\cap\Q$.
\end{list}
In particular, $(\alpha)-(\gamma)$ depend only on $s$, not on $B$.
If they hold, then $s\in(-2,2)$ and the number of branches of the leaf
and of the function $f_{mult}(.,s,B)$ is 
$\min(l\in\Z_{>0}\, |, \lambda_+(s)^{2l}=1\}$.
\end{remarks}

The two algebraic manifolds $\C\times V^{mat}$ (with quotient
$\C\times V^{mat}/\langle m_{[1]}^2\rangle\cong M_{3FN}^{mon}$) and $M_{3FN}^{ini}$ are
related by the coordinates $(x,f_0,g_0)$ considered as holomorphic functions
in $(\beta,s,B)\in\C\times V^{mat}$,
i.e., by $x=\tfrac12{e^{-\beta/2}}$ and $f_{univ}$ and $g_{univ}$.

The following theorem shows that the dependence of $f_{mult}$ and $g_{mult}$
on $x$ is governed by $P_{III}(0,0,4,-4)$. 
It contains the Painlev\'e property for $P_{III}(0,0,4,-4)$.
It controls also the initial conditions of solutions at singular points.
The theorem is essentially due to 
\cite{FN80} and is also used in \cite{IN86}, \cite{FIKN06}, \cite{Ni09}, though the manifolds are treated here more carefully than in
these references.

\begin{theorem}\label{t10.3}
(a) $f_{mult}, g_{mult}$ and $\www g_{mult}$ restrict for
any $(s,B)\in V^{mat}$ to multi-valued meromorphic functions on $\C^*$
with coordinate $x$. The function $f_{mult}(.,s,B)$ satisfies
$P_{III}(0,0,4,-4)$, and together with $g_{mult}(.,s,B)$ it satisfies
\begin{eqnarray}\label{10.12}
2f_{mult}g_{mult} = \xdx f_{mult}.
\end{eqnarray}

(b) Let us call a local solution $f$ of $P_{III}(0,0,4,-4)$ 
{\it singular} at $x_0\in\C^*$ 
if $f$ has at $x_0$ a zero or a pole, and {\it regular} otherwise,
i.e.\  if $f(x_0)\in\C^*$.

For any leaf in $M_{3FN}^{mon}$ with distinguished branch, 
a branch (not necessarily the distinguished one) 
of the corresponding solution $f_{mult}(.,s,B)$ 
is regular at $x_0\in\C^*$ if the corresponding intersection point of the
leaf with $M_{3FN}^{mon}(x_0)$ is in $M_{3FN}^{reg}(x_0)$, and singular otherwise.
More precisely, the following holds if the intersection point of the 
leaf with $M_{3FN}^{mon}(x_0)$ is in $M_{3FN}^{[k]}(x_0)$
($k$ and $(\varepsilon_1,\varepsilon_2)$ are related as in theorem \ref{t8.2} (b)):
\begin{eqnarray}
f_{mult}^{\varepsilon_1}(x_0,s,B)&=&0,\nonumber\\
(\varepsilon_2 \paa_x f_{mult}^{\varepsilon_1})(x_0,s,B)&=&-2,\nonumber\\ 
(\varepsilon_2 \paa_x^2f_{mult}^{\varepsilon_1})(x_0,s,B)&=&\frac{-2}{x_0},
\label{10.13}\\ 
(\varepsilon_2 \paa_x^3f_{mult}^{\varepsilon_1})(x_0,s,B)&=&
\frac{2}{x_0^2}+\frac{8}{x_0}\MGcdot\www g_{k,mult}(x_0,s,B).\nonumber
\end{eqnarray}
So, $f_{mult}(.,s,B)$ has a simple zero at $x_0$ if $k=0,2$ 
($\iff \varepsilon_1=1$), 
a simple pole at $x_0$ if $k=1,3$ ($\iff \varepsilon_1=-1$), 
and the two types of simple zeros (respectively, simple poles)
are distinguished by the sign of $\paa_xf_{mult}(x_0,s,B)=-2\varepsilon_2$
(respectively, $(\paa_xf_{mult}^{-1})(x_0,s,B)=-2\varepsilon_2$).

(c) The leaves in $M_{3FN}^{mon}$ (without/with distinguished branches) 
parameterize the global multi-valued solutions of $P_{III}(0,0,4,-4)$
(without/with distinguished branches).

(d) (The Painlev\'e property) Any local solution of $P_{III}(0,0,4,-4)$
extends to a global multi-valued meromorphic solution on $\C^*$ with only
simple zeros and simple poles.

(e) The space $M_{3FN}^{reg}(x_0)\cong \C^*\times \C$ with the coordinates
$(f_0,g_0)$ is the space of initial conditions at $x_0$
\begin{eqnarray}\label{10.14}
(f(x_0),\paa_xf(x_0))=(f_0,\frac{2f_0}{x_0}\MGcdot g_0)
\end{eqnarray}
for at $x_0$ regular local solutions $f$.
For $k=0,1,2,3$, the space $M_{3FN}^{[k]}(x_0)\cong \C$ with the coordinate
$\www g_k$ is the space of the initial condition
\begin{eqnarray}\label{10.15}
\varepsilon_2\paa_x^3 f^{\varepsilon_1}(x_0)
=\frac{2}{x_0^2}+\frac{8}{x_0}\MGcdot \www g_k
\end{eqnarray}
for the at $x_0$ singular local solutions $f$ with 
$f^{\varepsilon_1}(x_0)=0,\varepsilon_2\paa_xf^{\varepsilon_1}(x_0)=-2.$
This initial condition determines such a singular local solution.

(f) The algebraic manifolds $M_{3FN}^{ini}$ and $M_{3FN}^{reg}$ and the 
foliation on $M_{3FN}^{mon}$ had been considered in a different way in
\cite{Ok79}, in fact for all $P_{III}(D_6)$ equations.
The spaces $M_{3FN}^{ini}(x_0)$, $x_0\in\C^*$, are Okamoto's spaces
of initial conditions.
Descriptions of $M_{3FN}^{ini}$ by four charts are given in
\cite{MMT99}, \cite{NTY02}, \cite{Te07}, but part (b) and the four types of 
singular initial conditions are not made explicit there.
\end{theorem}

{\bf Proof:}
(a) By lemma \ref{t8.5}, any leaf intersects 
$M_{3FN}^{sing}:=\cup_{k=0}^3M_{3FN}^{[k]}$ only in a discrete set of points.
On the intersection of the leaf with $M_{3FN}^{reg}$, the restrictions of 
$f_0,g_0$ and $\www g_0$ to the leaf are multi-valued holomorphic functions
in $x$, and $f_0$ takes there values only in $\C^*$.
Part (b) will show that the restrictions 
$f_{mult}(.,s,B),g_{mult}(.,s,B)$ and $\www g_{mult}(.,s,B)$ 
of $f_0,g_0$ and $\www g_0$ 
to the leaf are meromorphic at the intersection points
with $M_{3FN}^{sing}$, with only simple zeros and simple poles.
It remains to establish the differential equations for $f_{mult}(.,s,B)$.

Let $U\subset (\textup{a leaf})\cap M_{3FN}^{reg}$ be a small open subset
in a leaf, and let $U^\prime:=pr_x(U)\subset\C^*$ be the (isomorphic to $U$)
open subset in $\C^*$. 
Theorem \ref{t8.2} (b) provides a basis $\uuuu\sigma_0$
of $\Gamma(\P^1,\OO(H))$ for any $P_{3D6}$-TEJPA bundle in $U$,
with \eqref{8.15}-\eqref{8.18} with $u^1_0=u^1_\iiii=x$.
The proof of theorem \ref{t8.2} (b), remark \ref{t8.1} (ii)+(iii),
theorem \ref{t7.3} (c) and the correspondence in \eqref{2.19} provide
a $4$-tuple of bases $\uuuu e^{\pm}_0,\uuuu e^\pm_\iiii$ with
\eqref{6.8}-\eqref{6.11}, \eqref{7.10}-\eqref{7.12},
and matrices $A^\pm_0\in GL(2,\AAA_{|I^\pm_0}),
A^\pm_\iiii\in GL(2,(\rho_1^*\AAA)_{|I^\pm_\iiii})$ with
$\widehat A^+_{0/\iiii}=\widehat A^-_{0/\iiii}, 
\widehat A^\pm_0(0)=\widehat A^\pm_\iiii(\iiii)={\bf 1}_2$ and
\begin{eqnarray}\label{10.16}
\uuuu\sigma_{0|\whhh I^\pm_0} &=& \uuuu e^\pm_0\MGcdot
\begin{pmatrix}e^{-x/z} &0 \\0 & e^{x/z}\end{pmatrix}\MGcdot A^\pm_0\MGcdot C
\ \textup{ near }0,\\ \label{10.17}
\uuuu\sigma_{0|\whhh I^\pm_\iiii} &=& \uuuu e^\pm_\iiii\MGcdot
\begin{pmatrix}e^{-xz} &0 \\0 & e^{xz}\end{pmatrix}\MGcdot A^\pm_\iiii\MGcdot C
\MGcdot\begin{pmatrix}f_0 & 0\\ 0 &f_0^{-1}\end{pmatrix} \ \textup{near }\iiii.
\end{eqnarray}
With respect to the isomonodromic extension of the connection $\nnn$ of the 
single $P_{3D6}$-TEJPA bundles in $U^\prime$, the $\uuuu e^\pm_0,\uuuu e^\pm_\iiii$
are flat families (in $x$) of flat sections (in $z$).
The basis $\uuuu\sigma_0$ and the matrices $A^\pm_0,A^\pm_\iiii$
depend holomorphically on $x$ and $z$.

By \eqref{8.15}, the pole parts at $0$ and $\iiii$ of $\nnn_\zdz\uuuu\sigma_0$
are $\uuuu\sigma_0\MGcdot\frac{x}{z}\begin{pmatrix}0&1\\1&0\end{pmatrix}$
and $\uuuu\sigma_0\MGcdot(-xz)\begin{pmatrix}0&f_0^{-2}\\f_0^2&0\end{pmatrix}$.
Because of \eqref{10.16} and \eqref{10.17} and
\[
\xdx e^{-x/z}=-\frac{x}{z}=-\zdz e^{-x/z},\quad
\xdx e^{-xz}=-xz e^{-xz} = \zdz e^{-xz},
\]
the pole parts at $0$ and $\iiii$ of $\nnn_\xdx\uuuu\sigma_0$ are
$\uuuu\sigma_0\MGcdot\frac{-x}{z}\begin{pmatrix}0&1\\1&0\end{pmatrix}$
and $\uuuu\sigma_0\MGcdot(-xz)\begin{pmatrix}0&f_0^{-2}\\f_0^2&0\end{pmatrix}$.
By calculations with $\xdx$ analogously to those with $\zdz$ in
\eqref{8.12} and \eqref{8.14}, one obtains that the part between the
pole parts has the form 
$\uuuu\sigma_0\MGcdot a(x,z)\begin{pmatrix}1&0\\0&-1\end{pmatrix}$
for some holomorphic function $a(x,z)$. Thus
\begin{eqnarray}\label{10.18}
\nnn_\xdx \uuuu\sigma_0 = \uuuu\sigma_0\left[
\frac{-x}{z}\begin{pmatrix}0&1\\1&0\end{pmatrix}
+a\begin{pmatrix}1&0\\0&-1\end{pmatrix}
-xz\begin{pmatrix}0&f_0^{-2}\\f_0^2&0\end{pmatrix}\right].
\end{eqnarray}
The flatness of $\nnn$ gives together with \eqref{8.15} and \eqref{10.18}
\begin{eqnarray}
0&=&(\nnn_\zdz\nnn_\xdx-\nnn_\xdx\nnn_\zdz)\uuuu\sigma_0\nonumber\\
&=& \uuuu\sigma_0\left[2\frac{x}{z}(a-g_0)\begin{pmatrix}0&-1\\1&0\end{pmatrix}
+(-2x^2(f_0^2-f_0^{-2})+\xdx g_0)\begin{pmatrix}1&0\\0&-1\end{pmatrix}\right. 
\nonumber\\
&+&\left. 2(g_0+a)xz\begin{pmatrix}0&f_0^{-2}\\-f_0^2&0\end{pmatrix} 
+xz\MGcdot\xdx \begin{pmatrix}0&f_0^{-2}\\f_0^2&0\end{pmatrix}\right].
\label{10.19}
\end{eqnarray}
Thus with $f_0=e^{\varphi/2}$ for a chosen branch $\varphi=2\log(f_0)$
\begin{eqnarray}\label{10.20}
a&=& g_0,\\
\xdx g_0 &=& 2x^2(f_0^2-f_0^{-2}) =4x^2\sinh \varphi,\label{10.21}\\
\xdx f_0^2&=& 4g_0f_0^2,\ \xdx f_0^{-2}=-4g_0f_0^{-2},\ \xdx\varphi=4g_0,
\label{10.22}\\
(\xdx)^2\varphi&=& \xdx(4g_0)=16x^2\sinh \varphi,\nonumber
\end{eqnarray}
which is the radial sinh-Gordon equation \eqref{9.5} and which is 
equivalent to $P_{III}(0,0,4,-4)$ for $f_0=e^{\varphi/2}$.

(b) Consider a leaf in $M_{3FN}^{mon}$ and the corresponding solution $f_0$ of $P_{III}(0,0,4,-4)$. By theorem \ref{t8.2} (g) and \eqref{8.20}
$f_0(x_0)\in\C^*$ holds (for one branch of the multi-valued function $f_0$)
if and only if the corresponding intersection point of the leaf with 
$M_{3FN}^{ini}(x_0)$ is in $M_{3FN}^{reg}(x_0)$.

Consider an intersection point of a leaf in $M_{3FN}^{mon}$ with
$M_{3FN}^{[k]}$, a neighbourhood $U$ in the leaf of this point and the
functions $(f_0,g_0)_{|U}$. Then 
$(\varepsilon_2f_0^{\varepsilon_1},g_k)_{|U}=(f_k,g_k)_{|U}$ by \eqref{8.20},
and this pair is mapped by $R_k$ (in remark \ref{t8.3}) to
$(f_0,g_0)_{|R_k(U)}$, and the intersection point is mapped to an intersection
point of the leaf $R_k\ (\textup{old leaf})$ with $M_{3FN}^{[0]}$.
Therefore it is sufficient to prove \eqref{10.13} for $k=0$.

Suppose now $k=0$. Then the corresponding branch of $f_0$ satisfies
$f_0(x_0)=0$ by theorem \ref{t8.2} (g).
\eqref{10.7} extends in the form
$$f_0g_0=-x+\frac{1}{2}f_0+\frac{1}{2}f_0^2 \www g_0$$ from 
$M_{3FN}^{reg}$ to $M_{3FN}^{[0]}$.
Together with \eqref{10.12} this shows
\begin{eqnarray*}
\paa_x f_0&=& -2+\frac{f_0}{x}+\frac{f_0^2}{x}\www g_0,\\
\paa_x f_0(x_0)&=&-2,\\
\paa_x^2f_0&=& -\frac{f_0}{x^2}+\frac{\paa_xf_0}{x}-\frac{f_0^2}{x^2}\www g_0
+\frac{2f_0\paa_xf_0}{x}\www g_0+\frac{f_0^2}{x}\www \paa_xg_0,\\
\paa_x^2 f_0(x_0)&=& \frac{-2}{x_0},\\
\paa_x^3 f_0(x_0)&=& \frac{2}{x_0^2} + \frac{8}{x_0}\www g_0(x_0).
\end{eqnarray*}

(c)+(d)+(e) 
The initial conditions $(f(x_0),\paa_xf(x_0))$ for a local
regular solution $f$ give rise to a leaf through the point
$(f_0,g_0)=(f(x_0),\frac{x_0}{2 f(x_0)}\paa_x f(x_0))$.
Then the restriction of $f_0$ to this leaf is the extension of the local
solution $f$ to a global multi-valued solution in $\C^*$. It has 
only simple zeros and simple poles by (b).
This establishes the Painlev\'e property.
The other statements follow directly from (a) and (b) as well.

(f) On the one hand, the equality of $M_{3FN}^{ini}$, $M_{3FN}^{reg}$
and the foliation on $M_{3FN}^{mon}$ with data of Okamoto follows
from the uniqueness of Okamoto's data, which is rather obvious.
On the other hand, one can directly and easily compare the four charts
for $M_{3FN}^{ini}$ with those in \cite{MMT99} or \cite{Te07},
who derive their charts from Okamoto's description of his data.
\hfill$\Box$

\begin{remarks}\label{t10.4}
(i) The construction of $M_{3FN}^{ini}$, $M_{3FN}^{reg}$, $M_{3FN}^{mon}$
and the foliation on $M_{3FN}^{mon}$ here is independent of \cite{Ok79},
in contrast to \cite{MMT99}, \cite{NTY02}, \cite{Te07}.
But \cite{Ok79} provides for any space $M_{3FN}^{ini}(x_0)$ of initial conditions
a natural compactification $S$ and a divisor $Y\subset S$ of type $\www D_6$
such that $M_{3FN}^{ini}(x_0)\cong S-Y$. It would be interesting to recover
$S$ and $Y$ using generalizations of $P_{3D6}$-TEJPA bundles,
and possibly derive new results on the solutions of $P_{III}(0,0,4,-4)$.

(ii) For any $P_{III}(D_6)$-equation, all solutions $f$ satisfy a generalization
of \eqref{10.13}.  This means that they have only simple zeros and simple poles,
and there are two types of simple zeros and two types of simple poles,
distinguished by the sign of $\paa_x f(x_0)$ or
$\paa_xf^{-1}(x_0)$, and the initial condition at a simple zero 
(respectively, pole) is $\paa_x^3f(x_0)$ (respectively, $\paa_x^3f^{-1}(x_0)$).
This can be proved by a power series ansatz.

(iii) By theorem \ref{t10.3}, the algebraic manifold $M_{3FN}^{ini}$ has
four algebraic hypersurfaces, one for each type of simple zeros or simple poles
of the solutions, and these are precisely the hypersurfaces for 
$(1,-1)$-twistors (cf.\  remark \ref{t4.3} (ii)).

The union of all spaces of initial conditions of all $P_{III}(D_6)$-equations
is an algebraic \cite{MMT99}, \cite{Te07} manifold and has by (ii) four analytic 
(in fact, algebraic) hypersurfaces, one for 
each type of simple zeros or simple poles of the solutions.
But in the linear system of \cite{JM81} for the $P_{III}(D_6)$-equations,
only one of the four hypersurfaces corresponds to $(1,-1)$-twistors.
For the linear system for $P_{III}(0,4,4,-4)$ this follows from the 
$4{:}1$ folding transformation which is mentioned in remark \ref{t9.1} (ii)
and which connects the $P_{3D6}$-TEJPA bundles for $P_{III}(0,0,4,-4)$
with the linear systems in \cite{JM81} for $P_{III}(0,4,4,-4)$.
The general case will be proved elsewhere, together with the claims in remark \ref{t9.1}.

(iv) In this paper we do not make use of the Hamiltonian description
of the $P_{III}(D_6)$-equations. For $P_{III}(0,0,4,-4)$
it is implicit in the coordinates $\www g_k$ and the equations \eqref{10.7}
and \eqref{10.12}.

(v) The matrices $A^\pm_0(x,z)$ and $A^\pm_\iiii(x,z)$ in \eqref{10.16}
and \eqref{10.17} deserve more careful consideration.
The properties of the basis $\uuuu\sigma_0$ and the bases $\uuuu e^\pm_0,
\uuuu e^\pm_\iiii$ with respect to $A, P$ and $J$ imply the following.
$A$ gives for $A^\pm_0$
\begin{eqnarray}\label{10.23}
A^\mp_0(x,-z)&=& \begin{pmatrix}0&-1\\1&0\end{pmatrix}
A^\pm_0(x,z)\begin{pmatrix}0&1\\-1&0\end{pmatrix}.
\end{eqnarray}
$P$ gives for $A^\pm_0$
\begin{eqnarray}\label{10.24}
A^\mp_0(x,-z)&=& (A^\pm_0(x,z)^t)^{-1}.
\end{eqnarray}
\eqref{10.23} and \eqref{10.24} are equivalent to \eqref{10.23}
and $\det A^\pm_0(x,z)=1$. $J$ gives for $A^\pm_0$ and $A^\pm_\iiii$
\begin{eqnarray}\label{10.25}
A^\pm_\iiii(x,z)&=& \begin{pmatrix}1&0\\0&-1\end{pmatrix}
A^\pm_0(x,\frac{1}{z})\begin{pmatrix}1&0\\0&-1\end{pmatrix}
\end{eqnarray}
(here $c=1$ because of $u^1_0=u^1_\iiii=x$).
Defining the matrices $B^\pm_{0/\iiii}(x,z)$ by 
\begin{eqnarray}\label{10.26}
C\MGcdot B^\pm_{0/\iiii}(x,z) &=& A^\pm_{0/\iiii}(x,z)\MGcdot C,
\end{eqnarray}
\eqref{10.23}-\eqref{10.25} are equivalent to
\begin{eqnarray}\label{10.27}
B^\mp_0(x,-z)&=& \begin{pmatrix}1&0\\0&-1\end{pmatrix}
B^\pm_0(x,z)\begin{pmatrix}1&0\\0&-1\end{pmatrix},\\
\det B^\pm_0(x,z)&=& 1,\label{10.28}\\
B^\pm_\iiii(x,z)&=& \begin{pmatrix}0&1\\1&0\end{pmatrix}
B^\pm_0(x,\frac{1}{z})\begin{pmatrix}0&1\\1&0\end{pmatrix}.\label{10.29}
\end{eqnarray}
In particular
\begin{eqnarray}\label{10.30}
\whhh B^\pm_0(x,z)={\bf 1}_2+z\MGcdot
\begin{pmatrix}0&\beta_1\\ \beta_2&0\end{pmatrix}
+z^2\MGcdot\begin{pmatrix}\beta_3&0\\0&\beta_4\end{pmatrix}+ \dots
\end{eqnarray}
where $\beta_1,\beta_2,\beta_3,\beta_4$ are holomorphic functions in $x$ with 
$\beta_3+\beta_4-\beta_1\beta_2=0$, and 
\begin{eqnarray}\label{10.31}
(\whhh B^\pm_0(x,z))^{-1}={\bf 1}_2-z\MGcdot
\begin{pmatrix}0&\beta_1\\ \beta_2&0\end{pmatrix}
+z^2\MGcdot\begin{pmatrix}\beta_4&0\\0&\beta_3\end{pmatrix}+ \dots
\end{eqnarray}
Taking the derivatives of \eqref{10.16} and \eqref{10.17}
by $\zdz$ and $\xdx$ gives \eqref{8.15} and \eqref{10.18} with
\begin{eqnarray}\label{10.32}
x(\beta_1-\beta_2)&\!=\!&g_0=a=-x(\beta_1-\beta_2)+f_0^{-1}\xdx f_0,\\
-2\beta_3+\beta_2^2-\frac{\beta_2}{x} &\!=\!& f_0^2 = 2\beta_3-\beta_2^2-\paa_x\beta_2,
\label{10.33}\\
-2\beta_4+\beta_1^2-\frac{\beta_1}{x} &\!=\!& f_0^{-2} 
= 2\beta_4-\beta_1^2-\paa_x\beta_1.
\label{10.34}
\end{eqnarray}
One recovers $2f_0g_0=\xdx f_0$ from \eqref{10.32} and 
$\xdx g_0=2x^2(f_0^2-f_0^{-2})$ from \eqref{10.32}-\eqref{10.34}.
This is an alternative proof of most of theorem \ref{t10.3} (a).

(vi) Any control of the dependence of the (multi-valued) solutions 
$f_{mult}(x,s,B)$ in $x$ of $P_{III}(0,0,4,-4)$ on the parameters $s$ and $B$ 
of the monodromy tuple would be very welcome. For example, are
there differential equations governing the dependence on $s$ or $B$?
We do not see any and wonder whether there might be reasons (such as Umemura's results on irreducibility) which prevent their existence.
\end{remarks}

\begin{remarks}\label{t10.5}
For each $u^1_0,u^1_\iiii$ there are two completely reducible 
$P_{3D6}$-TEP bundles (theorem \ref{t6.3} (c)), and above
each of them are two $P_{3D6}$-TEJPA bundles (theorem \ref{t7.6} (b)).
The two $P_{3D6}$-TEJPA bundles over the $P_{3D6}$-TEP bundle with
$(C_{11}),(C_{22})$ have monodromy data $(s,B)=(0,\pm{\bf 1}_2)$,
and they are fixed points of the action of $R_1$ on 
$M_{3TJ}^{mon}(u^1_0,u^1_\iiii)$. 
The two $P_{3D6}$-TEJPA bundles over the $P_{3D6}$-TEP bundle with
$(C_{12}),(C_{21})$ have monodromy data 
$(s,B)=(0,\pm 
\bsp
0&1\\-1&0
\esp
)$,
and they are fixed points of the action of $R_3$ on 
$M_{3TJ}^{mon}(u^1_0,u^1_\iiii)$.

As $R_1, R_2$ and $R_3$ map the coordinate $f_0$ on $M^{reg}_{3TJ}$ 
to $f_0^{-1}, -f_0$ and $-f_0^{-1}$, respectively (remark \eqref{t8.3}), 
we find 
\begin{eqnarray}\label{10.35}
f_{mult}(x,0,\pm{\bf 1}_2)&=& \varepsilon_1\MGcdot (\pm 1),\\
f_{mult}(x,0, \pm \begin{pmatrix}0&1\\-1&0\end{pmatrix}) &=& 
\varepsilon_2\MGcdot (\pm i).\label{10.36}
\end{eqnarray}
for any $x$, with some yet to be determined signs
$\varepsilon_1,\varepsilon_2\in\{\pm 1\}$.

Consider the case $(s,B)=(0,{\bf 1}_2)$. Then the four flat bases
$\uuuu e^\pm_0,\uuuu e^\pm_\iiii$ are restrictions of one global
flat basis $\uuuu e$. It satisfies 
\[
J(\uuuu e(z))=\uuuu e(\rho_c (z))\begin{pmatrix}1&0\\0&-1\end{pmatrix}.
\]
The basis $\uuuu\sigma_0$ from theorem \ref{t8.2} (b) is
\begin{eqnarray}\label{10.37}
\uuuu\sigma_0(z)=\uuuu e(z)\begin{pmatrix}e^{-u^1_0/z-u^1_\iiii z} & 0 \\
0 & e^{u^1_0/z+u^1_\iiii z} \end{pmatrix}\MGcdot C.
\end{eqnarray}
With
\[
u^1_0/z=u^1_\iiii\MGcdot\rho_c(z)\textup{ and }
u^1_\iiii \MGcdot z=u^1_0/\rho_c(z)
\]
one calculates
\begin{align}
\uuuu\sigma_0(\rho_c(&z))\begin{pmatrix}0& f_0^{-1}\\f_0&0\end{pmatrix}
= J(\uuuu\sigma_0(z))
\nonumber
\\
&= \uuuu e(\rho_c (z))\begin{pmatrix}1&0\\0&-1\end{pmatrix}
\begin{pmatrix}e^{-u^1_0/z-u^1_\iiii z} & 0 
\\
0 & e^{u^1_0/z+u^1_\iiii z} \end{pmatrix}\MGcdot C
\nonumber 
\\
&= \uuuu\sigma_0(\rho_c(z))\begin{pmatrix}0& 1\\1&0\end{pmatrix}.
\label{10.38}
\end{align}
This shows $\varepsilon_1=1$. 

Consider the case $(s,B)=(0,\pm 
\bsp
0&1\\-1&0
\esp
)$. 
Then the two bases $\uuuu e^\pm_0$ are restrictions of a global flat basis 
$\uuuu e_0$, the two bases $\uuuu e^\pm_\iiii$ are restrictions of a global
flat basis $\uuuu e_\iiii$, and 
$$\uuuu e_\iiii= \uuuu e_0\begin{pmatrix}0&1\\-1&0\end{pmatrix}.$$
Again
$$J(\uuuu e_0(z))=\uuuu e_\iiii(\rho_c(z))
\begin{pmatrix}1&0\\0&-1\end{pmatrix}.$$
The basis $\uuuu\sigma_0$ from theorem \ref{t8.2} (b) is
\begin{eqnarray}\label{10.39}
\uuuu\sigma_0(z)=\uuuu e_0(z)\begin{pmatrix}e^{-u^1_0/z+u^1_\iiii z} & 0 \\
0 & e^{u^1_0/z-u^1_\iiii z} \end{pmatrix}\MGcdot C.
\end{eqnarray}
Similarly as above, one calculates
\begin{eqnarray}\label{10.40}
\uuuu\sigma_0(\rho_c(z))\begin{pmatrix}0& f_0^{-1}\\f_0&0\end{pmatrix}
= J(\uuuu\sigma_0(z))
= \uuuu\sigma_0(\rho_c(z))\begin{pmatrix}0& -i\\i&0\end{pmatrix}.
\end{eqnarray}
This shows $\varepsilon_2=1$. 
\end{remarks}

\chapter[Comparison with the setting of Its et al.]
{Comparison with the setting of Its, Novokshenov, and Niles}\label{s11}
\setcounter{equation}{0}

\noindent
Theorem \ref{t10.3} fixes the relation between the solutions of 
$P_{III}(0,0,4,-4)$ and the isomonodromic families of $P_{3D6}$-TEJPA bundles.
As said above, this relation is due to \cite{FN80}, though without the
vector bundle language and the structural data $P,A,J$ which incorporate
the symmetries. 

Good use of this relation has been made in 
\cite{IN86}, \cite{FIKN06}, \cite{Ni09}.  In particular, 
\cite{IN86} connects for many, but not all, solutions of $P_{III}(0,0,4,-4)$
their asymptotic behaviour as $x\to 0$ and $x\to\iiii$ 
with the Stokes data and the central connection matrix 
of the corresponding isomonodromic families of $P_{3D6}$-TEJPA bundles.
\cite{FIKN06} rewrites and extends \cite{IN86}.

\cite{Ni09} builds on \cite{IN86}, \cite{FIKN06} and provides asymptotic formulae as $x\to 0$
for all solutions of $P_{III}(0,0,4,-4)$.
These formulae and a proof of theorem \ref{t4.2} for the case of 
$P_{III}(0,0,4,-4)$ are the two main results of \cite{Ni09}.
In chapter \ref{s12} we shall rewrite Niles' asymptotic formulae for 
$x\to 0$ and simplify and slightly extend them.

As preparation for that, and in order to facilitate access to the rich results
in \cite{IN86}, \cite{FIKN06}, we shall now connect our setting
with that in \cite{Ni09} (and hence \cite{IN86}, \cite{FIKN06}).

As we have seen,  equation $P_{III}(0,0,4,-4)$ \eqref{9.1} for $f$ can be rewritten, using $f=e^{\varphi/2}$, as the radial sinh-Gordon equation 
$(\xdx)^2\varphi=16 x^2\sinh\varphi$ \eqref{9.5}. 
In \cite[(1.1)]{Ni09} the sine-Gordon equation
\begin{eqnarray}\label{11.1}
(\paa_{x^{NI}}^2+\frac{1}{x^{NI}}\paa_{x^{NI}})u^{NI}(x^{NI})
=-\sin u^{NI}(x^{NI})
\end{eqnarray}
is used. Recall that $(\xdx)^2=x^2(\paa_x^2+\frac{1}{x}\paa_x).$
We have:

\begin{lemma}\label{t11.1}
\eqref{9.5} and \eqref{11.1} are related by
\begin{eqnarray}\label{11.2}
x^{NI}&=&4x,\\
u^{NI}(x^{NI})=u^{NI}(4x)&=&u(x)=i\varphi(x)+\pi,\label{11.3}\\
\varphi(x)&=&-iu(x)+i\pi. \nonumber
\end{eqnarray}
Here $u(x)$ satisfies the sine-Gordon equation
\begin{eqnarray}\label{11.4}
(\paa_x^2+\tfrac{1}{x}\paa_x)u=-16\sin u,\quad\textup{i.e.\ }
(\xdx)^2u=-16x^2\sin u.
\end{eqnarray}
and $f=e^{\varphi/2}$ satisfies
\begin{eqnarray}\label{11.5}
f(x)=e^{(-iu(x)+i\pi)/2}=ie^{-iu(x)/2}.
\end{eqnarray}
\end{lemma}

\begin{remarks}\label{t11.2}
(i) NI stands for Niles-Novokshenov-Its.

(ii) $u(x)$ is introduced in lemma \ref{t11.1} so that we can rewrite
below the formulae in \cite{Ni09} with $x$ and $u(x)$ instead of 
$x^{NI}$ and $u^{NI}(x^{NI})$.
%
%(iii) \cite{IN86} and \cite{FIKN06} use
%$$x^{IN}=2x=\frac{1}{2}x^{NI}
%\textup{ and }u^{IN}(x^{IN})=u(x)=u^{NI}(x^{NI}).$$
\end{remarks}

\cite{Ni09} considers in formulae (1.25)-(1.28) the first order matrix
differential equation
\begin{eqnarray}\label{11.6}
\zdz\Psi(z) &=& z A(z)\Psi(z)\quad
\textup{ with }\\
z A(z)&=& \frac{i}{z}\begin{pmatrix}\cos u & i\sin u\\-i\sin u&-\cos u\end{pmatrix}
+ g_0\begin{pmatrix}0&1\\1&0\end{pmatrix}
+ ix^2\MGcdot z\begin{pmatrix}1&0\\0&-1\end{pmatrix}.\nonumber
\end{eqnarray}
Here we have called $x,u(x),\lambda,w$ in \cite{Ni09} 
$x^{NI},u^{NI}(x^{NI}),z,g_0$ and replaced $x^{NI}$ and $u^{NI}(x^{NI})$
by $4x$ and $u(x)$.

We associate with \eqref{11.6} the trivial holomorphic vector bundle $H$
of rank 2 on $\P^1$ with global basis $\uuuu\rho=(\rho_1,\rho_2)$
and meromorphic connection $\nnn$ given by
\begin{eqnarray}\label{11.7}
\nnn_\zdz\uuuu\rho =\uuuu\rho\MGcdot (-zA(z)).
\end{eqnarray}
Then $\Psi$ is a solution of \eqref{11.6} in a sector of $\C^*$ if and only if
$\uuuu\rho\MGcdot\Psi$ is a flat basis of $(H,\nnn)$ in this sector.

\begin{lemma}\label{t11.3}
(a) The pair $(H,\nnn)$ takes the normal form \eqref{8.15}
\begin{eqnarray}\label{11.8}
\nnn_\zdz\uuuu\sigma_0 = \uuuu\sigma_0\MGcdot\left[
\frac{u^1_0}{z}\begin{pmatrix}0&1\\1&0\end{pmatrix} 
-g_0\begin{pmatrix}1&0\\0&-1\end{pmatrix}
-u^1_\iiii\MGcdot z\begin{pmatrix}0&f_0^{-2}\\f_0^2&0\end{pmatrix}\right]
\end{eqnarray}
of a pure $P_{3D6}$-TEJPA bundle with
\begin{eqnarray}\label{11.9}
f_0&=& e^{\varphi/2}=ie^{-iu/2},\quad \varphi=-iu+i\pi,\\
u^1_0&=&-u^2_0=-i,\ u^1_\iiii=-u^2_\iiii=ix^2,\label{11.10}\\
\uuuu\sigma_0 &=& \uuuu\rho\MGcdot P_0\MGcdot \begin{pmatrix}1&0\\0&i\end{pmatrix}
\MGcdot C = 
\uuuu\rho\MGcdot (-i)\MGcdot 
\begin{pmatrix}e^{\varphi/2}& -e^{-\varphi/2}\\e^{\varphi/2}&e^{-\varphi/2}\end{pmatrix},
\label{11.11}\\
P_0&=& \begin{pmatrix}\cos\frac{u}{2} & -i\sin\frac{u}{2}\\
-i\sin\frac{u}{2} & \cos\frac{u}{2}\end{pmatrix}
=(-i)\begin{pmatrix}\sinh\frac{\varphi}{2} & \cosh\frac{\varphi}{2}\\
\cosh\frac{\varphi}{2} & \sinh\frac{\varphi}{2}\end{pmatrix}\label{11.12}
\end{eqnarray}
and $C$ as in \eqref{8.5}.

(b) (Definition) We enrich the pair $(H,\nnn)$ to the $P_{3D6}$-TEJPA bundle
in \eqref{8.15}-\eqref{8.18} with $u^1_0,u^1_\iiii,f_0,g_0$ as above.

Furthermore, the two points in $(c^{2{:}1})^*M_{3TJ}^{reg}$ above
this $P_{3D6}$-TEJPA bundle have the coordinates $(x,y,f_0,g_0)=(x,ix,f_0,g_0)$
and $(-x,-ix,f_0,g_0)$. We distinguish the point $(x,ix,f_0,g_0)$.

(Lemma) Then $f_0,\varphi$ and $u$ restricted to the leaf through this point
are solutions of $P_{III}(0,0,4,-4)$, \eqref{9.5} and \eqref{11.4},
respectively.
\end{lemma}

{\bf Proof:}
(a) This is an elementary calculation, which we omit.

(b) Only the last statement is not a definition.
It follows from theorem \ref{t10.3}.\hfill$\Box$

Lemma \ref{t11.3} connects the normal form \eqref{11.6} above from 
\cite[(1.25)-(1.28)]{Ni09}
with the normal form in theorem \ref{t8.2} (b) respectively with a 
point in $(c^{2{:}1})^*M_{3TJ}^{reg}$.
Now we shall connect the Stokes data and the central connection matrix in 
\cite{Ni09} with those in theorem \ref{t6.3} (e) and theorem \ref{t7.5} (c).

The following data are fixed in \cite[1.4.1]{Ni09}.
A small $\varepsilon>0$ is chosen. Only $x$ in the sector
\begin{eqnarray}\label{11.13}
S^{NI}:=\{x\in\C^*\, |\, \arg x \in(-\tfrac{\pi}{4}+\varepsilon,
\tfrac{\pi}{4}-\varepsilon)\ \textup{mod }2\pi \}
\end{eqnarray}
is considered. The four sectors for $k=1,2,$
\begin{equation}\label{11.14}
\Omega^{(0)}_k:= \{z\in\C^*\, |\, \arg z \in(\pi(k-2),\pi k)\ \textup{mod } 2\pi\}, 
\end{equation}
\begin{equation}\label{11.15}
\Omega^{(\iiii)}_k:= \{z\in\C^*\, |\, \arg z \in
(\pi(k\!-\!\tfrac{3}{2})\!-\!2\varepsilon,\pi (k\!-\!\tfrac{1}{2})+2\varepsilon)
\ \textup{mod } 2\pi \}, 
\end{equation}
are used, and in each of them a solution $\Psi_k^{(0)}$ or $\Psi_k^{(\iiii)}$
of \eqref{11.6} is considered. 

%Later 1 picture
\includegraphics[width=0.9\textwidth]{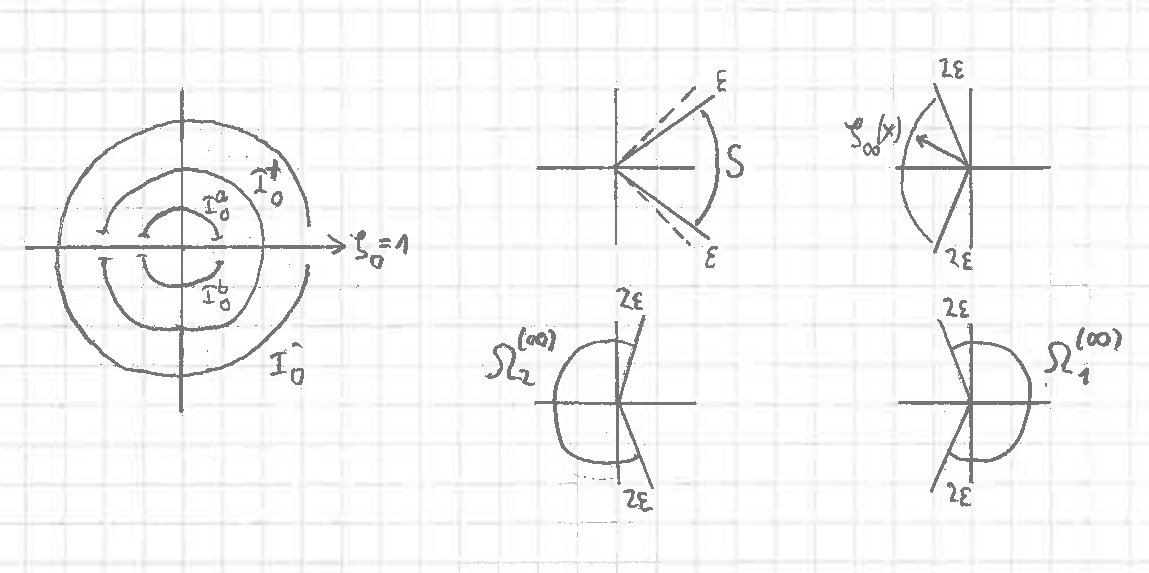} 

These solutions are fixed by their asymptotics
near $z=0$ respectively $z=\iiii$,
\begin{eqnarray}\label{11.16} 
\Psi_k^{(0)}(z)&:=& A^{(0)}_k(z)
\begin{pmatrix}e^{-i/z}&0\\0&e^{i/z}\end{pmatrix}\\ 
&&\textup{with }
A^{(0)}_k\in GL(2,\AAA|_{\Omega^{(0)}_k\cap S^1}),
A^{(0)}_k(0)=P_0,\nonumber \\
\Psi_k^{(\iiii)}(z)&:=& A^{(\iiii)}_k(z)
\begin{pmatrix}e^{-ix^2z}&0\\0&e^{ix^2z}\end{pmatrix}\label{11.17}\\ 
&&\textup{with }
A^{(\iiii)}_k\in GL(2,(\rho_1^*\AAA)|_{\Omega^{(\iiii)}_k\cap S^1}),
A^{(\iiii)}_k(\iiii)={\bf 1}_2.\nonumber
\end{eqnarray}
Then four Stokes matrices $S^{(0)}_k$ and $S^{(\iiii)}_k$, $k=1,2,$ are defined by
\begin{eqnarray}\label{11.18}
\begin{split}
S^{(0)}_k(z):=& \Psi_k^{(0)}(z)^{-1} \MGcdot \Psi_{k\pm 1}^{(0)}(z)
\quad \textup{for }z\textup{ with }\\ 
&\arg z\in(\pi(k-1),\pi k)\ \textup{mod } 2\pi,\\
S^{(\iiii)}_k(z):=& \Psi_k^{(\iiii)}(z)^{-1}
\MGcdot \Psi_{k\pm 1}^{(\iiii)}(z)\quad\textup{for }z\textup{ with } \\
& \arg z\in(\pi(k-\tfrac{1}{2})-2\varepsilon,
\pi (k-\tfrac{1}{2})+2\varepsilon)\  \textup{mod } 2\pi,
\end{split}
\end{eqnarray}
and a connection matrix by
\begin{eqnarray}\label{11.19}
E^{NI}: =\Psi_1^{(0)}(z)^{-1}
\MGcdot \Psi_1^{(\iiii)}(z)\quad\textup{for }z\in\Omega^{(\iiii)}_1
(\subset \Omega^{(0)}_1).
\end{eqnarray}
The Stokes matrices satisfy
\begin{eqnarray}\label{11.20}
S_1^{(\iiii)}=S_2^{(0)}=\begin{pmatrix}1&0\\ \zeta^{NI}&1\end{pmatrix},\quad
S_2^{(\iiii)}=S_1^{(0)}=\begin{pmatrix}1&\zeta^{NI}\\0&1\end{pmatrix}
\end{eqnarray}
for some $\zeta^{NI}\in\C$. 
For the connection matrix $E^{NI}$, two cases are 
distinguished, the special case
(also called the {\it separatrix case} in \cite{Ni09}):
\begin{eqnarray}\label{11.21}
E^{NI}=\pm i\begin{pmatrix}0&1\\1&0\end{pmatrix},
\end{eqnarray}
and the general case:
\begin{eqnarray}\label{11.22}
&&E^{NI}=\frac{\pm 1}{\sqrt{1+pq}}\begin{pmatrix}1&p\\-q&1\end{pmatrix}\\
&&\textup{for some }p,q\in\C\textup{ with }pq\neq -1\textup{ and }
\zeta^{NI}=p+q.\nonumber
\end{eqnarray}

\begin{lemma}\label{t11.4}
(a) For $u^1_0,u^2_0,u^1_\iiii(x),u^2_\iiii(x)$ as in \eqref{11.10}
and $x\in S^{NI}$, the data in \eqref{2.1} and \eqref{2.2}
are related as follows to the sectors $\Omega^{(0)}_k,\Omega^{(\iiii)}_k$
in \cite[1.4.1]{Ni09}.
\begin{eqnarray}\label{11.23}
&\zeta_0=1,\ \zeta_\iiii(x)=-{|x^2|}/{x^2},&\\
&\whhh I^+_0=\C^*-\R_{>0}(-\zeta_0)=\Omega^{(0)}_1,\ 
\whhh I^-_0=\C^*-\R_{>0}\zeta_0=\Omega^{(0)}_2,&\label{11.24}\\
& \whhh I^+_\iiii(x)=\C^*-\R_{>0}(-\zeta_\iiii(x)),\ 
\whhh I^-_\iiii(x)=\C^*-\R_{>0}\zeta_\iiii(x),&\label{11.25}\\
& \bigcap_{x\in S^{NI}}\whhh I^+_\iiii(x) = \oooo{\Omega_2^{(\iiii)}}-\{0\},\ 
\bigcap_{x\in S^{NI}}\whhh I^-_\iiii(x) = \oooo{\Omega_1^{(\iiii)}}-\{0\}.&
\label{11.26}
\end{eqnarray}

(b) The $P_{3D6}$-TEJPA bundle associated in lemma \ref{t11.3} (b)
to \eqref{11.6} comes equipped with the canonical $4$-tuple 
$\uuuu e^\pm_0,\uuuu e^\pm_\iiii$ of flat bases from theorem \ref{t7.3} (c).
For an isomonodromic family with $x\in S^{NI}$ the bases
$\uuuu e^+_0,\uuuu e^-_0,\uuuu e^+_\iiii,\uuuu e^-_\iiii$ are constant
(i.e.\  they depend flatly on $x$) on 
$\Omega^{(0)}_1,\Omega^{(0)}_2,\Omega^{(\iiii)}_2,\Omega^{(\iiii)}_1,$
respectively, in view of \eqref{11.24} and \eqref{11.26}.
We have
\begin{eqnarray}\label{11.27}
\uuuu e^+_0&=& \uuuu\rho\MGcdot \Psi_1^{(0)}\MGcdot
\begin{pmatrix}1&0\\0&i\end{pmatrix},\quad
\uuuu e^-_0= \uuuu\rho\MGcdot \Psi_2^{(0)}\MGcdot
\begin{pmatrix}1&0\\0&i\end{pmatrix},\\ \label{11.28}
\uuuu e^+_\iiii&=& \uuuu\rho\MGcdot \Psi_2^{(\iiii)}\MGcdot
\begin{pmatrix}0&1\\-i&0\end{pmatrix},\quad
\uuuu e^-_\iiii= \uuuu\rho\MGcdot \Psi_1^{(\iiii)}\MGcdot
\begin{pmatrix}0&1\\-i&0\end{pmatrix}.
\end{eqnarray}

Let $s$ and $B(\beta)$ be the Stokes parameter and the central connection
matrix from \eqref{6.13}  
(for any $\beta$ satisfying \eqref{2.23}, i.e.\  $e^{-\beta}=4x^2$)
for $\uuuu e^\pm_0,\uuuu e^\pm_\iiii$. We have
\begin{eqnarray}\label{11.29}
s=i\zeta^{NI},
\end{eqnarray}
and, for the unique $\beta$ with $\Im(\beta)\in
(-\frac{\pi}{2}+2\varepsilon,\frac{\pi}{2}-2\varepsilon)$,
\begin{eqnarray}\label{11.30}
B(\beta)&=&\begin{pmatrix}b_1&b_2\\-b_2&b_1+sb_2\end{pmatrix}
=\begin{pmatrix}1&0\\0&-i\end{pmatrix}E^{NI}
\begin{pmatrix}0&1\\-i&0\end{pmatrix},\\
E^{NI}&=&\begin{pmatrix}1&0\\0&i\end{pmatrix}B(\beta)
\begin{pmatrix}0&i\\1&0\end{pmatrix}
= \begin{pmatrix} b_2&ib_1\\i(b_1+sb_2)&b_2\end{pmatrix}.\label{11.31}
\end{eqnarray}

(c) For this $\beta$ the special case \eqref{11.21} is the case
\begin{eqnarray}\label{11.32}
b_2=0,\quad b_1\in\{\pm 1\},\quad E^{NI}=ib_1\begin{pmatrix}0&1\\1&0\end{pmatrix},
\end{eqnarray}
and the general case \eqref{11.22} is the case
\begin{eqnarray}\label{11.33}
b_2\neq 0,\quad p=\frac{ib_1}{b_2},\ q=\frac{-i(b_1+sb_2)}{b_2},\
E^{NI}=b_2\begin{pmatrix}1&p\\-q&1\end{pmatrix}.
\end{eqnarray}

For fixed $x$ and fixed $u^1_0,u^1_\iiii,\beta$ with 
$e^{-\beta}=4x^2,\Im(\beta)\in(-\frac{\pi}{2}+2\varepsilon,
\frac{\pi}{2}-2\varepsilon)$,
$M_{3TJ}^{mon}(-i,ix^2) $ and $V^{mat}$ can be identified with the manifold
obtained by glueing
two copies of the manifold 
\[
\{\zeta^{NI}\in\C\}\cong\C,
\]
(for the special
case \eqref{11.21}, taking both signs)
and the double cover of the manifold
\[
\{(\zeta^{NI},p,q)\in\C^3\, |\, pq\neq -1,\zeta^{NI}=p+q\},
\]
on which $\sqrt{1+pq}$ becomes a well defined function
(for the general case \eqref{11.22}, taking both signs).

For fixed $(\zeta^{NI},p,q)$, the two choices $(s,b_1,b_2)$
and $(s,-b_1,-b_2)$ correspond to two solutions $f_0$ and $-f_0$
of $P_{III}(0,0,4,-4)$.
\end{lemma}

{\bf Proof:} 
(a) This follows from the definitions of $\zeta_0,\zeta_\iiii$
and $\whhh I^\pm_0,\whhh I^\pm_\iiii(x),\Omega_k^{(0)},\Omega_k^{(\iiii)}$.

(b) First \eqref{11.27} will be proved. 
The global holomorphic bases $\uuuu\sigma_0$ and $\uuuu\rho$ are related
by \eqref{11.11}. The basis $\uuuu\sigma_0$ and the flat basis $\uuuu e^\pm_0$
on the sector $\whhh I^\pm_0$ are related by
\begin{eqnarray}\label{11.34}
\uuuu\sigma_0 = \uuuu e^\pm_0\MGcdot
\begin{pmatrix}e^{i/z}&0\\0&e^{-i/z}\end{pmatrix}\MGcdot A^\pm_0\MGcdot C
\end{eqnarray}
with $A^\pm_0\in GL(2,\AAA|_{I^\pm_0}),\ \whhh A^+_0=\whhh A^-_0,\
\whhh A^\pm_0(0)={\bf 1}_2$.
The basis $\uuuu\rho$ and the flat basis $\uuuu \rho\MGcdot\Psi^{(0)}_k$
on $\Omega^{(0)}_k$ are related by
\begin{eqnarray}\label{11.35}
\uuuu \rho=\uuuu\rho\MGcdot\Psi^{(0)}_k\MGcdot (\Psi^{(0)}_k)^{-1}
=\uuuu\rho\MGcdot\Psi^{(0)}_k\MGcdot
\begin{pmatrix}e^{i/z}&0\\0&e^{-i/z}\end{pmatrix}\MGcdot (A^{(0)}_k)^{-1}
\end{eqnarray}
with $A^{(0)}_k\in GL(2,\AAA|_{\Omega^{(0)}_k\cap S^1}),\ 
A^{(0)}_k(0)=P_0$. 
Therefore (with $+\leftrightarrow 1, -\leftrightarrow 2$) 
\begin{eqnarray*}
\uuuu e^\pm_0 &=& \uuuu\rho\MGcdot \Psi^{(0)}_{1/2}\MGcdot
\begin{pmatrix}e^{i/z}&0\\0&e^{-i/z}\end{pmatrix}\MGcdot (A^{(0)}_{1/2})^{-1}
\MGcdot P_0\MGcdot\begin{pmatrix}1&0\\0&i\end{pmatrix}\MGcdot C\\
&& \quad\MGtimes  C^{-1}
\MGcdot (A^\pm_0)^{-1}\MGcdot\begin{pmatrix}e^{-i/z}&0\\0&e^{i/z}\end{pmatrix}
= \uuuu \rho\MGcdot \Psi^{(0)}_{1/2}\MGcdot \begin{pmatrix}1&0\\0&i\end{pmatrix}.
\end{eqnarray*}
Now \eqref{11.28} will be proved. The basis $\uuuu\sigma_0$ and the flat 
basis $\uuuu e^\pm_\iiii$ on the sector $\whhh I^\pm_\iiii(x)$ are related by
\begin{eqnarray*}
J(\uuuu e^\pm_0(z)) &=& \uuuu e^\pm_\iiii(\rho_c(z))\MGcdot
\begin{pmatrix}1&0\\0&-1\end{pmatrix},\\
J(\uuuu\sigma_0(z)) &=& \uuuu\sigma_0(\rho_c(z))\MGcdot
\begin{pmatrix}0&f_0^{-1}\\f_0&0\end{pmatrix}
=\uuuu\sigma_0(\rho_c(z))\MGcdot
\begin{pmatrix}0&e^{-\varphi/2}\\e^{\varphi/2}&0\end{pmatrix},\\
\textup{ so} &&
\uuuu\sigma_0(z)\begin{pmatrix}0&e^{-\varphi/2}\\e^{\varphi/2}&0\end{pmatrix}
=J(\uuuu\sigma_0(\rho_c(z)))\\
&=& e^\pm_\iiii(z)\MGcdot \begin{pmatrix}1&0\\0&-1\end{pmatrix}
\begin{pmatrix}e^{i/\rho_c(z)}&0\\0&e^{-i/\rho_c(z)}\end{pmatrix}
\MGcdot A^\pm_0(\rho_c(z))\MGcdot C\\
&=& \uuuu e^\pm_\iiii(z)\begin{pmatrix}1&0\\0&-1\end{pmatrix}
\begin{pmatrix}e^{-ix^2z}&0\\0&e^{ix^2z}\end{pmatrix}\MGcdot 
A^\pm_0(-{1}/{x^2z})\MGcdot C.
\end{eqnarray*}
The basis $\uuuu\rho$ and the flat basis $\uuuu\rho\MGcdot\Psi^{(\iiii)}_k$
are related by
\begin{eqnarray*}
\uuuu\rho =\uuuu\rho\MGcdot\Psi^{(\iiii)}_k\MGcdot (\Psi^{(\iiii)}_k)^{-1}
=\uuuu\rho\MGcdot\Psi^{(\iiii)}_k\MGcdot
\begin{pmatrix}e^{ix^2z}&0\\0&e^{-ix^2z}\end{pmatrix}\MGcdot A^{(\iiii)}_k(z)^{-1}
\end{eqnarray*}
with $A^{(\iiii)}_k\in GL(2,\AAA|_{\Omega^{(\iiii)}_k\cap S^1}),\ 
A^{(\iiii)}_k(\iiii)={\bf 1}_2$. Therefore (with $+\leftrightarrow 2,
-\leftrightarrow 1$)
\begin{eqnarray*}
\uuuu e^\pm_\iiii &=& 
\uuuu \rho\MGcdot\Psi^{(k)}_{2/1}
\begin{pmatrix}e^{ix^2z}&0\\0&e^{-x^2z}\end{pmatrix}
A^{(\iiii)}_k(z)^{-1}\MGcdot (-i)
\begin{pmatrix}e^{\varphi/2}&-e^{-\varphi/2}\\e^{\varphi/2}&e^{-\varphi/2}
\end{pmatrix}\\
&& \quad\MGtimes\begin{pmatrix}0&e^{-\varphi/2}\\e^{\varphi/2}&0\end{pmatrix}
C^{-1}\MGcdot A^\pm_0(-{1}/{x^2z})^{-1}
\begin{pmatrix} e^{ix^2z}&0\\0&e^{-ix^2z}\end{pmatrix}
\begin{pmatrix} 1&0\\0&-1\end{pmatrix}\\
&=& 
\uuuu \rho\MGcdot\Psi^{(k)}_{2/1}
\begin{pmatrix}e^{ix^2z}&0\\0&e^{-x^2z}\end{pmatrix}
A^{(\iiii)}_k(z)^{-1}\MGcdot 
\begin{pmatrix}0&-1\\-i&0\end{pmatrix}\\
&& \quad\MGtimes A^\pm_0(-{1}/{x^2z})^{-1}
\begin{pmatrix} e^{ix^2z}&0\\0&e^{-ix^2z}\end{pmatrix}
\begin{pmatrix} 1&0\\0&-1\end{pmatrix}\\
&=& 
\uuuu \rho\MGcdot\Psi^{(k)}_{2/1}
\MGcdot \begin{pmatrix}0&-1\\-i&0\end{pmatrix}
\MGcdot \begin{pmatrix} 1&0\\0&-1\end{pmatrix}
=
\uuuu \rho\MGcdot\Psi^{(k)}_{2/1}
\MGcdot \begin{pmatrix}0&1\\-i&0\end{pmatrix}.
\end{eqnarray*}
Now, $S^{(0)}_1=\begin{pmatrix} 1 & \zeta^{NI}\\0&1\end{pmatrix}$
is the Stokes matrix with 
$$
(\uuuu\rho\MGcdot \Psi^{(0)}_1)|_{I^a_0}\MGcdot S^{(0)}_1 
= (\uuuu\rho\MGcdot \Psi^{(0)}_2)|_{I^a_0},
$$
and $S^a_0=\begin{pmatrix} 1 & s\\0&1\end{pmatrix}$
is the Stokes matrix with 
$$
\uuuu e^+_{0}|_{I^a_0}\MGcdot S^a_0 
= \uuuu e^-_{0}|_{I^a_0}.
$$
Thus
\[
\begin{pmatrix}1&-is\\0&1\end{pmatrix}
= \begin{pmatrix}1&0\\0&i\end{pmatrix}
S^a_0\begin{pmatrix} 1&0\\0&-i\end{pmatrix}
= S^{(0)}_1=\begin{pmatrix}1 &\zeta^{NI}\\0&1\end{pmatrix},
\textup{ so }s=i\zeta^{NI}.
\]

For $x\in S^{NI}$, a $\beta$ with $e^{-\beta}=4x^2$ \eqref{2.23}
and $\Im(\beta)\in(-\frac{\pi}{2}+2\varepsilon,\frac{\pi}{2}-2\varepsilon)$
encodes a path $[\beta]$ from 1 to $e^\beta$ which stays within
$\Omega^{(\iiii)}_1\subset \whhh I^+_0\cap\whhh I^-_\iiii(x)$.
Therefore $B(\beta)$ is the matrix with
$$\uuuu e^-_{\iiii}|_{\Omega^{(\iiii)}_1} = 
\uuuu e^+_{0}|_{\Omega^{(\iiii)}_1}\MGcdot B(\beta).$$
On the other hand
$$\uuuu\rho\MGcdot\Psi^{(\iiii)}_1 = \uuuu\rho\MGcdot \Psi^{(0)}_1\MGcdot E^{NI}.$$
This establishes \eqref{11.30} and \eqref{11.31}

(c) \eqref{11.32} and \eqref{11.33} follow immediately from \eqref{11.31}.
They induce a natural bijection between the two copies of $\{\zeta^{NI}\in\C\}$
and the double cover for $\sqrt{1+pq}$ of the manifold
$\{(\zeta^{NI},p,q)\in\C^3\, |\, pq\neq -1,\zeta^{NI}=p+q\}$ on one side 
and $M_{3TJ}^{mon}(-i,ix^2)$ on the other side. 
Here two points $(s,b_1,b_2)$ and $(s,-b_1,-b_2)$ are mapped to the same
point $(\zeta,p,q)$ by \eqref{11.33}.
The symmetry $R_2$ from \eqref{7.30} and remark \ref{t8.3} maps
$(s,b_1,b_2,f_0)$ to $(s,-b_1,-b_2,-f_0)$.
\hfill$\Box$

\begin{remarks}\label{t11.5}
(i) In \cite{IN86} and \cite{FIKN06} the special case \eqref{11.21} 
is mentioned, but its role is not emphasized.

(ii) Chapters 8 and 11 of \cite{IN86} contain rich results on the solutions on $\R_{>0}$ of 
\eqref{11.1}. Using lemma \ref{t11.1} they can be translated into results on the solutions
of $P_{III}(0,0,4,-4)$ on $\R_{>0}$.

(iii) Chapter 8 of \cite{IN86} studies first the asymptotics as $x\to\iiii$
of solutions $f$   of $P_{III}(0,0,4,-4)$ on $\R_{>0}$ which satisfy 
\begin{eqnarray}\label{11.36}
\log f=O(x^{-1/2}),\quad \paa_x(\log f)=O(x^{-1/2})\quad
\textup{for }x\to\iiii.
\end{eqnarray}
In particular we have $f\to 1$ as $x\to\iiii$ for these solutions.
\cite[theorem 8.1 and remark 8.2]{IN86} state that a solution 
$f_{mult}(.,s,B)|_{\R_{>0}}$ satisfies \eqref{11.36}  in the case
\eqref{11.22} with $pq>-1$. This is equivalent to $b_2\in \R^*$.
The case $b_2=0$ is the case \eqref{11.21} and is not considered in 
\cite[ch.\  8]{IN86}. We expect that \eqref{11.36} holds there, too.

In \cite[ch.\  8]{IN86} the condition $pq>-1$ is split into four cases,
and two subcases of two of the cases are also given.
The following table lists their conditions and notation
(the manifolds $\Mm_r,\Mm_i,\Mm_p,\Mm_q,\Mm_\R,\Mm_I$) and the
translation into conditions for $s,b_1,b_2$:
\begin{eqnarray*}
\begin{array}{l|l|l}
\Mm_r & pq>0 & b_2\in]-1,1[-\{0\}\\
\Mm_i & -1<pq<0 & b_2\in \R_{<-1}\cup\R_{>1} \\
\Mm_p & q=0 & b_2=\pm 1,b_1+sb_2=0\\
\Mm_q & p=0 & b_2=\pm 1, b_1=0\\
\Mm_\R\subset \Mm_r & pq>0,|p|=|q| & s\in i\R,b_1+\frac{s}{2}b_2\in\R,b_2\in\R^*\\
\Mm_I\subset \Mm_i & -1\!<\!pq\!<\!0,|p|\!=\!|q| & s\in (-2,2),b_1\!+\!\tfrac12{s}b_2\in i\R,
\pm b_2\in\R_{>1}
\end{array}
\end{eqnarray*}
The manifold $V^{mat,S^1}-\{(s,B)\, |\, b_2=0\}$ from \eqref{15.15} is a double cover 
of $\Mm_\R$, and $V^{mat,i\R_{>0}}$ from \eqref{15.12} is isomorphic to $\Mm_I$.
For $(s,B)\in V^{mat,S^1}$, $f_{mult}(.,s,B)|_{\R_{>0}}$ takes values in $S^1$
and $u^{NI}$ has real values. For $(s,B)\in V^{mat,i\R_{>0}}$, 
$f_{mult}(.,s,B)|_{\R_{>0}}$ takes values in $i\R_{>0}$ and $u^{NI}$ takes values
in $i\R_{>0}$; see lemma \ref{15.2} and theorem \ref{t15.5} below
and also \cite[ch.\  8]{IN86}. 

For the four cases $\Mm_r,\Mm_i,\Mm_p,\Mm_q$ \cite[theorem 8.1]{IN86} gives 
precise asymptotic formulae for $u^{NI}$ (and thus for $\log f$)
as $x\to\iiii$.

In the later part of \cite[ch.\  8]{IN86}, the monodromy data $(s,B)$
of a solution $f_{mult}(.,s,B)|_{\R_{>0}}$ which is smooth near $0$ is connected
with the asymptotics near $0$. This is generalized in \cite{Ni09},
and we reformulate it in chapter \ref{s12} and reprove it in chapter \ref{s13}.
Further comments on \cite[ch.\  8]{IN86} are made in remark \ref{t15.6} (iii)
below.

(iv) \cite[ch.\  11]{IN86} studies the singularities of real solutions $\varphi$
on $\R_{>0}$ of \eqref{9.5} and thus implicitly the zeros and poles
of real solutions  of $P_{III}(0,0,4,-4)$ on $\R_{>0}$. 
Again they restrict to the case \eqref{11.22}, which is precisely the case
of real solutions  of $P_{III}(0,0,4,-4)$ on $\R_{>0}$ which are
not smooth near $\iiii$. We take this up in chapters \ref{s15}
($V^{mat,\R}$ in \eqref{15.6}) and \ref{s18} (i.e.\ which solutions have which
sequences of zeros and poles).
However we do not reformulate formula \cite[(11.10)]{IN86} for the 
positions of zeros and poles of a solution for large $x$.
\end{remarks}

\chapter{Asymptotics of all solutions near $0$}\label{s12}
\setcounter{equation}{0}

\noindent
In this chapter we shall rewrite and extend one of the two main results
of \cite{Ni09}, the asymptotic formulae
as $x\to 0$ for all solutions of $P_{III}(0,0,4,-4)$.
In \cite[1.4.2]{Ni09},
Niles distinguishes
three cases $a$, $b$ and $c$ and makes implicitly the following finer
separation into five cases $a$, $b+$, $b-$, $c+$, $c-$. Let $\C^{[sto]}$ be the complex
plane with coordinate $s$, and define
\begin{eqnarray}\label{12.1}
\begin{split}
\C^{[sto,a]}:=& \ \C -(\R_{\leq -2}\cup\R_{\geq 2}),\\
\C^{[sto,b\pm]}:=& \ (\pm 1)\R_{>2},
\quad \C^{[sto,b]}:=\C^{[sto,b+]}\cup\C^{[sto,b-]},\\
\C^{[sto,c\pm]}:=& \ \{\pm 2\},
\quad \C^{[sto,c]}:=\C^{[sto,c+]}\cup\C^{[sto,c-]},\\
%\C^{[sto,J_1\cup \dots\cup J_n]}
%:=& \C^{[sto,J_1]}\cup\dots\cup\C^{[sto,J_n]}\\
%& \textup{ for } J_1,\dots,J_n\in\{a,b\pm,b,c\pm,c\},
\C^{[sto,J_1\cup J_2\cup\dots]}
:=& \ \C^{[sto,J_1]}\cup\C^{[sto,J_2]}\cup\dots
\ (J_1,J_2,\dots\!\in\!\{a,b\pm,b,c\pm,c\}),
\end{split}
\end{eqnarray}
and 
\begin{eqnarray}\label{12.2}
\begin{split}
V^{mat,J}:=& \{(s,B)\in V^{mat}\, |\, s\in \C^{[sto,J]}\}\\
&\textup{ for }J\in\{a,b\pm,b,c\pm,c,a\cup b+,\dots\}.
\end{split}
\end{eqnarray}
For each of cases $a$, $b$ and $c$, 
Niles has one formula for the special case
\eqref{11.21} and one formula for the general case \eqref{11.22}.
The formulae for the cases $b$ and $c$ contain a sign which distinguishes 
$b+,b-$ and $c+,c-$.

We shall give a refined version \eqref{12.20} of the formula for the case $a$.
Our versions \eqref{12.24} and \eqref{12.28} of Niles' formulae
for the cases $b+$ and $c+$ will follow from \eqref{12.20} 
by analytic continuation.
Our versions comprise the special case \eqref{11.21} and the general case
\eqref{11.22}. For the proof we only need Niles' formula for the case $a$
in the general case \eqref{11.22} and a basic property of his formula
for the case $b+$ in the general case \eqref{11.22}.

Remarks \ref{t12.1}-\ref{t12.3} prepare the formulae 
in theorem \ref{t12.4} and make them more transparent.

\begin{remarks}\label{t12.1}
Lemma \ref{t5.2} (b) gives analytic isomorphisms
(real analytic with respect to $s$ in the case $b$)
\begin{eqnarray}\label{12.3}
\begin{split}
V^{mat,a}\to \C^{[sto,a]}\times \C^*,\quad (s,B)\mapsto (s,b_-),\\
V^{mat,b}\to \C^{[sto,b]}\times \C^*,\quad (s,B)\mapsto (s,b_-),
\end{split}
\end{eqnarray}
so $(s,b_-)$ are global coordinates on $V^{mat,a}$ as an analytic manifold.
The fibration $V^{mat,a\cup b}\to \C^{[sto,a\cup b]}$ is an analytically locally
trivial fibre bundle with fibres isomorphic to $\C^*$. 
Recall the definition of $\sqrt{\tfrac14{s^2}-1}$ in chapter \ref{s5}.
The first row in the following table lists some holomorphic functions
on $V^{mat,a}$.  The functions on $V^{mat,a}$ in the second row are obtained from these by analytic continuation over $V^{mat,b+}$.
\begin{eqnarray}\label{12.4}
\renewcommand{\arraystretch}{1.9}
\begin{array}{c|c|c|c|c}
s & \sqrt{\tfrac14{s^2}-1} & \alpha_\pm & \lambda_\pm & b_\pm\\ \hline
s & -\sqrt{\tfrac14{s^2}-1} & \alpha_\mp\mp 1 & \lambda_\mp & b_\mp
\end{array}
\end{eqnarray}
In future, we denote a holomorphic function on $\C^{[sto,a]}$
by $\kappa(s,\sqrt{\tfrac14{s^2}-1})$, and then denote the function on 
$\C^{[sto,a]}$  obtained by analytic continuation over
$\C^{[sto,b+]}$ by $\kappa(s,-\sqrt{\tfrac14{s^2}-1})$
\end{remarks}

\begin{remarks}\label{t12.2}
(i) Recall the action $m_{[1]}$ on $\C\times V^{mat}$ (with coordinates
$(\xi,s,B)$ with $-\frac{\beta}{2}=\xi+\log 2$) in \eqref{10.11}:
\begin{eqnarray*}
m_{[1]}:\C\times V^{mat}\to\C\times V^{mat},\quad
(\xi,s,B)\mapsto (\xi-i\pi ,s,(\Mon_0^{mat})^{-1}\MGcdot B).
\end{eqnarray*}
The functions in the first row of the table \eqref{12.4} are with the exception of $b_\pm$
functions on $\C\times V^{mat,a}$ which are invariant under the action
of $m_{[1]}$. As $b_-$ is mapped under $m_{[1]}$ to $\lambda_-^{-1}\MGcdot b_-$,
one obtains the function
\begin{eqnarray}\label{12.5}
\www b_-:= e^{2(\xi-\log 2)\alpha_-}\MGcdot b_-
=\left(\frac{x}{2}\right)^{2\alpha_-}\MGcdot b_-
\end{eqnarray}
which is invariant under $m_{[1]}$ 
(here $x$ comes with the choice of $\xi=\log x$).
The summand $-\log 2$
is inserted in order to simplify the formulae in theorem \ref{t12.4}.
We obtain analytic isomorphisms (real analytic with respect to $s$ 
in case $b$)
\begin{eqnarray}\label{12.6}
\C\times V^{mat,a}/\langle m_{[1]}\rangle\to \C^*\times\C^{[sto,a]}\times \C^*,
\ [(\xi,s,B)]\mapsto (x^2,s,\www b_-),\\
\C\times V^{mat,b}/\langle m_{[1]}\rangle\to \C^*\times\C^{[sto,b]}\times \C^*,
\ [(\xi,s,B)]\mapsto (x^2,s,\www b_-),
\nonumber
\end{eqnarray}
so $(x^2,s,\www b_-)$ are global coordinates on 
$\C\times V^{mat,a}/\langle m_{[1]}\rangle$. 
The fibration $\C\times V^{mat,a\cup b}/\langle m_{[1]}\rangle 
\to \C^{[sto,a\cup b]}$ is an analytically locally
trivial fibre bundle with fibres isomorphic to $\C^*\times\C^*$. 

Let us erase in table \eqref{12.4} the column with $b_\pm$ and extend the table
by 
\begin{eqnarray}\label{12.7}
\renewcommand{\arraystretch}{1.7}
\begin{array}{c|c}
x^2 & \www b_-\\ \hline
x^2 & \left(\frac{x}{2}\right)^2\www b_-^{-1}
\end{array}
\end{eqnarray}
Then the first row lists some holomorphic functions
on $\C\times V^{mat,a}/\langle m_{[1]}\rangle$, and 
the functions on $\C\times V^{mat,a}/\langle m_{[1]}\rangle$ in the second row are obtained from these by analytic continuation over $\C\times V^{mat,b+}/\langle m_{[1]}\rangle$.

If now an open neighbourhood $U$ in 
$\C\times V^{mat,a\cup b_+}/\langle m_{[1]}\rangle$ of a point in 
$\C\times V^{mat,c+}/\langle m_{[1]}\rangle$
and a holomorphic function in this neighbourhood are given,
then the function has a convergent Laurent expansion
\begin{eqnarray}\label{12.8}
\sum_{(a_1,a_2)\in\Z^2}\kappa_{a_1,a_2}(s,\sqrt{\tfrac14{s^2}-1})\MGcdot
\left(\tfrac12{x}\right)^{2a_1}\MGcdot (\www b_-)^{a_2},
\end{eqnarray}
where $\kappa_{a_1,a_2}$ are coefficients holomorphic in 
$s\in \C^{[sto,a]}\cap \pr_s(U)$.
The Laurent expansion must stay the same after continuation over
$\C\times V^{mat,b+}/\langle m_{[1]}\rangle$. Therefore
the coefficients must satisfy 
\begin{eqnarray}\label{12.9}
\kappa_{a_1+a_2,-a_2}(s,\sqrt{\tfrac14{s^2}-1})
=\kappa_{a_1,a_2}(s,-\sqrt{\tfrac14{s^2}-1}).
\end{eqnarray}

(ii) We also need  coordinates on $\C\times V^{mat,c\pm}/\langle m_{[1]}\rangle$.
Suppose $s\in\{\pm2\}$. Then by lemma \ref{t5.2} (b), 
$\www b_1:=b_1+\tfrac12{s}b_2\in\{\pm 1\}$
and $B\in V^{mat}$ is determined by $(\www b_1,b_2)\in\{\pm 1\}\times\C$ 
and any such pair is realized by some matrix $B$.
Since
\begin{eqnarray}\label{12.10}
(\Mon_0^{mat})^{-1} \!\! \MGcdot \begin{pmatrix} b_1 & b_2\\ -b_2 & b_1+sb_2\end{pmatrix}
\!=\! \begin{pmatrix} \! b_1-s(b_2+sb_1) \! & \! b_2+sb_1 \!\\-(b_2+sb_1) & b_1\end{pmatrix}\!,
\end{eqnarray}
$m_{[1]}$ maps the functions in the first row to the functions in the second
row,
\begin{eqnarray}\label{12.11}
\renewcommand{\arraystretch}{1.7}
\begin{array}{c|c|c|c|c|c}
b_1 & b_2 & \www b_1 & \xi & x\www b_1 & 
\xi\!-\!\frac{i\pi}{4}s\www b_1 b_2=:\www b_2
\\ \hline
\!b_1\!-\!s(b_2+sb_1)\! & \!s\www b_1 \!-\! b_2\! & \!-\www b_1 \!& \!\xi\!-\!i\pi \!& \!x\www b_1 \!& 
\!\xi\!-\!\frac{i\pi}{4}s\www b_1 b_2=\www b_2\!
\end{array}
\end{eqnarray}
In particular, $x\www b_1$ and $\www b_2$ are invariant functions
and serve as global coordinates on 
$\C\times V^{mat,c\pm}/\langle m_{[1]}\rangle$. The map 
\begin{eqnarray}\label{12.12}
\C\times V^{mat,c\pm}/\langle m_{[1]}\rangle\to \C^*\times\C,\quad
[(\xi,s,B)]\mapsto (x\www b_1,\www b_2)
\end{eqnarray}
is an analytic isomorphism.
\end{remarks}

\begin{remark}\label{t12.3}
The group $G^{mon}=\{\id,R_1,R_2,R_3\}$ 
defined in \eqref{7.30}
(with $R_3=R_1\circ R_2=R_2\circ R_1$)
acts on $V^{mat}\cong V^{mon}$, $\C\times V^{mat}$,
$\C\times V^{mat}/\langle m_{[1]}\rangle$ and on $M_{3FN}^{mon}$
and respects the foliation there. Remark \ref{t8.3} fixes how it acts
on the functions $f_0,g_0$ on $M_{3FN}^{reg}$. 
The action of $m_{[1]}$ on $\C\times V^{mat}$ is fixed in remark \ref{t10.2} (iv).
The following table lists some meromorphic functions on 
$\C\times V^{mat,a}$ and the actions of $R_1$, $R_2$ and $m_{[1]}$ 
on them.
\begin{eqnarray}\label{12.13} 
\renewcommand{\arraystretch}{1.3}
\begin{array}{c|c|c|c}
 & R_1 & R_2 & m_{[1]} \\ \hline
s & -s & s & s \\
\lambda_\pm & \lambda_\mp & \lambda_\pm & \lambda_\pm \\
\alpha_\pm & -\alpha_\pm & \alpha_\pm & \alpha_\pm \\
B & \left(\begin{smallmatrix}1&0\\0&-1\end{smallmatrix}\right)B
\left(\begin{smallmatrix}1&0\\0&-1\end{smallmatrix}\right)&
-B & \Mon_0^{mat}(s)^{-1}\MGcdot B \\
b_1 & b_1 & -b_1 & (1-s^2)b_1-sb_2 \\
b_2 & -b_2 & -b_2 & b_2+sb_1 \\
b_1+\frac{s}{2}b_2 & b_1+\frac{s}{2}b_2 & -b_1-\frac{s}{2}b_2 & 
(1-\frac{s^2}{2})b_1-\frac{s}{2}b_2 \\
b_\pm & b_\mp & -b_\pm & \lambda_\pm^{-1}b_\pm \\
\xi & \xi & \xi & \xi -i\pi \\
x & x & x & -x \\
f_{univ} & f_{univ}^{-1} & -f_{univ} & f_{univ} \\
g_{univ} & -g_{univ} & g_{univ} & g_{univ}
\end{array}
\end{eqnarray}
\end{remark}

The core of theorem \ref{t12.4} is
the asymptotic formulae in \cite[1.4.2 and ch. 3]{Ni09}
for $x$ near $0$ of the solutions of the $P_{III}(0,0,4,-4)$ equation.
In theorem \ref{t12.4} they are extended, simplified and made more transparent.

\begin{theorem}\label{t12.4}
There exist two continuous $m_{[1]}$-invariant and $G^{mon}$-invariant maps
\begin{eqnarray}\label{12.14}
B_1:i\R\times V^{mat}\to \R_{>0},\quad
B_2:i\R\times V^{mat,a}\to \R_{>0}
\end{eqnarray}
such that $B_2\leq B_1$ (where both are defined)
and such that the $m_{[1]}$-invariant and $G^{mon}$-invariant open subsets
\begin{eqnarray}\label{12.15}
\begin{split}
U_1:=&\ \{(\xi,s,B)\in \C\times V^{mat}\, |\, |e^\xi|<B_1(i\Im(\xi),s,B)\}\\
U_2:=&\ \{(\xi,s,B)\in \C\times V^{mat,a}\, |\, |e^\xi|<B_2(i\Im(\xi),s,B)\}
\end{split}
\end{eqnarray}
satisfy the following:

(a) $f_{univ}$ is holomorphic and invertible on $\C\times V^{mat,a}\cap U_2$
and holomorphic on $\{(\xi,s,B)\in\C\times V^{mat}\cap U_1\, |\, \Re(s)>-1\}$.
$f^{-1}_{univ}$ is holomorphic on 
$\{(\xi,s,B)\in\C\times V^{mat}\cap U_1\, |\, \Re(s)<1\}$.
The restriction of $f_{univ}$ to 
$\C\times V^{mat,c}\cap U_1$ is holomorphic and invertible.

(b) $f_{mult}(.,s,B)$ is holomorphic near $0$ for 
$(s,B)\in V^{mat,a\cup b+\cup c}$ and invertible near $0$ for 
$(s,B)\in V^{mat,a\cup c}$. $f_{mult}^{-1}(.,s,B)$ is holomorphic near $0$
for $(s,B)\in V^{mat,a\cup b-\cup c}$
((b) rewrites (a) with $f_{mult}(.,s,B)$, but is not precise about 
$U_1$ and $U_2$).

(c) $f_{univ}$ is on $\C\times V^{mat,a\cup b+}\cap U_1$ a convergent sum
\begin{eqnarray}\label{12.16}
\sum_{(a_1,a_2)\in L}\kappa_{a_1,a_2}(s,\sqrt{\tfrac14{s^2}-1})\MGcdot
\left(\tfrac12{x}\right)^{2a_1}\MGcdot (\www b_-)^{a_2},
\end{eqnarray}
where $L\subset\Z^2$ is defined by the three inequalities
\begin{equation}\label{12.17}
2a_1-a_2\geq -1,\ 
2a_1+a_2\geq 1,\ 
2a_1+3a_2\geq -1,
\end{equation}
and where $\kappa_{a_1,a_2}$ are coefficients holomorphic in 
$s\in \C^{[sto,a]}\cap U$ which satisfy \eqref{12.9}.
In particular
\begin{eqnarray}\label{12.18}
\kappa_{0,1}=\frac{\Gamma(\frac{1}{2}-\alpha_-)}{\Gamma(\frac{1}{2}+\alpha_-)},\quad
\kappa_{1,-1}=\frac{\Gamma(\frac{1}{2}-(\alpha_++1))}
{\Gamma(\frac{1}{2}+(\alpha_++1))}.
\end{eqnarray}

(d) Recall that as multi-valued functions in $x$
\begin{eqnarray}\label{12.19}
\begin{split}
\www b_-&=\left(\tfrac12{x}\right)^{2\alpha_-}\MGcdot b_-,\\
\left(\tfrac12{x}\right)^2\MGcdot (\www b_-)^{-1}&=
\left(\tfrac12{x}\right)^{2-2\alpha_-}\MGcdot b_-^{-1}.
\end{split}
\end{eqnarray}
If $(s,B)\in V^{mat,a}$ with $\Re(s)\geq 0$ 
($\iff \Re(\alpha_-)\in[0,\tfrac{1}{2})$) then for small $x$
\begin{eqnarray}\label{12.20}
f_{mult}(x,s,B)&=& 
\kappa_{0,1}\MGcdot\www b_-
+ \kappa_{1,-1}\MGcdot \left(\tfrac12{x}\right)^2\MGcdot (\www b_-)^{-1}
+O(|x|^2)\\
&=& \kappa_{0,1}\MGcdot\www b_-
+O(|x|^{2-2\Re(\alpha_-)}).\label{12.21}
\end{eqnarray}
If $(s,B)\in V^{mat,a}$ with $\Re(s)\leq 0$ 
($\iff \Re(\alpha_-)\in(-\frac{1}{2},0]$) then for small $x$
\begin{eqnarray}\label{12.22}
f_{mult}(x,s,B)&=& 
\kappa_{0,1}\MGcdot\www b_-
+O(|x|^{2+6\Re(\alpha_-)}).
\end{eqnarray}
If $(s,B)\in V^{mat,b+}$ ($\iff \Re(\alpha_-)=\frac{1}{2},s\neq 2$) then
for small $x$
\begin{eqnarray}\label{12.23}
f_{mult}(x,s,B)&=& 
\kappa_{0,1}\MGcdot\www b_-
+ \kappa_{1,-1}\MGcdot \left(\frac{x}{2}\right)^2\MGcdot (\www b_-)^{-1}
+O(|x|^2)\\
&=& -\frac{x}{t^{NI}}\MGcdot
\sin\left(2t^{NI}\log\frac{x}{2}-2\arg\Gamma(1+it^{NI})+\delta^{NI}\right)
\nonumber\\
&&+O(|x|^2),\label{12.24}
\end{eqnarray}
where $e^{i\delta^{NI}}=b_-$ and $\Re(\delta^{NI})\in[0,2\pi)$,
and $t^{NI}$ is as in \eqref{5.17}.

(e) Part (a) and 
formula \eqref{12.24} show that $f_{mult}(.,s,B)$ has zeros arbitrarily
close to $0$ within $x$ with bounded argument $\arg(x)=\Im(\xi)$ if and only
if $s\in \R_{>2}$, and that then the only such zeros are approximately
given by 
\begin{eqnarray}\label{12.25}
\begin{split}
x_k=& 2 \exp\left(\frac{1}{2t^{NI}}(2\arg\Gamma(1+it^{NI})-\delta^{NI})\right)
\MGcdot \exp\left(\frac{-k\pi}{2t^{NI}}\right)\\
&\quad\MGtimes (1+ O\left(\exp\left(\frac{-k\pi}{2t^{NI}}\right)\right))
\quad\textup{ for }k\in\Z_{>\textup{a bound}}.
\end{split}
\end{eqnarray}
In particular, they all have approximately the same argument, which is
\begin{eqnarray}\label{12.26}
\frac{\Im(-\delta^{NI})}{2t^{NI}}=\frac{\log |b_-|}{2 t^{NI}} .
\end{eqnarray}

(f) Using $R_1$ in remark \ref{12.3}, one can derive from (d) and (g) 
asymptotic formulae for $f_{mult}^{-1}(.,s,B)$ for 
$(s,B)\in V^{mat,a\cup b-\cup c-}$. Using $R_1$ and (e), one finds
that $f_{mult}(.,s,B)$ has poles near $0$ if and only if $s\in \R_{<-2}$
and that then there is only the sequence of poles tending to $0$ in \eqref{12.25},
where now $t^{NI}=t^{NI}(|s|)$ and $\delta^{NI}=\delta^{NI}(-s,b_1,-b_2)
=-\delta^{NI}(s,b_1,b_2)$.

(g) $f_{univ}$ is on $\C\times V^{mat,c+}\cap U_1$ holomorphic and invertible
with leading term (leading with respect to $x\to 0$ with constant argument)
\begin{eqnarray}\label{12.27}
-2e^\xi \www b_1\left(\www b_2 -\log 2 + \gamma_{Euler}\right).
\end{eqnarray}
Consider $(s,B)\in V^{mat,c+}$, i.e.\ $s=2$.
Then 
\begin{equation}\label{12.28}
f_{mult}(x,s,B) \!=\!  -2x\www b_1 \!\left( \log\frac{x}{2}-
\frac{i\pi}{2}\www b_1 b_2 +\gamma_{Euler}\right)
\!+\! O(|x|^2).
\end{equation}
\end{theorem}

{\bf Proof:}
Niles has six formulae: for each of cases $a$, $b$ and $c$, one in the 
special case \eqref{11.21} and one in the general case \eqref{11.22}.
We shall use only the formula for case $a$ 
in the general case \eqref{11.22},
and a basic property of the formula for case $b+$
in the general case \eqref{11.22}. The basic property is that 
$f_{mult}(.,s,B)$ for $(s,B)\in V^{mat,b+}$ is holomorphic 
for $x\in S^{NI}$ near $0$.
First we shall treat the cases $a$ and $b$, that is, (a)-(f)
without the statements in (a),(b) on the case $c$, 
then the case $c$, that is, (g) and the rest of (a),(b).

Two more remarks on Niles' formulae are appropriate. 
First, he restricts attention to $x\in S^{NI}$, 
but remarks at the end of \cite{Ni09}
that similar formulae should hold for all $x\in \C^*$ near 0.
In fact, the same formulae hold for all $x\in\C^*$ near 0. 
This is contained in the proof below.
Second, Niles' formulae start with the monodromy data $(\zeta^{NI},p,q)$.
Lemma \eqref{11.4} (c) shows that one such triple corresponds to two triples
$(s,b_1,b_2)$ and $(s,-b_1,-b_2)=R_2(s,b_1,b_2)$ and to two solutions 
$f_0$ and $-f_0$ of $P_{III}(0,0,4,-4)$. Thus his formulae are not specific
about a global sign. We shall fix the sign using remark \ref{t10.5}.

Niles' formula for case $a$ in the general case \eqref{11.22} is as follows
\cite[1.4.2 and 3.1 theorem 7]{Ni09}: for $(s,B)\in V^{mat,a}$ and 
for $x\in S^{NI}$, any branch of the multi-valued function
$u(x)$ with \eqref{11.5} $f_{mult}(x,s,B)=ie^{-iu(x)/2}$ 
(here a priori the distinguished branch of $f_{mult}(.,s,B)$
should be taken) satisfies for $x\to 0$ (with constant argument)
\begin{eqnarray}\label{12.29}
u(x)&=&r^{NI}\log (4x)+s^{NI}+O(|x|^{2-|\Im(r^{NI})|}).
\end{eqnarray}
The constants $r^{NI}$ and $s^{NI}$ are as follows. $r^{NI}$ is given by 
\begin{eqnarray}
r^{NI}&=& 4i\alpha^{NI},\quad\textup{with }
\alpha^{NI}\textup{ given by}\nonumber\\
\zeta^{NI}&=& -2i\sin(\pi\alpha^{NI}),\ 
\Re(\alpha^{NI})\in(-\tfrac{1}{2},\tfrac{1}{2}).\nonumber 
\end{eqnarray}
\eqref{11.29} $s=i\zeta^{NI}$ and \eqref{5.14} $\sin(\pi \alpha_-)=
\frac12{s}$
show $\alpha^{NI}=\alpha_-$. Thus 
\begin{eqnarray}\label{12.30}
r^{NI}&=& 4i\alpha_-.
\end{eqnarray}
 $s^{NI}$ in \eqref{12.29} is given in \cite[1.4.2 and 3.1 theorem 7]{Ni09} 
by the first line of the following
calculation, which make use of formulae in lemma \ref{t5.2} (b)
and lemma \ref{t11.4} (c),
\begin{eqnarray}
e^{is^{NI}}&=&2^{-3ir^{NI}}\frac{\Gamma^2(\frac{1}{2}-\frac{ir^{NI}}{4})}
{\Gamma^2(\frac{1}{2}+\frac{ir^{NI}}{4})}
\frac{1+p\zeta^{NI}-p^2}{(e^{-\pi r^{NI}/4}-p)^{-2}}\\
&=& 2^{12\alpha_-}\frac{\Gamma^2(\frac{1}{2}+\alpha_-)}
{\Gamma^2(\frac{1}{2}-\alpha_-)}
\frac{1+pq}{(e^{-\pi i\alpha_-}-p)^2}\nonumber \\
&=& 2^{12\alpha_-}\frac{\Gamma^2(\frac{1}{2}+\alpha_-)}
{\Gamma^2(\frac{1}{2}-\alpha_-)}
\frac{1}{b_2^2(e^{-\pi i\alpha_-}-\frac{i b_1}{b_2})^2}\nonumber\\
&=& 2^{12\alpha_-}\frac{\Gamma^2(\frac{1}{2}+\alpha_-)}
{\Gamma^2(\frac{1}{2}-\alpha_-)}
\frac{-1}{b_-^2}\nonumber
\end{eqnarray}
We obtain for $(s,B)\in V^{mat,a}$ and $x\in S^{NI}$ and the 
distinguished branch of $f_{mult}(.,s,B)$ for $x\to 0$
the asymptotic formula
\begin{eqnarray}
f_{mult}(x,s,B) &=& ie^{-iu(x)/2}\nonumber \\
&=& (4x)^{2\alpha_-}\MGcdot i\MGcdot e^{-is^{NI}/2} + O(|x|^{2-4|\Re(\alpha_-)|+2\Re(\alpha_-)})
\nonumber \\
&=&\varepsilon\MGcdot \frac{\Gamma(\frac{1}{2}-\alpha_-)}
{\Gamma(\frac{1}{2}+\alpha_-)}
\left(\tfrac12{x}\right)^{2\alpha_-}b_-\label{12.32}\\
&& + O(|x|^{2-4|\Re(\alpha_-)|+2\Re(\alpha_-)})\nonumber
\end{eqnarray}
with a sign $\varepsilon\in\{\pm 1\}$, which has yet to be determined 
and which is not fixed by the formulae in \cite{Ni09}.

Remark \ref{t10.5} gives $f_{mult}(x,0,{\bf 1}_2)=1$.
In the case $(s,B)=(0,{\bf 1}_2)$ one has $\alpha_-=0$, $b_-=1$.
Together with \eqref{12.32} this shows that the sign is $\varepsilon=1$.

A priori the asymptotic formula \eqref{12.32} for $f_{univ}$ is established only 
in the general case \eqref{11.22}, but by continuity it holds also
in the special case \eqref{11.21}.

Now \eqref{12.21} and \eqref{12.22} with $\kappa_{0,1}$ as in \eqref{12.18} 
are established for $x\in S^{NI}\subset \C^*$.
For some $U_2$ as in \eqref{12.15}, $f_{univ}$ restricted to 
$\{(\xi,s,B)\in \C\times V^{mat,a}\cap U_2\, |\, e^\xi\in S^{NI}\}$
is a convergent sum 
\begin{eqnarray}\label{12.33}
\sum_{(a_1,a_2)\in\www L}\kappa_{a_1,a_2}(s,\sqrt{\tfrac14{s^2}-1})
\MGcdot \left(\tfrac12{x}\right)^{2a_1}\MGcdot (\www b_-)^{a_2},
\end{eqnarray}
where $\www L\subset\Z^2$ is determined by the condition that
\begin{eqnarray}
&&\textup{for }\Re(\alpha_-)\in(-\tfrac{1}{2},\tfrac{1}{2}),
\ 2a_1+2a_2\MGcdot\Re(\alpha_-)\geq 2\Re(\alpha_-),\nonumber \\
\textup{equivalently:}&& \textup{for }r\in \{-\tfrac{1}{2},\tfrac{1}{2}\},
\ 2a_1+2a_2r\geq 2r,\nonumber\\
\textup{equivalently:}&& 2a_1-a_2\geq -1,\ 2a_1+a_2\geq 1.
\label{12.34}
\end{eqnarray}
But then it is convergent on the whole intersection
$\C\times V^{mat,a}\cap U_2$ for some $U_2$ as in \eqref{12.15}.

We do not need Niles' precise formula in \cite[1.4.2 and 3.2 theorem 8]{Ni09}
for the case $b+$ in the general case \eqref{11.22}, but just the following
basic property, which is a consequence of that formula: 
$f_{univ}$ is holomorphic (without poles) (or real analytic
with respect to $s\in\C^{[sto,b+]}$) on
$S^{NI}\times V^{mat,b+}\cap U_1$ for an open set $U_1$ as in \eqref{12.15}. 

Then the Laurent expansion is valid there, i.e.\  the Laurent
expansion \eqref{12.33} must stay constant after continuation
over $S^{NI}\times V^{mat,b+}\cap U_1$. Thus it satisfies \eqref{12.9}.
Therefore if $(a_1,a_2)\in\www L$ then also $(a_1+a_2,-a_2)\in\www L$.
Replacing $(a_1,a_2)$ by $(a_1+a_2,-a_2)$ 
in the two inequalities in \eqref{12.34} 
gives only the one new inequality $2a_1+3a_2\geq -1$.
Thus one can replace
$\www L$ by $L$ 
in the Laurent expansion \eqref{12.33}
to obtain the Laurent expansion \eqref{12.16}.
Now the convergence includes also $(\C-S^{NI})\times V^{mat,b+}\cap U_1$
for a suitable $U_1$ as in \eqref{12.15}.  Part (c) is proved.

Thus $f_{univ}$ is holomorphic on a neighbourhood of 
$\C\times V^{mat,b+}\cap U_1$ in $\C\times V^{mat,a\cup b+}\cap U_1$.
This neighbourhood can be chosen for example as 
$\{(\xi,s,B)\in \C\times V^{mat,a\cup b+}\cap U_1\, |\, \Re(s)>-1\}$,
because \eqref{12.21} and \eqref{12.22} show that $f_{univ}$
is holomorphic and invertible on $\C\times V^{mat,a}\cap U_2$
for a set $U_2$ as in \eqref{12.15}.
Part (a),(b) are proved up to the statements on the case $c$.

The following two pictures show the most relevant part of $L$
and the graphs of the maps 
\[
(-\tfrac{1}{2},\tfrac{1}{2})\to\R,\quad r\mapsto a_1+a_2r
\]
for the most important points in $L$. 

%Later 2 pictures
\includegraphics[width=0.45\textwidth]{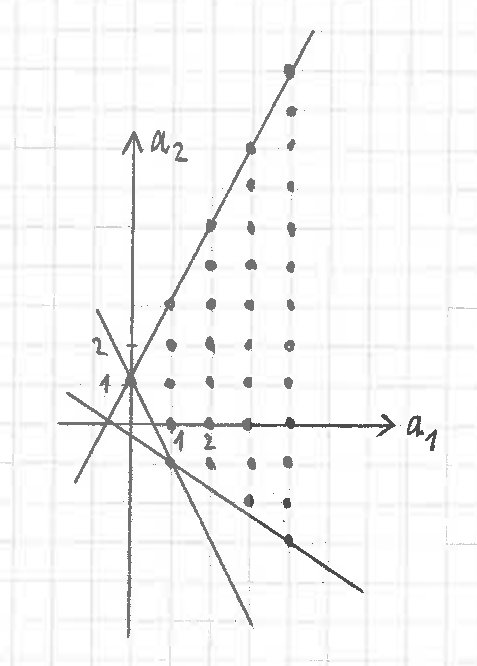} 
\includegraphics[width=0.5\textwidth]{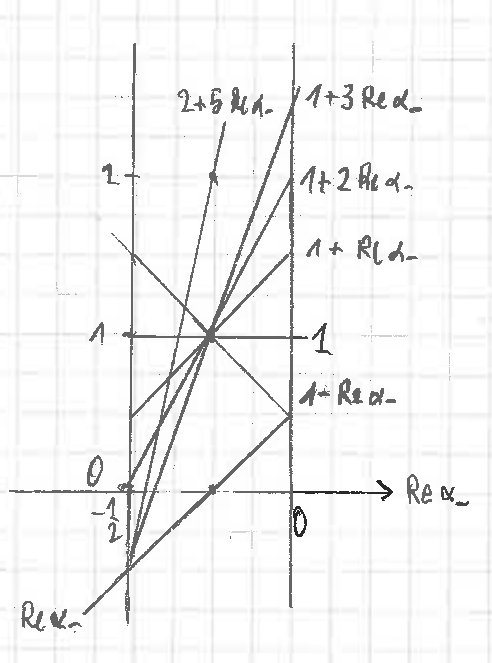} 

\eqref{12.20} and \eqref{12.23} follow immediately.
The next calculation derives \eqref{12.24} from \eqref{12.23}.
But of course \eqref{12.24} is equivalent to the formula
for the case $b+$ in the general case \eqref{11.22} in 
\cite[1.4.2 and 3.2 theorem 8]{Ni09}, up to the global sign, which is 
not fixed in \cite{Ni09}.
The calculation uses lemma \ref{t5.1} (b).
\begin{eqnarray}
&& \kappa_{0,1}\MGcdot\www b_-
+ \kappa_{1,-1}\MGcdot \left(\tfrac12{x}\right)^2\MGcdot (\www b_-)^{-1}\nonumber\\
&=& \frac{\Gamma(-it^{NI})}{\Gamma(1+it^{NI})}
\left(\tfrac12{x}\right)^{1+2it^{NI}}\MGcdot e^{i\delta^{NI}}
+ \frac{\Gamma(it^{NI})}{\Gamma(1-it^{NI})}
\left(\tfrac12{x}\right)^{1-2it^{NI}}\MGcdot e^{-i\delta^{NI}}\nonumber\\
&=& \frac{x}{2it^{NI}}\MGcdot\left(
-\frac{\Gamma(1-it^{NI})}{\Gamma(1+it^{NI})}
\MGcdot e^{2it^{NI}\log(x/2)}\MGcdot e^{i\delta^{NI}}\right. \nonumber\\
&&\hspace*{2cm} + \left. \frac{\Gamma(1+it^{NI})}{\Gamma(1-it^{NI})}
e^{-2it^{NI}\log(x/2)}\MGcdot e^{-i\delta^{NI}}\right)\nonumber\\
&=& -\frac{x}{t^{NI}}\MGcdot
\sin\left(2t^{NI}\log(x/2)-2\arg\Gamma(1+it^{NI})+\delta^{NI}\right) .
\label{12.35}
\end{eqnarray}
Part (d) is proved.

Now we turn to part (e). Denote by $x\MGcdot h(x,s,B)$ the remainder term of order
$O(|x|^2)$ in \eqref{12.24}, so $h(.,s,B)$ has order $O(|x|)$.
Then for 
$k$ sufficiently large,
$x_k$ is a zero of $f_{mult}(.,s,B)$ if 
\[
2t^{NI}\log (x_k/2)-2\arg\Gamma(1+it^{NI})+\delta^{NI}+k\pi
\]
is small and 
\[
0\!=\!\sin\left(2t^{NI}\log (x_k/2)\!-\!2\arg\Gamma(1\!+\!it^{NI})\!+\!\delta^{NI}\!+\!k\pi\right)
\!+\! t^{NI}\MGcdot h(x_k,s,B).
\]
This establishes \eqref{12.25}. \eqref{12.26} is obvious. Part (e) is proved.

Part (f) follows immediately from the statements on $R_1$ in remark \ref{t12.3}.

We turn to part (g) and the case $c+$.
$f_{univ}$ must be holomorphic on $\C\times V^{mat,c+}\cap U_1$:
by lemma \ref{t8.5} a hypersurface of poles in 
$\C\times V^{mat}$ would intersect $\C\times V^{mat,c+}\cap U_1$
in a codimension 1 set. But as there are no poles in 
$\C\times V^{mat,a\cup b+}\cap U_1$
for $\Re(s)>-1$, the hypersurface of poles would have codimension 2
which is absurd. 

An asymptotic formula for $f_{mult}(.,s,B)$ for 
$(s,B)\in V^{mat,c+}$ can be determined from \eqref{12.20} or \eqref{12.24}
by analytic (respectively, real analytic) continuation. 
We choose \eqref{12.24}, as the calculations are slightly shorter.
We just need formulae which connect some functions with the 
coordinates $(s,b_2)$ on $V^{mat}$ near $V^{mat,c+}$.

%Near a point $(s^0,b_2^0)=(2,b_2^0)\in V^{mat,c+}$, $(s,b_2)$ are coordinates
%on $V^{mat}$ because
%$$\frac{\paa b_1^2+b_2^2+sb_1b_2-1}{\paa {b_1}}=2b_1+sb_2
%\equiv 2\www b_1^0=2\MGcdot(\pm 1)\neq 0.$$
If $s\in\C^{[sto,b+]}=\R_{>2}$ is close to 2, then $t^{NI}$ is close to $0$,
and the left hand sides in the following table are approximated by the
right hand sides. Recall that $\Gamma'(1)=-\gamma_{Euler}$.
\begin{eqnarray}
\renewcommand{\arraystretch}{1.3}
\begin{array}{c|c}
\Gamma(1+it^{NI}) & 1-\gamma_{Euler}\MGcdot it^{NI}\\
\arg\Gamma(1+it^{NI}) & -\gamma_{Euler}\MGcdot t^{NI}\\ 
\sqrt{\tfrac14{s^2}-1} & \sqrt{s-2}\\
\lambda_- & -1 -2\sqrt{s-2}\\
t^{NI}=\frac{1}{2\pi}\log|\lambda_-| & \frac{1}{\pi}\sqrt{s-2}\\
b_- & \www b_1(1+\www b_1\sqrt{\tfrac14{s^2}-1}\,b_2)\\
b_- & \www b_1(1+\www b_1b_2\pi t^{NI})\\
\delta^{NI} & (1-\www b_1)\frac{\pi}{2}-i\www b_1b_2\pi t^{NI}
\end{array}\nonumber
\end{eqnarray}
Then \eqref{12.24} is close to 
\begin{eqnarray}
&&-\frac{x}{t^{NI}}\sin\left(2t^{NI}\log\frac{x}{2}+2t^{NI}\gamma_{Euler}
+(1-\www b_1)\frac{\pi}{2}-i\www b_1b_2\pi t^{NI}\right)\nonumber\\
&= &
-\frac{x\www b_1}{t^{NI}}
\sin\left(t^{NI}\left(2\log\frac{x}{2}+2\gamma_{Euler}
-i\pi\www b_1b_2\right) \right)\nonumber\\
&\approx & 
-2x\www b_1\left(\log\frac{x}{2}+\gamma_{Euler}
-\frac{i\pi}{2}\www b_1b_2\right).\label{12.36}
\end{eqnarray}
This proves \eqref{12.27} and part (g). 
Of course, \eqref{12.27} is equivalent to the formula
for the case $c+$ in the general case \eqref{11.22} in 
\cite[1.4.2 and 3.3 theorem 9]{Ni09}, up to the global sign, which is 
not fixed in \cite{Ni09}.
\eqref{12.27} shows that $U_1$
can be chosen such that also the statements on the case $c+$ in (a),(b) hold.
\hfill$\Box$

\begin{remarks}\label{t12.5}
(i) By theorem \ref{t10.3}, $f_0$ on $M_{3FN}^{mon}$ has simple
zeros along the two hypersurfaces $M_{3FN}^{[0]}$ and $M_{3FN}^{[2]}$
and simple poles along the two hypersurfaces $M_{3FN}^{[1]}$ and $M_{3FN}^{[3]}$
and no other zeros or poles. The two types of zeros are distinguished by the
value $\mp 2$ of the first derivative with respect to $x$.

(ii) By theorem \ref{t12.4} (a), $f_{univ}$ is holomorphic on 
$\{(\xi,s,B)\in\C\times V^{mat}\, |\, \Re(s)>-1\}\cap U_1$.
It has zeros of both types.
More precisely, 
for any $(s,B)\in V^{mat,b+}$, the zero $x_k$ of $f_{mult}(.,s,B)$
satisfies $\paa_x f_{mult}(x_k,s,B)=(-1)^{k+1}\MGcdot 2$ because then the
argument of the sine in \eqref{12.24} is approximately 
$k\pi\ \textup{mod }2\pi$.
So, $x_k$ is a zero of type $M_{3FN}^{[0]}$ for even $k$,
and a zero of type $M_{3FN}^{[2]}$ for odd $k$.
Hence the two hypersurfaces $M_{3FN}^{[0]}$ and $M_{3FN}^{[2]}$
intersect the image in $M_{3FN}$ of the set
$\left(\{(\xi,s,B)\in\C\times V^{mat}\, |\, \Re(s)>-1\}\cap U_1\right)
/\langle m_{[1]}\rangle$.

(iii) The two hypersurfaces $M_{3FN}^{[0]}$ and $M_{3FN}^{[2]}$
cannot restrict to the image in $M_{3FN}$ of the set
$\left(\{(\xi,s,B)\in\C\times V^{mat,b+}\, |\, \Re(s)>-1\}\cap U_1\right)
/\langle m_{[1]}\rangle$,
as then they would have real codimension 3. As for any $(s^0,B^0)\in V^{mat,b+}$
the $x_k$ tend to zero for $k\to\iiii$, also the functions $f_{mult}(.,s,B)$
for $(s,B)\in V^{mat,a}$ close to $(s^0,B^0)$ must have zeros rather close
to 0. This implies the following:
if $(s,B)\in V^{mat,a}$ approaches a point $(s^0,B^0)\in V^{mat,b+}$, then
$B_2(s,B)$ must tend to $0$, in contrast to $B_1(s,B)$ which tends to $B_1(s^0,B^0)>0$.

(iv) Consider now a point $(2,B^0)\in V^{mat,c+}$ and $(s,B)\in V^{mat,b+}$ which tends to this point. Then $t^{NI}$ tends to $0$,
\begin{eqnarray*}
x_k=2 \exp\left(\frac{1}{2t^{NI}}(2\arg\Gamma(1+it^{NI})-\delta^{NI})\right)
\MGcdot \exp\left(\frac{-k\pi}{2t^{NI}}\right)
\end{eqnarray*}
tends by the proof of part (g) to 
\begin{eqnarray*}
2 \exp\left(-\gamma_{Euler}+\frac{i\pi \www b_1b_2}{2}\right)
\MGcdot \exp\left(\frac{-(k+(1-\www b_1)/2)\pi}{2t^{NI}}\right) ,
\end{eqnarray*}
so its argument tends to
\[
\tfrac12{\pi \www b_1b_2},
\]
and $x_k$ itself tends rapidly to $0$.
This helps to understand why $f_{mult}(.,2,B^0)$ has no zeros near 0,
but $f_{mult}(.,s,B)$ for $(s,B)\in V^{mat,b+}$ close to $(2,B^0)\in V^{mat,c+}$
has zeros arbitrarily close to $0$.
\end{remarks}

\chapter{Rank 2 TEPA bundles with a logarithmic pole}\label{s13}
\setcounter{equation}{0}

\noindent
This chapter offers an independent proof of the results in \cite{Ni09}
which are used above to prove theorem \ref{t12.4}, so together
with the arguments in chapter \ref{s12} it reproves theorem \ref{t12.4}.
It uses the language of this paper.
But some crucial arguments are close to those in \cite{Ni09}
and \cite[ch.\  8]{IN86}: the approximation of (sections of)
the $P_{3D6}$-TEJPA bundles for small $x$ by (sections of)
a closely related bundle, and the explicit control of sections
by Hankel functions.

Still, the details are quite different. We believe that the proof
offered here is useful, as it shows clearly the origins
of the different ingredients in the formulae in theorem \ref{t12.4}.

The function $f_0$ on $M_{3FN}^{ini}$ (lemma \ref{t10.1} (c))
induces the function $f_{univ}:\C\times V^{mat}\to\C$ 
which is the subject of theorem \ref{t12.4}.
The function $f_0$ has simple zeros along the smooth 
hypersurfaces $M_{3FN}^{[0]}$ and $M_{3FN}^{[2]}$ and simple poles
along $M_{3FN}^{[1]}$ and $M_{3FN}^{[3]}$.
On $M_{3FN}^{reg}$ the bases $\uuuu\sigma_0$ and $ \uuuu\sigma_2$
of the normal forms in theorem \ref{t8.2} (b) coincide,
see \eqref{8.19}. Here we shall call this basis
$\uuuu\sigma=(\sigma_1,\sigma_2)$. 

Formula \eqref{8.30} in theorem \ref{t8.2} (g) shows:
the holomorphic family of sections $\sigma_1$ for the pure
$P_{3D6}$-TEJPA bundles in $M_{3FN}^{reg}$ extends holomorphically
to the hypersurfaces $M_{3FN}^{[0]}$ and $M_{3FN}^{[2]}$,
so it is a holomorphic family of sections $\sigma_1$ 
on the union of the chart with $k=0$ and the chart with $k=2$ of $M_{3FN}^{ini}$.
On the union of these two charts $f_0$ can now be calculated
using $\sigma_1$ via 
\begin{eqnarray}\label{13.1}
2f_0 = P(\sigma_1(z),J(\sigma_1(-\rho_1(z)))
\end{eqnarray}
because $J(\sigma_1(-\rho_1(z))))=f_0\MGcdot \sigma_2(-z)$
by \eqref{8.18} and
$P(\sigma_1(z),\sigma_2(-z))=2$ by \eqref{8.16}.

The idea of the proof of the crucial formulae \eqref{12.21}
and \eqref{12.22} in theorem \ref{t12.4} is to show that $\sigma_1$ is 
holomorphic on a suitable neighbourhood of $x=0$ in $M_{3FN}$ and to 
replace in \eqref{13.1} $\sigma_1$ by an approximation for $x$ near $0$ and 
calculate the right hand side of \eqref{13.1} for this approximation.
The approximation will be a linear combination of 
{\it elementary sections}, which are defined after the following
choices.

We shall restrict most of the time to $x\in\C-\R_{\leq 0}$.
This is not a serious restriction, but it has the advantage that
we can fix branches of $\log x$ on $\C-\R_{\leq 0}$ and of
$\log z$ on $\whhh I^+_0$.
We choose the branch $\xi$ of $\log x$ with 
$\arg x=\Im(\xi)\in (-\pi,\pi)$ and $\beta\in\C$ with 
$\tfrac12 e^{-\beta/2}=x$ and $\Im(\beta)\in(-2\pi,2\pi)$.

The choice $u^1_0=u^1_\iiii=x=-u^2_0=-u^2_\iiii$ gives 
$c={u^1_\iiii}/{u^1_0}=1$,
\begin{eqnarray}\label{13.2}
\begin{split}
\zeta_0(x) ={ix}/{|x|},& \quad I^+_0(x)=S^1-\{-\zeta_0\},\\
\zeta_\iiii(x) ={|x|}/{ix},& \quad I^+_\iiii(x)=S^1-\{-\zeta_\iiii\}.
\end{split}
\end{eqnarray}
We choose on $\whhh I^+_0$ the branch of $\log z$ with
\[
\arg z=\Im(\log z)\in (-\tfrac12{\pi}+\arg x,\tfrac32{\pi}+\arg x).
\]

Now let $(H,\nnn,x,x,P,A,J)$ be a $P_{3D6}$-TEJPA bundle with
$x\in \C-\R_{\leq 0}$. 
The $\Phi^{mon}$ in \eqref{10.4} associates to it a unique
tuple $(\beta,s,B)$ with $\beta$ as above, namely
$\tfrac12 e^{-\beta/2}=x,\Im(\beta)\in(-2\pi,2\pi)$.
Throughout this chapter we shall assume that 
$s\in \C^{[sto,a\cup b]}=\C-\{\pm 2\}$. 
This is equivalent to the monodromy being semisimple.

Let $\uuuu e^\pm_0,\uuuu e^\pm_\iiii$ be the (up to a global sign)
unique $4$-tuple of bases in theorem \ref{t7.3} (c). Then
\eqref{2.22}, \eqref{5.1} and \eqref{6.13} show
\begin{eqnarray}\label{13.3}
\Mon(\uuuu e^+_0)&=&\uuuu e^+_0\MGcdot \Mon_0^{mat}(s)\\
\textup{ with }\Mon_0^{mat}(s)&=&S^t\MGcdot S^{-1}, \quad 
S=\begin{pmatrix} 1& s\\ 0 & 1\end{pmatrix}.\nonumber
\end{eqnarray}
\eqref{5.4} and \eqref{5.15} give
\begin{eqnarray}\label{13.4}
\Mon_0^{mat}(s)\MGcdot v_\pm &=& \lambda_\pm\MGcdot v_\pm\\
\textup{ for }v_\pm &=& 
\begin{pmatrix}
\vphantom{\dfrac12}
1
\\
\mp \sqrt{\tfrac14{s^2}-1}+\frac{s}{2}\end{pmatrix}
=\begin{pmatrix}
\vphantom{\dfrac12}
1
\\
\mp i\MGcdot e^{-\pi i\alpha_\pm}\end{pmatrix},\nonumber
\end{eqnarray}
so
\begin{eqnarray}\label{13.5}
\Mon(\uuuu f^+_0)=\uuuu f^+_0\MGcdot 
\begin{pmatrix}\lambda_+&0\\0&\lambda_-\end{pmatrix}
\quad\textup{for}\quad\uuuu f^+_0:=\uuuu e^+_0\MGcdot
(v_+\, v_-).
\end{eqnarray}
$\uuuu f^+_0$ is a flat basis on $\whhh I^+_0$ of eigenvectors
of the monodromy.
Here $z\in \whhh I^+_0$ if 
$\arg z\in(-\tfrac12{\pi}+\arg x,\tfrac32{\pi}+\arg x)$. 
For $z$ with different argument $\arg z$, $\uuuu f^+_0(z)$ is the
flat extension of $\uuuu f^+_0$ on $\whhh I^+_0$. 

Now we define for $k\in\Z$ the {\it elementary sections} on $\whhh I^+_0$
\begin{eqnarray}\label{13.6}
\begin{split}
es_1^{(k)}(z)&:= f^{+1}_0(z)\MGcdot z^{\alpha_++k}
= f^{+1}_0(z)\MGcdot exp((\alpha_++k)\log z),\\
es_2^{(k)}(z)&:= f^{+2}_0(z)\MGcdot z^{\alpha_-+k}
= f^{+2}_0(z)\MGcdot exp((\alpha_-+k)\log z),
\end{split}
\end{eqnarray}
and extend them to (a priori) multi-valued sections on $\C^*$ 
by extending $f_0^{+1/2}$ flatly and $z^{\alpha_\pm+k}$ 
holomorphically 
(on $\whhh I^+_0$ $\log z$ is the branch chosen above with
$\arg z=\Im(\log z) \in(-\tfrac12{\pi}+\arg x,\tfrac32{\pi}+\arg x)$).
The elementary sections are single-valued holomorphic sections on
$H|_{\C^*}$ because
\begin{eqnarray}\label{13.7}
\begin{split}
es_{1/2}^{(k)}(z\MGcdot e^{2\pi i}) 
&= f^{+1/2}_0(z\MGcdot e^{2\pi i})\MGcdot (z\MGcdot e^{2\pi i})^{\alpha_\pm+k}\\
&= \Mon(f^{+1/2}_0(z))\MGcdot \lambda_\mp\MGcdot z^{\alpha_\pm+k}
=es_{1/2}^{(k)}(z).
\end{split}
\end{eqnarray}
They form a holomorphic basis 
$\uuuu{es}^{(k)} =(es_1^{(k)},es_2^{(k)})$ of $H|_{\C^*}$ with
\begin{eqnarray}\label{13.8}
\nnn_\zdz \uuuu{es}^{(k)}(z)=\uuuu{es}^{(k)}(z)\MGcdot
\begin{pmatrix}\alpha_++k&0\\0&\alpha_-+k\end{pmatrix}.
\end{eqnarray}
The following lemma says how the elementary sections behave 
with respect to $A,J,P$.

\begin{lemma}\label{t13.1}
For $z\in \whhh I^+_0$ the following holds.

(a)
\begin{eqnarray}\label{13.9}
A(\uuuu f^+_0(z))&=& \uuuu f^+_0(z\MGcdot e^{\pi i})\MGcdot
\begin{pmatrix} i\MGcdot e^{-\pi i\alpha_-}&0\\
0&-i\MGcdot e^{-\pi i\alpha_+}\end{pmatrix},\\
A(\uuuu{es}^{(k)}(z)) &=& \uuuu{es}^{(k)}(-z)\MGcdot
\begin{pmatrix}i & 0\\ 0&-i\end{pmatrix}\MGcdot (-1)^k.
\label{13.10}
\end{eqnarray}

(b) In \eqref{13.11} $\uuuu f^+_0({1}/{z})$ is the flat
extension of $\uuuu f^+_0$ on $\whhh I^+_0$ with 
$\arg{1}/{z}=-\arg z$. 
\begin{eqnarray}\label{13.11}
J(\uuuu f^+_0(z))&=& \uuuu f^+_0({1}/{z})\MGcdot
\begin{pmatrix} 0&b_+\\b_-&0\end{pmatrix},\\
J(\uuuu{es}^{(k)}(z)) &=& \uuuu{es}^{(-k)}({1}/{z})\MGcdot
\begin{pmatrix}0&b_+\\b_-&0\end{pmatrix}.
\label{13.12}
\end{eqnarray}

(c) 
\begin{equation}\label{13.13}
P(\uuuu f^+_0(z)^t,\uuuu f^+_0(z\MGcdot e^{\pi i}))=
\begin{pmatrix} 0&\lambda_-+1 \\ \lambda_++1&0\end{pmatrix},
\end{equation}
\begin{equation}\label{13.14}
P(\uuuu{es}^{(k)}(z)^t,\uuuu{es}^{(l)}(-z)) = 
(-1)^l z^{k+l}\MGcdot 2\cos(\pi \alpha_+)
\begin{pmatrix}0&1\\1&0\end{pmatrix}.
\end{equation}
\end{lemma}

{\bf Proof:} We can assume $x=\frac{1}{2}$, because 
for any $x\in \C-\R_{\leq 0}$ a $P_{3D6}$-TEJPA bundle
is connected by an isomonodromic family with a 
$P_{3D6}$-TEJPA bundle with $x=\frac{1}{2}$,
and $\uuuu f^+_0(z)$, $\uuuu{es}^{(k)}(z)$, $A$, $J$ and $P$ 
(for fixed $z$) vary flatly with respect to $x$ 
in this isomonodromic family.

(a) 
\begin{eqnarray*}
A(\uuuu f^+_0(z))&=& 
A(\uuuu e^+_0(z))\MGcdot(v_+\, v_-)\\
&\stackrel{\eqref{7.11}}{=}&
\uuuu e^-_0(-z)\MGcdot\begin{pmatrix}0&-1\\1&0\end{pmatrix}\MGcdot 
(v_+\, v_-)\\
&\stackrel{\eqref{2.21}}{=}&
\uuuu e^+_0(z\MGcdot e^{\pi i})\MGcdot\begin{pmatrix}1&s\\0&1\end{pmatrix}
\MGcdot\begin{pmatrix}0&-1\\1&0\end{pmatrix}
\MGcdot(v_+\, v_-)\\
&=& \uuuu e^+_0(z\MGcdot e^{\pi i}) \MGcdot T(s)^{-1}\MGcdot
(v_+\, v_-)\\
&\stackrel{\eqref{5.29}}{=}&
\uuuu e^+_0(z\MGcdot e^{\pi i})\MGcdot
(v_+\, v_-)\MGcdot
\begin{pmatrix}ie^{-\pi i\alpha_-}&0\\ 
0 & -ie^{-\pi i\alpha_+}\end{pmatrix}\\
&=&
\uuuu f^+_0(z\MGcdot e^{\pi i})\MGcdot
\begin{pmatrix}ie^{-\pi i\alpha_-}&0\\ 
0 & -ie^{-\pi i\alpha_+}\end{pmatrix}.
\end{eqnarray*}

\begin{eqnarray*}
A(es_{1/2}^{(k)}(z)) &=& 
A(f^{+1/2}_0(z)\MGcdot z^{\alpha_\pm+k})\\
&\stackrel{\eqref{13.9}}{=}&
f^{+1/2}_0(z\MGcdot e^{\pi i})\MGcdot (\pm i)\MGcdot e^{-\pi i\alpha_\mp}
\MGcdot (z\MGcdot e^{\pi i})^{\alpha_\pm+k}
\MGcdot (-1)^k\MGcdot e^{\pi i\alpha_\mp}\\
&=& es_{1/2}^{(k)}(z\MGcdot e^{\pi i})\MGcdot (\pm i)\MGcdot (-1)^k.
\end{eqnarray*}

(b) $x=\frac{1}{2}$ implies $\zeta_0=i=-\zeta_\iiii$, 
$\whhh I^+_0=\whhh I^-_\iiii$, and $z\in \whhh I^+_0$
satisfies $\arg(z)\in (-\tfrac12{\pi},\tfrac32{\pi})$,
$\arg({1}/{z})\in (-\tfrac32{\pi},\tfrac12{\pi})$.
Thus $\uuuu e^-_\iiii({1}/{z})$ is the flat extension
of $\uuuu e^-_\iiii$ on $\whhh I^-_\iiii$ clockwise, and
$$\uuuu e^-_\iiii({1}/{z})=\uuuu e^+_\iiii({1}/{z})\MGcdot S^b_\iiii
=\uuuu e^+_\iiii({1}/{z})\MGcdot (S^{-1})^t$$
by \eqref{2.21}.
Furthermore, $x=\frac{1}{2}$ implies $\beta=0$, and \eqref{2.24} gives
\[
\uuuu e^-_\iiii({1}/{z})=\uuuu e^+_0({1}/{z})\MGcdot B.
\]
Finally, observe that
\begin{align*}
S^t \begin{pmatrix}1&0\\ 0&-1\end{pmatrix}
&(v_+\, v_-)
\! = \!  
\begin{pmatrix}1&0\\s&-1\end{pmatrix}
\!\! 
\begin{pmatrix}
\vphantom{\dfrac12}
1&1\\ \!  -\sqrt{\tfrac14{s^2}-1}+\frac{s}{2} \! &
\! \sqrt{\frac{s^2}{4}-1}+\frac{s}{2} \! \end{pmatrix} 
\\
&\! =\!  \begin{pmatrix}
\vphantom{\dfrac12}
1&1\\ \sqrt{\tfrac14{s^2}-1}+\tfrac12{s} &
-\sqrt{\tfrac14{s^2}-1}+\tfrac12{s}\end{pmatrix}
\! =\!  (v_-\, v_+).
\end{align*}
Now one calculates
\begin{eqnarray*}
J(\uuuu f^+_0(z)) &=& J(\uuuu e^+_0(z))\MGcdot (v_+\, v_-)\\
&\stackrel{\eqref{7.12}}{=}&
\uuuu e^+_\iiii({1}/{z})\MGcdot\begin{pmatrix}1&0\\0&-1\end{pmatrix}
\MGcdot (v_+\, v_-) \\
&=& \uuuu e^-_\iiii({1}/{z}) \MGcdot S^t \MGcdot 
\begin{pmatrix}1&0\\0&-1\end{pmatrix} \MGcdot (v_+\, v_-)\\
&=& \uuuu e^-_\iiii({1}/{z}) \MGcdot(v_-\, v_+)\\
&=& \uuuu e^+_0({1}/{z}) \MGcdot B\MGcdot (v_-\, v_+)\\
&=& \uuuu e^+_0({1}/{z}) \MGcdot \begin{pmatrix}b_-v_-&b_+v_+\end{pmatrix}\\
&=& \uuuu f^+_0({1}/{z})\MGcdot 
\begin{pmatrix}0&b_+\\b_-&0\end{pmatrix}.
\end{eqnarray*}
\begin{eqnarray*}
J(es_{1/2}^{(k)}(z)) &=& J(f^{+1/2}_0(z)\MGcdot z^{\alpha_\pm+k}) \\
&\stackrel{\eqref{13.11}}{=}& f^{+2/1}_0({1}/{z})\MGcdot 
b_\mp\MGcdot ({1}/{z})^{\alpha_\mp-k}\\
&=& b_\mp\MGcdot es^{(-k)}_{2/1}({1}/{z}).
\end{eqnarray*}

(c) 
\begin{eqnarray*}
&&P(\uuuu f^+_0(z)^t,\uuuu f^+_0(z\MGcdot e^{\pi i})) \\
&=&  (v_+\, v_-)^t\MGcdot
P(\uuuu e^+_0(z)^t,\uuuu e^+_0(z\MGcdot e^{\pi i})) \MGcdot
(v_+\, v_-) \\
&\stackrel{\eqref{2.21}}{=} &
 (v_+\, v_-)^t \MGcdot 
P(\uuuu e^+_0(z)^t,\uuuu e^-_0(-z))\MGcdot
S^{-1}\MGcdot (v_+\, v_-)\\
&\stackrel{\eqref{7.10}}{=} &
 (v_+\, v_-)^t \MGcdot {\bf 1}_2\MGcdot
S^{-1}\MGcdot (v_+\, v_-)\\
&=& \begin{pmatrix}1& -i\MGcdot e^{-\pi i\alpha_+}\\
1 & i\MGcdot e^{-\pi i \alpha_-}\end{pmatrix}
\begin{pmatrix}1&-s\\0&1\end{pmatrix}
\begin{pmatrix}1& 1\\-i\MGcdot e^{-\pi i\alpha_+}& 
i\MGcdot e^{-\pi i \alpha_-}\end{pmatrix}\\
&=& \begin{pmatrix} 
1+is\MGcdot e^{-\pi i\alpha_+}-e^{-2\pi i\alpha_+}&1-is\MGcdot e^{-\pi i\alpha_-}+1\\
1+is\MGcdot e^{-\pi i\alpha_+}+1 & 1-is\MGcdot e^{-\pi i\alpha_-}-e^{-2\pi i\alpha_-}
\end{pmatrix}\\
&=& \begin{pmatrix} 0 &\lambda_-+1\\ \lambda_++1&0\end{pmatrix},
\end{eqnarray*}
\begin{eqnarray*}
&& P(es_1^{(k)}(z),es_2^{(l)}(-z)) \\
&=& P(f^{+1}_0(z)\MGcdot z^{\alpha_++k},f^{+2}_0(z\MGcdot e^{\pi i})\MGcdot
(z\MGcdot e^{\pi i})^{\alpha_-+l})\\
&\stackrel{\eqref{13.13}}{=}& 
(\lambda_-+1)\MGcdot z^{k+l}\MGcdot (-1)^l\MGcdot e^{\pi i\alpha_-}\\
&=& (-1)^l\MGcdot z^{k+l}\MGcdot (e^{-\pi i\alpha_-}+e^{\pi i\alpha_-})\\
&=& (-1)^l\MGcdot z^{k+l}\MGcdot 2\cos(\pi\alpha_-)\\
&=& (-1)^l\MGcdot z^{k+l}\MGcdot 2\cos(\pi\alpha_+).
\end{eqnarray*}
\hfill $\Box$

We keep the $P_{3D6}$-TEJPA bundle $(H,\nnn,x,x,P,A,J)$ which was chosen above with
$x\in\C-\R_{\leq 0}$ and its triple $(\beta,s,B)$.

\begin{definition/lemma}\label{t13.2}
In the following, 
\begin{eqnarray*}
\textup{either }(\www \alpha_+,\www \alpha_-)&:=&(\alpha_+(s),\alpha_-(s))\\
\textup{or }(\www \alpha_+,\www \alpha_-)&:=&(\alpha_-(s)-1,\alpha_+(s)+1)
\end{eqnarray*}
(remark \ref{t13.4} motivates the second choice).
In any case $\www\alpha_++\www\alpha_-=0$.

(a) (Definition) Let $\www H\to\P^1$ be a flat holomorphic pure rank
2 vector bundle on $\P^1$ with global basis $\uuuu\chi=(\chi_1,\chi_2)$.
Define a flat meromorphic connection $\www\nnn$, a pairing $\www P$
and an automorphism $\www A$ on $\www H$ by
\begin{eqnarray}\label{13.15}
\www\nnn_\zdz \uuuu\chi(z) &=& 
\uuuu\chi(z)\MGcdot\left[\frac{x}{z}\begin{pmatrix}0&1\\1&0\end{pmatrix}
+\begin{pmatrix}\www\alpha_+&0\\0&\www\alpha_-\end{pmatrix}\right],\\
\www P(\uuuu\chi(z)^t,\uuuu\chi(-z)) &=& 
\begin{pmatrix}0&2\\2&0\end{pmatrix},\label{13.16}\\
\www A(\uuuu\chi(z)) &=& \uuuu\chi(-z)\MGcdot
\begin{pmatrix}i&0\\0&-i\end{pmatrix}.\label{13.17}
\end{eqnarray}

(b) (Lemma) $(\www H\to\P^1,\www\nnn,\www P,\www A)$ is a 
TEPA bundle (definition \ref{t6.1} (a) and definition \ref{t7.1} (a)).
\end{definition/lemma}

{\bf Proof:} (a) Definition.
(b) Obviously $\www P$ is symmetric and nondegenerate,
${\www A}^2=-\id$, and $\www P(\www A \, a,\www A \, b)=\www P(a,b)$.
The flatness of $\www P$ and $\www A$ follows from
\begin{eqnarray*}
&& \www P(\nnn_\zdz\uuuu\chi(z)^t,\uuuu\chi(-z)) 
+  \www P(\uuuu\chi(z)^t,\nnn_\zdz\uuuu\chi(-z)) \\
&=& \left[\frac{x}{z}\begin{pmatrix}0&1\\1&0\end{pmatrix}
+\begin{pmatrix}\www\alpha_+&0\\0&\www\alpha_-\end{pmatrix}\right]
\begin{pmatrix}0&2\\2&0\end{pmatrix}\\
&&+\begin{pmatrix}0&2\\2&0\end{pmatrix}
\left[\frac{x}{-z}\begin{pmatrix}0&1\\1&0\end{pmatrix}
+\begin{pmatrix}\www\alpha_+&0\\0&\www\alpha_-\end{pmatrix}\right]\\
&=& 0 = \zdz \www P(\uuuu\chi(z)^t,\uuuu\chi(-z)) 
\end{eqnarray*}
and
\begin{eqnarray*}
\nnn_\zdz (\www A(\uuuu\chi(z))) 
&=& \uuuu\chi(-z)\left[\frac{x}{-z}\begin{pmatrix}0&1\\1&0\end{pmatrix}
+\begin{pmatrix}\www\alpha_+&0\\0&\www\alpha_-\end{pmatrix}\right]
\begin{pmatrix}i&0\\0&-i\end{pmatrix}\\
&=&  \uuuu\chi(-z)\begin{pmatrix}i&0\\0&-i\end{pmatrix}
\left[\frac{x}{z}\begin{pmatrix}0&1\\1&0\end{pmatrix}
+\begin{pmatrix}\www\alpha_+&0\\0&\www\alpha_-\end{pmatrix}\right]\\
&=& \www A(\nnn_\zdz \uuuu\chi(z)).
\end{eqnarray*}
\hfill $\Box$ 

\begin{lemma}\label{t13.3}
Consider the above $P_{3D6}$-TEJPA bundle $(H,\nnn,x,x,P,A,J)$
and the TEPA bundle $(\www H,\www\nnn,\www P,\www A)$
in definition/lemma \ref{t13.2}.
Their restrictions to $\C$ are isomorphic, and the isomorphism
\begin{eqnarray}\label{13.18}
(\www H,\www\nnn,\www P,\www A)|_{\C}\to (H,\nnn,P,A)|_{\C}
\end{eqnarray}
is unique up to a sign.

The restrictions to $\C$ of the bundles will be identified via this
isomorphism (i.e.\ by fixing one of the two isomorphisms which differ
only by a sign).
\end{lemma}

{\bf Proof:}
The tuple $(\www H,\www\nnn,\www P,\www A)|_\C$ satisfies all properties
of a $P_{3D6}$-TEPA bundle which do not concern the pole at $\iiii$,
but only the restriction to $\C$.
Therefore all parts of theorem \ref{t2.3}, theorem \ref{t6.3}
and theorem \ref{t7.3} hold which concern only the restriction to $\C$.
In particular, for $u^1_0=u^1_\iiii=-u^2_0=-u^2_\iiii=x$, 
$\alpha^1_0=\alpha^2_0=0$, $\zeta_0={ix}/{|x|}$ and
$I^a_0,I^b_0,I^\pm_0$ as in \eqref{2.1}, and $\www L:=\www H|_{\C^*}$,
the pole at $0$ induces flat rank 1 subbundles 
$\www L^{\pm j}_0$ ($j=1,2$) as after the remarks \ref{t2.4}.
Theorem \ref{t7.3} (c) provides 2 up to a global sign unique bases
$\uuuu{\www e}^\pm_0=(\www e^{\pm 1}_0,\www e^{\pm 2}_0)$ of $\www L|_{\whhh I^\pm_0}$
such that
\begin{eqnarray}\label{13.19}
\begin{split}
& \www e^{\pm j}_0\textup{ are flat generating sections of }\www L^{\pm j}_0,\\
& \textup{the first half of \eqref{2.21} holds with }s^a_0=s^b_0=\www s
\textup{ for some }\www s\in\C,\\
& \www P(\uuuu{\www e}^\pm_0(z)^t,\uuuu{\www e}^\pm_0(-z))={\bf 1}_2,\\
& \www A(\uuuu{\www e}^\pm_0(z))=\uuuu{\www e}^\mp_0(-z)\MGcdot
\begin{pmatrix}0&-1\\1&0\end{pmatrix}.
\end{split}
\end{eqnarray}

By \eqref{2.22} the monodromy matrix with respect to $\uuuu{\www e}^\pm_0$
is $\Mon_0^{mat}(\www s)$, so its eigenvalues are $\lambda_\pm(\www s)$.

On the other hand, the logarithmic pole at $\iiii$ of the bundle
$(\www H,\www\nnn,\www P,\www A)$ shows that the eigenvalues of the
monodromy are $e^{-2\pi i\www\alpha_\pm}$ which is either $\lambda_\pm(s)$ or
$\lambda_\mp(s)$. 
Thus $\{\lambda_+(s),\lambda_-(s)\}
=\{\lambda_+(\www s),\lambda_-(\www s)\},$ so $\www s=s$ or $\www s=-s$.

We have to show $\www s =s$. Then the bases $\uuuu{\www e}^\pm_0$
coincide up to a global sign with the analogous bases
$\uuuu e^\pm_0$ in theorem \ref{t7.3} (c) of the $P_{3D6}$-TEJPA bundle
$(H,\nnn,x,x,P,A,J)$ with tuple $(\beta,s,B)$.
This gives then the (unique up to a sign) isomorphism in \eqref{13.18}.

Define 
\begin{eqnarray}\label{13.20}
\begin{split}
& \left.
\begin{array}{ll}
\uuuu{\www f}^+_0(z)&:= \uuuu{\www e}^+_0(z)\MGcdot 
\begin{pmatrix}v_+(\www s)&v_-(\www s)\end{pmatrix}\\
(\www\alpha_+(\www s),\www\alpha_-(\www s))&:= (\alpha_+(\www s),\alpha_-(\www s))
\end{array} 
\right\} \textup{ if }\www\alpha_\pm=\alpha_\pm(s),\\
& \left.
\begin{array}{ll}
\uuuu{\www f}^+_0(z)&:= \uuuu{\www e}^+_0(z)\MGcdot 
\begin{pmatrix}v_-(\www s)&v_+(\www s)\end{pmatrix}\\
(\www\alpha_+(\www s),\www\alpha_-(\www s))&:= (\alpha_-(\www s)-1,\alpha_+(\www s)+1)
\end{array} 
\right\} \textup{ if }\www\alpha_\pm=\alpha_\mp(s)\mp 1,
\end{split}
\end{eqnarray}
and in both cases
$$\www{es}_{1/2}^{(k)}(z):= \www f^{+1/2}_0(z)\MGcdot z^{\www\alpha_\pm(\www s)}
\quad (\textup{with }1\leftrightarrow +,\ 2\leftrightarrow -).$$
The formulae \eqref{13.9}, \eqref{13.10}, \eqref{13.13} and \eqref{13.14}
hold in both cases if one replaces $A,P,\uuuu f^+_0,\uuuu{es}^{(k)},\alpha_\pm,
\lambda_\pm$ by 
$$\www A,\www P,\uuuu{\www f}^+_0,
\uuuu{\www{es}}^{(k)}:=(\www{es}^{(k)}_1,\www{es}^{(k)}_2),
\www\alpha_\pm, \www\lambda_\pm:=e^{-2\pi i\www\alpha_\pm(\www s)}.$$
\eqref{13.10} and \eqref{13.17} give for $\chi_1$ and $\chi_2$ the 
elementary section expansions in 
\begin{eqnarray}\label{13.21}
\begin{split}
\chi_1 &= \sum_{k\in2\Z}a_{k,1}^+\MGcdot\www{es}_1^{(k)} + 
\sum_{k\in2\Z-1}a_{k,1}^-\MGcdot\www{es}_2^{(k)},\\
\chi_2 &= \sum_{k\in2\Z-1}a_{k,2}^+\MGcdot\www{es}_1^{(k)} + 
\sum_{k\in2\Z}a_{k,2}^-\MGcdot\www{es}_2^{(k)}
\end{split}
\end{eqnarray}
with suitable coefficients $a_{k,1/2}^\pm\in\C$.
Now the logarithmic pole at $\iiii$ in \eqref{13.15} shows
\begin{eqnarray}\label{13.22}
\begin{split}
&\www\alpha_\pm(\www s)=\www\alpha_\pm,\quad\textup{so }
\alpha_\pm(\www s)=\alpha_\pm(s),\\
&\textup{so }\www s=s\textup{ and }\uuuu{\www e}^\pm_0=\uuuu e^\pm_0,
\end{split}
\end{eqnarray}
which is what we had to prove, and furthermore
\begin{eqnarray}\label{13.23}
&&a_{k,1}^+=a_{k,1}^-=a_{k,2}^+=a_{k,2}^-=0\quad\textup{for }k\geq 1,\\
&&a_{0,1}^+\neq 0,\quad a_{0,2}^-\neq 0,\label{13.24}
\end{eqnarray}
and
\begin{eqnarray}\label{13.25}
\begin{split}
\uuuu{\www{es}}^{(k)}&= \uuuu{es}^{(k)}(z)\quad
\textup{in the case}\quad\www\alpha_\pm=\alpha_\pm(s),\\
\uuuu{\www{es}}^{(k)}&= \uuuu{es}^{(k)}(z)\begin{pmatrix}0&z\\z^{-1}&0\end{pmatrix}
\quad \textup{in the case}\quad\www\alpha_\pm=\alpha_\mp(s)\mp 1.
\end{split}
\end{eqnarray}
\hfill$\Box$

\begin{remark}\label{t13.4}
In definition/lemma \ref{t13.2} and lemma \ref{t13.3} we are primarily interested
in the case $\www\alpha_\pm=\alpha_\pm(s)$. We consider also the other case
$\www\alpha_\pm=\alpha_\mp(s)\mp 1$ because it arises by analytic continuation
of $s$ over $\C^{[sto,b+]}$ from the case $\www\alpha_\pm=\alpha_\pm(s)$.

We chose $\www\alpha_\pm,\ (\www H,\www\nnn,\www P,\www A),\ \uuuu{\www f}^+_0$
(in \eqref{13.20}) and $\uuuu{\www{es}}^{(k)}$ such that they fit together
in the two cases by analytic continuation of $s$ over $\C^{[sto,b+]}$.
This will allow us to prove statements on the case $s\in\C^{[sto,b+]}$ by
holomorphic extension from both cases for $s\in \C^{[sto,a]}$.

It also allows us to consider in the proof of lemma \ref{t13.5} after \eqref{13.40}
only the case $\www\alpha_\pm=\alpha_\pm(s)$.
In fact, we need that in lemma \ref{t13.5} the case $\www\alpha_\pm=\alpha_\mp(s)\mp 1$
follows by analytic continuation from the case $\www\alpha_\pm=\alpha_\pm(s)$,
because the leading parts of $\chi_1$ and $\chi_2$ for $x$ near $\iiii$
in the case $\www\alpha_\pm=\alpha_\mp(s)\mp 1$ determine the coefficients
$a_{-1,1}^-(1)$ and $a_{0,2}^-(1)$, but not the coefficient $a_{0,1}^+(1)$.
\end{remark}

\begin{lemma}\label{t13.5}
In \eqref{13.21} the coefficients $a_{-2k,1}^+$, $a_{-2k-1,1}^+$, 
$a_{-2k,2}^-$, $a_{-2k-1,2}^-$
for $k\in\Z_{\geq 0}$ (all others are $0$ because of \eqref{13.23}) are
\begin{eqnarray}\label{13.26}
\begin{split}
a_{-2k,1}^+(x)&= \Gamma(\www\alpha_+-\tfrac12{(2k-1)})\MGcdot
\frac{2^{\www\alpha_+-2k}\MGcdot(-1)^k}{\sqrt{\pi}\MGcdot k!}\MGcdot x^{-\www\alpha_++2k},\\
a_{-2k-1,2}^+(x)&= \Gamma(\www\alpha_+-\tfrac12{(2k+1)})\MGcdot
\frac{2^{\www\alpha_+-2k-1}\MGcdot(-1)^k}{\sqrt{\pi}\MGcdot k!}\MGcdot 
x^{-\www\alpha_++2k+1},\\
a_{-2k,2}^-(x)&= \Gamma(\www\alpha_--\tfrac12{(2k-1)})\MGcdot
\frac{2^{\www\alpha_--2k}\MGcdot(-1)^k}{\sqrt{\pi}\MGcdot k!}\MGcdot x^{-\www\alpha_-+2k},\\
a_{-2k-1,1}^-(x)&= \Gamma(\www\alpha_--\tfrac12{(2k+1)})\MGcdot
\frac{2^{\www\alpha_--2k-1}\MGcdot(-1)^k}{\sqrt{\pi}\MGcdot k!}\MGcdot x^{-\www\alpha_-+2k+1}.
\end{split}
\end{eqnarray}
\end{lemma}

{\bf Proof:}
We make the dependence on $x$ and $z$
explicit
in the basis $\uuuu\chi$ of sections of $\www H$  by writing $\uuuu\chi(x,z)$, and we use the fact that we have an isomonodromic
family in $x$. Then the shape of \eqref{13.15} shows
\begin{eqnarray}\label{13.27}
\uuuu\chi(x,z)=\uuuu\chi(1,{z}/{x}),
\end{eqnarray}
so in \eqref{13.21}
\begin{eqnarray}\label{13.28}
a_{k,1}^+(x)\MGcdot \www{es}_1^{(k)}(z) = 
a_{k,1}^+(1)\MGcdot \www{es}_1^{(k)}({z}/{x}),
\end{eqnarray}
and analogously for the other terms. 
This gives the dependence on $x$ in \eqref{13.26}.
From now on we restrict to the case $x=1$.

In that case \eqref{13.15} together with \eqref{13.21}, \eqref{13.23}, 
\eqref{13.24} and
\[
\nnn_\zdz \www{es}_{1/2}^{(k)}(z) 
= (\www\alpha_\pm+k)\MGcdot \www{es}_{1/2}^{(k)}(z)
\]
gives the recursive relations
\begin{eqnarray}\label{13.29}
\begin{split}
(2\www\alpha_++k)\MGcdot a_{k,2}^+(1) 
&= a_{k+1,1}^+(1)\quad\textup{for }k\in 2\Z_{\leq 0}-1,\\
k\MGcdot a_{k,1}^+(1) 
&= a_{k+1,2}^+(1)\quad\textup{for }k\in 2\Z_{\leq 0}-2,\\
(2\www\alpha_-+k)\MGcdot a_{k,1}^-(1) 
&= a_{k+1,2}^-(1)\quad\textup{for }k\in 2\Z_{\leq 0}-1,\\
k\MGcdot a_{k,2}^-(1) 
&= a_{k+1,1}^-(1)\quad\textup{for }k\in 2\Z_{\leq 0}-2.
\end{split}
\end{eqnarray}
With these recursive relations, the equations in \eqref{13.26} for $x=1$
reduce to the special cases for $a_{0,1}^+(1)$ and $a_{0,2}^-(1)$,
\begin{eqnarray}\label{13.30}
\begin{split}
a_{0,1}^+(1) &= \Gamma(\www\alpha_++\frac{1}{2}) 
\MGcdot \frac{2^{\www\alpha_+}}{\sqrt{\pi}},\\
a_{0,2}^-(1) &= \Gamma(\www\alpha_-+\frac{1}{2}) 
\MGcdot \frac{2^{\www\alpha_-}}{\sqrt{\pi}}.
\end{split}
\end{eqnarray}

\eqref{13.30} will now be proved as follows. Equation \eqref{13.15}
for the connection will be rewritten as two second order linear differential equations,
which turn out to be Bessel equations and are solved by Hankel functions.
From the known asymptotics of the Hankel functions near $0$ and near $\iiii$
we can read off \eqref{13.30}.

\eqref{13.15} in the case $x=1$ has at $0$ a semisimple pole of order 2 with eigenvalues
$\pm 1$. As in the proof of theorem \ref{t10.3} (a), formula \eqref{10.16}, 
one sees that the solution $\uuuu\chi$ takes the form
\begin{eqnarray}\label{13.31}
\uuuu\chi(z) = \uuuu e^+_0(z)\MGcdot
\begin{pmatrix} e^{-1/z}&0\\ 0 & e^{1/z}\end{pmatrix}\MGcdot A^+_0(z)\MGcdot C
\end{eqnarray}
with $C=
\bsp
1&1\\-i&i
\esp$, $A^+_0\in GL(2,\AAA|_{I^\pm_0})$,
$\whhh{A}^+_0={\bf 1}_2$. 

We write
\[
\uuuu\chi(z) = \uuuu e^+_0(z)\MGcdot h(z)=\uuuu e^+_0(z)\MGcdot
\begin{pmatrix}h_{11}(z)& h_{12}(z)\\ h_{21}(z)&h_{22}(z)\end{pmatrix}.
\]
Then the leading parts of $h_{jk}(z)$ near $0$ are
\begin{eqnarray}\label{13.32}
\begin{pmatrix}h_{11}(z)& h_{12}(z)\\ h_{21}(z)&h_{22}(z)\end{pmatrix}
\sim \begin{pmatrix} e^{-1/z}&0\\ 0 & e^{1/z}\end{pmatrix}\MGcdot {\bf 1}_2\MGcdot C.
\end{eqnarray}
The second order linear differential equations satisfied by the functions $h_{jk}(z)$
are calculated as follows.
\begin{eqnarray*}
0&=& \nnn_\zdz \uuuu\chi -\uuuu\chi\MGcdot
\left[\frac{1}{z}\begin{pmatrix}0&1\\1&0\end{pmatrix}
+\begin{pmatrix}\www\alpha_+&0\\0&\www\alpha_-\end{pmatrix}\right]\\
&=& \uuuu e^+_0(z)\left[\zdz h-h\MGcdot \frac{1}{z}\begin{pmatrix}0&1\\1&0\end{pmatrix}
-h\MGcdot \begin{pmatrix}\www\alpha_+&0\\0&\www\alpha_-\end{pmatrix}\right]\\
&=& \uuuu e^+_0(z)\MGcdot 
\begin{pmatrix}\zdz h_{11}-\frac{1}{z}h_{12}-\www\alpha_+ h_{11} & 
\zdz h_{12}-\frac{1}{z}h_{11}-\www\alpha_- h_{12} \\
\zdz h_{21}-\frac{1}{z}h_{22}-\www\alpha_+ h_{21} & 
\zdz h_{22}-\frac{1}{z}h_{21}-\www\alpha_- h_{22}\end{pmatrix},
\end{eqnarray*}

\begin{eqnarray*}
0&=& (\zdz)^2h_{j1}-\zdz (\tfrac{1}{z}h_{j2}+\www\alpha_+h_{j1})\\
&=& (\zdz)^2h_{j1}+\tfrac{1}{z}h_{j2}-\tfrac{1}{z}\zdz 
h_{j2}-\www\alpha_+\zdz h_{j1}\\
&=& (\zdz)^2h_{j1}+\tfrac{1}{z}h_{j2}
-\tfrac{1}{z}(\tfrac{1}{z}h_{j1}+\www\alpha_- h_{j2})-\www\alpha_+\zdz h_{j1}\\
&=& (\zdz)^2h_{j1}-\www\alpha_+\zdz h_{j1}
-\tfrac{1}{z^2}h_{j1}+(1-\www\alpha_-)(\zdz h_{j1}-\www\alpha_+ h_{j1})\\
&=& (\zdz)^2h_{j1} + \zdz h_{j1}
- (\tfrac{1}{z^2}+(1-\www\alpha_-)\www\alpha_+) h_{j1}.
\end{eqnarray*}
Write 
\begin{eqnarray}\label{13.33}
z={i}/{t},\quad 
h_{jk}(z)=h_{jk}({i}/{t})=\sqrt{t}\MGcdot \www h_{jk}(t).
\end{eqnarray}
Then
\begin{eqnarray*}
0&=& (\tdt)^2(\sqrt{t}\MGcdot\www h_{j1})-\tdt(\sqrt{t}\MGcdot\www h_{j1})
+(t^2-(1-\www\alpha_-)\www\alpha_+)\sqrt{t}\MGcdot\www h_{j1}\\
&=& \left(\tfrac{1}{4}\sqrt{t}\MGcdot\www h_{j1} +\sqrt{t}\MGcdot\tdt \www h_{j1}+
\sqrt{t}\MGcdot(\tdt)^2\www h_{j1}\right)\\
&& -\left(\tfrac{1}{2}\sqrt{t}\MGcdot\www h_{j1} 
+\sqrt{t}\MGcdot\tdt \www h_{j1}\right)
+ (t^2-(1-\www\alpha_-)\www\alpha_+)\sqrt{t}\MGcdot\www h_{j1}\\
&=& \sqrt{t}\left[(\tdt)^2\www h_{j1}
+(t^2-(1-\www\alpha_-)\www\alpha_+-\tfrac{1}{4})
\www h_{j1}\right] .
\end{eqnarray*}
Thus
\begin{eqnarray}\label{13.34}
0&=& (\tdt)^2\www h_{j1}+(t^2-(\tfrac{1}{2}+\www\alpha_+)^2 )\MGcdot\www h_{j1},
\end{eqnarray}
and analogously
\begin{eqnarray}\label{13.35}
0&=& (\tdt)^2\www h_{j2}+(t^2-(\tfrac{1}{2}+\www\alpha_-)^2 )\MGcdot\www h_{j2}.
\end{eqnarray}
\eqref{13.34} and \eqref{13.35} are the cases $\nu=\tfrac{1}{2}+\www\alpha_\pm$
of the Bessel equation 
(one of many possible references is \cite[(9.1)]{Te96})
\begin{eqnarray}\label{13.36}
0=(\tdt)^2 y(t) +(t^2-\nu^2)\MGcdot y(t).
\end{eqnarray}
A basis of solutions are the Hankel functions $H_\nu^{(1)}(t)$ and $H_\nu^{(2)}(t)$
\cite[(9.3)]{Te96}, which are multi-valued on $\C^*$. 
Another basis is 
\begin{eqnarray}\label{13.37}
H_{-\nu}^{(1)}(t)=e^{i\nu\pi}H_\nu^{(1)}(t),\quad
H_{-\nu}^{(2)}(t)=e^{-i\nu\pi}H_\nu^{(2)}(t).
\end{eqnarray}
Their asymptotics near the logarithmic pole at $0$ and the irregular pole at $\iiii$
are known. If $\Re(\nu)>0$ the leading part near $0$ is 
\cite[(9.13)]{Te96}
\begin{eqnarray}\label{13.38}
H_\nu^{(1)}(t)\sim \frac{1}{\pi i}\Gamma(\nu)\left({2}/{t}\right)^\nu, \quad
H_\nu^{(2)}(t)\sim \frac{-1}{\pi i}\Gamma(\nu)\left({2}/{t}\right)^\nu.
\end{eqnarray}
The leading part near $\iiii$ is \cite[(9.50) and before (9.52)]{Te96}
\begin{eqnarray}\label{13.39}
\begin{split}
H_\nu^{(1)}(t)&\sim \sqrt{{2}/{\pi t}}\MGcdot e^{i(t-\frac{1}{2}\nu\pi-\frac{1}{4}\pi)}
\quad\textup{for }-\pi<\arg t<2\pi,\\
H_\nu^{(2)}(t)&\sim \sqrt{{2}/{\pi t}}\MGcdot e^{-i(t-\frac{1}{2}\nu\pi-\frac{1}{4}\pi)}
\quad\textup{for }-2\pi<\arg t<\pi.
\end{split}
\end{eqnarray}
Then \eqref{13.32}, \eqref{13.34}, \eqref{13.35} and \eqref{13.39} show that
\begin{eqnarray}\label{13.40}
\begin{split}
h_{11}(z) &=\sqrt{{\pi}/{2}}\MGcdot e^{i\frac{1}{2}\pi(\www\alpha_++1)}
\MGcdot H_{\www\alpha_++\frac{1}{2}}^{(1)}({i}/{z})\MGcdot\sqrt{{i}/{z}},\\
h_{21}(z) &=\sqrt{{\pi}/{2}}\MGcdot e^{-i\frac{1}{2}\pi(\www\alpha_++1)}
\MGcdot H_{\www\alpha_++\frac{1}{2}}^{(2)}({i}/{z})\MGcdot(- i)\MGcdot\sqrt{{i}/{z}},\\
h_{12}(z) &=\sqrt{{\pi}/{2}}\MGcdot e^{i\frac{1}{2}\pi(\www\alpha_-+1)}
\MGcdot H_{\www\alpha_-+\frac{1}{2}}^{(1)}({i}/{z})\MGcdot\sqrt{{i}/{z}},\\
h_{22}(z) &=\sqrt{{\pi}/{2}}\MGcdot e^{-i\frac{1}{2}\pi(\www\alpha_-+1)}
\MGcdot H_{\www\alpha_-+\frac{1}{2}}^{(2)}({i}/{z})\MGcdot i\MGcdot\sqrt{{i}/{z}}.
\end{split}
\end{eqnarray}

Now for a moment we restrict to the case $\www\alpha_\pm=\alpha_\pm(s)$ and to
$s\in\C^{[sto,a]}$, so we exclude the case $\www\alpha_\pm=\alpha_\mp(s)\mp 1$ and
also $s\in \C^{[st,b]}$.
Then $\nu=\frac{1}{2}+\www\alpha_\pm=\frac{1}{2}+\alpha_\pm(s)$ and $\Re(\nu)>0$.
By \eqref{13.38} and \eqref{13.40} the leading parts of 
$h_{11}(z),h_{21}(z),h_{12}(z),h_{22}(z)$ for $z$ near $0$ are
\begin{eqnarray}\label{13.41}
\begin{split}
h_{11}(z)&\sim \sqrt{{\pi}/{2}}\MGcdot e^{i\frac{1}{2}\pi(\alpha_++1)}
\MGcdot \frac{1}{\pi i}\MGcdot\Gamma(\alpha_++\tfrac{1}{2})\MGcdot
\left({2z}/{i}\right)^{\alpha_++\frac{1}{2}}\MGcdot\sqrt{{i}/{z}}\\
&= \Gamma(\alpha_++\tfrac{1}{2})\MGcdot\frac{1}{\sqrt{\pi}}\MGcdot (2z)^{\alpha_+},\\
h_{21}(z)&\sim \Gamma(\alpha_++\tfrac{1}{2})\MGcdot\frac{1}{\sqrt{\pi}}
\MGcdot (-i)e^{-i\pi\alpha_+}\MGcdot(2z)^{\alpha_+},\\
h_{12}(z)&\sim \Gamma(\alpha_-+\tfrac{1}{2})\MGcdot\frac{1}{\sqrt{\pi}}
\MGcdot(2z)^{\alpha_-},\\
h_{22}(z)&\sim \Gamma(\alpha_-+\tfrac{1}{2})\MGcdot\frac{1}{\sqrt{\pi}}
\MGcdot ie^{-i\pi\alpha_-}\MGcdot(2z)^{\alpha_-},
\end{split}
\end{eqnarray}
so the leading parts of $\chi_1(z)$ and $\chi_2(z)$ for $z$ near $\iiii$ are
\begin{eqnarray}\label{13.42}
\begin{split}
\chi_1(z) &= \uuuu e^+_0\MGcdot v_+\MGcdot \Gamma(\alpha_++\tfrac{1}{2})
\MGcdot\frac{1}{\sqrt{\pi}}\MGcdot (2z)^{\alpha_+},\\
\chi_2(z) &= \uuuu e^+_0\MGcdot v_-\MGcdot \Gamma(\alpha_-+\tfrac{1}{2})
\MGcdot\frac{1}{\sqrt{\pi}}\MGcdot (2z)^{\alpha_-}.
\end{split}
\end{eqnarray}
This establishes \eqref{13.30} in the case $\www\alpha_\pm=\alpha_\pm(s)$ and
$s\in \C^{[sto,a]}$. 

In the case $\www\alpha_\pm =\alpha_\mp(s)\mp 1$ and/or $s\in \C^{[sto,b]}$
\eqref{13.30} follows by analytic continuation over $\C^{[sto,b+]}$
from the case $\www\alpha_\pm=\alpha_\pm(s)$ and $s\in\C^{[sto,a]}$,
see remark \ref{t13.4}. \hfill$\Box$   

\begin{lemma}\label{t13.6}
There exists a continuous $m_{[1]}$-invariant and $G^{mon}$-invariant map
\begin{eqnarray}\label{13.43}
B_3:i\R\times V^{mat,a\cup b+}\to \R_{>0}
\end{eqnarray}
such that the set
\begin{eqnarray}\label{13.44}
U_3:= \{(\xi,s,B)\in \C\times V^{mat,a\cup b+}\, |\, |e^\xi|<B_3(i\Im(\xi),s,B)\}
\end{eqnarray}
and the induced subset $V_3:= (\Phi^{mon})^{-1}(U_3/\langle m_{[1]}^2\rangle)$ 
of $M_{3FN}^{mon}$ (see lemma \ref{t10.1} (b)) satisfy:
$\sigma_1$ is a holomorphic family of sections for the $P_{3D6}$-TEJPA bundles
in $V_3$, and for those bundles with $s\in \C^{[sto,a]}$ 
an approximation of $\sigma_1$
for $x$ near $0$ is the leading part of $\chi_1$ for the case 
$\www\alpha_\pm=\alpha_\pm(s)$, i.e.
\begin{eqnarray}\label{13.45}
a_{0,1}^+(1)\MGcdot x^{-\alpha_+}\MGcdot es_1^{(0)}(z).
\end{eqnarray}
\end{lemma}

Before we prove this lemma, we show that it gives those results of \cite{Ni09}
which are used in chapter \ref{s12} together with the other arguments there to
prove theorem \ref{t12.4}. $\sigma_1$ is holomorphic on $V_3\subset M_{3FN}$,
therefore formula \eqref{13.1} $2f_0=P(\sigma_1(z),J(\sigma_1(-\rho_1(z))))$ shows that
$f_0$ is holomorphic on $V_3$. This implies the basic property for the case
$b+$ that $f_{univ}$ is holomorphic (real analytic with respect to $s\in \C^{[sto,b]}$)
on $S^{NI}\times V^{mat,b+}\cap U_1$ for an open set $U_1$ as in \eqref{12.15}.
For $s\in \C^{[sto,a]}$, \eqref{13.1} and \eqref{13.45} give the leading term of
$f_{mult}$ for $x$ near 0,
\begin{eqnarray}
&& f_{mult}(x,s,B)\nonumber \\ &\sim& 
\tfrac{1}{2}P(a_{0,1}^+(1)\MGcdot x^{-\alpha_+}\MGcdot es_1^{(0)}(z),
J(a_{0,1}^+\MGcdot x^{-\alpha_+}\MGcdot es_1^{(0)}(-{1}/{z})))\nonumber \\
&\stackrel{\eqref{13.12}}{=}& \tfrac{1}{2}\MGcdot(a_{0,1}^+)^2\MGcdot 
x^{-2\alpha_+}\MGcdot P(es_1^{(0)}(z),b_-\MGcdot es_2^{(0)}(-z))\nonumber \\
&\stackrel{\eqref{13.14}}{=}& \tfrac{1}{2}\MGcdot\tfrac{1}{2}
\MGcdot(a_{0,1}^+)^2\MGcdot x^{-2\alpha_+}\MGcdot b_-\MGcdot 2\cos(\pi\alpha_+)\nonumber \\
&\stackrel{\eqref{13.29}}{=}& \Gamma(\alpha_++\tfrac{1}{2})^2\MGcdot\frac{1}{\pi}
\MGcdot \left({x}/{2}\right)^{-2\alpha_+}\MGcdot b_-\MGcdot 2\cos(\pi\alpha_+)
\nonumber \\
&=& \frac{\Gamma(\frac{1}{2}+\alpha_+)}{\Gamma(\frac{1}{2}-\alpha_-)}
\MGcdot \left({x}/{2}\right)^{-2\alpha_+}\MGcdot b_- .\label{13.46}
\end{eqnarray}
The last equality uses the well known formula
$
\Gamma(\tfrac{1}{2}+x)\Gamma(\tfrac{1}{2}-x)={\pi}/{\cos(\pi x)}.
$

\eqref{13.46} gives the leading parts in \eqref{12.21} and \eqref{12.22}.
That the next parts are of order $O(|x|^{2-2\Re(\alpha_-)})$ (respectively,
$O(|x|^{2+6\Re(\alpha_-)})$), follows from the arguments in the proof of 
theorem \ref{t12.4},
(see in particular \eqref{12.16} and the two pictures in the proof of theorem \ref{t12.4}).
This completes our approach to theorem \ref{t12.4}.
It remains to prove lemma \ref{t13.6}.

{\bf Proof of lemma \ref{t13.6}:}
A $P_{3D6}$-TEJPA bundle is given by its restrictions to $\C$ and $\P^1-\{0\}$ 
together with $P,A,J$. By lemma \ref{t13.3} its restriction to $\C$ is
$\OO_\C\MGcdot \chi_1+\OO_\C\MGcdot\chi_2$, 
and its restriction to $\P^1-\{0\}$ is 
$\OO_{\P^1-\{0\}}\MGcdot J(\chi_1)+ \OO_{\P^1-\{0\}}\MGcdot J(\chi_2)$.
Therefore the 2-dimensional space $\Gamma(\P^1,\OO(H))$ of global sections 
of $H$ is
$$\Gamma(\C,\OO_\C\MGcdot \chi_1+\OO_\C\MGcdot\chi_2)\cap
\Gamma(\P^1-\{0\},\OO_{\P^1-\{0\}}\MGcdot J(\chi_1)+ \OO_{\P^1-\{0\}}\MGcdot J(\chi_2)).$$
It splits into the 1-dimensional eigenspaces of $A$ with eigenvalues $\pm i$.
We have to show that a generating section of the eigenspace with eigenvalue
$i$ is nonzero at $z=0$. Then it is (up to a scalar) $\sigma_1$.

We consider the sections $\chi_1$ and $\chi_2$ in definition/lemma \ref{t13.2} in the
case $\www\alpha_\pm=\alpha_\pm(s)$. Write
\begin{eqnarray}\label{13.47}
\begin{split}
f_1^+(t)&:=\sum_{k\in 2\Z_{\geq 0}} a_{-k,1}^+(1)\MGcdot t^k,\\
f_1^-(t)&:=\sum_{k\in 2\Z_{\geq 0}} a_{-k-1,1}^-(1)\MGcdot t^k,\\
f_2^+(t)&:=\sum_{k\in 2\Z_{\geq 0}} a_{-k-1,2}^+(1)\MGcdot t^k,\\
f_2^-(t)&:=\sum_{k\in 2\Z_{\geq 0}} a_{-k,2}^-(1)\MGcdot t^k
\end{split}
\end{eqnarray}
so that
\begin{eqnarray}\label{13.48}
f_1^+,f_1^-,f_2^+,f_2^-\in\C\{t^2\}\cap \Gamma(\C,\OO_\C)
\end{eqnarray}
and
\begin{eqnarray}\label{13.49}
\begin{split}
\chi_1(x,z)&= f_1^+(\frac{x}{z})\MGcdot es_1^{(0)}(\frac{z}{x})
+f_1^-(\frac{x}{z})\MGcdot es_2^{(-1)}(\frac{z}{x}),\\
\chi_2(x,z)&= f_2^+(\frac{x}{z})\MGcdot es_1^{(-1)}(\frac{z}{x})
+f_2^-(\frac{x}{z})\MGcdot es_2^{(0)}(\frac{z}{x}).
\end{split}
\end{eqnarray}
Then
\begin{eqnarray}\label{13.50}
\begin{split}
J(\chi_1)(x,z)&= b_-\MGcdot f_1^+(xz)\MGcdot es_2^{(0)}(xz)
+b_+\MGcdot f_1^-(xz)\MGcdot es_1^{(1)}(xz),\\
J(\chi_2)(x,z)&= b_-\MGcdot f_2^+(xz)\MGcdot es_2^{(1)}(xz)
+b_+\MGcdot f_2^-(xz)\MGcdot es_1^{(0)}(xz).
\end{split}
\end{eqnarray}
For small $x$  we ask about the existence of
\begin{eqnarray}\label{13.51}
g_1(x,z),g_2(x,z)\in\C\{z^2\},\quad h_1(x,z),h_2(x,z)\in\C\{z^{-2}\}
\end{eqnarray}
with
\begin{eqnarray}\label{13.52}
\begin{split}
&g_1\MGcdot x^{\alpha_+}\MGcdot \chi_1 + g_2\MGcdot x^{\alpha_-}\MGcdot z\MGcdot \chi_2\\
= &\ h_1\MGcdot x^{-\alpha_-}\MGcdot z^{-1}\MGcdot J(\chi_1) 
+ h_2\MGcdot x^{-\alpha_+}\MGcdot J(\chi_2).
\end{split}
\end{eqnarray}
Define 
\begin{eqnarray}\label{13.53}
\gamma:= \alpha_+-(\alpha_--1)=2\alpha_++1=2-((\alpha_-+1)-\alpha_+).
\end{eqnarray}
In view of \eqref{13.49},
\eqref{13.52} is equivalent to the following two equations,
which collect the coefficients of $es_1^{(0)}(z)$ and $es_2^{(-1)}(z)$:
\begin{eqnarray}\label{13.54}
\begin{split}
&g_1\MGcdot f_1^+(\frac{x}{z}) + g_2\MGcdot x^{2-\gamma}\MGcdot f_2^+(\frac{x}{z}) \\
= &\ h_1\MGcdot x^\gamma\MGcdot b_+\MGcdot f_1^-(xz) 
+ h_2\MGcdot b_+\MGcdot f_2^-(xz),
\end{split}\\
\begin{split}
&g_1\MGcdot x^\gamma\MGcdot f_1^-(\frac{x}{z}) 
+ g_2\MGcdot z^2\MGcdot f_2^-(\frac{x}{z}) \\
= &\ h_1\MGcdot b_-\MGcdot f_1^+(xz) 
+ h_2\MGcdot x^{2-\gamma}\MGcdot z^2\MGcdot b_-\MGcdot f_2^+(xz).
\end{split} \label{13.55}
\end{eqnarray}
The system \eqref{13.54} \& \eqref{13.55} can be solved inductively in powers of $x^2$
if we fix the normalization
\begin{eqnarray}\label{13.56}
\textup{coefficient of }z^0\textup{ in }g_1=1.
\end{eqnarray}
Each step in the induction can be separated into four substeps:
\begin{list}{}{}
\item[(1)]
\eqref{13.54} determines the coefficient in $\C[x^{\pm \gamma},z^2]$
of a power of $x^2$ in $g_1$.
\item[(2)]
\eqref{13.55} determines the coefficient in $\C[x^{\pm \gamma},z^{-2}]$
of the same power of $x^2$ in $h_1$.
\item[(3)]
\eqref{13.54} determines the coefficient in $\C[x^{\pm \gamma},z^{-2}]$
of the same power of $x^2$ in $h_2$.
\item[(4)]
\eqref{13.55} determines the coefficient in $\C[x^{\pm \gamma},z^2]$
of the same power of $x^2$ in $g_2$.
\end{list}

As a result, at most the following monomials turn up in $g_1,g_2,h_1,h_2$:
\begin{eqnarray}\label{13.57}
\begin{split}
g_1:&\ (xz)^{2k}\MGcdot x^{2l+m\gamma}\textup{ with }&
k\in\Z_{\geq 0},l\in 2\Z_{\geq 0},\\&&m\in\{-l,-l+2,\dots,l+2\},\\
g_2:&\ (xz)^{2k}\MGcdot x^{2l+m\gamma}\textup{ with }&
k\in\Z_{\geq 0},l\in 2\Z_{\geq 0}+1,\\&&m\in\{-l,-l+2,\dots,l+2\},\\
h_1:&\ \left(\frac{x}{z}\right)^{2k}\MGcdot x^{2l+m\gamma}\textup{ with }&
k\in\Z_{\geq 0},l\in 2\Z_{\geq 0},\\&&m\in\{-l+1,-l+3,\dots,l+1\},\\
h_2:&\ \left(\frac{x}{z}\right)^{2k}\MGcdot x^{2l+m\gamma}\textup{ with }&
k\in\Z_{\geq 0},l\in 2\Z_{\geq 0},\\&&m\in\{-l,-l+2,\dots,l+2\}.
\end{split}
\end{eqnarray}
We leave the details of the proof of \eqref{13.57} to the reader.

The number $\gamma$ satisfies $0\leq \Re(\gamma)<2$.
Therefore the coefficients of all powers of $z^2$ in $g_1,g_2,h_1,h_2$ tend to $0$, 
if $x\to 0$, except the coefficient $1$ of $z^0$ in $g_1$, and the 
coefficient of $z^0$ in $h_2$, which is also bounded if $x\to 0$.
This indicates that the power series $g_1,g_2,h_1,h_2$ are convergent if $x\to 0$.
Also the details of the proof of this convergence is left to the reader.
The convergence implies the first part of lemma \ref{t13.6}.

For $x$ near 0, $\sigma_1$ and $\chi_1$ are dominated by
\begin{eqnarray}\label{13.58}
x^{-\alpha_+}\MGcdot \left(a_{0,1}^+(1)\MGcdot es_1^{(0)}(z) + 
x^\gamma\MGcdot a_{-1,1}^-(1)\MGcdot es_2^{(-1)}(z)\right).
\end{eqnarray}
If $s\in \C^{[sto,a]}$ then $x^\gamma\to 0$ if $x\to 0$. This establishes the 
second part \eqref{13.45} of lemma \ref{t13.6}.
\hfill$\Box$

\begin{remark}\label{t13.7}
If $s\in\C^{[sto,b+]}$ then $\Re(\gamma)=0$, and $\sigma_1$ and $\chi_1$ are
dominated by \eqref{13.58}. If one replaces $\sigma_1$ in \eqref{13.1} by \eqref{13.58}
one obtains directly the leading part in \eqref{12.24}.
\end{remark}

\chapter[Symmetries of the set of all solutions of $P_{III}(0,0,4,-4)$]
{Symmetries of the universal family of solutions
of $P_{III}(0,0,4,-4)$}\label{s14}
\setcounter{equation}{0}

\noindent
In remark \ref{t10.2} the holomorphic function 
$f_{univ}:\C\times V^{mat}\to \C$ was introduced.
By theorem \ref{t10.3}, the associated multi-valued functions on $\C^*$
(with distinguished branch, see remark \ref{t10.2} (ii))
\[
f_{mult}(.,s,B), \quad (s,B)\in V^{mat}, 
\]
give all solutions of $P_{III}(0,0,4,-4)$.
We have 
$f_{mult}(x,s,B)=f_{univ}(\log x,s,B)$.

Certain symmetries of the function $f_{univ}$ are given by the group
$G^{mon}\times\langle m_{[1]}\rangle \cong (\Z_2\times\Z_2)\times \Z$
(see remark \ref{t12.3} and remark \ref{t10.2} (iv)):
\begin{eqnarray}\label{14.1}
f_{univ}\circ R_1&=&f_{univ}^{-1},\\
\textup{i.e.\ }
f_{univ}(\xi,-s,\left(\begin{smallmatrix}1&0\\0&-1\end{smallmatrix}\right) B
\left(\begin{smallmatrix}1&0\\0&-1\end{smallmatrix}\right)) 
&=& f_{univ}^{-1}(\xi,s,B),\nonumber \\
f_{univ}\circ R_2&=&-f_{univ},\label{14.2}\\
\textup{i.e.\ }
f_{univ}(\xi,s,-B) &=& -f_{univ}(\xi,s,B),\nonumber\\
f_{univ}\circ R_3&=&-f_{univ}^{-1},\label{14.3}\\
\textup{i.e.\ }
f_{univ}(\xi,-s,-\left(\begin{smallmatrix}1&0\\0&-1\end{smallmatrix}\right) B
\left(\begin{smallmatrix}1&0\\0&-1\end{smallmatrix}\right)) 
&=& -f_{univ}^{-1}(\xi,s,B),\nonumber
\end{eqnarray}
and
\begin{eqnarray}\label{14.4}
f_{univ}\circ m_{[1]}&=&f_{univ},\\
\textup{i.e.\ }
f_{univ}(\xi-i\pi,s,B) &=& f_{univ}(\xi,s,\Mon_0^{mat}(s)\MGcdot B).\nonumber
\end{eqnarray}

In theorem \ref{t14.1} and theorem \ref{t14.2} two more symmetries of $f_{univ}$
will be given. The proofs use theorem \ref{12.4}, which allows us to recover the monodromy data 
$(s,B)$ from the asymptotic behaviour of $f_{mult}(.,s,B)$ as $x\to 0$.
The symmetry in theorem \ref{t14.1} will not be used in the rest of the
paper, but it is too beautiful to be omitted.
The symmetry in theorem \ref{t14.2} will be used in chapter \ref{s15}.

\begin{theorem}\label{t14.1}
(a) Let
\begin{eqnarray}\label{14.5}
\begin{split}
R_4:\C\times V^{mat}&\to \C\times V^{mat},\\
(\xi,s,B)&\mapsto (\xi+\tfrac12{\pi i},s,
\bsp 0&1\\-1&s\esp
\MGcdot B).
\end{split}
\end{eqnarray}
Then
\begin{eqnarray}\label{14.6}
f_{univ}\circ R_4 &=& i\MGcdot f_{univ},\\
\textup{ i.e.\ }
(-i)f_{univ}(\xi+\tfrac12{\pi i},s,B) 
&=&f_{univ}(\xi,s,
\bsp s&-1\\1&0\esp
\MGcdot B),\nonumber
\end{eqnarray}
and
\begin{eqnarray}\label{14.7}
R_4^2 = R_2\circ m_{[1]}^{-1}=m_{[1]}^{-1}\circ R_2,\ R_4\circ R_2=R_2\circ R_4\\
R_4\circ R_1=R_2\circ R_1\circ R_4,\ 
R_4\circ m_{[1]}=m_{[1]}\circ R_4.\label{14.8}
\end{eqnarray}

(b) For $(s,B)\in V^{mat}$
\begin{eqnarray}\label{14.9}
(-i)\MGcdot f_{mult}(x\MGcdot e^{\pi i/2},s,B)=f_{mult}(x,s,
\bsp s&-1\\1&0\esp \MGcdot B).
\end{eqnarray}
This means that the multi-valued function 
$(x\mapsto (-i)f_{mult}(x\MGcdot e^{\pi i/2},s,B)$ is the solution
$f_{mult}(.,s,\left(\begin{smallmatrix}s&-1\\1&0\end{smallmatrix}
\right)\MGcdot B)$ 
of $P_{III}(0,0,4,-4)$.
\end{theorem}

{\bf Proof:} 
(a) and (b) The form  \eqref{9.1} of the equation $P_{III}(0,0,4,-4)$ shows
immediately that for any $(s,B)\in V^{mat}$ the multi-valued function
$(x\mapsto (-i)f_{mult}(x\MGcdot e^{\pi i/2},s,B))$ is a solution of 
$P_{III}(0,0,4,-4)$. By theorem \ref{t10.3} (c) there exist unique monodromy data
$(\www s,\www B)\in V^{mat}$  such that this function is
$f_{mult}(.,\www s,\www B)$. The map
\[
V^{mat}\to V^{mat},\quad (s,B)\mapsto (\www s,\www B),
\]
is an automorphism of $V^{mat}$ which is
at least analytic, and possibly algebraic.
We shall show that $(\www s,\www B)=(s,
\bsp  s&-1\\1&0\esp
\MGcdot B)$.
This will establish \eqref{14.6} and \eqref{14.9}.
The identities in \eqref{14.7} and \eqref{14.8} are easily checked.

It is sufficient to show
$(\www s,\www B)=(s,
\bsp  s&-1\\1&0\esp
\MGcdot B)$
for $(s,B)\in V^{mat,a}$.  This is equivalent to 
$(\www s,\www b_-)=(s,(-i)e^{\pi i\alpha_-}\MGcdot b_-)$, by 
lemma \ref{t5.2} (b) and lemma \ref{t5.3}. 
The asymptotic formula \eqref{12.21} looks as follows
for the multi-valued function 
$(x\mapsto (-i)f_{mult}(x\MGcdot e^{\pi i/2},s,B))$:
\begin{align*}
(-i)&f_{mult}(x\MGcdot e^{\pi i/2},s,B)
\\ 
&= (-i)
\frac{\Gamma(\frac{1}{2}-\alpha_-)}{\Gamma(\frac{1}{2}+\alpha_-)}
\left(\tfrac12{x}\MGcdot e^{\pi i/2}\right)^{2\alpha_-}\MGcdot b_-
+O(|x|^{2-2\Re(\alpha_-)})
\\
&= \frac{\Gamma(\frac{1}{2}-\alpha_-)}{\Gamma(\frac{1}{2}+\alpha_-)}
\left(\tfrac12{x}\right)^{2\alpha_-}\MGcdot ((-i)e^{\pi i\alpha_-}\MGcdot b_-)
+O(|x|^{2-2\Re(\alpha_-)}) .
\end{align*}
This gives $(\www s,\www b_-)=(s,(-i)e^{\pi i\alpha_-}\MGcdot b_-)$.
\hfill$\Box$

\begin{theorem}\label{t14.2}
(a) Let
\begin{equation}\label{14.10}
R_5:\C\times V^{mat}\to \C\times V^{mat},\ 
(\xi,s,B)\mapsto (\oooo\xi,\oooo s,\oooo B^{-1}).
\end{equation}
Then
\begin{eqnarray}\label{14.11}
f_{univ}\circ R_5 &=& \oooo{f_{univ}},\\ 
\textup{ i.e.\ } \oooo{f_{univ}(\oooo \xi,s,B)} 
&=&f_{univ}(\xi,\oooo s,\oooo B^{-1}),\nonumber
\end{eqnarray}
and
\begin{eqnarray}\label{14.12}
R_5^2=\id,\ 
R_5\circ R_2=R_2\circ R_5,\ R_5\circ R_1=R_1\circ R_5,\\ 
R_5\circ R_4=R_2\circ m_{[1]}\circ R_4\circ R_5.\label{14.13}
\end{eqnarray}

(b) For $(s,B)\in V^{mat}$
\begin{eqnarray}\label{14.14}
\oooo{f_{mult}(\oooo x,s,B)}=f_{mult}(x,\oooo s,\oooo B^{-1}).
\end{eqnarray}
This means that the multi-valued function 
$(x\mapsto \oooo{f_{mult}(\oooo x,s,B)})$ is the solution
$f_{mult}(.,\oooo s,\oooo B^{-1})$ of
$P_{III}(0,0,4,-4)$.

(c) (Addendum to lemma \ref{t5.1} and lemma \ref{t5.2}) 
If $s\in\C$ with $s^2\notin \R_{\geq 4}$ then 
\begin{eqnarray}\label{14.15}
\oooo{\sqrt{\tfrac14{s^2}-1}} = -\sqrt{\tfrac14{\oooo s^2}-1},\ 
\lambda_\pm(\oooo s)=\oooo{\lambda_\mp(s)}, \ 
v_\pm(\oooo s)=\oooo{v_\mp(s)},\\ 
\alpha_\pm(\oooo s)=\oooo{\alpha_\pm(s)},\ 
b_\pm(\oooo s,\oooo B^{-1})=b_\mp(\oooo s,\oooo B)=\oooo{b_\pm(s,B)}.
\label{14.16}
\end{eqnarray}
If $s\in \C$ with $s^2\in \R_{\geq 4}$, then $\oooo s=s$ and 
\begin{eqnarray}\label{14.17}
\sqrt{\tfrac14{s^2}-1}\in\R_{\geq 0},\ 
\lambda_\pm(s)\in \R_{<0},\ 
v_\pm(s)\textup{ real},\\
\oooo{\alpha_\pm(s)}=\alpha_\mp(s)\mp 1,\ 
b_\pm(s,\oooo B^{-1})=b_\mp(s,\oooo B)=\oooo{b_\mp(s,B)}.
\label{14.18}
\end{eqnarray}
If $s=\pm 2$ then 
\begin{eqnarray}
v_{1/2}(s)=\oooo{v_{1/2}(s)},\ 
\www b_1 (s,\oooo B^{-1})=\oooo{\www b_1(s,B)},\label{14.19}\\
b_2(s,\oooo B^{-1})=-b_2(s,\oooo B)= -\oooo{b_2(s,B)}.
\label{14.20}
\end{eqnarray}
\end{theorem}

{\bf Proof:} (c) The proof consists of elementary calculations and is omitted.

(a) and (b) The form \eqref{9.1} of the equation $P_{III}(0,0,4,-4)$  shows
immediately that for any $(s,B)\in V^{mat}$ the multi-valued function
$x\mapsto \oooo{f_{mult}(\oooo x,s,B)}$ is a solution of 
$P_{III}(0,0,4,-4)$. By theorem \ref{t10.3} (c) there exist unique monodromy data
$(\www s,\www B)\in V^{mat}$ such that this function is
$f_{mult}(.,\www s,\www B)$. The map
\[
V^{mat}\to V^{mat},\quad (s,B)\mapsto (\www s,\www B),
\]
is an automorphism of $V^{mat}$ which is
at least real analytic, possibly real algebraic.
We shall show that $(\www s,\www B)=(\oooo s,\oooo B^{-1})$.
This will establish \eqref{14.11} and \eqref{14.14}.
The identities in \eqref{14.12} and \eqref{14.13} are easily checked.
It is sufficient to show that
$(\www s,\www B)=(\oooo s,\oooo B^{-1})$
for $(s,B)\in V^{mat,a}$. This is equivalent to 
$(\www s,\www b_-)=(\oooo s,\oooo{b_-})$, by \eqref{14.16}.
The asymptotic formula \eqref{12.21} looks as follows for the 
multi-valued function 
$\oooo{f_{mult}(\oooo x,s,B)}$:
\begin{eqnarray*}
\oooo{f_{mult}(\oooo x,s,B)} &=& 
\frac{\Gamma(\frac{1}{2}-\oooo{\alpha_-})}{\Gamma(\frac{1}{2}+\oooo{\alpha_-})}
\left(\tfrac12{x}\right)^{2\oooo{\alpha_-}}\MGcdot \oooo{b_-}
+O(|x|^{2-2\Re(\oooo{\alpha_-})}).
\end{eqnarray*}
This gives $(\www s,\www b_-)=(\oooo s,\oooo b_-)$.
\hfill$\Box$

\begin{remark}\label{t14.3}
The following table extends table \eqref{12.13} and lists some meromorphic
functions on $\C\times V^{mat,a}$ (so $s^2\notin\R_{\geq 4}$) and their images under
the actions of $R_4, R_5,R_1\circ R_5$ and $R_2\circ R_5$:
\begin{eqnarray*}
\renewcommand{\arraystretch}{1.5}
\begin{array}{c|c|c|c|c}
 & R_4 & R_5 & R_1\circ R_5 & R_2\circ R_5 \\ \hline
s & s & \oooo s & -\oooo s & \oooo s\\
\!\sqrt{\tfrac14{s^2}-1}\! & \sqrt{\tfrac14{s^2}-1} & 
\!-\oooo{\sqrt{\tfrac14{s^2}-1}} \!& 
-\oooo{\sqrt{\tfrac14{s^2}-1}} & 
\!\!-\oooo{\sqrt{\tfrac14{s^2}-1}} \!\! \\
\lambda_\pm & \lambda_\pm & \oooo{\lambda_\mp} & \oooo{\lambda_\pm} 
& \oooo{\lambda_\mp} \\
\alpha_\pm & \alpha_\pm & \oooo{\alpha_\pm} & -\oooo{\alpha_\pm} 
& \oooo{\alpha_\pm} \\
B & \left(\begin{smallmatrix}0&1\\-1&s\end{smallmatrix}\right) B & 
\oooo B^{-1} & 
\!\left(\begin{smallmatrix}1&0\\0&-1\end{smallmatrix}\right)\!\oooo B^{-1} \!\!
\left(\begin{smallmatrix}1&0\\0&-1\end{smallmatrix}\right)\! & -\oooo B^{-1} \\
b_1 & -b_2& \oooo{b_1} +\oooo s\oooo{b_2} & 
\oooo{b_1} +\oooo s\oooo{b_2} & 
-\oooo{b_1} -\oooo s\oooo{b_2} \\
b_2 & b_1+sb_2 & -\oooo{b_2} & \oooo{b_2} & \oooo{b_2} \\
b_1+\tfrac12{s}b_2 & 
\!-(1\!-\!\tfrac12{s^2})b_2\!+\!\tfrac12{s}b_1\! & 
\oooo{b_1}+\tfrac12{\oooo s}\oooo{b_2} & 
\oooo{b_1}+\tfrac12{\oooo s}\oooo{b_2} & 
\!\!  -\oooo{b_1}-\tfrac12{\oooo s}\oooo{b_2} \\
b_\pm & \mp ie^{-\pi i\alpha_\pm}\MGcdot b_\pm & 
\oooo{b_\pm} & \oooo{b_\mp} & -\oooo{b_\pm}\\
\xi & \xi +\tfrac12{\pi i}& \oooo\xi & \oooo\xi & \oooo\xi \\
x & x\MGcdot e^{\pi i/2} & \oooo x & \oooo x & \oooo x \\
f_{univ} & i\MGcdot f_{univ} & \oooo{f_{univ}} & \oooo{f_{univ}}^{-1} & 
-\oooo{f_{univ}}
\end{array}
\end{eqnarray*}
Beware that for $s^2\in \R_{\geq 4}$ the images of 
$\sqrt{\tfrac14{s^2}-1},\lambda_\pm,\alpha_\pm$ and $b_\pm$ under 
$R_5,R_1\circ R_5,R_2\circ R_5$ are different, in particular we have:
\begin{eqnarray*}
\renewcommand{\arraystretch}{1.3}
\begin{array}{c|c|c|c}
 & R_5 & R_1\circ R_5 & R_2\circ R_5 \\ \hline
\lambda_\pm & \lambda_\pm=\oooo{\lambda_\pm} & \lambda_\mp & \lambda_\pm \\
\alpha_\pm & \alpha_\pm & -\alpha_\pm & \alpha_\pm \\
b_\pm & \oooo{b_\mp} & \oooo{b_\pm} & -\oooo{b_\mp}
\end{array}
\end{eqnarray*}
\end{remark}

\chapter{Three families of solutions on $\R_{>0}$}\label{s15}
\setcounter{equation}{0}

\noindent
In this chapter we are interested in restrictions to $\R_{>0}$ of solutions 
of $P_{III}(0,0,4,-4)$, which are related to real solutions (possibly with
singularities) on $\R_{>0}$ of one of the following three differential 
equations:
\begin{eqnarray}\label{15.1}
\textup{radial sinh-Gordon}\oplus:&& 
(\xdx)^2\varphi = 16 x^2\sinh\varphi,\\
\textup{radial sinh-Gordon}\ominus:&& 
(\xdx)^2\psi = -16 x^2\sinh\psi,\label{15.2}\\
\textup{radial sine-Gordon}:&& 
(\xdx)^2 u = -16 x^2\sin u.\label{15.3}
\end{eqnarray}
Recall that $(\xdx)^2=x^2(\paa_x^2+\frac{1}{x}\paa_x)$.

\begin{remark}\label{t15.1}
Solutions of \eqref{15.2} lead to CMC surfaces in $\R^3$ with rotationally symmetric metric
 (cf.\ \cite{FPT94}).
Solutions of \eqref{15.1} lead to CMC surfaces in $\R^{2,1}$ with rotationally symmetric metric
 (cf.\ \cite{DGR10}). 
Solutions of \eqref{15.3} lead to timelike surfaces
of constant negative Gauss 
curvature in Minkowski space with rotationally symmetric metric
(cf.\ \cite{Ko11}).
\end{remark}

The following lemma makes precise how real solutions of 
\eqref{15.1}-\eqref{15.3} are related to those solutions of 
$P_{III}(0,0,4,-4)$ which are either real or purely imaginary or have
values in $S^1$. 

\begin{lemma}\label{t15.2}
(a) (i) $f$ is a local solution of $P_{III}(0,0,4,-4)$ on $\R_{>0}$
with values in $\R_{>0}$ (respectively, $\R_{<0}$) $\iff$ 
$\varphi=2\log f$ is a local solution of \eqref{15.1} on $\R_{>0}$
with values in $\R$ (mod $4\pi i$) 
(respectively, $2\pi i + \R$ (mod $4\pi i$)).

(ii) By theorem \ref{t10.3} (b), a real solution of $P_{III}(0,0,4,-4)$
on $\R_{>0}$ can have four types of singularities, which are indexed by
$(\varepsilon_1,\varepsilon_2)\in\{\pm 1\}^2$ 
(or equivalenty by $k\in\{0,1,2,3\},$ see theorem \ref{t8.2} (b)).
These are listed in \eqref{10.13}.
If $x_0\in \R_{>0}$ is a singularity of $f$ of type 
$(\varepsilon_1,\varepsilon_2)$, then the corresponding singularity of
$\varphi=2\log f$ is as follows: for some $\www g_k\in \R$
\begin{eqnarray}
\varphi(x) &=& (1-\varepsilon_2)\pi i \varepsilon_1 
+ 2\varepsilon_1\log ((-2)(x-x_0)) 
+\varepsilon_1\MGcdot \frac{x-x_0}{x_0}\nonumber \\
&&-\varepsilon_1\left(\frac{7}{12 x_0^2}+\frac{4}{3x_0}\www g_k\right) 
(x-x_0)^2+O((x-x_0)^3).   \label{15.4}
\end{eqnarray}

(b) (i) $f$ is a local solution of $P_{III}(0,0,4,-4)$ on $\R_{>0}$ 
with values in $i\R_{>0}$ $\iff$ 
$\psi=2\log f-\pi i$ is a local solution of \eqref{15.2} on $\R_{>0}$
with values in $\R$ (mod $4\pi i$).

(ii) Any local solutions $f$ and $\psi$ as in (i) extend to real analytic
global solutions without singularities on $\R_{>0}$.

(c) (i) $f$ is a local solution of $P_{III}(0,0,4,-4)$ on $\R_{>0}$ 
with values in $S^1$ $\iff$ 
$u=2i\log f+\pi$ is a local solution of \eqref{15.3} on $\R_{>0}$
with values in $\R$. $u$ is unique up to addition of multiples
of $4\pi$.

(ii) Any local solutions $f$ and $u$ as in (i) extend to real analytic
global solutions without singularities on $\R_{>0}$.
\end{lemma}

{\bf Proof:}
(a) (i) This was stated up to the restriction to $\R_{>0}$ 
and the reality condition already in \eqref{9.1}, \eqref{9.2} and \eqref{9.5}.

(ii) \eqref{10.13} gives the Taylor expansion of 
$\varepsilon_2\MGcdot f^{\varepsilon_1}$ at a singularity $x_0$ of type
$(\varepsilon_1,\varepsilon_2)$:
\begin{align*}
\varepsilon_2\MGcdot f^{\varepsilon_1}
&= (-2)(x-x_0)+\frac{-1}{x_0}(x-x_0)^2 
\\
&\quad\quad\quad+ 
\left(\frac{1}{3x_0^2}+\frac{4}{3x_0}\www g_k\right)(x-x_0)^3
+ O((x-x_0)^4)
\\
&= (-2)(x-x_0)\left[ 1+\frac{1}{2x_0}(x-x_0)\right.
\\
&\quad\quad\quad\left.
-\left(\frac{1}{6x_0^2}+\frac{2}{3x_0}\www g_k\right)
(x-x_0)^2+O((x-x_0)^3)\right] .
\end{align*}
Using $\log(1+x)=x-\frac{x^2}{2}+O(x^3)$ and 
$$
2\log f= 2\varepsilon_1\log \epsilon_2 + 
2\varepsilon_1\log (\varepsilon_2 \MGcdot f^{\varepsilon_1}) 
= (1-\varepsilon_2)\pi i \varepsilon_1 + 
2\varepsilon_1\log (\varepsilon_2\MGcdot f^{\varepsilon_1}),
$$
we obtain \eqref{15.4}.

(b) (i) is obvious. (ii) $f$ as in (i) extends uniquely to a real
analytic global solution on $\R_{>0}$ with values in $i\R$ 
which might have singularities a priori. 
But for any $(\varepsilon_1,\varepsilon_2)\in\{\pm 1\}^2$ one has 
$\varepsilon_2\MGcdot \paa_x f^{\varepsilon_1}(x_0)\in i\R$, thus
$\neq -2$. This and \eqref{10.13} show that $f$ has no singularities.

(c) (i) is lemma \ref{t11.1} and is obvious.
(ii) is also obvious, as any local solution with values in $S^1$
extends to a global solution with values in $S^1$.
\hfill$\Box$

Lemma \ref{t15.4} will describe the restrictions of the spaces  
$M_{3FN}^{ini}(x_0)$ of initial conditions to the three
families of solutions on $\R_{>0}$ in lemma \ref{t15.2}.
This will be an elementary consequence of lemma \ref{t10.1} (c)
and theorem \ref{t10.3} (e).

Theorem \ref{t15.5} will describe the restrictions of the spaces 
$M_{3FN}^{mon}(x_0)$ of monodromy data to the three families.
Theorem \ref{t15.5} will follow from the fact that the three families of solutions
in lemma \ref{t15.2} are the fixed points of the symmetries
$R_5$, $R_2\circ R_5$ and $R_1\circ R_5$, respectively, in the 
table in remark \ref{t14.3}.

\begin{definition}\label{t15.3}
For $x_0\in\R_{>0}$, the sets 
\[
\text{$M_{3FN,\R}(x_0)$, $M_{3FN,i\R_{>0}}(x_0)$
and $M_{3FN,S^1}(x_0)$}
\]
are the subsets of $M_{3FN}(x_0)$ of those
$P_{3D6}$-TEJPA bundles which correspond by theorem \ref{t10.3} (e) to the
regular or singular initial conditions at $x_0$ of solutions of
$P_{III}(0,0,4,-4)$ on $\R_{>0}$ with values in $\R$, $i\R_{>0}$ or $S^1$,
respectively.
Their unions over all $x_0\in\R_{>0}$ are called 
\[
\text{
$M_{3FN,\R}$, $M_{3FN,i\R_{>0}}$ and $M_{3FN,S^1}$}
\] 
respectively.

The upper index $\ini$ in 
$M_{3FN,\R}^{ini}(x_0)$, $M_{3FN,i\R_{>0}}^{ini}(x_0)$, $M_{3FN,S^1}^{ini}(x_0)$, 
$M_{3FN,\R}^{ini}$, $M_{3FN,i\R_{>0}}^{ini}$ and $M_{3FN,S^1}^{ini}$
denotes the additional structure (real algebraic or semi-algebraic manifold,
charts, natural functions) which the corresponding set inherits from
$M_{3FN}^{ini}(x_0)$ or $M_{3FN}^{ini}$.

The upper index $\mon$ in 
$M_{3FN,\R}^{mon}(x_0)$, $M_{3FN,i\R_{>0}}^{mon}(x_0)$, $M_{3FN,S^1}^{mon}(x_0)$, 
$M_{3FN,\R}^{mon}$, $M_{3FN,i\R_{>0}}^{mon}$ and $M_{3FN,S^1}^{mon}$
denotes the additional structure (real algebraic or semi-algebraic 
or real-analytic manifold, foliation)
which the corresponding set inherits from
$M_{3FN}^{mon}(x_0)$ or $M_{3FN}^{mon}$.
\end{definition}

\begin{lemma}\label{t15.4}
(a) (i) $M_{3FN,\R}^{ini}(x_0)$ is a real algebraic manifold with four
natural charts. The charts have coordinates $(f_k,\www g_k)$ for $k=0,1,2,3$
and are isomorphic to $\R\times\R$. Each chart intersects each other chart
in $\R^*\times\R$. The coordinate changes are given by \eqref{8.20}
and \eqref{10.7}.
The intersection of all four charts is called $M_{3FN,\R}^{reg}(x_0)$.
Each chart consists of $M_{3FN,\R}^{reg}(x_0)$ and a hyperplane
$M_{3FN,\R}^{[k]}(x_0)$ isomorphic to $\{0\}\times \R$.

(ii) $M_{3FN,\R}^{reg}(x_0)$ with the coordinates $(f_0,g_0)$
is the space of initial conditions \eqref{10.14} at $x_0$ for at $x_0$
regular solutions $f$  of $P_{III}(0,0,4,-4)$ on $\R_{>0}$ with values
in $\R$.
For $k=0,1,2,3$, the space $M_{3FN,\R}^{[k]}(x_0)\cong\R$ with the coordinate
$\www g_k$ is the space of the initial condition \eqref{10.15}
for the at $x_0$ singular solutions $f$ with 
$f^{\varepsilon_1}(x_0)=0$, $\varepsilon_2\paa_x f^{\varepsilon_1}(x_0)=-2$.

(iii) Because $x_0\in\R_{>0}$, $M_{3FN,\R}^{ini}$ is only a real semi-algebraic
manifold with four charts, which are given by putting together the charts for
all $x_0\in\R_{>0}$. 

(b) $M_{3FN,i\R_{>0}}^{ini}(x_0)$ is a real semi-algebraic manifold isomorphic
to $i\R_{>0}\times \R$ with coordinates $(f_0,g_0)$. 
It is the space of initial conditions \eqref{10.14} at $x_0$ for solutions
$f$  of $P_{III}(0,0,4,-4)$ on $\R_{>0}$ with values in $i\R_{>0}$.

(c) $M_{3FN,S^1}^{ini}(x_0)$ is a real algebraic manifold isomorphic
to $S^1\times i\R$ with coordinates $(f_0,g_0)$. 
It is the space of initial conditions \eqref{10.14} at $x_0$ for solutions
$f$  of $P_{III}(0,0,4,-4)$ on $\R_{>0}$ with values in $S^1$.
\end{lemma}

{\bf Proof:}
Everything follows easily from lemma \ref{t15.2}, lemma \ref{t10.1} (c) 
and theorem \ref{t10.3}, in particular formula \eqref{10.14} for the 
initial conditions at regular points.
For (c), observe that if $f$ takes values in $S^1$, then 
$\paa_x f(x_0)\in f(x_0)\MGcdot i\MGcdot \R$.
Conversely, if $f(x_0)\in S^1$ and $\paa_x f(x_0)\in f(x_0)\MGcdot i\MGcdot \R$,
then $u=2i\log f+\pi$ satisfies
$u(x_0)\in\R$, $\paa_x u(x_0)=2i\frac{\paa_x f(x_0)}{f(x_0)}\in\R$,
so $u$ is real on $\R_{>0}$, and $f$ takes  values in $S^1$ on $\R_{>0}$.
\hfill$\Box$ 

Recall that by lemma \ref{t10.1} (b), there is a canonical isomorphism
for $\beta\in\R$ with $\frac12{e^{-\beta/2}}=x_0$,
\begin{eqnarray}\label{15.5}
\begin{split}
M_{3FN}^{mon}(x_0)&\cong V^{mat}\\
(H,\nnn,x_0,x_0,P,A,J)&\mapsto (s,B(\beta)).
\end{split}
\end{eqnarray}
Here $s$ and $B(\beta)$ are associated to $(H,\nnn,x_0,x_0,P,A,J)$ as
in theorem \ref{t7.5} (b).

\begin{theorem}\label{t15.5}
(a) (i) \eqref{15.5} restricts to an isomorphism of real algebraic manifolds
\begin{eqnarray}\label{15.6}
\begin{split}
M_{3FN,\R}^{mon}(x_0&)= V^{mat,\R}\\
&:=\{(s,B)\in V^{mat}\, |\, s\in\R,b_1+\tfrac12{s} b_2\in\R,b_2\in i\R\}\\
&=\{(s,B)\in V^{mat}\, |\, (s,B)=(\oooo s,\oooo B^{-1})\} \\
&\cong \{(s,b_5,b_6)\in\R^3\, |\, b_5^2+(\tfrac14{s^2}-1)b_6^2=1\}\\
\textup{with }&b_5=b_1+\tfrac12{s} b_2,b_6=ib_2. 
\end{split}
\end{eqnarray}

(ii) $V^{mat,\R}$ decomposes into 
$V^{mat,J,\R}:=V^{mat,\R}\cap V^{mat,J}$ for $J\in\{a,b+,b-,c+,c-\}$.
As real analytic manifolds
\begin{eqnarray}\label{15.7}
V^{mat,a,\R}&\stackrel{\cong}{\longrightarrow}&
(-2,2)\times\R^*,\ (s,B)\mapsto (s,b_-),\\
V^{mat,b\pm,\R}&\stackrel{\cong}{\longrightarrow}&
(\pm 1)\MGcdot\R_{> 2}\times S^1,\ (s,B)\mapsto (s,b_-),\label{15.8}\\
V^{mat,c\pm,\R}&\stackrel{\cong}{\longrightarrow}&
\{\pm 2\}\times\{\pm 1\}\times\R,\ (s,B)\mapsto (s,\www b_1,ib_2).\label{15.9}
\end{eqnarray}

(iii) As a $C^\iiii$-manifold $V^{mat,\R}$ is a sphere with four holes.
The symmetries $R_1$ and $R_2$ (and $R_3$) act on it ---
see the table in remark \ref{t12.3}.
The following schematical picture shows $V^{mat,\R}$ with these symmetries.

%Later 1 picture
\includegraphics[width=1.0\textwidth]{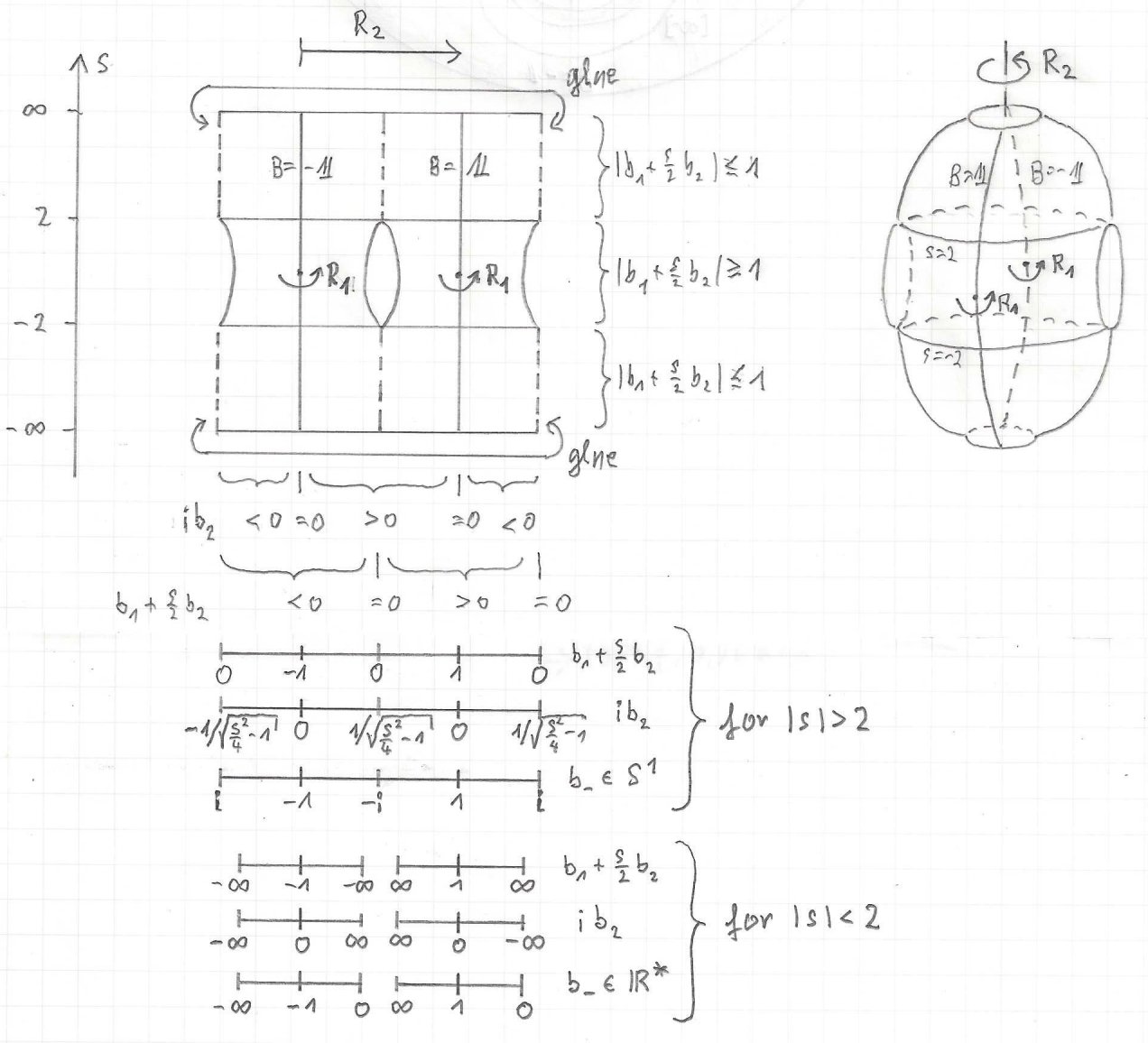} 

(iv) $M_{3FN,\R}^{mon}$ is a real semi-algebraic manifold, isomorphic to
$\R_{>0}\times V^{mat,\R}$. The foliation which $M_{3FN,\R}^{mon}$ inherits
from $M_{3FN}^{mon}$, corresponds to the foliation of 
$\R_{>0}\times V^{mat,\R}$ with leaves $\R_{>0}\times\{(s,B)\}$.

(v) The map 
\begin{eqnarray}\label{15.10}
V^{mat,\R}&\to& \{\textup{solutions of }P_{III}(0,0,4,-4)\textup{ on }\R_{>0}\\
&&\textup{ with values in }\R\}\nonumber \\
(s,B)&\mapsto& f_{mult}(.,s,B)|_{\R_{>0}}
\quad(\textup{i.e.\  }x\textup{ with }\arg x=0)\nonumber
\end{eqnarray}
is an isomorphism.
$f_{mult}(.,s,B)|_{\R_{>0}}$ has a singularity (zero or pole)
at $x_0$ if and only if $(x_0,s,B)$ belongs to the image in $V^{mat,\R}(x_0)$
of $M_{3FN,\R}^{[k]}(x_0)$ for some $k\in\{0,1,2,3\}$.
A pair $(\varepsilon_1,\varepsilon_2)\in\{\pm1\}^2$ is associated with $k$
as in theorem \ref{t8.2} (b).
The function has a zero at $x_0$ if $\varepsilon_1=1$ ($\iff k=0\textup{ or }2$)
and a pole at $x_0$ if $\varepsilon_1=-1$ ($\iff k=1\textup{ or }3$),
and in any case
\begin{eqnarray}\label{15.11}
\paa_x f_{mult}(.,s,B)_{|\R_{>0}}^{\varepsilon_1}(x_0)=\varepsilon_2\MGcdot (-2).
\end{eqnarray}

(b) (i) 
\eqref{15.5} restricts to an isomorphism of real semi-algebraic manifolds
\begin{eqnarray}
&&M_{3FN,i\R_{>0}}^{mon}(x_0)= V^{mat,i\R_{>0}}\nonumber \\
&:=& \{(s,B)\in V^{mat}\, |\, s\in(-2,2), b_1+\tfrac12{s} b_2\in i\R,\nonumber \\
&& b_2=\sqrt{(1-(b_1+\tfrac12{s} b_2)^2)/(1-\tfrac14{s^2})}\in\R_{>1}\} \nonumber\\
&=& \{(s,B)\in V^{mat}\, |\, (s,B)=(\oooo s,-\oooo B^{-1}), b_2\in \R_{>1}\}
\label{15.12}\\
&\cong& \{(s,b_7,b_2)\in(-2,2)\times\R\times \R_{>1}\, |\, 
-b_7^2+(1-\tfrac14{s^2})b_2^2=1\}
\nonumber\\
\textup{with}&&b_7=i(b_1+\tfrac12{s} b_2). \nonumber
\end{eqnarray}

(ii) As a real analytic manifold
\begin{eqnarray}\label{15.13}
V^{mat,i\R_{>0}}&\stackrel{\cong}{\longrightarrow}&
(-2,2)\times\R_{>0},\ (s,B)\mapsto (s,-ib_-).
\end{eqnarray}

(iii) 
The symmetry $R_3=R_1\circ R_2$ acts on $V^{mat,i\R_{>0}}$ --- see the table in remark \ref{t12.3}.

(iv) $M_{3FN,i\R_{>0}}^{mon}$ is a real semi-algebraic manifold, isomorphic to
$\R_{>0}\times V^{mat,i\R_{>0}}$. The foliation which 
$M_{3FN,i\R_{>0}}^{mon}$ inherits from $M_{3FN}^{mon}$, 
corresponds to the foliation of 
$\R_{>0}\times V^{mat,i\R_{>0}}$ with leaves $\R_{>0}\times\{(s,B)\}$.

(v) The map 
\begin{eqnarray}\label{15.14}
V^{mat,i\R_{>0}}&\to& 
\{\textup{solutions of }P_{III}(0,0,4,-4)\textup{ on }\R_{>0}\\
&&\textup{ with values in }i\R_{>0}\}\nonumber \\
(s,B)&\mapsto& f_{mult}(.,s,B)|_{\R_{>0}}
\quad(\textup{i.e.\ }x\textup{ with }\arg x=0)\nonumber
\end{eqnarray}
is an isomorphism.

(c) (i) 
\eqref{15.5} restricts to an isomorphism of real algebraic manifolds
\begin{eqnarray}
&&M_{3FN,S^1}^{mon}(x_0)= V^{mat,S^1}\nonumber \\
&:=& \{(s,B)\in V^{mat}\, |\, s\in i\R,b_1+\tfrac12{s} b_2\in \R,
b_2\in\R \} \nonumber\\
&=& \{(s,B)\in V^{mat}\, |\, (s,B)=(\oooo s,
\left(\begin{smallmatrix}1&0\\0&-1\end{smallmatrix}\right)
\oooo B^{-1}
\left(\begin{smallmatrix}1&0\\0&-1\end{smallmatrix}\right) )  \}\label{15.15} \\
&\cong& \{(s,b_5,b_2)\in i\R\times\R\times \R \, |\, 
b_5^2+(1-\tfrac14{s^2})b_2^2=1\}
\nonumber\\
&&\textup{with}\quad b_5=b_1+\tfrac12{s}b_2. \nonumber
\end{eqnarray}

(ii) As a real analytic manifold
\begin{eqnarray}\label{15.16}
V^{mat,S^1}&\stackrel{\cong}{\longrightarrow}&
\R\times S^1,\ (s,B)\mapsto (is,b_-).
\end{eqnarray}

(iii) 
The symmetries $R_1$ and $R_2$ act on $V^{mat,S^1}$ --- see the table in remark \ref{t12.3}.

(iv) $M_{3FN,S^1}^{mon}$ is a real semi-algebraic manifold, isomorphic to
$\R_{>0}\times V^{mat,S^1}$. The foliation which 
$M_{3FN,S^1}^{mon}$ inherits from $M_{3FN}^{mon}$, 
corresponds to the foliation of 
$\R_{>0}\times V^{mat,S^1}$ with leaves $\R_{>0}\times\{(s,B)\}$.

(v) The map 
\begin{eqnarray}\label{15.17}
V^{mat,S^1}&\to& 
\{\textup{solutions of }P_{III}(0,0,4,-4)\textup{ on }\R_{>0}\\
&&\textup{ with values in }S^1\}\nonumber\\
(s,B)&\mapsto& f_{mult}(.,s,B)|_{\R_{>0}}
\quad(\textup{i.e.\ }x\textup{ with }\arg x=0)\nonumber
\end{eqnarray}
is an isomorphism.
\end{theorem}

{\bf Proof:}
Observe that for $(s,B)\in V^{mat}$
\begin{eqnarray}\label{15.18}
1=b_1^2+b_2^2+sb_1b_2=(b_1+\tfrac12{s} b_2)^2+(1-\tfrac14{s^2})b_2^2.
\end{eqnarray}

(a) $f_{mult}(.,s,B)|_{\R_{>0}}$ takes values in $\R$ if and only if
it is a fixed point of the symmetry $R_5$ in theorem \ref{t14.2}
and remark \ref{t14.3}.
And that holds if and only if $(s,B)$ is a fixed point of the action
of $R_5$ on $V^{mat}$, i.e.\ 
\begin{eqnarray}\label{15.19}
&&(s,B)=(\oooo s,\oooo B^{-1}),\\
\textup{equivalently:}&& s\in\R, b_2\in i\R,b_1+\tfrac12{s}b_2\in \R.
\nonumber
\end{eqnarray}
This establishes the isomorphisms \eqref{15.6} and \eqref{15.10}.
The table in remark \ref{t14.3} shows furthermore that
\begin{eqnarray*}
(s,B)\in V^{mat,a}&\Rightarrow& (B=\oooo B^{-1}\iff b_-\in\R),\\
(s,B)\in V^{mat,b\pm}&\Rightarrow& (B=\oooo B^{-1}\iff b_-\in S^1),\\
(s,B)\in V^{mat,c\pm}&\Rightarrow& (B=\oooo B^{-1}\iff b_2\in i\R).
\end{eqnarray*}
This gives the isomorphisms \eqref{15.7}-\eqref{15.9}.
The rest is clear from lemma \ref{t10.1} and theorem \ref{t10.3}.

(b) $f_{mult}(.,s,B)|_{\R_{>0}}$ takes values in $i\R$ if and only if
it is a fixed point of the symmetry $R_2\circ R_5$ in remark \ref{t14.3}.
And that holds if and only if $(s,B)$ is a fixed point of the action
of $R_2\circ R_5$ on $V^{mat}$, i.e.\
\begin{eqnarray}\label{15.20}
&&(s,B)=(\oooo s,-\oooo B^{-1}),\\
\textup{equivalently:}&& s\in\R, b_2\in \R,b_1+\tfrac12{s}b_2\in i\R,
\nonumber\\
\textup{equivalently:}&& s\in(-2,2), b_2\in \R,b_1+\tfrac12{s}b_2\in i\R
\textup{ by }\eqref{15.18},
\nonumber\\
\textup{equivalently:}&& s\in(-2,2), b_1+\tfrac12{s}b_2\in i\R,\nonumber\\
&& b_2=\pm\sqrt{(1-(b_1+\tfrac12{s}b_2)^2)/(1-\tfrac14{s^2})}\in\R^* .
\nonumber
\end{eqnarray}
The set of these $(s,b_1,b_2)$ has two components. 
By lemma \ref{t15.2} (b) (ii) and the symmetry $R_2$ (or $R_1$),
$f_{mult}(.,s,B)|_{\R_{>0}}$ takes values in $i\R_{>0}$ on one component
and values in $i\R_{<0}$ on the other component.

Because of \eqref{10.36} with $\varepsilon_2=1$,
$f_{mult}(x,0,\pm\left(\begin{smallmatrix}0&1\\-1&0\end{smallmatrix}\right))
= \pm i$, so
$f_{mult}(.,s,B)|_{\R_{>0}}$ takes values in $i\R_{>0}$ on the component
$V^{mat,i\R_{>0}}$. 
This establishes the isomorphisms \eqref{15.12} and \eqref{15.14}.
The formulae in lemma \ref{t5.2} (b)
\[
b_\pm =(b_1+\tfrac12{s}b_2)\mp\sqrt{\tfrac14{s^2}-1}b_2\quad\textup{and}
\quad b_+\MGcdot b_-=1
\]
give $b_-\in i\R_{>0},b_+\in i\R_{<0}$ and the isomorphism \eqref{15.13}.
The rest is clear from lemma \ref{t10.1} and theorem \ref{t10.3}.

(c) 
$f_{mult}(.,s,B)|_{\R_{>0}}$ takes values in $S^1$ if and only if
it is a fixed point of the symmetry $R_1\circ R_5$ in remark \ref{t14.3}.
And that holds if and only if $(s,B)$ is a fixed point of the action
of $R_1\circ R_5$ on $V^{mat}$, i.e.\ 
\begin{eqnarray}\label{15.21}
&&(s,B)=(-\oooo s,
\left(\begin{smallmatrix}1&0\\0&-1\end{smallmatrix}\right)\oooo B^{-1}
\left(\begin{smallmatrix}1&0\\0&-1\end{smallmatrix}\right)),\\
\textup{equivalently:}&& s\in i\R, b_2\in \R,b_1+\tfrac12{s}b_2\in \R.
\nonumber
\end{eqnarray}
This establishes the isomorphisms \eqref{15.15} and \eqref{15.17}.
The table in remark \ref{t14.3} as well as \eqref{15.18} 
give $b_-\in S^1$ and the isomorphism \eqref{15.16}.
The rest is clear from lemma \ref{t10.1} and theorem \ref{t10.3}.
\hfill$\Box$

\begin{remarks}\label{t15.6}
(i) The proof above of theorem \ref{t15.5} uses the symmetry $R_5$ and 
theorem \ref{t14.2}. Theorem \ref{t15.5} can also be proved with the 
additional structures which are put onto the 
$P_{3D6}$-TEJPA bundles (for the three families) in chapter \ref{s16}.

(ii) $f_{mult}(.,s,B)|_{\R_{>0}}$ has singularities (i.e.\ zeros or poles)
arbitrarily close to $0$ if and only if 
$(s,B)\in V^{mat,b,\R}$. This follows from \eqref{12.26}
and the symmetry $R_1$.

(iii) \cite{IN86} and \cite{FIKN06} study primarily the restriction of 
solutions of $P_{III}(0,0,4,-4)$ to $\R_{>0}$. In \cite[ch.\  11]{IN86}
the real solutions with asymptotic behaviour \eqref{12.21} as $x\to 0$ 
and with $b_2\neq 0$,
that means the solutions with $(s,B)\in V^{mat,a,\R}$ and $b_2\neq 0$, 
and the distribution of their singularities (zeros or poles) as $x\to\iiii$
are studied. We shall discuss this in chapters \ref{s16} and \ref{s18}.

In \cite[ch.\  8]{IN86} the solutions $f_{mult}(.,s,B)|_{\R_{>0}}$ for
$(s,B)\in V^{mat}$ with $b_2\in \R^*$ and their asymptotic behaviour as
$x\to\iiii$ and $x\to 0$ are studied. These solutions have no singularities
(zeros or poles) near $x=0$ because of remark (ii) above and 
$\R^*\cap i\R=\emptyset$. Theorem 8.1 in \cite{IN86} says that the solutions
$f_{mult}(.,s,B)$ with $b_2\in\R^*$ have no singularities near $x=\iiii$,
and it gives asymptotic formulae as $x\to\iiii$ for these solutions.

As these solutions are real analytic and have no singularities near
$0$ or near $\iiii$, each of them can have only finitely many singularities.
It is an interesting open question \cite[p. 105]{IN86} whether these
solutions have any singularities at all. The solutions in the subfamilies
$V^{mat,i\R_{>0}}$ and $V^{mat,S^1}-\{(s,B)\, |\, b_2=0\}$ have no singularities
on $\R_{>0}$ by lemma \ref{t15.2} (b)(ii) +(c)(ii).
\end{remarks}

\chapter{TERP structures and $P_{3D6}$-TEP bundles}\label{s16}
\setcounter{equation}{0}

\noindent
The solutions of $P_{III}(0,0,4,-4)$ on $\R_{>0}$ which take values in $\R$ or
in $S^1$ are related to the TERP structures which the second author had
defined in \cite{He03}, motivated by \cite{CV91}, \cite{CV93}, \cite{Du93},
and which were studied subsequently in
\cite{HS07}, \cite{HS10}, \cite{HS11} 
%\cite{KKP09} 
\cite{Mo11b}, \cite{Sa05a}, \cite{Sa05b}
and other papers.
They generalize variations of (polarized) Hodge structures.
The concept of TERP(0) bundle is defined below in definition \ref{t16.1}.
It is a TEP bundle with an additional real structure.
It can be pure or not, and if it is pure, it can be polarized or not.
A pure polarized TERP(0) bundle generalizes a polarized Hodge structure.

A variation of polarized Hodge structures on a punctured disk 
can be approximated by
a so-called nilpotent orbit of Hodge structures, and this gives rise to a
limit polarized mixed Hodge structure \cite{Sch73}, \cite{CKS86}.
This story was generalized to TERP structures in \cite{HS07}
(for the regular singular case, building on \cite{Mo03})
and in \cite{Mo11b} (for the irregular case, building on \cite{Mo11a}).

The case of semisimple TERP structures of rank 2 is the first case which
is far from the regular singular setting and the classical theory of 
variations of Hodge structures.
It is closely related to the solutions of $P_{III}(0,0,4,-4)$ on $\R_{>0}$ 
with values in $\R$ or in $S^1$. 
A good understanding of these solutions gives a good understanding 
of this case and vice versa.
In particular, singularity freeness or existence of singularities 
for real solutions on $\R_{>0}$ near $x=0$ or $x=\iiii$
can be derived from general results on certain 1-parameter orbits
of TERP structures and limit mixed TERP structures.
In chapter \ref{s17} we shall recall these results in the case of semisimple
TERP structures.

In this chapter, first TERP(0) bundles will be defined.
Then it will be shown that, of the $P_{3D6}$-TEJPA bundles with
$u^1_0=u^1_\iiii=x\in\R_{>0}$, only those with $(s,B)\in V^{mat,\R}\cup V^{mat,S^1}$
can be enriched to TERP(0) bundles.
Finally some general results on TERP structures will be cited and their implications
for singularities of real solutions on $\R_{>0}$ of $P_{III}(0,0,4,-4)$
will be explained. The $P_{3D6}$-TEJPA bundles with 
$u^1_0=u^1_\iiii=x\in\R_{>0}$ and $(s,B)\in V^{mat,i\R_{>0}}$ can be enriched
by a quaternionic structure instead of a real structure.
This is explained below.

Recall the antiholomorphic involution of $\P^1$ in \eqref{6.3},
\[
\gamma:\P^1\to\P^1,\quad z\mapsto 1/\oooo z.
\]
\begin{definition/lemma}\label{t16.1}
(a) A TERP(0) bundle is a TEP bundle $(H\to\P^1,\nnn,P)$ together with
a $\C$-linear flat (on $H|_{\C^*}$) isomorphism
\begin{eqnarray}\label{16.1}
\textup{(pointwise)}&& \tau:H_z\to \oooo{H_{\gamma(z)}},\quad
\textup{for all }z\in\P^1,\\
\textup{(for sections)}&& \tau:\OO(H)\to \oooo{\gamma^*\OO(H)},
\quad\OO_{\P^1}\textup{-linear}\nonumber
\end{eqnarray}
with
\begin{eqnarray}\label{16.2}
\tau^2&=&\id,\\
P(\tau(a),\tau(b))&=&\oooo{P(a,b)}\ \textup{for }a\in H_z,b\in H_{-z}.\label{16.3}
\end{eqnarray}
Here $\oooo{H_{\gamma(z)}}$ is $H_{\gamma(z)}$ with the complex conjugate
$\C$-linear structure, so $\tau:H_z\to H_{\gamma(z)}$ is $\C$-antilinear, and  $\oooo{\gamma^*\OO(H)}$ is indeed a free $\OO_{\P^1}$-module
whose rank is $\rank H$. 

(b) A TERP(0) bundle is {\it pure} if $H$ is pure, i.e., if $H$ is a trivial
holomorphic vector bundle (remark \ref{t4.1} (iv)).

(c) (Lemma \cite[lemma 2.5]{HS10}) 
The pairing $P(.,\tau .):H_z\times H_{-1/\oooo z}\to \C$
is sesquilinear (linear $\times$ antilinear) and hermitian, in the sense that
\begin{eqnarray}\label{16.4}
P(b,\tau(a))=\oooo{P(a,\tau(b))}\quad \textup{for }a\in H_z,b\in H_{-1/\oooo z}.
\end{eqnarray}
This follows from \eqref{16.2}, \eqref{16.3} and the symmetry of $P$.
$P(.,\tau .)$ has constant values on global sections of $H$.
The induced pairing 
\begin{eqnarray}\label{16.5}
S_{TERP}:=P(.,\tau .):\Gamma(\P^1,\OO(H))\times \Gamma(\P^1,\OO(H))\to\C
\end{eqnarray}
is sesquilinear and hermitian. It is nondegenerate if and only if 
the TERP(0) bundle is pure.

(d) A pure TERP(0) bundle is {\it polarized} if $S_{TERP}$ is positive definite.
\end{definition/lemma}

\begin{remarks}\label{t16.2}
(i) $\tau$ restricted to $H|_{S^1}$ is fibrewise a $\C$-antilinar involution
and defines a flat real subbundle 
$\ker(\tau-\id:H|_{S^1}\to H|_{S^1})$ which extends to a flat
real subbundle $H'_\R$ of $H':=H|_{\C^*}$ with 
$H_z=H'_{\R,z}\oplus iH'_{\R,z}$ for $z\in\C^*$.

(ii) One can recover $(H,\nnn,P,\tau)$ from $((H,\nnn,P)|_{\C},H'_\R)$:
$\tau$ restricted to $H|_{S^1}$ is the complex conjugation and is extended
flatly to $H'$. Then $H|_\C$ and $\oooo{\gamma^*H|_\C}$ are glued using
$\tau$ to $H$. 
Now $P$ extends holomorphically to $\iiii$, $\nnn$ extends meromorphically to 
$\iiii$. 
In \cite{He03} and \cite{HS07} this observation was taken as starting point,
there a TERP structure (in the case without parameters) was defined by 
the data $((H,\nnn,P)|_\C,H'_\R)$.

(iii) This has the advantage that a TERP structure with parameters,
i.e.\  a variation of TERP(0) bundles, can be defined by the \lq\lq almost
entirely holomorphic\rq\rq\  object $((H,\nnn,P)|_{\C\times M},H'_\R)$.
Here $M$ is a complex manifold, $H|_{\C\times M}$ is a holomorphic vector bundle,
$\nnn$ is a flat connection on $H':=H|_{\C^*\times M}$ with pole of
Poincar\'e rank $\leq 1$ along $\{0\}\times M$, $P$ is a flat holomorphic 
symmetric nondegenerate pairing, and $H'_\R$ is a flat real subbundle
with $H_{(z,t)}=H'_{\R,(z,t)}\oplus iH'_{\R,(z,t)}$ for 
$(z,t)\in\C^*\times M$ and such that $P$ takes real values
on $H'_\R$. The only non-holomorphic ingredient is $H'_\R$. 
Then $\tau$ is defined on $H'$ as above, and
$H|_{\C\times M}$ and $\oooo{\gamma^* H|_{\C\times M}}$ are glued using
$\tau$ to the bundle $H\to\P^1\times M$. This bundle is holomorphic
in $z\in\P^1$, but only real analytic with respect to the parameters
$t\in M$.

(iv) This procedure, passing from the almost holomorphic object
$((H,\nnn,P)|_{\C\times M},H'_\R)$ to the analytic $\times$ (real analytic)
object $(H\to\P^1\times M,\nnn,P,\tau)$, is related to the 
topological-antitopological fusion in \cite{CV91}, \cite{CV93} and to the
Dorfmeister-Pedit-Wu method \cite{DPW98}, \cite{Do08} in the construction
of CMC surfaces.

(v) In \cite{He03}, \cite{HS07} TERP structures of weight $w\in\Z$
are defined. These occur naturally, for example in singularity theory.
They can be reduced (without losing information) to TERP structures
of weight $0$. TERP(0) bundles are TERP structures (without parameters)
of weight $0$.

(vi) Because of the involution $\tau$, the pole at $\iiii$
of a TERP(0) bundle is a \lq\lq twin\rq\rq\  of the pole at $0$.
In particular, if $u^1_0,\dots,u^n_0$ are the eigenvalues of the pole part
$[z\nnn_\zdz]:H_0\to H_0$, then $\oooo{u^1_0},\dots,\oooo{u^n_0}$
are the eigenvalues of the pole part $[-\nnn_{\paa_z}]:H_\iiii\to H_\iiii$ 
at $\iiii$, because $\oooo{\gamma^*(-u^1_0/z)}=-\oooo{u^1_0}\MGcdot z$.

(vii) If $(H,\nnn,P,\tau)$ is a semisimple TERP(0) bundle of rank 2
with eigenvalues $u^1_0\neq u^2_0$ of the pole part at $0$,
the eigenvalues of the pole part at $\iiii$ are
$u^1_\iiii=\oooo{u^1_0},u^2_\iiii=\oooo{u^2_0}$.
Tensoring the bundle with 
$$\OO_{\P^1}\MGcdot e^{(u^1_0+u^2_0)/(2z)+(\oooo{u^1_0}+\oooo{u^2_0})\MGcdot z/2}$$
is a mild twist and leads to a semisimple rank 2 TERP(0) bundle with
eigenvalues $\www{u^1_0}=\frac{u^1_0-u^2_0}{2}$ and $\www{u^2_0}=-\www{u^1_0}$.
\end{remarks}

\begin{theorem}\label{t16.3}
Let $(H,\nnn,x,x,P,A,J)$ be a $P_{3D6}$-TEJPA bundle with $x\in \R_{>0}$
and with canonically associated monodromy data $(s,B)\in V^{mat}$
(lemma \ref{t10.1} (b), $B=B(\beta)$ for $\beta\in\R$ with
$\frac12{e^{-\beta/2}}=x$). Then
\begin{eqnarray}\label{16.6}
&& \zeta_0=i,\quad \zeta_\iiii=-i,\quad c=
{u^1_\iiii}/{u^1_0}={x}/{x}=1,\\
&& I^+_0=I^-_\iiii=S^1-\{-i\},\ I^-_0=I^+_\iiii=S^1-\{i\},\
I^a_0=I^a_\iiii,\ I^b_0=I^b_\iiii.\nonumber
\end{eqnarray}

(a) The underlying TEP bundle can be enriched to a TERP(0) bundle 
if and only if $(s,B)\in V^{mat,\R}\cup V^{mat,S^1}$
(see theorem \ref{t15.5} (a)(i)+(c)(i) for the definition of these subsets
of $V^{mat}$).

(b) (i) Suppose $(s,B)\in V^{mat,\R}$. 
Consider for all $x\in \R_{>0}$ the $P_{3D6}$-TEJPA bundle 
$(H,\nnn,x,x,P,A,J)$ with monodromy data $(s,B)$, and consider the associated 
real solution $f:=f_{mult}(.,s,B)|_{\R_{>0}}$  of $P_{III}(0,0,4,-4)$ on $\R_{>0}$.
Then there are only two involutions which enrich the underlying TEP bundles
to a TERP(0) bundles. One is the $\tau$ defined in \eqref{16.7},
the other is $-\tau$. In \eqref{16.7} $\uuuu e^\pm_0,\uuuu e^\pm_\iiii$
is the (unique up to a global sign) $4$-tuple of bases associated
to the $P_{3D6}$-TEJPA bundle in theorem \ref{t7.3} (c).
\begin{eqnarray}\label{16.7}
\tau(\uuuu e^\pm_0(z))=e^\mp_\iiii(1/\oooo z),\quad 
\tau(\uuuu e^\pm_\iiii(z))=e^\mp_0(1/\oooo z).
\end{eqnarray}

(ii) $(H,\nnn,P,\tau)$ is a pure TERP(0) bundle if and only if 
$x$ is not a zero or a pole of $f$.

(iii) Suppose that the $P_{3D6}$-TEJPA bundle is pure.
Then it has the normal form \eqref{8.15}--\eqref{8.18} for $k=0$ with
$f_0=f(x)$. Then $\tau$ is the map with
\begin{eqnarray}\label{16.8}
\uuuu\tau(\sigma_0(z))=\uuuu\sigma_0(1/\oooo z)\MGcdot 
\begin{pmatrix}0&f(x)^{-1}\\f(x)&0\end{pmatrix}=J(\uuuu \sigma_0(\oooo z)).
\end{eqnarray}
The pairing $S=P(.,\tau\, .)$ on $\Gamma(\P^1,\OO(H))$ has the matrix
\begin{eqnarray}\label{16.9}
S(\uuuu \sigma_0(z)^t,\uuuu\sigma_0(-1/\oooo z))=
\begin{pmatrix} 2f(x)& 0 \\ 0&2f(x)^{-1}\end{pmatrix}.
\end{eqnarray}
So $S$ is positive definite and $(H,\nnn,P,\tau)$ is a pure and polarized
TERP(0) bundle if $f(x)>0$, and $S$ is negative definite and 
$(H,\nnn,P,\tau)$ is a pure, but not polarized TERP(0) bundle if $f(x)<0$.

(iv)
\begin{eqnarray}\label{16.10}
A\circ \tau &=& \tau\circ A,\\
J\circ \tau &=& \tau\circ J.\label{16.11}
\end{eqnarray}

(c) (i) Suppose $(s,B)\in V^{mat,S^1}$, and let $f:=f_{mult}(.,s,B)|_{\R_{>0}}$
be the corresponding solution  of $P_{III}(0,0,4,-4)$ on $\R_{>0}$
with values in $S^1$.
Then there are only two involutions which enrich the underlying TEP bundle
to a TERP(0) bundle. One is the $\tau$ defined in \eqref{16.12},
the other is $-\tau$. In \eqref{16.12} $\uuuu e^\pm_0,\uuuu e^\pm_\iiii$
is the (unique up to a global sign) $4$-tuple of bases associated
to the $P_{3D6}$-TEJPA bundle in theorem \ref{t7.3} (c).
\begin{eqnarray}\label{16.12}
\begin{split}
\tau(\uuuu e^\pm_0(z))&=e^\mp_\iiii(1/\oooo z)
\MGcdot\begin{pmatrix}1&0\\0&-1\end{pmatrix}, \\ 
\tau(\uuuu e^\pm_\iiii(z))&=e^\mp_0(1/\oooo z)
\MGcdot\begin{pmatrix}1&0\\0&-1\end{pmatrix}.
\end{split}
\end{eqnarray}

(ii) $(H,\nnn,P,\tau)$ is a pure TERP(0) bundle.

(iii) It has the normal form \eqref{8.15}--\eqref{8.18} for $k=0$ with
$f_0=f(x)\in S^1$. Then $\tau$ is the map with
\begin{eqnarray}\label{16.13}
\uuuu\tau(\sigma_0(z))=\uuuu\sigma_0(1/\oooo z)\MGcdot 
\begin{pmatrix}f(x)^{-1}&0\\0&f(x)\end{pmatrix}.
\end{eqnarray}
The pairing $S=P(.,\tau\, .)$ on $\Gamma(\P^1,\OO(H))$ has the matrix
\begin{eqnarray}\label{16.14}
S(\uuuu \sigma_0(z)^t,\uuuu\sigma_0(-1/\oooo z))=
\begin{pmatrix} 0&2f(x) \\ 2f(x)^{-1}&0\end{pmatrix}.
\end{eqnarray}
So $S$ is nondegenerate, but indefinite, 
and $(H,\nnn,P,\tau)$ is a pure, but not polarized
TERP(0) bundle.

(iv)
\begin{eqnarray}\label{16.15}
A\circ \tau &=& -\tau\circ A,\\
J\circ \tau &=& \tau\circ J.\label{16.16}
\end{eqnarray}
\end{theorem}

{\bf Proof:}
\eqref{16.6} is obvious. Let $\uuuu e^\pm_0,\uuuu e^\pm_\iiii$ be the
$4$-tuple of bases associated to the $P_{3D6}$-TEJPA bundle in 
theorem \ref{t7.3} (c). Then
\begin{eqnarray}\label{16.17}
\uuuu e^-_\iiii &=& \uuuu e^+_0\MGcdot B.
\end{eqnarray}

(a) A map $\tau$ enriches the underlying TEP bundle to a TERP(0) bundle
if and only if it is an isomorphism between $(H,\nnn,P)$ and
$\oooo{\gamma^*(H,\nnn,P)}$ with $\tau^2=\id$.

We want to describe $\tau$ via its action on the $4$-tuple of bases
$\uuuu e^\pm_0,\uuuu e^\pm_\iiii$. These bases generate the restrictions
of $H$ to the sectors $\whhh I^\pm_0,\whhh I^\pm_\iiii$.
We have to take care that the images under $\tau$ of these restrictions
of $H$ glue in the same way as the restrictions of $H$ themselves,
that $\tau$ respects $P$ in the sense of \eqref{16.3} and that $\tau^2=\id$
holds.
We reorder these conditions as follows:

\begin{list}{}{}
\item[($\alpha$)]
$\tau$ respects the splittings of $H$ on the sectors  
$\whhh I^\pm_0,\whhh I^\pm_\iiii$ from the Stokes structure, and $\tau$
respects $P$.
\item[($\beta$)]
$\tau^2=\id$.
\item[($\gamma$)]
$\tau$ respects the glueing of the restrictions to $\whhh I^+_0$ and
$\whhh I^-_0$, that means, the Stokes structure at $0$,
and similarly the Stokes structure at $\iiii$.
\item[($\delta$)]
$\tau$ respects the glueing via the connection matrix of the bundles
on $\C$ and on $\P^1-\{0\}$.
\end{list}

($\alpha$) says that $\tau(\uuuu e^\mp_\iiii),\tau(\uuuu e^\mp_0)$
is one of the eight $4$-tuples of bases in \eqref{6.12},
\begin{eqnarray}
\varepsilon_0\MGcdot (e^{\pm 1}_0,\varepsilon_1 e^{\pm 2}_0),
\varepsilon_0\varepsilon_2\MGcdot (e^{\pm 1}_\iiii,\varepsilon_1 e^{\pm 2}_\iiii)
\quad \textup{for}\quad \varepsilon_0,\varepsilon_1,\varepsilon_2
\in\{\pm\},\nonumber\\
\textup{so}\quad 
\tau(\uuuu e^\mp_\iiii)=\varepsilon_0\MGcdot \uuuu e^\pm_0
\MGcdot\begin{pmatrix}1&0\\0&\varepsilon_1\end{pmatrix},\ 
\tau(\uuuu e^\mp_0)=\varepsilon_0\varepsilon_2\MGcdot \uuuu e^\pm_\iiii
\MGcdot\begin{pmatrix}1&0\\0&\varepsilon_1\end{pmatrix}.\label{16.18}
\end{eqnarray}

($\beta$) $\tau^2=\id$ says that $\varepsilon_2=1$.

Concerning ($\gamma$): recall from \eqref{6.13} that the $4$-tuple of
bases $\uuuu e^\pm_0,\uuuu e^\pm_\iiii$ satisfies \eqref{2.21}
with 
$$S=S^a_0\quad\textup{and}\quad S^a_\iiii=S^{-1}, 
S^b_0=S^t, S^b_\iiii=(S^t)^{-1},$$
and that any of the eight $4$-tuples of bases in \eqref{6.12} satisfies 
\eqref{2.21} with $S$ replaced by 
$\begin{pmatrix}1&0\\0&\varepsilon_1\end{pmatrix} \MGcdot S\MGcdot  
\begin{pmatrix}1&0\\0&\varepsilon_1\end{pmatrix}$.
Now the following four calculations show that the $4$-tuple
$\tau(\uuuu e^\mp_\iiii),\tau(\uuuu e^\mp_0)$ satisfies \eqref{2.21}
with $S$ replaced by $\oooo S$. We give the details only for the first
calculation. It uses \eqref{2.21} and \eqref{16.6}.
\begin{eqnarray*}
\tau(\uuuu e^+_\iiii)|_{I^a_0}&=&\tau({\uuuu e^+_\iiii}|_{I^a_\iiii})
=\tau({\uuuu e^-_\iiii}|_{I^a_\iiii}\MGcdot (S^a_\iiii)^{-1})\\
&=& \tau(\uuuu e^-_\iiii)|_{I^a_0}\MGcdot (\oooo S^a_\iiii)^{-1}
=\tau(\uuuu e^-_\iiii)|_{I^a_0}\MGcdot \oooo S,\\
\tau(\uuuu e^+_\iiii)|_{I^b_0}&=& \dots
=\tau(\uuuu e^-_\iiii)|_{I^b_0}\MGcdot \oooo S^t,\\
\tau(\uuuu e^+_0)|_{I^a_\iiii}&=& \dots
=\tau(\uuuu e^-_0)|_{I^a_\iiii}\MGcdot \oooo S^{-1},\\
\tau(\uuuu e^+_0)|_{I^b_\iiii}&=& \dots
=\tau(\uuuu e^-_0)|_{I^b_\iiii}\MGcdot (\oooo S^t)^{-1}.
\end{eqnarray*}
($\gamma$) is equivalent to 
\begin{eqnarray}\label{16.19}
\begin{pmatrix}1 &0\\0&\varepsilon_1\end{pmatrix} \MGcdot S\MGcdot  
\begin{pmatrix}1 &0\\0&\varepsilon_1\end{pmatrix} =\oooo S,\quad
\textup{so}\quad \varepsilon_1 s=\oooo s.
\end{eqnarray}

Concerning ($\delta$): the analogue of \eqref{16.17} for 
$\varepsilon_0 \uuuu e^-_\iiii
\bsp 1&0\\0&\varepsilon_1\esp$
and $\varepsilon_0 \uuuu e^+_0
\bsp 1&0\\0&\varepsilon_1\esp$ 
is
\[
\varepsilon_0 \uuuu e^-_\iiii\begin{pmatrix}1&0\\0&\varepsilon_1\end{pmatrix}
= \varepsilon_0 \uuuu e^+_0\begin{pmatrix}1&0\\0&\varepsilon_1\end{pmatrix}
\MGcdot 
\left[
\begin{pmatrix}1&0\\0&\varepsilon_1\end{pmatrix}\MGcdot B\MGcdot 
\begin{pmatrix}1&0\\0&\varepsilon_1\end{pmatrix}
\right].
\]
The analogue of \eqref{16.17} for $\tau(\uuuu e^+_0)$ and $\tau(\uuuu e^-_\iiii)$ is 
$$\tau(\uuuu e^+_0)=\tau(\uuuu e^-_\iiii)\MGcdot \oooo B^{-1}.
$$
($\delta$) is equivalent to 
\begin{eqnarray}\label{16.20}
\begin{pmatrix}1&0\\0&\varepsilon_1\end{pmatrix}\MGcdot B\MGcdot
\begin{pmatrix}1&0\\0&\varepsilon_1\end{pmatrix}=\oooo B^{-1}.
\end{eqnarray}

\eqref{16.19} and \eqref{16.20} give \eqref{15.18} for $\varepsilon_1=1$
and \eqref{15.21} for $\varepsilon_1=-1$.
This shows that $\tau$ with all desired properties exists if and only if
$(s,B)\in V^{mat,\R}\cup V^{mat,S^1}$. It also shows that $\tau$ is unique
up to  sign and that it is given by \eqref{16.7} for $\varepsilon_1=1$
and by \eqref{16.12} for $\varepsilon_1=-1$. This proves
(b) (i) and (c) (i).

(b) (i) is proved above.

(ii) This follows from theorem \ref{t10.3} (b).

(iii) Recall remark \ref{t8.1} (ii): by the correspondence \eqref{2.19},
the bases $\uuuu e^\pm_0$ correspond to a basis $\uuuu v_0$ of $H_0$,
the bases $\uuuu e^\pm_\iiii$ correspond to a basis $\uuuu v_\iiii$ of
$H_\iiii$. Then $\tau(\uuuu v_0)=\uuuu v_\iiii$. The construction of the
basis $\sigma_0$ in the proof of theorem \ref{t8.2} (b) shows that
\[
\uuuu\sigma_0(0)=\uuuu v_0\MGcdot C,
\quad \uuuu\sigma_0(\iiii)=\uuuu v_\iiii \MGcdot
C\MGcdot \begin{pmatrix} f(x)&0\\0&f(x)^{-1}\end{pmatrix}.
\]
Then 
\begin{eqnarray*}
\tau(\uuuu\sigma_0(0))&=&\uuuu v_\iiii\MGcdot\oooo C
=\uuuu\sigma_0(\iiii)\MGcdot \begin{pmatrix} f(x)^{-1}&0\\0&f(x)\end{pmatrix}
\MGcdot C^{-1}\MGcdot \oooo C\\
&=&\uuuu\sigma_0(\iiii)\MGcdot \begin{pmatrix} 0&f(x)^{-1}\\f(x)&0\end{pmatrix}.
\end{eqnarray*}
as $C^{-1}\MGcdot\oooo C=\begin{pmatrix}0&1\\1&0\end{pmatrix}$.

As $\tau$ acts on $\Gamma(\P^1,\OO(H))$, this establishes \eqref{16.8}.
\eqref{16.9} and the conclusions for $f(x)>0$ and for $f(x)<0$ 
follow from \eqref{16.8} and \eqref{8.16}.

(iv) Compare \eqref{16.7} with \eqref{7.11} and \eqref{7.12}.

(c) (i) is proved above.

(ii) This follows from theorem \ref{t10.3} (b)
and lemma \ref{t15.2} (c) (ii).

(iii) This is analogous to the proof of (b) (iii). The only change is
\begin{eqnarray*}
\tau(\uuuu v_0)&=&\uuuu v_\iiii\MGcdot\begin{pmatrix}1&0\\0&-1\end{pmatrix},\\
\tau(\uuuu\sigma_0(0))
&=&\uuuu v_\iiii\MGcdot\begin{pmatrix}1&0\\0&-1\end{pmatrix}\MGcdot\oooo C\\
&=&\uuuu\sigma_0(\iiii)\MGcdot \begin{pmatrix} f(x)^{-1}&0\\0&f(x)\end{pmatrix}
\MGcdot C^{-1}\MGcdot\begin{pmatrix}1&0\\0&-1\end{pmatrix}\MGcdot \oooo C\\
&=&\uuuu\sigma_0(\iiii)\MGcdot \begin{pmatrix} f(x)^{-1}&0\\0&f(x)\end{pmatrix}.
\end{eqnarray*}

(iv) Compare \eqref{16.12} with \eqref{7.11} and \eqref{7.12}.
\hfill$\Box$

\begin{remark}\label{t16.4}
Theorem \ref{t16.3} (c) (ii) says that $M_{3FN,S^1}\subset M_{3FN}^{reg}$.
The proof above uses the fact that points in $M_{3FN}^{sing}$ correspond to zeros or
poles of solutions $f$ of $P_{III}(0,0,4,-4)$ and that the solutions
on $\R_{>0}$ with values in $S^1$ have no poles or zeros.

The following is an alternative proof which does not use
isomonodromic families and solutions of $P_{III}(0,0,4,-4)$, but \eqref{16.15}.

Suppose that $(H,\nnn,x_0,x_0,P,A,J)$ is a $P_{3D6}$-TEJPA bundle in 
$M_{3FN,S^1}(x_0)\cap M_{3FN}^{[k]}(x_0)$ for some $k\in\{0,1,2,3,\}$
(and the corresponding $(\varepsilon_1,\varepsilon_2)\in\{\pm 1\}^2$).
Consider its normal form in \eqref{8.21}--\eqref{8.24}. Then
\begin{eqnarray*}
\Gamma(\P^1,\OO(H))&=& \C\MGcdot\psi_1\oplus \C\MGcdot z\psi_1,\\
A(\psi_1)&=& i\varepsilon_1\MGcdot \psi_1,\quad 
A(z\psi_1)= -i\varepsilon_1\MGcdot z\psi_1,\\
\tau(\psi_1)&=& \kappa\MGcdot z\psi_1\quad\textup{for some }\kappa\in\C^*.
\end{eqnarray*}
This and \eqref{16.15} $A\circ \tau=-\tau\circ A$ give a contradiction
by the following calculation,
\begin{eqnarray*}
-i\varepsilon_1\MGcdot \tau(\psi_1)=A(\tau(\psi_1)) = -\tau(A(\psi_1))
=i\varepsilon_1\MGcdot \tau(\psi_1).
\end{eqnarray*}
\end{remark}

One can ask about variants of the notion of TERP(0) bundles.
In particular, one can ask: are there $P_{3D6}$-TEP bundles which can be 
enriched by a $\tau$ which satisfies all properties of a TERP(0) bundle
except $\tau^2=\id$?  The question is justified by the positive answer
which singles out the $P_{3D6}$-TEJPA bundles in 
$M_{3FN,i\R_{>0}}\cup M_{3FN,i\R_{<0}}$, where
$M_{3FN,i\R_{<0}}:=R_1(M_{3FN,i\R_{>0}})=R_2(M_{3FN,i\R_{>0}})$.
This is made precise in the following theorem.

\begin{theorem}\label{t16.5}
(i) A $P_{3D6}$-TEJPA bundle $(H,\nnn,x,x,P,A,J)$ with $x\in \R_{>0}$
can be enriched by a $\tau$ which satisfies all properties in 
definition \ref{t16.1} (a) except that now $\tau^2=\id$ is replaced
by $\tau^2\neq\id$, if and only if the $P_{3D6}$-TEJPA bundle is in 
$M_{3FN,i\R_{>0}}\cup M_{3FN,i\R_{<0}}$.
Then $\tau$ is unique up to a sign. One choice is $\tau$ as in \eqref{16.21},
\begin{eqnarray}\label{16.21}
\tau(\uuuu e^\pm_0(z))=e^\mp_\iiii(1/\oooo z)\MGcdot, \quad 
\tau(\uuuu e^\pm_\iiii (z))=-e^\mp_0(1/\oooo z).
\end{eqnarray}
Here $\uuuu e^\pm_0,\uuuu e^\pm_\iiii$ is the $4$-tuple of bases in
theorem \ref{t7.3} (c).

(ii) $H$ is a pure twistor.

(iii) The $P_{3D6}$-TEJPA bundle has the normal form \eqref{8.15}--\eqref{8.18}
for $k=0$ with $f_0=f(x)$. Then $\tau$ is the map with
\begin{eqnarray}\label{16.22}
\uuuu\tau(\sigma_0(z))=\uuuu\sigma_0(1/\oooo z)\MGcdot 
\begin{pmatrix}0&f(x)^{-1}\\f(x)&0\end{pmatrix}=J(\uuuu \sigma_0(\oooo z)).
\end{eqnarray}

(iv)
\begin{eqnarray}\label{16.23}
\tau^2 &=& -\id,\\
A\circ \tau &=& \tau\circ A,\label{16.24}\\
J\circ \tau &=& \tau\circ J.\label{16.25}
\end{eqnarray}
$\tau^2=-\id$ says that $\tau$ enriches the complex structure on each
fibre $H_z$ to a quaternionic structure.
\end{theorem}

{\bf Proof:}
One can follow the proof of theorem \ref{t16.3} (a)+(b).
($\alpha$) and \eqref{16.18} are unchanged. $\tau^2\neq \id$ requires
$\varepsilon_2=-1$. This implies $\tau^2=-\id$.

($\gamma$) and \eqref{16.19} are unchanged. But ($\delta$)  changes to
\begin{eqnarray}\label{16.26}
-\begin{pmatrix}1&0\\0&\varepsilon_1\end{pmatrix} \MGcdot B \MGcdot
\begin{pmatrix}1&0\\0&\varepsilon_1\end{pmatrix} =\oooo B^{-1},
\quad\textup{together with }\varepsilon_1\MGcdot s=\oooo s.
\end{eqnarray}
In the case $\varepsilon_1=-1$ this yields $s,b_2,b_1+\frac{s}{2}b_2\in i\R$, 
which is not possible because
\[
1=(b_1+\tfrac12{s}b_2)^2+(1-\tfrac14{s^2})\MGcdot b_2^2.
\]
Thus $\varepsilon_1=1$. In this case \eqref{16.28} is \eqref{15.20},
so $(s,B)\in M_{3FN,i\R_{>0}}\cup M_{3FN,i\R_{<0}}$. \eqref{16.21} is proved.

(ii) This follows from theorem \ref{t10.3} (b)
and lemma \ref{t15.2} (b) (ii).

(iii) This is analogous to the proof of theorem \ref{t16.3} (b) (iii). 

(iv) Compare \eqref{16.21} with \eqref{7.11} and \eqref{7.12}.
\hfill$\Box$

\begin{remarks}\label{t16.6}
(i) The pairing $iP(.,\tau .):H_z\times H_{-1/\oooo z}\to \C$ is hermitian.
It induces a hermitian, nondegenerate and indefinite pairing on 
$\Gamma(\P^1,\OO(H))$ with matrix
\begin{eqnarray}\label{16.27}
iP(\uuuu\sigma^t_0(z),\tau(\uuuu \sigma_0(-1/\oooo z))) = 
\begin{pmatrix} 2if(x) & 0 \\ 0 & 2if(x)^{-1}\end{pmatrix}.
\end{eqnarray}
The pairing $P(.,\tau\circ A(.)):H_z\times H_{1/\oooo z}\to\C$
is hermitian. Its restriction to $H|_{S^1}$ and the induced pairing on 
$\Gamma(\P^1,\OO(H))$ are hermitian, nondegenerate and positive definite
in the case $M_{3FN,i\R_{>0}}(x_0)$, with matrix
\begin{eqnarray}\label{16.28}
iP(\uuuu\sigma^t_0(z),\tau\circ A(\uuuu\sigma_0(1/\oooo z))) = 
\begin{pmatrix} -2if(x) & 0 \\ 0 & 2if(x)^{-1}\end{pmatrix}.
\end{eqnarray}

(ii) Theorem \ref{t16.5} (ii) says that $M_{3FN,i\R_{>0}}\subset M_{3FN}^{reg}$.
Its proof is similar to the proof of theorem \ref{t16.3} (c) (ii).
As in remark \ref{16.4}, we can give an alternative proof which does not use
isomonodromic families and solutions of $P_{III}(0,0,4,-4)$.
But it is very different from the proof in remark \ref{16.4}.
It uses the pairing $P(.,\tau\circ A(.))$ on $H|_{S^1}$.
The fact that this pairing is hermitian and positive definite permits application of the
Iwasawa decomposition for loop groups to conclude that $H$ is a pure twistor.
We shall not give details here.

(iii) The appearance of a quaternionic structure on the vector bundle $H$ in theorem
\ref{t16.5} is not a surprise. An isomonodromic family of such bundles can be 
related via equation \eqref{15.2} or more directly (via the Sym-Bobenko formula)
to a CMC surface in $\R^3$ with rotationally symmetric metric.
The appearance of a quaternionic structure in the construction of CMC surfaces
in $\R^3$ from vector bundles with connections is well established
\cite{FLPP01}.
\end{remarks}

\chapter{Orbits of TERP structures and mixed TERP structures}\label{s17}
\setcounter{equation}{0}

\noindent
The real solutions (possibly with zeros and/or poles) of $P_{III}(0,0,4,-4)$
 on $\R_{>0}$ 
are by \eqref{15.10} the functions $f=f_{mult}(.,s,B)|_{\R_{>0}}$ for $(s,B)\in V^{mat,\R}$.
They will be studied in detail in chapter \ref{s18}.

By theorem \ref{t10.3} any solution corresponds to an isomonodromic family 
of $P_{3D6}$-TEJPA bundles $(H(x),\nnn,x,x,P,A,J)$ for $x\in \R_{>0}$ with 
monodromy data $(s,B)$. By theorem \ref{t16.3} (b) with the $\tau$ in \eqref{16.9},
this can be enriched to an isomonodromic family of TERP(0) bundles.
Then $H(x)$ is pure if and only if $x$ is not a zero or a pole of $f$,
and then $H(x)$ is polarized if $f(x)>0$, and it has negative definite pairing $S$ if 
$f(x)<0$. 

The 1-parameter isomonodromic families $(H(x),\nnn,P,\tau)$ for $x\in \R_{>0}$ 
turn out to be special cases of 1-parameter isomonodromic families of 
TERP(0) bundles which had been studied in \cite{HS07}, \cite{Mo11b} and which
are called Euler orbits in definition \ref{t17.1} below.  
The situation where all members with large $x$ or when all members with small $x$
are pure and polarized was investigated in In \cite{HS07}, \cite{Mo11b}.
In our case this corresponds to smoothness and positivity
of $f$ near $\iiii$ or $0$. 

The characterization for large $x$ will be stated below
in the semisimple case, and the characterisation for small $x$ will be discussed 
informally. This will be preceded by the definition of Euler orbits of TERP(0) bundles
and a discussion of semisimple TERP(0) bundles.
At the end of the chapter a result from \cite{HS11} will be explained in the special
case of semisimple TERP(0) bundles. It provides semisimple TERP(0) bundles such that
they and all their semisimple isomonodromic deformations are pure and polarized.
This will give all solutions $f$ of $P_{III}(0,0,4,-4)$ on $\R_{>0}$ which are 
smooth and positive on $\R_{>0}$.

\begin{definition}\label{t17.1}
Let $(H,\nnn,P,\tau)$ be a TERP(0) bundle. Let
\begin{eqnarray}\label{17.1}
\pi_{orb}:\C\times\C^*\to\C,\quad (z,r)\mapsto z\MGcdot r^{-1}.
\end{eqnarray}
Then $(\pi_{orb}^*(H|_\C,\nnn,P),\pi_{orb}^*(H'_\R))$ is the almost entirely holomorphic
object of a variation on $M:=\C^*$ of TERP(0) bundles mentioned in remark \ref{t16.2} (iii).
$\tau$ on $\pi_{orb}^*(H|_\C)|_{S^1\times M}$ 
is defined as complex conjugation with respect to the 
real structure $\pi_{orb}^*(H'_\R)|_{S^1\times M}$, and it is extended flatly to
$\C^*\times M$. Then $\pi_{orb}^*(H|_\C)$ and $\oooo{\gamma^*(\pi_{orb}^*(H|_\C))}$
are glued by $\tau$ to a bundle $G\to\P^1\times M$, which is holomorphic
in $z\in\P^1$, but only real analytic with respect to $r\in M=\C^*$.

$(G,\nnn,P,\tau)$ is a special variation of TERP(0) bundles with parameter
space $M=\C^*$ 
which is called here an {\it Euler orbit}.
(See \cite[ch.\  2]{He03} or remark \ref{t16.2} (iii) for the general
notion of a variation of TERP(0) bundles --- in \cite{He03} this is called a TERP structure.)

The single TERP(0)  bundle for $r\in M=\C^*$ is $(G(r),\nnn,P,\tau)$ with
$G(r)=G|_{\P^1\times\{r\}}$. It is obtained by glueing $\pi_r^*(H|_\C)$ and $\oooo{\gamma^*(\pi_r^*(H|_\C))}$
via $\tau$, here
\begin{eqnarray}\label{17.2}
\pi_r:\C\to\C,\quad z\to z\MGcdot r^{-1}.
\end{eqnarray}
Of course $(G(1),\nnn,P,\tau)=(H,\nnn,P,\tau)$.
\end{definition}

\begin{remarks}\label{t17.2}
(i) If $u^1_0,\dots,u^n_0$ are the eigenvalues of the pole part 
$[z\nnn_{\zdz}]:H_0\to H_0$ at $0$, then $r u^1_0,\dots,r u^n_0$ are the 
eigenvalues of the pole part of $(G(r),\nnn)$ at $0$.

(ii) In the semisimple case the Euler orbit of a TERP(0) bundle is simply the 
isomonodromic 1-parameter family 
of TERP(0) bundles $(G(r),\nnn,P,\tau)$ for $r\in\C^*$ 
such that for fixed $r$ the eigenvalues of the pole part
at $0$ are $r u^1_0,\dots,ru^n_0$ and the eigenvalues of the pole part 
$[-\nnn_{\paa_z}]:G(r)_\iiii\to G(r)_\iiii$ are $\oooo{ru^1_0},\dots,\oooo{ru^n_0}$.

(iii) Any member $(G(r),\nnn,P,\tau)$ of the Euler orbit of a TERP(0) bundle
$(H,\nnn,P,\tau)$ induces up to a rescaling of the parameter space the same Euler orbit.

(iv) The definition of the pull-back $\pi_{orb}^*(H|_\C)=G|_{\C\times \C^*}$ implies the existence of
canonical isomorphisms
\begin{eqnarray*}
G(r_1)_{z_1}\cong H_{z_1r_1^{-1}}= H_{z_2r_2^{-1}}\cong G(r_2)_{z_2}
\quad \textup{if }z_1r_1^{-1}=z_2r_2^{-1}.
\end{eqnarray*}
It is shown in \cite[lemma 4.4]{HS07} that these isomorphisms extend
to $z_1=z_2=\iiii$ if $|r_1|=|r_2|$ and induce a bundle isomorphism
$G(r_1)\cong G(r_2)$ over the automorphism 
$\P^1\to\P^1,\ z\mapsto z\frac{r_2}{r_1}$.

Therefore the restriction of the variation $(G,\nnn,P,\tau)$ on $\C^*$
to any circle $\{r\in\C^*\, |\, |r|=r_0\}$, $r_0\in\R_{>0}$,
is a (rather trivial) variation of TERP(0) bundles.
All bundles are pure/polarized if one bundle is pure/polarized.
The interesting part of the variation $(G,\nnn,P,\tau)$
on $\C^*$ is its restriction to $\R_{>0}\subset\C^*$.

(v) Because of (iv) and remark \ref{t16.2} (vi) and (vii),
for the study of semisimple rank 2 TERP(0) bundles we can restrict
to the case $u^1_0=u^1_\iiii=x\in\R_{>0}$,
$u^2_0=u^2_\iiii=-x$.
Such a TERP(0) bundle is a $P_{3D6}$-TEP bundle $(H,\nnn,x,x,P)$ with additional
real structure given by $\tau$. 
Theorem \ref{t16.3} relates such a TERP(0) bundle to a value $f(x)$ of a
solution $f=f_{mult}(.,s,B)|_{\R_{>0}}$  of $P_{III}(0,0,4,-4)$ on $\R_{>0}$ with
$(s,B)\in V^{mat,\R\cup S^1}$.

By theorem \ref{t16.3}, 
in the case of $(s,B)\in V^{mat,S^1}$ the TERP(0) bundle is pure and has indefinite
pairing $S$. In the case of $(s,B)\in V^{mat,\R}$ it is pure if and only if $x$ is not a 
zero or pole of $f$, otherwise $\OO(H)\cong\OO_{\P^1}(1)\oplus \OO_{\P^1}(-1)$.
If $f(x)>0$ it is pure and polarized, if $f(x)<0$ it is pure with negative definite
pairing $S$.

Theorem \ref{t16.3} and the results in chapter \ref{s18} give a complete picture
of the semisimple rank 2 TERP(0) bundles and their Euler orbits.

(vi) By \eqref{15.10} and theorem \ref{t16.3} and remark \ref{t17.2} (ii),
the solutions $f$  of $P_{III}(0,0,4,-4)$ on $\R_{>0}$ with values in $\R$ or 
in $S^1$ correspond to the restrictions to $\R_{>0}$ of the Euler orbits
of the TERP(0) bundles with $u^1_0=u^1_\iiii=x\in\R_{>0}, \ 
u^2_0=u^2_\iiii=-x$. Results on their Euler orbits are equivalent to results
on solutions $f$.

(vii) The TERP structures in \cite{He03} give a framework for the data studied in
\cite{CV91}, \cite{CV93}, \cite{Du93}, which are essentially certain variations of
TERP bundles. 
The Euler orbits appear in \cite{CV91}, \cite{CV93} from the renormalisation group flow.
There the limits $r\to\iiii$ and $r\to 0$ are called {\it infrared limit} and 
{\it ultraviolet limit}.
\end{remarks}

\begin{definition}\label{t17.3} \cite[def. 4.1]{HS07}
Let $(H,\nnn,P,\tau)$ be a TERP(0) bundle.
Its Euler orbit $(G(r),\nnn,P,\tau)$, $r\in\C^*$, is called a {\it nilpotent orbit}
of TERP(0) bundles if, for all large $|r|$, $(G(r),\nnn,P,\tau)$ is a pure and polarized
TERP(0) bundle. 
Its Euler orbit is called a {\it Sabbah orbit} of TERP(0) bundles 
if, for all small $|r|$, $(G(r),\nnn,P,\tau)$ is a pure and polarized TERP(0) bundle.
\end{definition}

\begin{remarks}\label{t17.4}
(i) A real solution (possibly with zeros and/or poles) on $\R_{>0}$ of
$P_{III}(0,0,4,_4)$ is positive for large $x$ (respectively, small $x$)
if and only if its Euler orbit of TERP(0) bundles is a nilpotent orbit 
(respectively, a Sabbah orbit).

(ii) \cite[theorem 7.3]{HS07} gives a precise characterisation of Sabbah orbits:

{\it The Euler orbit of a TERP(0) bundle is a Sabbah orbit if a certain candidate for
a sum of two polarized mixed Hodge structures (briefly: PMHS) is indeed a sum
of two PMHS.}

The definition of the candidate starting from the TERP(0) bundle is lengthy.
It uses a variant due to Sabbah of the Kashiwara-Malgrange V-filtration, and the 
definition of the polarizing form needs special care. Also the definition of a PMHS
is nontrivial. All definitions can be found in \cite{HS07}, so we shall not reproduce them here. 

(iii) Another reason for omitting the details in the characterization of 
Sabbah orbits is that the application here would be a characterization
of those $(s,B)\in V^{mat,\R}$ for which $f_{mult}(.,s,B)|_{\R_{>0}}$ is positive
for small $x$. But chapter \ref{s12} settles this completely and provides richer 
information, which we cannot extract easily from \cite[theorem 7.3]{HS07}.
See theorem \ref{t18.2} (c)+(d) for the details.

(iv) In fact, for $(s,B)\in V^{mat,\R}$ $f_{mult}(.,s,B)$ is positive for small $x$
if and only if $|s|\leq 2$ and $b_1+\frac{s}{2}b_2\in\R_{\geq 1}$.
If $|s|<2$ then $\Mon(s)$ is semisimple and the candidate in \cite[theorem 7.3]{HS07}
is a pure polarized Hodge structure. But if $|s|=2$ then $\Mon(s)$ has a $2\times 2$
Jordan block and the candidate is a polarized (and truly) {\it mixed} Hodge structure.
\end{remarks}

Now we shall review the notion of semisimple TERP(0) bundles and their monodromy data,
following \cite[ch.\  8,10]{HS07}. Then we shall define semisimple mixed TERP structures
and formulate the characterization of nilpotent orbits. Finally a result of \cite{HS11}
will show that for certain semisimple mixed TERP(0)  bundles they and all their
semisimple isomonodromic deformations are pure and polarized TERP(0) bundles.

First we cite and explain a lemma from \cite{HS07} on semisimple TEP bundles.
Recall (remark \ref{t6.2} (i)) that the TEP structures of weight $0$ in \cite{HS07}
are the restrictions to $\C$ of the TEP bundles in definition \ref{t6.1} (a)
(which are defined on $\P^1$).

Let us call two matrices $T$ and $T'$ in $M(n\times n,\C)$ {\it sign equivalent} 
if there is a matrix $B=\textup{diag}(\varepsilon_1,\dots,\varepsilon_n)$ with
$\varepsilon_1,\dots,\varepsilon_n\in\{\pm 1\}$ such that $BTB=T'$.

\begin{lemma}\label{t17.5} \cite[lemma 10.1 (1)]{HS07}
Fix pairwise distinct values $u_1,\dots,u_n\in\C$ and $\xi\in S^1$ with 
$\Re(\tfrac{u_i-u_j}{\xi})<0$ for $i<j$. 
There is a natural $1{:}1$ correspondence between the set of restrictions to $\C$
of semisimple TEP bundles with pole part at $0$ having eigenvalues $u_1,\dots,u_n$,
and the set of sign equivalence classes of upper triangular matrices 
$T\in Gl(n,\C)$ with diagonal entries equal to $1$. 
The matrices $T$ are the Stokes matrices of the pole at $0$ of the TEP bundle.
\end{lemma}

A similar statement for the case without pairing holds and is a special case
of a Riemann-Hilbert correspondence between germs of holomorphic vector
bundles with merormorphic connections and their Stokes data, see for example
\cite{Si67}, \cite{Si90}, \cite{BJL79}, \cite{Ma83a}, \cite[ch.\  II 5,6]{Sa02}, 
\cite{Bo01}, \cite[theorem 4.3.1]{Mo11a}. 
There one has two Stokes matrices and $n$ {\it exponents} in $\C$ 
determining the regular singular
rank 1 pieces in the formal isomorphism class. In the case of a TEP bundle,
the $n$ exponents are all equal to $0$, and the second Stokes matrix is the transpose
of the first. 

Let us explain the $1{:}1$ correspondence in more detail.
The proof can be found in \cite{HS07}.

The formal decomposition of Turrittin works in the semisimple case without ramification
(e.g. \cite[II theorem 5.7]{Sa02}) and gives a formal isomorphism
\begin{align}\label{17.3} 
\Psi: (\OO(H)_0&[z^{-1}],\nnn)\otimes_{\C\{z\}[z^{-1}]}\C[[z]][z^{-1}]
\\ \nonumber
&\to\oplus_{j=1}^n e^{-u_j/z}\otimes 
(\RR_j,\nnn_j)\otimes _{\C\{z\}[z^{-1}]}\C[[z]][z^{-1}] 
\end{align}
where $\RR_j$ is a $\C\{z\}[z^{-1}]$ vector space of dimension 1 and $\nnn_j$ is a 
meromorphic connection on it with a regular singular pole at $0$.

With the general results in \cite[II 5. and III 2.1]{Sa02} it follows easily
that $\Psi$ extends to a formal isomorphism
\begin{align}\label{17.4} 
\Psi: (\OO(H)_0,\nnn,&P)\otimes_{\C\{z\}}\C[[z]]
\\ \nonumber
&\to\oplus_{j=1}^n e^{-u_j/z}\otimes 
(\OO(H_j)_0,\nnn_j,P_j)\otimes _{\C\{z\}}\C[[z]] 
\end{align}
of germs at $0$ of TEP bundles \cite[lemma 8.2]{HS07}.
Here $(\OO(H_j)_0,\nnn_j,P_j)$ is the germ at $0$ of a rank 1 TEP bundle with regular
singular pole at $0$ and with $\RR_j=\OO(H_j)_0[z^{-1}]$, 
and $e^{-u_j/z}\otimes \OO(H_j)_0$ is the germ at $0$ whose 
holomorphic sections are obtained by multiplying those in $\OO(H_j)_0$ with
$e^{-u_j/z}$. Then $e^{-u_j/z}\otimes (\OO(H_j)_0,\nnn_,P_j)$ is the germ
at $0$ of a rank 1 TEP bundle with eigenvalue $u_j$ of the pole part at $0$.

We remark that up to isomorphism there is only one germ at $0$ of a rank 1 TEP bundle
with regular singularity at $0$: the regular singularity leads to a generating section
$z^\alpha\MGcdot e_0$ where $e_0$ is a flat multi-valued section on $\C^*$,
the pairing $P$ implies $\alpha=0$ and the single-valuedness of $e_0$ and
$P(e_0(z),e_0(-z))=\textup{constant}\in\C^*$.

Denote by $f_j$ a flat generating section with $P(f_j(z),f_j(-z))=1$ of the
germ $(\OO(H_j)_0,\nnn_j,P_j)$. It is unique up to sign.

For any $n$ distinct values $u_1,\dots,u_n\in\C$ the set
\[
\Sigma:=\{\xi\in\C\, |\, \exists\ i,j\textup{ with }i\neq j\textup{ and }
\Re(\tfrac{u_i-u_j}{\xi})=0\}
\]
of {\it Stokes directions} is finite. 
For any $\xi\in S^1-\Sigma$ one can renumber the $n$ distinct values
such that then $\Re(\tfrac{u_i-u_j}{\xi})<0$ for $i<j$ holds.
The choice of a $\xi$ and a numbering of $u_1,\dots,u_n$ with this property
is assumed in lemma \ref{t17.5} and in the following.
Let $I^a_0$ (respectively, $I^b_0$) be the component of $S^1-\Sigma$ which 
contains $\xi$ (respectively, $-\xi$), and denote
\[
I^\pm_0:= I^a_0\cup I^b_0\cup\{z\in S^1\, |\, 
\pm \Im({z}/{\xi})\leq 0\}.
\]
This generalizes the notation in chapter \ref{s2} for the rank 2 case.
Each of the sets $I^+_0$ and $I^-_0$ contains exactly one of the two 
Stokes directions $\pm\xi'$ for any $\xi'\in \Sigma$.

Denote by $\AAA$ the sheaf on $S^1$ of holomorphic functions in neighbourhoods
of $0$ in sectors which have an asymptotic expansion in $\C[[z]]$ in the 
sense of \cite[ch.\  3]{Ma83a}. 

A result of Hukuhara and many others (e.g. \cite[ch.\  3-5]{Ma83a},
\cite[II 5.12]{Sa02}) says in our case that the formal isomorphism $\Psi$
in \eqref{17.4} lifts in the sectors $\widehat I^\pm_0$ to unique
isomorphisms $\Psi^\pm$ with coefficients in $\AAA|_{I^\pm_0}$,
\begin{align}\label{17.5} 
\Psi^\pm: (\OO(H)_0&,\nnn)\otimes_{\C\{z\}}\AAA|_{I^\pm_0}
\\
\nonumber
&\to\oplus_{j=1}^n e^{-u_j/z}\otimes 
(\OO(H_j)_0,\nnn_j)\otimes _{\C\{z\}}\AAA|_{I^\pm_0}
\end{align}
which together respect the pairing in the following sense:
The preimages $e^\pm_j:=(\Psi^\pm)^{-1}(f_j)$ are together for $j=1,\dots,n$
flat bases $\uuuu e^\pm=(e^{\pm 1}_0,\dots,e^{\pm n}_0)$ of $H|_{\whhh I^\pm_0}$ with
\begin{eqnarray}\label{17.6}
P(\uuuu e^+_0(z)^t, e^-_0(-z))={\bf 1}_n
=P(\uuuu e^-_0(-z)^t, e^+_0(z))\quad \textup{for }z\in \whhh I^+_0
\end{eqnarray}
\cite[lemma 8.4]{HS07}. The base change matrix which is defined by 
\eqref{17.7} is upper triangular with $1$'s on the diagonal
and, by \eqref{17.6} and \eqref{17.8}, satisfies 
\begin{eqnarray}\label{17.7}
\uuuu e^-_0|_{I^a_0} &=& \uuuu e^+_0|_{I^a_0}\MGcdot T,\\
\uuuu e^-_0|_{I^b_0} &=& \uuuu e^+_0|_{I^b_0}\MGcdot T^t \label{17.8}
\end{eqnarray}
\cite[lemma 8.3]{HS07}. The sign equivalence class of $T$ is unique.

For any $T$ which is upper triangular with $1$'s on the diagonal
a unique germ $(\OO(H)_0,\nnn,P)$ exists \cite[lemma 10.1 (1)]{HS07}.

Now we turn to TERP(0) bundles.

\begin{definition}\label{t17.6}
Let $(H,\nnn,P,\tau)$ be a a semisimple TERP(0) bundle with
pairwise distinct eigenvalues $u_1,\dots,u_n\in\C$ of the pole part at $0$
and with a $\xi\in S^1$ such that $\Re(\tfrac{u_i-u_j}{\xi})<0$ for $i<j$.

(a) {\it The real structure and Stokes structure are compatible} if 
$\lambda_1,\dots,\lambda_n\in S^1$ exist such that 
$(\lambda_1 e^{\pm 1}_0,\dots,\lambda_n e^{\pm n}_0)$ is a flat basis of 
$H'_{\R}|_{\whhh I^\pm _0}$, i.e., if the splitting
$\oplus_{j=1}^n \C\MGcdot e^{\pm j}_0$ of the flat bundle on $\whhh I^\pm_0$
is compatible with the real structure.

(b) Remark: Then $\lambda_j\in\{\pm 1,\pm i\}$ because by \eqref{16.3} 
for $z\in I^+_0$
\begin{eqnarray}\label{17.9}
\oooo{\lambda_j^2}=\oooo{P(\lambda_j e^{+ j}_0(z),\lambda_j e^{- j}_0(-z))}
= P(\lambda_j e^{+ j}_0(z),\lambda_j e^{- j}_0(-z)) =\lambda_j^2.
\end{eqnarray}

(c) The TERP(0) bundle is a {\it mixed TERP structure} if in (a) all 
$\lambda_j\in\{\pm 1\}$.
\end{definition}

\begin{remarks}\label{t17.7}
(i) The Stokes structure is encoded in the splittings
$\oplus_{j=1}^n \C\MGcdot e^{\pm j}_0$ of the flat
bundle $H|_{\whhh I^\pm_0}$ and the Stokes matrix $T$,
which depend on the choice of $\xi$.
One can reformulate definition \ref{t17.6} (a) in a way which does not depend
on this choice \cite[8.1.1]{Mo11b}. It says that the {\it Stokes filtration}
and the real structure are compatible.

(ii) An equivalent condition for the Stokes structure and the real structure to be compatible is that
the real structure of the TERP(0) bundle induces on the pieces
$e^{-u_j/z}\otimes (\OO(H_j)_0,\nnn_j,P_j)$ a natural real structure 
such that these pieces become TERP(0) bundles. Then $\lambda_jf_j$ is
a generating flat real section on $\C^*$. 
The TERP(0) bundle $(H_j,\nnn_j,P_j,\tau_j)$ is automatically pure.
It is polarized if $P(\lambda_jf_j,\lambda_jf_j)=1$, and that holds
if $\lambda_j\in\{\pm 1\}$. 
The alternative is $P(\lambda_jf_j,\lambda_j,f_j)=-1$ and
$\lambda_j\in\{\pm  i\}$.
This motivates the definition of a mixed TERP structure in part (c) above.

(iii) In the non-semsimple case \eqref{17.4}--\eqref{17.8} hold
mutatis mutandis. Then the pieces of the decomposition do not all have rank 1.
Part (a) of definition \ref{t17.6} goes through, and if it holds the pieces
become TERP(0) bundles. But part (c) has to be replaced by the condition that
for each regular singular piece a certain candidate for a sum of two PMHS
is indeed a sum of two PMHS. See \cite{HS07} for the definition
and explanations.

(iv) Unfortunately, in the proof of \cite[lemma 10.1 (2)]{HS07},
the case $\lambda_j\in\{\pm i\}$ and the fact that it is excluded in a mixed TERP structure because of the polarization of the PMHS, had been missed.
(Furthermore, in (10.3) and (10.4) $(A^-)^{tr}$ has to be replaced
by $(A^-)^{-1}$;  they follow from (10.1) instead of (10.2), and in
(10.4) the sign $(-1)^w$ has to be deleted.)

(v) The pole at $\iiii$ of a TEP bundle or a TERP(0) bundle can be 
described analogously to the explanations after lemma \ref{t17.5}.
We consider only the case of a TERP(0) bundle. Then the eigenvalues of
the pole part $[-\nnn_{\paa_z}]:H_\iiii\to H_\iiii$ are 
$\oooo{u_1},\dots,\oooo{u_n}$ (remark \ref{t16.2} (vi)).
Then for $u_1,\dots,u_n$ and $\xi$ as above one can define
\begin{eqnarray}\label{17.10}
I^a_\iiii:= I^a_0,\ I^b_\iiii:= I^b_0,\ I^+_\iiii:= I^-_0,\ I^-_\iiii:=I^+_0.
\end{eqnarray}
Because of the involution $\tau$, the pole at $\iiii$ is a twin of
the pole at $0$, and the analogously defined flat bases 
$\uuuu e^\pm_\iiii$ on $\whhh I^\pm_\iiii$ are
\begin{eqnarray}\label{17.11}
\uuuu e^\pm_\iiii = \tau(\uuuu e^\mp_0)
\end{eqnarray}
with base change matrices
\begin{eqnarray}\label{17.12}
\uuuu e^-_\iiii|_{I^a_\iiii} &=& \uuuu e^+_\iiii|_{I^a_\iiii}\MGcdot
\oooo T^{-1},\\
\uuuu e^-_\iiii|_{I^b_\iiii} &=& \uuuu e^+_\iiii|_{I^b_\iiii}\MGcdot
(\oooo T^{-1})^t. \label{17.13}
\end{eqnarray}
Define the connection matrix $B$ by
\begin{eqnarray}\label{17.14}
\tau(\uuuu e^+_0)=\uuuu e^-_\iiii = \uuuu e^+_0\MGcdot B.
\end{eqnarray}
Then definition \ref{t17.6} says that 
$B=\textup{diag}(\varepsilon_1,\dots,\varepsilon_n)$ with 
$\varepsilon_1,\dots,\varepsilon_n\in\{\pm 1\}$ if and only if 
Stokes structure and real structure are compatible, and that 
\begin{align}\label{17.15}
B={\bf 1}_n \iff  \ &\text{the TERP(0) bundle}
\\
\nonumber
&\textup{is a mixed TERP structure.}
\end{align}
\end{remarks}

\begin{corollary}\label{t17.8}
\cite[lemma 10.1 (2)]{HS07}
Fix pairwise distinct values $u_1,\dots,u_n\in\C$ and $\xi\in S^1$ with 
$\Re(\tfrac{u_i-u_j}{\xi})<0$ for $i<j$. 
The $1{:}1$ correspondence in lemma \ref{t17.5} restricts to a $1{:}1$ correspondence
between the set of semisimple mixed TERP structures with pole part at $0$ 
having eigenvalues $u_1,\dots,u_n$, and the set of sign equivalence classes
of upper triangular matrices $T\in GL(n,\R)$ with diagonal entries equal to $1$.
\end{corollary}

{\bf Proof:} In the case of a semisimple mixed TERP structure the flat bases
$\uuuu e^+_0=\uuuu e^-_\iiii$ and $\uuuu e^-_0=\uuuu e^+_\iiii$
of $H|_{\whhh I^\pm_0}= H|_{\whhh I^\mp_\iiii}$ are real bases,
thus $T$ in \eqref{17.7} has real entries.

Conversely, if one starts with the restriction to $\C$ of a semisimple TEP bundle
with matrix $T$ with real entries, the real structure with real bases
$\uuuu e^+_0$ and $\uuuu e^-_0$ is well defined and compatible with the
Stokes structure, and satisfies $B={\bf 1}_n$.\hfill$\Box$

The following theorem generalizes to all TERP(0) bundles whose formal 
decomposition at $0$ is valid without ramification. This generalization was
conjectured in \cite[conjecture 9.2]{HS07}.
The simpler direction $\Leftarrow$ and the regular singular case of the
direction $\Rightarrow$ were proved in \cite[theorem 9.3]{HS07}, building on \cite{Mo03}.
The general case of the more difficult direction $\Rightarrow$
was proved by T. Mochizuki \cite[corollary 8.15]{Mo11b}, building on \cite{Mo11a}.

\begin{theorem}\label{t17.9}
A semisimple TERP(0) bundle induces a nilpotent orbit $\iff$
it is a mixed TERP structure.
\end{theorem}

\begin{corollary}\label{t17.10}
A real solution $f_{mult}(.,s,B)|_{\R_{>0}}$ for $(s,B)\in V^{mat,\R}$ 
(possibly with zeros and/or poles)  of $P_{III}(0,0,4,-4)$ on $\R_{>0}$
is smooth and positive for large $x$ if and only if $B={\bf 1}_2$.
\end{corollary}

{\bf Proof:} This follows from remark \ref{t17.4} (i), theorem \ref{t17.9}
and the formulae \eqref{17.14} and \eqref{17.15}.\hfill$\Box$

The following theorem is a special case of a result in \cite{HS11}.
It is relevant for globally smooth real solutions  of
$P_{III}(0,0,4,-4)$ on $\R_{>0}$.

\begin{theorem}\label{t17.11} \cite[theorem 5.9]{HS11}
Let $(H,\nnn,P,\tau)$ be a semisimple TERP(0) bundle with pairwise different
eigenvalues $u_1,\dots,u_n\in\C$ of the pole part at $0$ and 
with $\xi\in S^1$ such that $\Re(\tfrac{u_i-u_j}{\xi})<0$ for 
$i<j$, and suppose that the Stokes matrix $T$ has real entries,
so that the TERP(0) bundle is also a mixed TERP structure.

Then all TERP(0) bundles in its Euler orbit are pure and polarized
if $T+T^t$ is positive semidefinite.
\end{theorem}

\begin{corollary}\label{t17.12}
A real solution $f_{mult}(.,s,B)|_{\R_{>0}}$ for $(s,B)\in V^{mat,\R}$
 of $P_{III}(0,0,4,-4)$ on $\R_{>0}$ is everywhere smooth and positive
if $|s|\leq 2$ and $B={\bf 1}_2$.
\end{corollary}

{\bf Proof:} Suppose that $|s|\leq 2$ and $B={\bf 1}_2$. 
Then the TERP(0) bundles in the Euler orbit for the solution 
$f_{mult}(.,s,{\bf 1}_2)|_{\R_{>0}}$
are mixed TERP structures because of $B={\bf 1}_2$, and 
their Stokes matrix ist 
\[
T=S=\begin{pmatrix}1&s\\0&1\end{pmatrix}.
\]
Then
\[
T+T^t=\begin{pmatrix}2&s\\s&2\end{pmatrix}
\]
is positive semidefinite.
By theorem \ref{t17.11} all TERP(0) bundles are pure and polarized.
By theorem \ref{t16.3} (b) (iii) $f_{mult}(x,s,{\bf 1}_2)>0$ for all $x\in\R_{>0}$.
\hfill$\Box$

\begin{remarks}\label{t17.13}
(i) $B={\bf 1}_2$ is necessary for smoothness and positivity for large $x$
(corollary \ref{t17.10}), $|s|\leq 2$ is necessary for smoothness near 0
(theorem \ref{t18.2} (c)), therefore $|s|\leq 2$ and $B={\bf 1}_2$ are
also necessary conditions for global smoothness and positivity
(theorem \ref{t18.2} (e)). 

(ii) \cite[theorem 5.9]{HS11} is more general in several aspects,
but it claims only that the original TERP(0) bundle is pure and polarized.
But one can go easily to all TERP(0) bundles in the Euler orbit,
because the eigenvalues of their pole parts at $0$ are of the form
$ru^1_0,\dots,ru^n_0$ for $r\in\R_{>0}$ and the Stokes matrix $T$ and the
connection matrix $B$ do not change (remark \ref{t16.2} (iii)).

(iii) In fact, theorem \ref{t17.11} together with the well known behaviour
of the Stokes matrix under the change of the tuple $(u^1_0,\dots,u^n_0,\xi)$
show that all semisimple TERP structures which are in the isomonodromic
family of one TERP structure in theorem \ref{t17.11} are pure and polarized.
If the tuple $(u^1_0,\dots,u^n_0,\xi)$ moves, then in general the Stokes directions change,
$\xi$ crosses some Stokes directions, and the values $u^1_0,\dots,u^n_0$
have to be renumbered. Then the Stokes matrix $T$ changes by some braid group
action to a new Stokes matrix $\www T$ \cite[theorem 4.6]{Du99}.
But $\www T$ inherits from $T$ the property that $\www T+\www T^t$ 
is positive semidefinite.
As we do not need this for the Euler orbits, we shall not give details here.

(iv) The data in \cite[theorem 5.9]{HS11} differ in several ways from
the TERP(0) bundles. First, they are more general:  instead of a real structure
and a $\C$-bilinear pairing one sesquilinear pairing is considered.
Second, the structures there do not have weight 0, but weight 1.
Therefore the second Stokes matrix and the formal eigenvalues there differ 
by a sign from the second Stokes matrix and the formal eigenvalues here. 
Third, not only the semisimple case is considered.
Finally, there a {\it minimality condition} is assumed. But in the semisimple case
it is trivial, then $K_c=0$ holds for all $c$ ($=u^1_0,\dots,u^n_0$ here).
$K_c\neq 0$ can arise only if the monodromy of the corresponding piece of the formal 
decomposition of a TERP(0) bundle has $-1$ as an eigenvalue. 
In the semsimple case all pieces have rank 1, and all formal eigenvalues are
equal to 1.
\end{remarks}

\chapter{Real solutions  of $P_{III}(0,0,4,-4)$ on $\R_{>0}$}\label{s18}
\setcounter{equation}{0}

\noindent
The real solutions (possibly with zeros and/or poles) of $P_{III}(0,0,4,-4)$ on $\R_{>0}$ 
are by \eqref{15.10} the functions $f_{mult}(.,s,B)|_{\R_{>0}}$ for 
$(s,B)\in V^{mat,\R}$. In this chapter we shall obtain complete results about the
sequences of zeros and/or poles of these solutions. 
These global results will be derived by combining, 
first, the local behaviour for small $x$ and for large $x$, and second, the geometry
of the spaces $M_{3FN,\R}^{mon}$ and $M_{3FN,\R}^{ini}=M_{3FN,\R}^{reg}\cup
M_{3FN,\R}^{sing}$, which was described in lemma \ref{t15.4} (a) and theorem
\ref{t15.5} (a).

This kind of argument, from the local behaviour of the solutions near $0$ and $\iiii$ and 
from the geometry of certain moduli spaces to the global behaviour of the solutions,
seems to be new in the theory of the Painlev\'e III equations.

Theorem \ref{t18.2} collects the known results on the behaviour of individual solutions
near $0$ and $\iiii$. Theorem \ref{t18.3} develops how the zeros and/or poles behave
in certain families of solutions. Theorem \ref{t18.4} derives from this the global
results on the sequences of zeros and/or poles of all solutions.

The sources for the local results in theorem \ref{t18.2} are:
\cite{MTW77}, \cite[ch.\  11]{IN86} and \cite[ch.\  15]{FIKN06}, \cite{Ni09} and 
chapter \ref{s12}, and chapter \ref{s17} which builds on 
\cite{HS07}, \cite{Mo11b}, \cite{HS11}.
Some of the local results can be derived from several sources. 
The proof of theorem \ref{t18.2}
will explain this. Theorem \ref{t18.2} has a qualitative character and does not rewrite
the precise asymptotic formulae in \cite[ch.\  11]{IN86}, \cite{Ni09} and chapter \ref{s12}.
It is a feature of the argument leading to theorem \ref{t18.4} that it does not
depend on precise asymptotic formulae.

\begin{notation}\label{t18.1}
The following notation allows a concise formulation of the local and global
results about sequences of zeros and/or poles of the real solutions 
of $P_{III}(0,0,4,-4)$ on $\R_{>0}$. Recall that there are two types of zeros and two types of poles.
If $x_0$ is a zero (respectively, a pole) of  a solution $f$ then $\paa_x f(x_0)=\pm 2$
(respectively, $\paa_x(f^{-1})(x_0)=\pm 2$). A zero $x_0$ with $\paa_x f(x_0)=\pm 2$
is denoted $[0\pm]$, a pole $x_0$ with $\paa_x (f^{-1})(x_0)=\pm 2$ is denoted
$[\iiii \pm]$.

For a pair of subsequent zeros or poles of a solution $f$, the value of $f$ is positive between
them if they are of one of the types
$$[0+][0-],\ [0+][\iiii-],\ [\iiii+][0-],\ [\iiii+][\iiii-]$$
and negative if they are of one of the types
$$[0-][0+],\ [0-][\iiii+],\ [\iiii-][0+],\ [\iiii-][\iiii+].$$
Other pairs are not possible.

Let $f$ be a real solution (possibly with zeros and/or poles) of
$P_{III}(0,0,4,-4)$ on $\R_{>0}$. Then $f$ is of type $\olll{>0}$ (respectively, $\olll{<0}$)
if a $y_0\in \R_{>0}$ exists such that $f(y)>0$ (respectively, $f(y)<0$) for all
$y\in (0,y_0)$. Then also $f|_{(0,y_0)}$ is called of type $ \olll{>0}$ (respectively, $\olll{<0}$).

$f$ is of type $\olll{[0+][0-]}$ if a $y_0\in\R_{>0}$ exists such that $f(y_0)\neq 0$,
$f$ has infinitely many zeros or poles in $(0,y_0)$ which are denoted 
$x_1,x_2,x_3,\dots$ with $x_1>x_2>x_3>\dots$, and $x_k$ is of type $[0-]$ for odd $k$ 
and of type $[0+]$ for even $k$. 
Then also $f|_{(0,y_0]}$ is called of type $\olll{[0+][0-]}$.
The types $\olll{[0-][0+]}$, $\olll{[\iiii+][\iiii-]}$  and $\olll{[\iiii-][\iiii+]}$
 are defined analogously.
Of course $f$ is of type $\olll{[0+][0-]}$ if it is of type $\olll{[0-][0+]}$.
But for $f|_{(0,y_0]}$ the type of the largest zero is important.

For large $x$ the types $\orrr{>0}$, $\orrr{<0}$, $\orrr{[0+][\iiii-]}$, 
$\orrr{[\iiii-][0+]}$, $\orrr{[0-][\iiii+]}$ and $\orrr{[\iiii+][0-]}$ are defined analogously.
For example $f|_{[y_0,\iiii)}$ and $f$ are of type $\orrr{[0+][\iiii-]}$ if $f(y_0)\neq 0$,
if $f$ has infinitely many zeros and poles in $(y_0,\iiii)$ which are denoted 
$x_1,x_2,x_3,\dots$ with $x_1<x_2<x_3<\dots$ and if $x_k$ is of type $[0+]$ for odd $k$
and of type $[\iiii-]$ for even $k$.

It will turn out that only the types defined here are relevant.
\end{notation}

Theorem \ref{t18.2} collects known results on the zeros and poles of single real solutions
$f$  of $P_{III}(0,0,4,-4)$ on $\R_{>0}$. 

\begin{theorem}\label{t18.2}
Fix a real solution $f=f_{mult}(.,s,B)|_{\R_{>0}}$  of $P_{III}(0,0,4,-4)$ on $\R_{>0}$
and its monodromy data $(s,B)\in V^{mat,\R}$.

(a) $f$ has no zeros or poles near $\iiii \iff B=\pm{\bf 1}_2$.

More precisely, $f$ is of type $\orrr{>0}$  $\iff B={\bf 1}_2$,
and $f$ is of type $\orrr{<0}$  $\iff B=-{\bf 1}_2$.

(b) If $B\neq \pm{\bf 1}_2$ then for any $y_0\in\R_{>0}$ 
$f|_{(y_0,\iiii)}$ has infinitely many zeros and/or poles.  This case
will be described more precisely in theorem \ref{t18.4}(a).

(c) $f$ has no zeros or poles near $0$ $\iff |s|\leq 2$.

More precisely, $f$ is of type $\olll{>0}$ $\iff |s|\leq 2$ and 
$b_1+\frac12{s}b_2\geq 1$, and $f$ is of type $\olll{<0}$  $\iff |s|\leq 2$ and 
$b_1+\frac12{s}b_2\leq -1$.

(d) If $s>2$ then $f$ is of type $\olll{[0+][0-]}$
(thus also of type $\olll{[0-][0+]}$).

If $s<-2$ then $f$ is of type $\olll{[\iiii+][\iiii-]}$
(thus also of type $\olll{[\iiii-][\iiii+]}$).

(e) If $|s|\leq 2$ and $B={\bf 1}_2$ then $f$ has no zeros or poles on $\R_{>0}$ and is 
positive on $\R_{>0}$. 

If $|s|\leq 2$ and $B=-{\bf 1}_2$ then $f$ has no zeros or poles on $\R_{>0}$ and is 
negative on $\R_{>0}$. 
\end{theorem}

{\bf Proof:} 
(a)+(b) Of course, part (a) implies part (b). 
Part (a) has at least two different sources, which will be described in the following
parts (I) and (II).

(I) One source which gives the complete result is the equivalence between nilpotent 
orbits of TERP structures and mixed TERP structures from 
\cite[corollary 8.15]{Mo11b}, \cite[theorem 9.3]{HS07} 
(the relevant special case is reformulated 
in theorem \ref{t17.11}) together with the relation between 
real solutions of $P_{III}(0,0,4,-4)$ on $\R_{>0}$  and Euler orbits of TERP(0) bundles
in chapter \ref{s17} (remarks \ref{t17.2} (vi) and \ref{t17.4} (i)). 
Corollary \ref{t17.10} gives the final result
\begin{eqnarray}\label{18.1} 
f_{mult}(.,s,B)\textup{ is of type }\orrr{>0}\iff B={\bf 1}_2.
\end{eqnarray}
The symmetry $R_2$ gives
$$f_{mult}(.,s,B)=-f_{mult}(.,R_2(s,B))=-f_{mult}(.,s,-B),$$
thus 
\begin{eqnarray}\label{18.2} 
f_{mult}(.,s,B)\textup{ is of type }\orrr{<0}\iff B=-{\bf 1}_2.
\end{eqnarray}
This establishes part (a)

(II) Another source is the combination of \cite{MTW77} and 
\cite[ch.\  2,11]{IN86}
(or \cite[ch.\  15]{FIKN06}).  However,  a priori it gives the result 
in part (a) only for
solutions $f_{mult}(.,s,B)|_{\R_{>0}}$ with $(s,B)\in V^{mat,a\cup c,\R}$ because in 
\cite[ch.\  11]{IN86} only such solutions are considered. 

In \cite{MTW77}, only those real solutions of $P_{III}(0,0,4,-4)$ 
are studied which
are of type $\orrr{>0}$ or type $\orrr{<0}$. 
In \cite[ch.\ 2]{IN86} isomonodromic families of $P_{3D6}$ bundles are associated
to all complex solutions  of $P_{III}(0,0,4,-4)$ on $\R_{>0}$. 
In \cite[ch.\ 11]{IN86} it is stated that for $(s,B)\in V^{mat,a\cup c,\R}$ 
a solution $f_{mult}(.,s,B)$ is considered in \cite{MTW77} if and only if 
$B=\pm {\bf 1}_2$, and the zeros and poles for large $x$ of the other solutions 
with $(s,B)\in V^{mat,a\cup c,\R}$ are studied.
If the arguments in \cite[ch.\  11]{IN86} work without the restriction to 
$(s,B)\in V^{mat,a\cup c,\R}$, which they probably do, they give also a proof of part (a).

(c)+(d) For (c) and (d) the best source consists of the asymptotic formulae from 
\cite{Ni09}, which are reformulated in theorem \ref{t12.4} in the formulae
\eqref{12.21}, \eqref{12.22}, \eqref{12.24} and \eqref{12.28}. 
They cover $(s,B)\in V^{mat, a\cup b+\cup c+,\R}$. For $(s,B)\in V^{mat,b-\cup c-}$ 
one needs \eqref{12.24} and \eqref{12.28} and the symmetry $R_1$ with
\[
f_{mult}^{-1}(.,s,B)=f_{mult}(.,R_1(s,B))
=f_{mult}(.,-s,
\bsp 1&0\\0&-1\esp
B
\bsp  1&0\\0&-1\esp).
\]
The formulae \eqref{12.24} and \eqref{12.28} give the type of $f_{mult}(.,s,B)|_{\R_{>0}}$ near $0$
immediately. In the case of the formulae \eqref{12.21} and \eqref{12.22} one has $s\in(-2,2)$
and thus $\alpha_-\in(-\frac{1}{2},\frac{1}{2})$, so $\kappa_{0,1}>0$, and for
$\delta_2\in\{\pm 1\}$ one has
\[
b_1+\tfrac12{s}b_2\in\delta_2\MGcdot\R_{\geq 1}\iff b_-\in \delta_2\MGcdot \R_{>0}.
\]
This establishes parts (c) and (d).

An alternative would be to apply the equivalence
between Sabbah orbits of TERP structures and certain polarized mixed Hodge structures 
in \cite[theorem 7.3]{HS07}, 
see remarks \ref{t17.4}. But this would require quite some additional work.
For $|s|\leq 2$ one would have to show that the TERP(0) bundles which are associated
to $f_{mult}(.,s,B)|_{\R_{>0}}$ induce a certain PMHS. For $|s|>2$ they do not because the
monodromy $\Mon$ does not have eigenvalues in $S^1$, and thus cannot be an automorphism
of a PMHS. But for $|s|>2$ the finer information that $f_{mult(.,s,B)]}$ is of
type $\olll{[0+][0-]}$ or of type $\olll{[\iiii+][\iiii-]}$ does not follow from
\cite[theorem 7.3]{HS07}. For these reasons we do not carry out the alternative method.

(e) The necessity of the conditions $|s|\leq 2$ and
$B=\pm{\bf 1}_2$ for the smoothness of a solution $f_{mult}(.,s,B)|_{\R_{>0}}$ 
follows from (a) and (c). 
For the sufficiency there are several sources. Three will be discussed in the following
points (I) to (III).

(I) In \cite{MTW77} the smooth solutions $f_{mult}(.s,B)$ for $|s|\leq 2$ 
and $B=\pm{\bf 1}_2$ are established, though without identifying $(s,B)$.
A simpler way to establish them is given in \cite{Wi00}.
\cite[ch.\  11]{IN86} makes statements about their monodromy data which amount to 
$|s|\leq 2$ and $B=\pm{\bf 1}_2$.

(II) \cite[theorem 5.9]{HS11} constructs certain pure and polarized TERP(0) bundles.
A special case is reformulated in theorem \ref{t17.11}.
In particular, the TERP(0) bundles for $(s,B)$ with $|s|\leq 2$ and $B={\bf 1}_2$ are
pure and polarized. This implies (corollary \ref{t17.12}) that for such $(s,B)$
the solution $f_{mult}(.,s,B)$ is smooth and positive on $\R_{>0}$.
From the symmetry $R_2$ one obtains that the solution is smooth and negative in the
case $|s|\leq 2$ and $B=-{\bf 1}_2$.

(III) Consider the projection $pr_{mat}:M_{3FN,\R}\to V^{mat,\R}$. 
For $V=V^{mat,a,\R}(B=\pm {\bf 1}_2)$ (and later also for 
$V=V^{mat,b,\R}(B=\pm {\bf 1}_2)$) we want to show that
the set $M_{3FN,\R}^{sing}\cap pr_{mat}^{-1}(V)$ is empty. 
If it were not empty,  parts (a) and (c) and the arguments in the proof of theorem
\ref{t18.3} (b) would show that the restricted projection 
\[
pr_{mat}:M_{3FN,\R}^{sing}\cap pr_{mat}^{-1}(V)\to V
\]
is a covering with finitely many sheets. 

But the existence of the two special smooth solutions $f_{mult}(.,0,\pm {\bf 1}_2)=\pm 1$
from remark \ref{t10.5} excludes this. Therefore the intersection above is empty.

For $V=V^{mat,b,\R}(B=\pm {\bf 1}_2)$, $pr_{mat}^{-1}(V)$ is in the closure of 
$pr_{mat}^{-1}(V^{mat,a,\R}(B=\pm {\bf 1}_2))$, therefore also 
 $M_{3FN,\R}^{sing}\cap pr_{mat}^{-1}(V)$ is empty.
\hfill$\Box$

The isomorphism in theorem \ref{t15.5} (a) (iv) between semi-algebraic manifolds
is now called 
\[
\psi_{mat}: M_{3FN,\R}^{mon}\to \R_{>0}\times V^{mat,\R},
\]
and the projection to $V^{mat,\R}$ is called
\[
pr_{mat}:M_{3FN,\R}^{mon}\to V^{mat,\R}.
\]
Part (a) of the following theorem \ref{t18.3} gives a stratification of $V^{mat,\R}$,
part (b) states that $pr_{mat}$ restricts above the strata to coverings from
$M_{3FN,\R}^{sing}$ and numbers the sheets, 
part (d) connects different strata,
the parts (b)+(c)+(e)+(f) say how the sheets of the coverings glue. 

\begin{theorem}\label{t18.3}
(a) The decomposition of $V^{mat,\R}$ into the three subsets
\[
V^{mat,\R}=V^{mat,a,\R}\cup V^{mat,b,\R}\cup V^{mat,c,\R}
\]
in theorem \ref{t15.5} (a) (ii)
is refined by the conditions $B=\pm{\bf 1}_2$ (respectively, $B\neq \pm{\bf 1}_2$)
into a decomposition into the six subsets
\begin{eqnarray}\label{18.4}
V^{mat,\R}&=& 
V^{mat,a,\R}(B=\pm{\bf 1}_2)\cup V^{mat,a,\R}(B\neq \pm{\bf 1}_2)\\  
&\cup& V^{mat,b,\R}(B=\pm{\bf 1}_2)\cup V^{mat,b,\R}(B\neq \pm{\bf 1}_2)\nonumber \\
&\cup& V^{mat,c,\R}(B=\pm{\bf 1}_2)\cup V^{mat,c,\R}(B\neq \pm{\bf 1}_2). \nonumber
\end{eqnarray}
They have 2, 4, 4, 4, 4, 8 components respectively
(cf.\ the two pictures in theorem \ref{t15.5} (a) (iii)).

(b) If $V$ is the set $V^{mat,a,\R}(B=\pm{\bf 1}_2)$ or the set
$V^{mat,b,\R}(B=\pm{\bf 1}_2)$, then $M_{3FN,\R}^{sing}\cap pr_{mat}^{-1}(V)$ is empty.
If $V$ is one of the other four sets in (a), then the restricted projection
\[
pr_{mat}:M_{3FN,\R}^{sing}\cap pr_{mat}^{-1}(V)\to V
\]
is a covering. 

The sheets of the coverings and the zeros and poles of any solution
$f_{mult}(.,s,B)|_{\R_{>0}}$ are indexed by $\Z$ or $\Z_{\geq 1}$ or $\Z_{\leq 0}$ such
that $x_k<x_l\iff k<l$ for zeros/poles $x_k$ and $x_l$. The second column in  the 
following table lists the zeros and poles of a solution 
$f_{mult}(.,s,B)|_{\R_>0}$ for $(s,B)$ in the set in the first column:
\begin{eqnarray}\label{18.3}
\begin{array}{c|c}
V^{mat,a,\R}(B\neq \pm{\bf 1}_2) & (x_k(s,B))_{k\in\Z_{\geq 1}}  \\
V^{mat,c,\R}(B\neq \pm{\bf 1}_2) & (x_k(s,B))_{k\in\Z_{\geq 1}}  \\
V^{mat,b,\R}(B\neq \pm{\bf 1}_2) & (x_k(s,B))_{k\in\Z}  \\
V^{mat,b,\R}(B=\pm{\bf 1}_2)     & (x_k(s,B))_{k\in\Z_{\leq 0}} 
\end{array}
\end{eqnarray}
The indexing is unique for the three sets with index sets $\Z_{\geq 1}$ or $\Z_{\leq 0}$.
It is unique up to a shift for each of the four components of the set 
$V^{mat,b,\R}(B\neq \pm{\bf 1}_2)$ with index set $\Z$. There it is fixed by the 
requirement that the sheets over $V^{mat,b,\R}(B= {\bf 1}_2)$ with indices 
$k\in\Z_{\leq 0}$ glue with the sheets with the same indices $k$ over each of the four 
components of the set $V^{mat,b,\R}(B\neq \pm{\bf 1}_2)$. 

(c) For $\delta_1,\delta_3\in\{\pm 1\}$ the sheet with index $k\in\Z_{\leq 0}$ over the
component $V^{mat,b\delta_1,\R}(B=-{\bf 1}_2)$ glues to the sheet with index 
$k-\delta_1\delta_3$ over the component $V^{mat,b\delta_1,\R}(ib_2\in\delta_3\R_{>0})$.

(d)There are continuous maps
\begin{eqnarray*}
&&\gamma_1: V^{mat,a,\R} \to \R_{>0},\\
&&\gamma_2: V^{mat,c,\R} \to \R_{>0},\\
&&\gamma_3: V^{mat,\R}(B=\pm{\bf 1}_2) \to \R_{>0},
\end{eqnarray*}
such that
\begin{equation}\label{18.5}
\{(x,s,B)\, |\, x<\gamma_1(s,B), (s,B)\in V^{mat,a,\R}\}\cap 
\psi_{mat}(M_{3FN,\R}^{sing})=\emptyset,
\end{equation}
\begin{equation}\label{18.6}
\{(x,s,B)\, |\, x<\gamma_2(s,B), (s,B)\in V^{mat,c,\R}\}\cap 
\psi_{mat}(M_{3FN,\R}^{sing})=\emptyset,
\end{equation}
\begin{equation}\label{18.7}
\{(x,s,B)\, |\, x>\gamma_3(s,B), B\in\{\pm 1\}\}\cap 
\psi_{mat}(M_{3FN,\R}^{sing})=\emptyset.
\end{equation}

(e) For $\delta_1,\delta_2,\delta_3\in\{\pm 1\}$ the sheet with index $k\in\Z_{\geq 1}$ 
over the component $V^{mat,c\delta_1,\R}(b_1+\frac{s}{2}b_2=\delta_2,ib_2\in\delta_3\R_{>0})$ 
glues to the sheet with index $k$ over the component 
$V^{mat,a,\R}(b_1+\frac12{s}b_2\in\delta_2\R_{\geq 1}, ib_2\in\delta_3\R_{>0})$.

(f) It glues to the sheet with index $k-\frac{1-\delta_2}{2}\delta_1\delta_3$
over the component 
$V^{mat,b\delta_1,\R}(ib_2\in\delta_3\R_{>0})$.

(g) Glueing all the sheets over the 26 components of the six sets in (a) gives precisely
the four smooth hypersurfaces $M_{3FN,\R}^{sing}=\cup_{k=0}^3 M_{3FN,\R}^{[k]}$
in $M_{3FN,\R}^{ini}$ of types $[0-], [\iiii-], [0+], [\iiii+]$ for $k=0,1,2,3$.
\end{theorem}

{\bf Proof:} 
(a) This is essentially a definition. The 26 components can be seen in the 
picture in theorem \ref{t15.5} (a) (iii).

(b) If $V=V^{mat,a,\R}(B=\pm{\bf 1}_2)$ or $V=V^{mat,c,\R}(B=\pm{\bf 1}_2)$,
then $M_{3FN,\R}^{sing}\cap pr_{mat}^{-1}(V)=\emptyset$ because of theorem \ref{t18.2} (e).

Let $V$ be one of the other four sets in (a).
The hypersurface $M_{3FN,\R}^{sing}\subset M_{3FN,\R}^{mon}$ is real analytic and smooth
with four components. Because the functions $f_{mult}(.,s,B)$ have only simple zeros
and poles, the hypersurface $M_{3FN,\R}^{sing}$ is everywhere transversal to the fibres
of the projection $pr_{mat}$.

Therefore, in order to show that $pr_{mat}:M_{3FN,\R}^{sing}\cap pr_{mat}^{-1}(V)\to V$
is a covering, it is sufficient to show that the following data do not exist:
\begin{align}\nonumber
&\textup{a }C^\iiii\textup{ path }\gamma:[0,1]\to V\textup{ and }
\textup{a }C^\iiii\textup{ lift } \textup{ of }\gamma:\vert_{[0,1)}
\\
\label{18.8}
&\www \xi=(\xi,\gamma):[0,1)\to 
\psi_{mat}(M_{3FN,\R}^{sing})\cap \R_{>0}\times V
\\  \nonumber
&\textup{such that for }r\to 1,\ 
\xi(r)\textup{ tends to $0$ or }\iiii.
\end{align}
%\begin{eqnarray}
%\textup{a }C^\iiii\textup{ path }\gamma:[0,1]\to V\textup{ and}\nonumber \\
%\label{18.8}
%\textup{a }C^\iiii\textup{ lift }\www \xi=(\xi,\gamma):[0,1[\to 
%\psi_{mat}(M_{3FN,\R}^{sing})\cap \R_{>0}\times V\\
%\textup{ of }\gamma:[0,1[\textup{ such that for }r\to 1\ 
%\xi(r)\textup{ tends to $0$ or }\iiii.\nonumber
%\end{eqnarray}

{\bf 1st case:} Suppose that such data exist with $\lim_{r\to 1}\xi(r)=0.$

{\bf 1st subcase:} $V=V^{mat,b,\R}(B=\pm{\bf 1}_2)$ 
or $V=V^{mat,b,\R}(B\neq\pm{\bf 1}_2)$.

By theorem \ref{t18.2} (d) for any $r\in[0,1]$ the set 
$\psi_{mat}(M_{3FN,\R}^{sing})\cap \R_{>0}\times \{\gamma(r)\}$ is discrete
and contains points with $x$-values arbitrarily close to $0$.

Choose any point $(\eta(1),\gamma(1))\in \psi_{mat}(M_{3FN,\R}^{sing})
\cap \R_{>0}\times \{\gamma(1)\}$. It extends to a $C^\iiii$ lift
$\www \eta=(\eta,\gamma):[1-\varepsilon,1]\to
\psi_{mat}(M_{3FN,\R}^{sing})\cap \R_{>0}\times \gamma([1-\varepsilon,1])$
of $\gamma:[1-\varepsilon,1]\to V$ for some $\varepsilon>0$.

For $r\in[1-\varepsilon,1)$, there can be only finitely many points in
$\psi_{mat}(M_{3FN,\R}^{sing})\cap \R_{>0}\times \{\gamma(r)\}$
between $\www \eta(r)$ and $\www\xi(r)$. 
But only their continuations to $r=1$ along the path $\gamma$ can be in
$\psi_{mat}(M_{3FN,\R}^{sing})\cap \R_{>0}\times \{\gamma(1)\}$
and lie with respect to their $x$-values under $\www\eta(1)=(\eta(1),\gamma(1))$.
This contradicts theorem \ref{t18.2} (d).

{\bf 2nd subcase:} $V=V^{mat,a,\R}(B\neq\pm{\bf 1}_2)$ 
or $V=V^{mat,c,\R}(B\neq\pm{\bf 1}_2)$.

The asymptotic formulae \eqref{12.21}, \eqref{12.22} and \eqref{12.28} in theorem
\ref{t12.4} for $f_{mult}(.,s,B)|_{\R_{>0}}$ for $x\to 0$ depend real analytically
on $(s,B)\in V$. Therefore the data \eqref{18.8} do not exist.
The 2nd subcase also leads to a contradiction.  Hence the 1st case is impossible.

{\bf 2nd case:} Suppose that the data \eqref{18.8} exist with 
$\lim_{r\to 1}\xi(r)=\iiii$.

{\bf 1st subcase:} $V=V^{mat,a,\R}(B\neq\pm{\bf 1}_2)$ 
or $V=V^{mat,c,\R}(B\neq\pm{\bf 1}_2)$. or $V=V^{mat,b,\R}(B\neq\pm{\bf 1}_2)$.

By theorem \ref{t18.2} (b) for any $r\in [0,1]$ the set 
$\psi_{mat}(M_{3FN,\R}^{sing})\cap \R_{>0}\times \{\gamma(r)\}$
is discrete and contains points with arbitrarily large $x$-values.
The same arguments as in the 1st subcase of the 1st case lead to a contradiction.

{\bf 2nd subcase:} $V=V^{mat,b,\R}(B=\pm{\bf 1}_2)$.

Compare the proof (i) of theorem \ref{t18.2} (a).
Choose $(s,B)=(s,{\bf 1}_2)\in V$ and consider the associated Euler orbit
of TERP(0) bundles. By \cite[theorem 9.3 (2)]{HS07} 
(see also the proof of corollary \ref{t17.10}) it is a nilpotent orbit,
so  there is a lower bound $x_{bound}(s)$ such that the TERP(0) bundles for 
$x>x_{bound}(s)$ are pure and polarized and thus $f_{mult}(x,s,{\bf 1}_2)>0$.

We claim that the bound $x_{bound}(s)$ can be chosen so that it depends 
continuously on $s$. This is not explicitly stated in \cite[theorem 9.3 (2)]{HS07},
but it follows from its proof. In the semisimple case, which is the only relevant case here,
that proof simplifies considerably. In our rank 2 case one just needs that for all $z\in\C^*$
with $\arg(z)\in(-\varepsilon,\varepsilon)$ for some $\varepsilon>0$ the matrix
\[
\begin{pmatrix} 
e^{-\frac{x}{z}-xz} & 
\\
 & e^{\frac{x}{z}+xz}
 \end{pmatrix} 
\MGcdot
\begin{pmatrix}1&s\\ 0 & 1\end{pmatrix}
\MGcdot
\begin{pmatrix} 
e^{-\frac{x}{z}-xz} & 
\\
 & e^{\frac{x}{z}+xz}
 \end{pmatrix} 
=
\begin{pmatrix}1&s e^{-\frac{2x}{z}-2xz}\\ 0 & 1\end{pmatrix}
\]
is sufficiently close to the unit matrix ${\bf 1}_2$. 
For this one can find a lower bound $x_{bound}(s)$ for $x$ which depends continuously on $s$.

As $f_{mult}(x,s,{\bf 1}_2)>0$ for $x>x_{bound}(s)$, data as in \eqref{18.8} with
$\lim_{r\to 1}\xi(r)_\iiii$ cannot exist. 

The case $B=-{\bf 1}_2$ can be reduced to the case $B={\bf 1}_2$ with the symmetry $R_2$.
The 2nd subcase also leads to a contradiction. So the 2nd case is impossible.

Therefore the data \eqref{18.8} do not exist. Therefore
$pr_{mat}:M_{3FN,\R}^{sing}\cap pr_{mat}^{-1}(V)\to V$ is a covering.

(c) Formula \eqref{12.25} in theorem \ref{12.4} gives approximations for 
the $x$-values of the sheets above $V^{mat,b\delta_1,\R}$ with indices $k\ll 0$.
But beware that the $k$ in \eqref{12.25} differs by the sign and a shift from 
the index $k$ of a sheet.

If $s$ is fixed and $b_-$ goes once around $0$ in the positive direction,
then $b_-$ must be replaced by $b_-\MGcdot e^{2\pi i}$, 
$\delta^{NI}$ must be replaced by $\delta^{NI}+2\pi$, 
and the sheet with index $k$ turns into the sheet with index $k-2$.
This establishes (c) for $\delta_1=1$. For $\delta_1=-1$ one applies the symmetry
$R_1$ with $R_1(f_{mult})=f_{mult}^{-1}$, $R_1(s)=-s$, $R_s(b_-)=b_-^{-1}$;
see table \eqref{12.13}.

(d) The existence of $\gamma_1$ with \eqref{18.5} and of $\gamma_2$ with 
\eqref{18.6} follows from the arguments in the 2nd subcase of the 1st case
in the proof of (b).

The existence of $\gamma_3$ with \eqref{18.7} follows from the arguments 
in the 2nd subcase of the 2nd case in the proof of (b).

(e) Choose any $x\in\R_{>0}$, consider the point 
\begin{eqnarray*}
(x,2\delta_1,\delta_2 {\bf 1}_2)&\in & 
\R_{>0}\times V^{mat,c\delta_1,\R}(B=\delta_2{\bf 1}_2)\\
&\subset& \psi_{mat}(M_{3FN,\R}^{reg})\quad (\Leftarrow\textup{ theorem \ref{t18.2} (e)}),
\end{eqnarray*}
and choose a small open neighbourhood $U\subset \psi_{mat}(M_{3FN,\R}^{reg})$
of this point. 
Because of $\R_{>0}\times V^{mat,a\cup c\delta_1,\R}(B=\delta_2{\bf 1}_2) \subset 
\psi_{mat}(M_{3FN,\R}^{sing})$  ($\Leftarrow$ theorem \ref{t18.2} (e))
and because of \eqref{18.5} and \eqref{18.6}  all values $x_k(.,s,B)$ for
$(s,B)\in V^{mat,a\cup c\delta_1,B}(b_1+\frac{s}{2}b_2\in \delta_2\R_{\geq 1},
ib_2\in\delta_3\R_{>0})$ tend along the lines 
$s=\textup{const}$ for $(b_1,b_2)\to (\delta_2,0)$ to $\iiii$.

Therefore the parts over $pr_{mat}(U)$ of all sheets over 
$V^{mat,c\delta_1,\R}(\www b_1=\delta_2,ib_2\in\delta_3\R_{>0})$ as well
as the parts over $pr_{mat}(U)$ of all sheets over 
$V^{mat,a,\R}(b_1\frac{s}{2}b_2\in\delta_2\R_{\geq 1},ib_2\in\delta_3\R_{>0})$
lie with respect to their $x$-values above $U$.

Therefore the parts over $\pr_{mat}(U)$ of the sheets with the same indices glue,
therefore the sheets with the same indices glue globally.

(f) By (b)+(c), the sheet over 
$V^{mat,b\delta_1,\R}(ib_2\in\delta_3\R_{>0})$ with index 
$k-\frac{1-\delta_2}{2}\delta_1\delta_3$ for $k\in\Z_{\leq 0}$ glues to the sheet over
$V^{mat,b\delta_1,\R}(B=\delta_2{\bf 1}_2)$ with index $k$.

As there are no sheets with indices $k\in \Z_{\geq 1}$ over 
$V^{mat,b\delta_1,\R}(B=\delta_2{\bf 1}_2)$, the values $x_k(s,B)$ of the sheets over
$V^{mat,b\delta_1,\R}(ib_2\in\delta_3\R_{>0})$ with indices 
$k\geq 1-\frac{1-\delta_2}{2}\delta_1\delta_3$ tend along all lines 
$s=\textup{const}$ for $(b_1,b_2)\to (\delta_2,0)$ to $\iiii$.

Therefore the parts over $pr_{mat}(U)$ of all the sheets over
$V^{mat,b\delta_1,\R}(ib_2\in\delta_3\R_{>0})$ with indices 
$k\geq 1-\frac{1-\delta_2}{2}\delta_1\delta_3$ lie with respect to their $x$-values
above $U$.

On the other hand, because of \eqref{18.7} and theorem \ref{t18.2} (e),
all values $x_k(s,\delta_2{\bf 1}_2)$ (automatically $k\le 0$)
for $(s,\delta_2{\bf 1}_2)\in V^{mat,b\delta_1,\R}(B=\delta_2{\bf 1}_2)$
tend to $0$ for $s\to \delta_2\MGcdot 2$. 
Therefore the parts over $pr_{mat}(U)$ of all the sheets over 
$V^{mat,b\delta_1,\R}(B=\delta_2{\bf 1}_2)$ lie with respect to their
$x$-values below $U$.

Because the sheets over $V^{mat,b\delta_1,\R}(B=\delta_2 {\bf 1}_2)$ with
indices $k$ (automatically $k\le0$) glue to the sheets over 
$V^{mat,b\delta_1,\R}(ib_2\in\delta_3\R_{>0})$ with indices 
$k-\frac{1-\delta_2}{2}\delta_1\delta_3$ by (b)+(c),
also the parts over $pr_{mat}(U)$ of the sheets over 
$V^{mat,b\delta_1,\R}(ib_2\in\delta_3\R_{>0})$ with indices 
$k-\frac{1-\delta_2}{2}\delta_1\delta_3$ 
for $k\leq 0$ lie with respect to their $x$-values below $U$.

The sheets over $V^{mat,b\delta_1,\R}(ib_2\in\delta_3\R_{>0})$
and over $V^{mat,c\delta_1,\R}(\www b_1=\delta_2,ib_2\in\delta_3\R_{>0})$
whose parts over $pr_{mat}(U)$ lie with respect to their $x$-values above $U$,
must glue pairwise. This establishes (f).

(g) This is clear as the sheets are the components of 
$M_{3FN,\R}^{sing}\cap pr_{mat}^{-1}(V)$.
It follows from \eqref{10.13}
that $M_{3FN,\R}^{[k]}$ for $k=0,1,2,3$
are, respectively,  of type 
$[0-], [\iiii -], [0+], [\iiii +].  
$\hfill$\Box$

\begin{theorem}\label{t18.4}
Fix a real solution $f=f_{mult}(.,s,B)|_{\R_{>0}}$  of $P_{III}(0,0,4,-4)$ on $\R_{>0}$
and its monodromy data $(s,B)\in V^{mat,\R}$.

(a) (A supplement to theorem \ref{t18.2} (b)) 
Suppose $B\neq \pm{\bf 1}_2$. 

If $ib_2<0$ then $f$ is of type
$\orrr{[0-][\iiii+]}$ (thus also of type $\orrr{[\iiii+][0-]}$).

If $ib_2>0$ then $f$ is of type
$\orrr{[0+][\iiii-]}$ (thus also of type $\orrr{[\iiii-][0+]}$).

(b) (Global results)
A $y_0\in \R_{>0}$ exists such that $f(y_0)\neq 0$ and such that 
$f|_{(0,y_0]}$ and $f|_{[y_0,\iiii)}$ are each of the type determined by the behaviour
near $0$ and $\iiii$. This means there is no intermediate mixed zone.

The following two tables list the 18 possible combinations of 
\[
\text{\em type of }f|_{(0,y_0]} \ \&\  \text{\em type of }f|_{[y_0,\iiii)}.
\]
In each of the four cases with $|s|> 2$ and $B\neq \pm{\bf 1}_2$ there are two ways 
to position $y_0$ and to split the sequence of zeros and poles into a type near $0$
and a type near $\iiii$. 
In the other 10 cases there is only one way to position $y_0$.

It follows that there are 14 possible sequences of zeros and poles. Each is realized above one
of the 14 components of 
\begin{eqnarray*}
V^{mat,a\cup c,\R}(B=\pm{\bf 1}_2) \quad \text{(2 components)},\\
V^{mat,a\cup c,\R}(B\neq \pm{\bf 1}_2) \quad \text{(4 components)}, \\
V^{mat,b,\R}(B=\pm{\bf 1}_2) \quad \text{(4 components)},\\
V^{mat,b,\R}(B\neq\pm{\bf 1}_2) \quad \text{(4 components)}.
\end{eqnarray*}

\begin{eqnarray}\label{18.9}
\renewcommand{\arraystretch}{1.7}
\begin{array}{ccc}
& |s|\leq 2 & \\ \hline
B={\bf 1}_2 & \vline & \olll{>0} \ \&\  \orrr{>0} \\
B=-{\bf 1}_2 & \vline & \olll{<0} \ \&\  \orrr{<0} \\
b_1+\frac12{s}b_2\geq 1,ib_2<0 & \vline & 
   \olll{>0} \ \&\  \orrr{[\iiii+][0-]} \\
b_1+\frac12{s}b_2\geq 1,ib_2>0 & \vline &  
   \olll{>0} \ \&\  \orrr{[0+][\iiii-]} \\
b_1+\frac12{s}b_2\leq -1,ib_2<0 & \vline & 
   \olll{<0} \ \&\  \orrr{[0-][\iiii+]} \\
b_1+\frac12{s}b_2\leq -1,ib_2>0& \vline & 
   \olll{<0} \ \&\  \orrr{[\iiii-][0+]}
\end{array}
\end{eqnarray}

\begin{eqnarray}\label{18.10}
\renewcommand{\arraystretch}{1.7}
\begin{array}{ccccc}
  & \vline & s> 2 & \vline  & s<-2 \\ \hline
B={\bf 1}_2 & \vline & \olll{[0-][0+]} \ \&\  \orrr{>0} 
& \vline & \olll{[\iiii-][\iiii+]} \ \&\  \orrr{>0} \\
B=-{\bf 1}_2 & \vline & \olll{[0+][0-]} \ \&\   \orrr{<0} 
& \vline & \olll{[\iiii+][\iiii-]} \ \&\  \orrr{<0} \\
ib_2<0 & \vline & \olll{[0-][0+]} \ \&\  \orrr{[0-][\iiii+]} 
  & \vline & \olll{[\iiii-][\iiii+]} \ \&\  \orrr{[0-][\iiii+]} \\
 & \vline &\textup{or } \olll{[0+][0-]} \ \&\  \orrr{[\iiii+][0-]} 
  & \vline &\textup{or } \olll{[\iiii+][\iiii-]} \ \&\  \orrr{[\iiii+][0-]} \\
ib_2>0 & \vline & \olll{[0-][0+]} \ \&\  \orrr{[\iiii-][0+]} 
  & \vline & \olll{[\iiii-][\iiii+]} \ \&\  \orrr{[\iiii-][0+]} \\
& \vline &\textup{or } \olll{[0+][0-]} \ \&\  \orrr{[0+][\iiii-]} 
  & \vline &\textup{or } \olll{[\iiii+][\iiii-]} \ \&\  \orrr{[0+][\iiii-]} 
\end{array}
\end{eqnarray}

\end{theorem}

{\bf Proof:} 
For each of the sheets over the 26 components of the
six sets in theorem \ref{t18.3} (a), it is crucial to determine to
which of the four hypersurfaces 
$M_{3FN,\R}^{[k]}$ it belongs.  
First we read off from the glueing information in theorem \ref{t18.3} (b)+(c)+(e)+(f)
which sheets glue together to one of these hypersurfaces, then we will find its type.

Theorem \ref{t18.3} (e) allows one to join for $\delta_2,\delta_3\in\{\pm 1\}$ the sheets over
$V^{mat,a,\R}(b_1+\frac{s}{2}b_2\in\delta_2\R_{\geq 1},ib_2\in\delta_3\R_{>0})$ 
and 
$V^{mat,c,\R}(\www b_1=\delta_2,ib_2\in\delta_3\R_{>0})$ 
to sheets over 
$V^{mat,a\cup c,\R}(b_1+\frac{s}{2}b_2\in\delta_2\R_{\geq 1},ib_2\in\delta_3\R_{>0})$. 
The sheets over 
$V^{mat,b,\R}(B=\pm{\bf 1}_2)$ will not be considered separately, as the sheets on both sides
of them will be glued via them. The sheets are denoted as follows:

\begin{list}{}{}
\item[$u_k^{\delta_2}$] sheet with index $k$ over
$V^{mat,b+,\R}(ib_2\in\delta_3\R_{>0})$, here $k\in\Z$,
\item[$d_k^{\delta_2}$] sheet with index $k$ over
$V^{mat,b-,\R}(ib_2\in\delta_3\R_{>0})$, here $k\in\Z$,
\item[$l_k^{\delta_2}$] sheet with index $k$ over
$V^{mat,a\cup c,\R}(b_1+\frac{s}{2}b_2\in\R_{\geq 1}, 
ib_2\in\delta_3\R_{>0})$, here $k\in\Z_{\leq 0}$,
\item[$r_k^{\delta_2}$] sheet with index $k$ over
$V^{mat,a\cup c,\R}(b_1+\frac{s}{2}b_2\in\R_{\leq -1}, 
ib_2\in\delta_3\R_{>0})$, here $k\in\Z_{\leq 0}.$
\end{list}

The glueing information in theorem \ref{t18.3} (b)+(c)+(f) is rewritten for $-3\leq k\leq 4$
in the following four chains of sheets. They extend in the obvious way to all $k$:
\begin{eqnarray}\label{18.11}
\begin{split}
.. u^-_{-3}\MGcdot u^+_{-3}\MGcdot u^-_{-1}\MGcdot u^+_{-1}\MGcdot
u^-_1\MGcdot r^-_1\MGcdot d^-_1\MGcdot l^-_2\MGcdot u^-_3\MGcdot r^-_3\MGcdot d^-_3\MGcdot l^-_4..\\  
.. u^-_{-2}\MGcdot u^+_{-2}\MGcdot u^-_0\MGcdot u^+_0\MGcdot
l^+_1\MGcdot d^+_2\MGcdot r^+_2\MGcdot u^+_2\MGcdot l^+_3\MGcdot d^+_4\MGcdot r^+_4\MGcdot u^+_4..\\  
.. d^+_{-3}\MGcdot d^-_{-3}\MGcdot d^+_{-1}\MGcdot d^-_{-1}\MGcdot
d^+_1\MGcdot r^+_1\MGcdot u^+_1\MGcdot l^+_2\MGcdot d^+_3\MGcdot r^+_3\MGcdot u^+_3\MGcdot l^+_4..\\  
.. d^+_{-2}\MGcdot d^-_{-2}\MGcdot d^+_0\MGcdot d^-_0\MGcdot
l^-_1\MGcdot u^-_2\MGcdot r^-_2\MGcdot d^-_2\MGcdot l^-_3\MGcdot u^-_4\MGcdot r^-_4\MGcdot d^-_4..
\end{split}
\end{eqnarray}
Each chain is one of the four hypersurfaces.

By theorem \ref{t18.2} (d), the first two chains are the hypersurfaces of types $[0+]$ and $[0-]$,
and the last two chains are the hypersurfaces of types $[\iiii+]$ and $[\iiii-]$.
But it remains to determine which is which.

The sheets $u^+_0$ and $u^-_0$ are glued via the highest sheet (it has index 0) over
$V^{mat,b+,\R}(B={\bf 1}_2)$. Above this sheet $f_{mult}(.,s,B)$ is of type $\orrr{>0}$.
Thus $f_{mult}(.,s,B)$ is positive right above $u^+_0$ and $u^-_0$.
Therefore the hypersurface which contains $u^+_0$ and $u^-_0$ is of type $[0+]$,
and the hypersurface which contains $u^-_1$ and $u^+_{-1}$ is of type $[0-]$.
Analogously, the hypersurface which contains $d^+_0$ and $d^-_0$ is of type $[\iiii+]$,
and the hypersurface which contains $d^-_{-1}$ and $d^+_1$ is of type $[\iiii-]$.

The following table lists for each of the 14 components in theorem \ref{t18.4} (b),
above (or under) which sheet the value $y_0$ has to be positioned so that 
$f|_{(0,y_0]}$ and $f|_{[y_0,\iiii)}$ are of one of the types in the table.
For the first 10 components there is only one possibility, for the last four components
there are two possibilities.
\[
\begin{array}{lcl}
V^{mat,a\cup c,\R}(B=\pm{\bf 1}_2) & \textup{(2 components):} & 
y_0\textup{ arbitrary in }(0,\iiii)\\
V^{mat,a\cup c,\R}(B\neq\pm{\bf 1}_2)
 & \textup{(4 components):} & 
y_0\textup{ under the lowest sheet }\\
V^{mat,b\delta_1,\R}(B=\delta_2\MGcdot{\bf 1}_2) & \textup{(4 components):} & 
y_0\textup{ above the highest sheet }\\
V^{mat,b\delta_1,\R}(ib_2\in\delta_3\R_{>0}) & \textup{(4 components):} & 
y_0\textup{ above the sheet }
\end{array} 
\]
\[
\begin{array}{r|r|r|r} 
u^+_0 & u^-_0 & d^+_0 & d^-_0 
\\ 
\textup{or }u^+_{-1} & \textup{or }u^-_1 & \textup{or }d^+_1 & \textup{or }d^-_{-1} 
\\
 \textup{for }(\delta_1,\delta_3)= (1,1)&(1,-1)&(-1,1)&(-1,-1) 
 \end{array}
 \]
This proves theorem \ref{t18.4} (b).
It also proves theorem \ref{t18.4} (a).
\hfill$\Box$ 

\begin{remarks}\label{t18.5}
(i) The following three pictures give an idea of the glueing of the sheets and of the 
sequences of zeros and/or poles of the real solutions  of
$P_{III}(0,0,4,-4)$ on $\R_{>0}$. 

The first picture is a sketch of the restriction to some fixed
value $s_0>2$ of $M_{3FN,\R}^{mon}$, that is
\[
M_{3FN,\R}^{mon}|_{s=s_0}\cong \R_{>0}\times V^{mat,\R}|_{s=s_0}
\cong \R_{>0}\times S^1\quad\text{(cf.\  \eqref{15.8})}.
\]
Here $V^{mat,\R}|_{s=s_0}\cong S^1$ is cut open at
$b_-=i$ (i.e.\ $b_1+\frac12{s_0}b_2=0$, $b_2=i/\sqrt{\frac14{s^2}-1}$ --- see \eqref{5.24}).
In the white region $f$ is positive, in the shaded region $f$ is negative.
The lines are the intersections of the glued sheets, i.e.\  of the four hypersurfaces,
with $M_{3FN,\R}^{mon}|_{s=s_0}$. 

%{\sc Later pic18-1.jpg}
\includegraphics[width=0.8\textwidth]{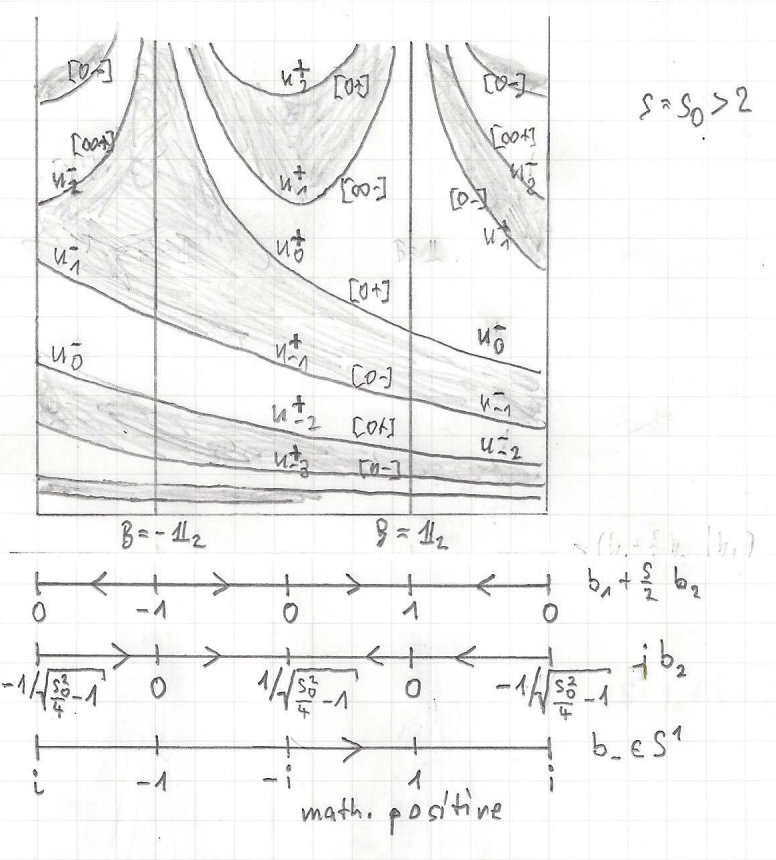}

The second picture is a sketch of the restriction to some value $s_1\in[-2,2]$ 
of $M_{3FN,\R}^{mon}$, that is 
\[
M_{3FN,\R}^{mon}|_{s=s_1}\cong \R_{>0}\times V^{mat,\R}|_{s=s_1}
\cong \R_{>0}\times \R^*\quad\text{(cf.\  \eqref{15.7})}.
\]
It has two components. The symmetry $R_2$ is a horizontal shift --- 
it maps the two components to one another and exchanges white and shaded regions
and the types $[0+]\leftrightarrow[0-]$ and $[\iiii+]\leftrightarrow[\iiii-]$. 

%{\sc Later pic18-2.jpg}
\includegraphics[width=0.8\textwidth]{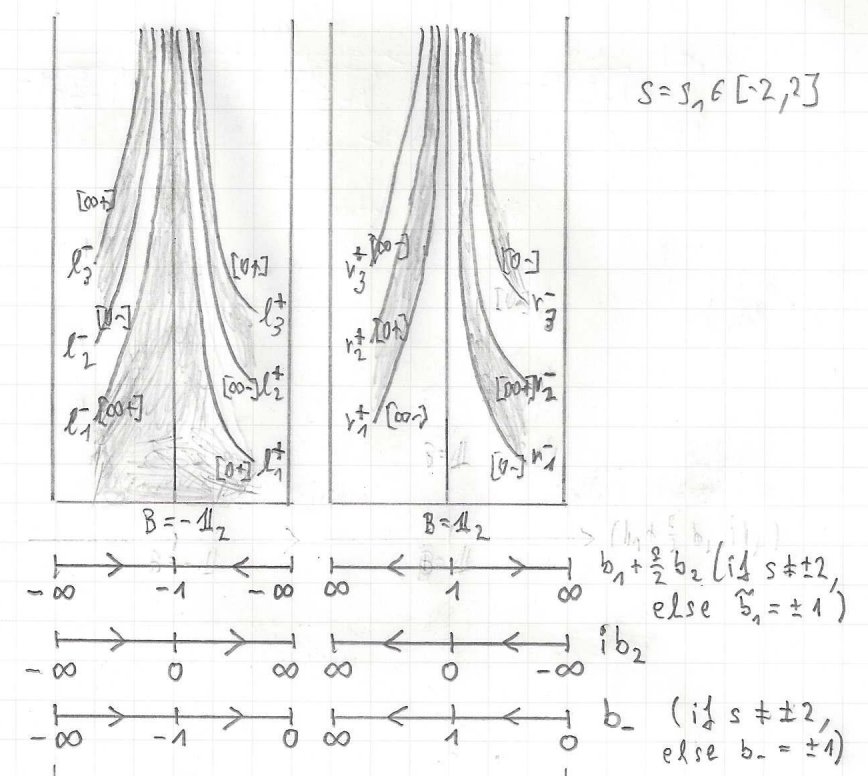} 

The third picture is a sketch of the restriction to the value 
$s=s_2=-s_0<-2$ of $M_{3FN}^{mon}$. It is obtained from the first
picture using the symmetry $R_1$, which acts as reflection along the line 
$B={\bf 1}_2$. It preserves white and shaded regions, but exchanges the types
$[0+]\leftrightarrow[\iiii+]$ and $[0-]\leftrightarrow[\iiii-]$.

The symmetry $R_2$ is visible in the first and the third pictures
as a horizontal shift of half the length
of the $S^1$, which exchanges white and shaded regions and the types
$[0+]\leftrightarrow[0-]$ and $[\iiii+]\leftrightarrow[\iiii-]$.

(ii) The pictures are only sketches. All the lines are drawn as convex graphs, which is what we expect. Conjecture \ref{t18.6} (a) makes this expectation
more precise and extends it.

(iii) The two regions in $M_{3FN,\R}^{mon}$ where $f_{mult}$ is positive (white region)
or negative (shaded region) are the two components of 
$M_{3FN,\R}^{reg}\cong \R_{>0}\times \R^*\times\R$, so they are connected and 
contractible. 

%{\sc Later pic18-3.jpg}
\includegraphics[width=0.8\textwidth]{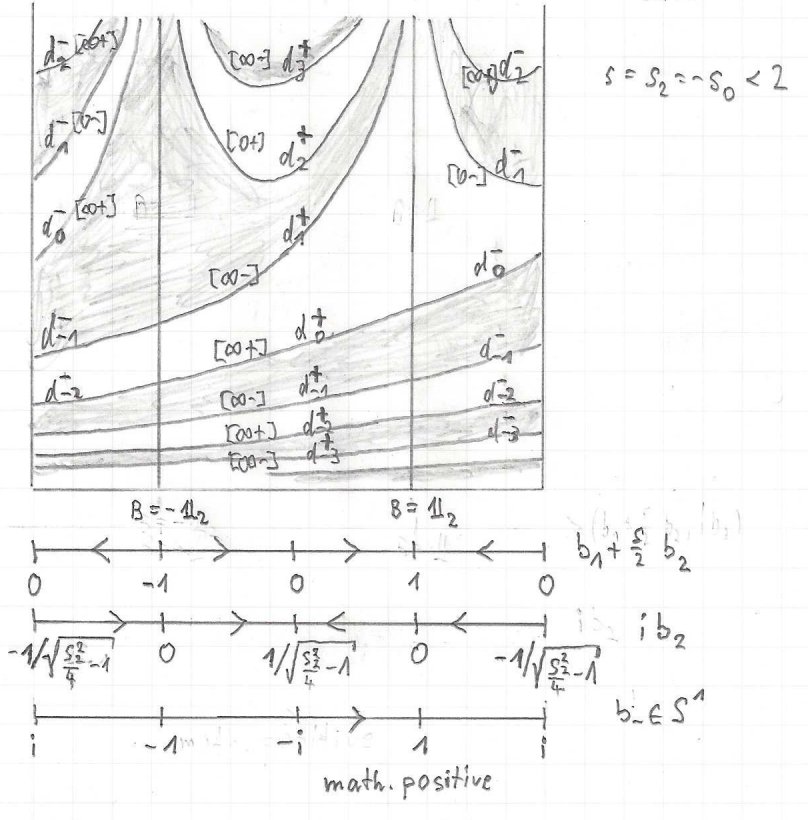} 

In the picture for $s=s_0>2$ the lower white and shaded strips, that is,
the strips between lines of types $[0+]$ and $[0-]$, are glued by glueing the two
vertical bounding lines where $b_-=i$. 
The movement to lower and lower strips in $M_{3FN,\R}^{mon}|_{s=s_0}$ 
within one region projects in $V^{mat,\R}$ to a positive 
winding around the hole with $s=\iiii$. The same holds for the
lower strips in the picture for $s=s_2=-s_0<-2$ and the hole with $s=-\iiii$.
Observe that, for $s=s_2$, $b_-\in S^1$ winds in the negative direction around 
the hole $s=-\iiii$.

But for the upper strips one needs all three pictures.
It is necessary to alternate two types of moves, moving $x_0$ appropriately so that
one stays in the white or the shaded region:
\begin{list}{}{}
\item[$(\alpha)$] 
for fixed $ib_2\in (-1/\sqrt{\tfrac14{s_0^2}-1},1/\sqrt{\tfrac14{s_0^2}-1})$ move between
$s=s_0$ and $s=s_2=-s_0$,
\item[$(\beta)$]
for fixed $s=\pm s_0$ move between points near $B={\bf 1}_2$ and points near $B=-{\bf 1}_2$.
\end{list}
Together this gives two movements in $M_{3FN,\R}^{mon}$ whose projections to $V^{mat,\R}$ 
are
windings
in the positive direction  around the hole with $ib_2=\iiii$ or the hole with $ib_2=-\iiii$.

Parts (b) and (c) of the following conjecture would allow us to make all these movements
within $\{x_0\}\times V^{mat,\R}$ for any fixed $x_0\in \R_{>0}$. 
\end{remarks}

\begin{conjecture}\label{t18.6}
(a) The components of 
$\psi_{mat}(M_{3FN,\R}^{sing})\cap \R_{>0}\times V^{mat,a\cup c,\R}$ and of 
$\psi_{mat}(M_{3FN,\R}^{sing})\cap \R_{>0}\times V^{mat,b,\R}$ are convex surfaces, where
$V^{mat,a\cup c,\R}$ is reparametrized by 
$(s,\delta_2,r)\in[-2,2]\times\{\pm 1\}\times (-1,1)$ with 
${b_-}/{|b_-|}=\delta_2,\ ib_2={r}/{(1-r^2)}$, 
and $V^{mat,b,\R}$ is reparametrized by 
$(s,r)\in(\R_{<-2}\cup\R_{>2})\times (\frac{-3\pi}{2},\frac{\pi}{2}]$ with 
$b_-=e^{ir}$.

(b) For fixed $b_-\in S^1$, and any $k\in\Z$ if $b_-\neq\pm 1$, or any
$k\in \Z_{\leq 0}$ if $b_-=\pm 1$, then $x_k(s,B)$ tends to $\iiii$ in both
limits $s\to\iiii$ and $s\to -\iiii$.

(c) For fixed $s_1\in[-2,2]$ and any $k\in\Z_{\geq 1}$, then $x_k(s_1,B)$ tends to $0$ in any of the four limits 
$b_-\to -\iiii, \ b_-\nearrow 0, \ b_-\searrow 0, \ b_-\to +\iiii$.
\end{conjecture}

\begin{remarks}\label{t18.7}
(i) In conjecture \ref{t18.6} (c), 
a priori some $x_k(s_1,B)$ could also have a finite limit or even tend to $+\iiii$.
Neither the asymptotic formulae in \cite[ch.\  11]{IN86} nor the results on 
nilpotent orbits of TERP structures in \cite{Mo11b} are formulated in a way which
would control the behaviour of $x_k(s_1,B)$ for the four limits of $b_-$.
Nor can conjecture \ref{t18.6} (b)  be answered from the asymptotic formulae
in \cite[ch.\  11]{IN86}.

(ii) If conjectures \ref{t18.6} (b) and (c) are correct, 
then also the following picture is essentially correct. 

It is a conjectural sketch of the restriction to some $x_0\in\R_{>0}$ of 
$M_{3FN,\R}^{mon}$, that is 
\[
M_{3FN,\R}^{mon}|_{x=x_0} \cong \{x_0\}\times V^{mat,\R},
\]
together with the lines which are the intersection of $M_{3FN,\R}^{mon}|_{x=x_0}$
with $M_{3FN,\R}^{sing}$.  The regions where $f$ is positive or negative are again white or shaded.
$V^{mat,\R}$ is a sphere with four holes. In the picture all four holes are shifted to the front.
The back of the sphere is not visible. It is part of the shaded region.

Conjecture \ref{t18.6} (b) gives the spirals around the two holes $s=\pm\iiii$. 
For large $x_0$ they are small, for small $x_0$ they are large, and the spiral
around $s=\varepsilon\MGcdot\iiii$ for $\varepsilon\in\{\pm 1\}$ approaches 
 the line $s=2\varepsilon$ near $B={\bf 1}_2$ and near $B=-{\bf 1}_2$ as $x_0\to 0$.

Conjecture \ref{t18.6} (c) gives the spirals around the two holes $ib_2=\pm\iiii$. 
For small $x_0$ they are small, for large $x_0$ they are large, and the spirals
around both holes approach  both lines $B=\pm{\bf 1}_2$ as $x_0\to\iiii$.
\end{remarks}

%{\sc Later pic18-4jpg}
\includegraphics[width=0.9\textwidth]{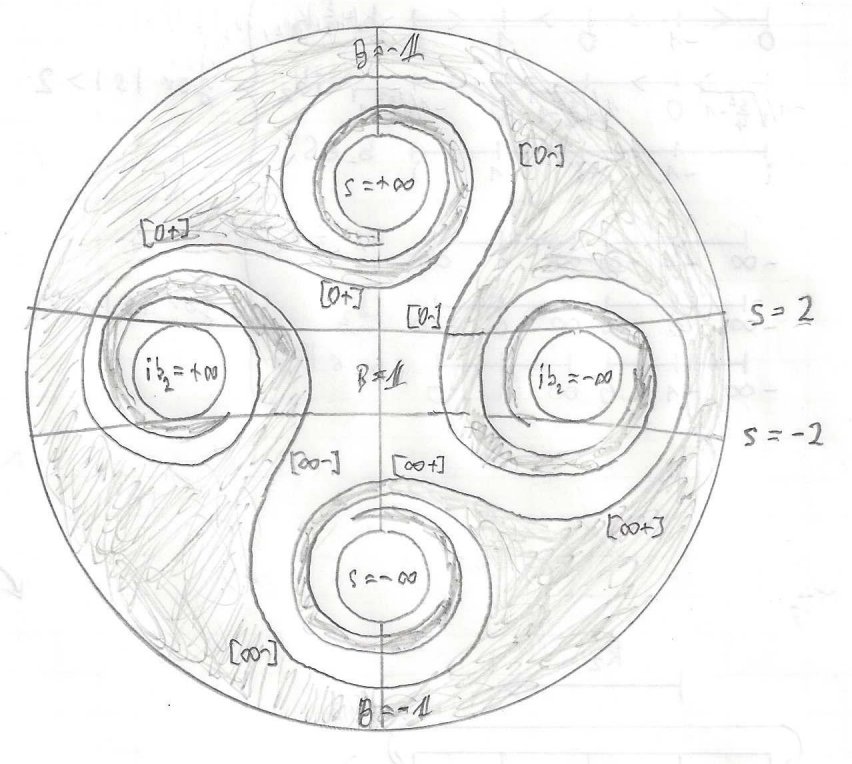}

\renewcommand{\indexname}{Index of notation}
\begin{theindex}

%\indexspace

\item $[1]$ \hfill \eqref{2.28}
\item $\orrr{>0},\orrr{<0},\orrr{[0+][\iiii-]}$ 
\hfill \ref{t1.6}, \ref{t18.1}, \ref{t18.2}, \ref{t18.4}
\item $\olll{>0},\olll{<0},\olll{[0+][0-]}$ 
\hfill \ref{t1.6}, \ref{t18.1}, \ref{t18.2}, \ref{t18.4}

\item $A$ \hfill \eqref{7.1}
\item $\alpha^j_0,\alpha^j_\iiii$ \hfill \ref{t2.2}
\item $\alpha_\pm$ \hfill \eqref{5.5}, \eqref{5.6}

\item $B$ \hfill \eqref{2.24}, \eqref{2.25}, \eqref{5.22}, \eqref{6.13}
\item $\www b_1$ \hfill \eqref{5.25}
\item $b_\pm$ \hfill \eqref{5.24}, \eqref{12.5}
\item $\www b_2$ \hfill \eqref{12.11}
\item $\www b_-$ \hfill \eqref{12.5}
\item Bessel equation \hfill \eqref{13.36}
\item $\beta$ \hfill \eqref{2.23}

\item $C$ \hfill \eqref{8.5}, \eqref{8.6}
\item $(C_{ij})$ \hfill \ref{t3.1}
\item $\C^{[sto,a]},\C^{[sto,J]}$ \hfill \eqref{12.1}\item
\item $c^{uc*}, c^{uc}, c^{beta}$ \hfill \eqref{4.1}, \eqref{4.4}
\item $\uuuu\chi$ \hfill \ref{t13.2}
\item compatible real and Stokes structure \hfill \ref{t17.6} (a)
\item completely reducible \hfill \ref{t3.1}

\item $\delta^{NI}$ \hfill \eqref{12.24}

\item $E^{NI}$ \hfill \eqref{11.19} 
\item elementary section \hfill \eqref{13.6}
\item $\uuuu{es}^{(k)}=(es_1^{(k)},es_2^{(k)})$ \hfill \eqref{13.6}
\item $\uuuu e^\pm_0,\uuuu e^\pm_\iiii$ \hfill \eqref{2.12}, \eqref{2.21}, 
\ref{t6.3}, \ref{t7.3}
\item $\varepsilon_0,\varepsilon_1,\varepsilon_2$ \hfill
\eqref{6.12}, \eqref{7.9}, \eqref{7.13}, \ref{t8.2} (b)
\item Euler orbit \hfill \ref{t17.1}

\item $\uuuu f^+_0$ \hfill \eqref{13.5}
\item $f_k$ \hfill \ref{t8.2} (b), \eqref{8.20}
\item folding transformation \hfill \ref{t9.2} (iii)+(iv)
\item $f_{mult}(.,s,B)$ \hfill \ref{t10.2} (iii), \ref{10.3}
\item $f_{univ}$ \hfill \eqref{10.9}

\item $g_k$ \hfill \ref{t8.2} (b), \eqref{8.20}
\item $\gamma$ \hfill \eqref{6.3}, \eqref{16.1}
\item $\gamma_1,\gamma_2,\gamma_3$ \hfill \eqref{18.5}--\eqref{18.7}
\item $G$ \hfill \eqref{6.15}, proof of \ref{t7.3}
\item $G^{mon}$ \hfill \eqref{7.30}
\item $g_{mult}$ \hfill \ref{t10.2} (iii), \ref{t10.3}

\item $H$ \hfill \ref{t2.1}
\item $H', H'_\R$ \hfill \ref{t16.2} (i)
\item Hankel function \hfill \eqref{13.38}, \eqref{13.39}
\item $H_\nu^{(1)}(t),H_\nu^{(2)}(t)$ \hfill \eqref{13.38}, \eqref{13.39}

\item $I^a_0,I^b_0,I^+_0,I^-_0$ \hfill \eqref{2.1}
\item $I^a_\iiii,I^b_\iiii,I^+_\iiii,I^-_\iiii$ \hfill \eqref{2.2}
\item $j$ \hfill \eqref{6.1}, \eqref{6.5}, \eqref{7.1}
\item $J$ \hfill \eqref{7.4}

\item $k$ \hfill \ref{t8.2} (b)
\item $\kappa_{a_1,a_2}$ \hfill \eqref{12.8}, \eqref{12.18}

\item $L,L^{\pm j}_0,L^{\pm j}_\iiii$ \hfill \ref{t2.2}
\item $\lambda_\pm$ \hfill \eqref{5.3}

\item $m_{[1]}$ \hfill \eqref{7.28}, \eqref{10.6}, \eqref{10.11}
\item $\Mon$ \hfill \eqref{2.22}
\item $\Mon_0^{mat}(s)$ \hfill \eqref{1.17}, \eqref{5.1}

\item $M_{3FN}^{ini},M_{3FN}^{reg},M_{3FN}^{sing},M_{3FN}^{[k]}$ \hfill
\ref{t1.3}, \ref{t10.1}
\item $M_{3FN}^{mon}$ \hfill \eqref{1.18}, \ref{t10.1} 
\item $M_{3FN,\R},M_{3FN,S^1},M_{3FN,i\R_{>0}}$  \hfill \ref{t15.3}
\item $M_{3FN,\R}^{ini},M_{3FN,S^1}^{ini},M_{3FN,i\R_{>0}}^{ini}$  \hfill \ref{t15.4}
\item $M_{3FN,\R}^{mon},M_{3FN,S^1}^{mon},M_{3FN,i\R_{>0}}^{mon}$  \hfill \ref{t15.5}
\item mixed TERP structure \hfill \ref{t17.6} (c), \ref{t17.9}

\item $M_{3T}$ \hfill \ref{t7.6} (b)
\item $M_{3TJ}^{ini}, M_{3TJ}^{reg}, M_{3TJ}^{[\varepsilon_1]}$ \hfill
\ref{t8.2} (c), \ref{t8.4} (b)
\item $M_{3TJ}^{mon}$ \hfill ch. \ref{s1.6}, \ref{t7.5}
\item $M_{3TJ}^{[\varepsilon_1,\varepsilon_2\sqrt{c}]}(u^1_0,u^1_\iiii)$ 
\hfill \ref{t8.2} (e)

\item $\nabla$ \hfill \ref{t2.1}
\item nilpotent orbit \hfill \ref{t17.3}, \ref{t17.9}

\item $\Omega_k^{(0)},\Omega_k^{(\iiii)}$ \hfill \eqref{11.14}, \eqref{11.15}

\item $P$ \hfill \eqref{6.5}
\item $p,q$ \hfill \eqref{11.22}, \eqref{11.33}
\item $P_I,\dots,P_{VI}$  \hfill ch. \ref{s9} 
\item $P_{III}(\alpha,\beta,\gamma,\delta)$ \hfill \eqref{1.1}, \eqref{9.1}
\item $P_{III}(D_6),P_{III}(D_7),P_{III}(D_8),P_{III}(Q)$ \hfill ch. \ref{s9}

\item $P_{3D6}$ bundle \hfill \ref{t2.1}
\item $P_{3D6}$ monodromy tuple \hfill \ref{t2.2}
\item $P_{3D6}$ numerical tuple \hfill \ref{t2.6}
\item $P_{3D6}$-TEP bundle \hfill \ref{t6.1}
\item $P_{3D6}$-TEPA bundle \hfill \ref{t7.1}
\item $P_{3D6}$-TEJPA bundle \hfill \ref{t7.1}

\item $pr_{3FN}$ \hfill \eqref{1.19}, \ref{t10.1} 
\item $pr_{3TJ}$ \hfill \eqref{1.39}
\item $pr_{mat}$ \hfill \eqref{1.20} 
\item pure twistor \hfill \ref{t4.1} (iv)

\item $\Phi_{3FN},\Phi_{3FN,\R}$ \hfill \eqref{1.21}, \eqref{1.25} 
\item $\pi_{orb}$ \hfill \eqref{17.1}
\item $\uuuu\psi$ \hfill \ref{t8.2} (d)
\item $\Psi_k^{(0)},\Psi_k^{(\iiii)}$ \hfill \eqref{11.16},  \eqref{11.17}

\item $R_1,R_2,R_3$ \eqref{7.30}, \eqref{8.38}, \eqref{12.13}
\item $R_4$ \hfill \ref{t14.1} (a), \ref{t14.3}
\item $R_5$ \hfill \ref{t14.2} (a), \ref{t14.3}
\item radial sinh-Gordon equation 
\hfill \eqref{9.5}, \eqref{15.1}--\eqref{15.2}
\item reducible \hfill \ref{t3.1}
\item $\rho_c$ \hfill \eqref{6.2}, \eqref{7.4}

\item $S$ \hfill \eqref{1.17}, \eqref{5.1}
\item $S^{NI}$ \hfill \eqref{11.13}
\item $S^a_0, S^b_0, S^a_\iiii, S^b_\iiii$ \hfill \eqref{2.20}, \eqref{2.21}
\item $S_k^{(0)},S_k^{(\iiii)}$ \hfill \eqref{11.18}
\item $\sigma$ \hfill \eqref{6.4}
\item $\uuuu \sigma_0,\uuuu \sigma_k$ \hfill \ref{t8.2} (b)
\item $\sqrt{\frac{s^2}{4}-1}$ \hfill ch. \ref{s5}
\item Sabbah orbit \hfill \ref{t17.3}
\item sine-Gordon equation \hfill \eqref{11.4}, \eqref{15.3}

\item $\tau$ \hfill \eqref{16.1}, \ref{t16.5}
\item $T$ \hfill \ref{t17.5}
\item $T_k$ \hfill \eqref{4.6}
\item $T(s)$ \hfill \ref{t5.3}
\item $t^{NI}$ \hfill \eqref{5.17}, \eqref{12.24}
\item TEP bundle \hfill \ref{t6.1}
\item TEPA bundle \hfill \ref{t7.1}
\item TEJPA bundle \hfill \ref{t7.1}

\item TERP(0) bundle \hfill ch. \ref{s1.7}, \ref{16.1}, ch. \ref{s17}
\item trace free \hfill \ref{t4.1} (ii)+(iii)
\item $(1,-1)$ twistor \hfill \ref{t4.1} (iv)

\item $u, u^{NI}$ \hfill \eqref{11.3}, \eqref{11.1}
\item $U, U^{beta}$ \hfill \eqref{4.4}
\item $u^1_0,u^1_\iiii,u^2_0,u^2_\iiii$ \hfill \ref{t2.1}
\item $u_k^{\delta_2}, d_k^{\delta_2}, l_k^{\delta_2}, r_k^{\delta_2}$
\hfill proof of \ref{t18.4}

\item $v_\pm$ \hfill \eqref{5.4}
\item $V^{mat},V^{mon}$ \hfill \eqref{1.15}, \eqref{7.26}, \eqref{7.25} 
\item $V^{mat,a},V^{mat,b\pm},V^{mat,c\pm}$ \hfill \eqref{1.30}, \eqref{12.2} 
\item $V^{mat,i\R_{>0}},V^{mat,S^1}$ \hfill \eqref{15.12}, \eqref{15.15} 
\item $V^{mat,\R}$ \hfill \eqref{1.23}, \eqref{15.6} 

\item $\uuuu\varphi_k$ \hfill \ref{t8.2} (f)
\item $\uuuu v_0,\uuuu v_\iiii$ \eqref{8.1}--\eqref{8.4}

\item $\uuuu w_0,\uuuu w_\iiii$ \hfill \eqref{8.7}--\eqref{8.10}

\item $x^{NI}$ \hfill \eqref{11.2}

\item $\zeta_0,\zeta_\iiii$ \hfill \eqref{2.1}, \eqref{2.2}
\item $\zeta^{NI}$ \hfill \eqref{11.20}, \eqref{11.29}

\end{theindex}

\end{document}